\documentclass[a5paper,11pt]{book}
\usepackage{amssymb}
\usepackage{amsmath}
\usepackage{theorem}
\usepackage{xspace}
\usepackage{cite}
\usepackage{endnotes}
\usepackage{mathrsfs}
\usepackage[russian,english]{babel}

\renewcommand{\notesname}{Comments}
\makeatletter
\renewcommand{\enoteheading}{\chapter*{\notesname
  \@mkboth{\notesname}{\notesname}
  \addcontentsline{toc}{chapter}{Comments} }     \leavevmode\par\vskip-\baselineskip}
\renewcommand{\@makeenmark}{\hbox{$^{\@theenmark )}$}}
\renewcommand{\enoteformat}{\rightskip\z@ \leftskip\z@ \parindent=1.8em
     \leavevmode\llap{\hbox{$^{\@theenmark )}$}}}
\makeatother

\renewcommand{\thechapter}{\Roman{chapter}}

\renewcommand{\chaptermark}[1]{\markboth{\chaptername\ \thechapter. \ #1}{}}

\makeatletter
\renewcommand{\@makechapterhead}[1]{  \vspace*{20\p@}  {\parindent \z@ \raggedright \normalfont
    \ifnum \c@secnumdepth >\m@ne
      \if@mainmatter
        \large \bfseries \@chapapp\space \thechapter        \par\nobreak
        \vskip 10\p@      \fi
    \fi
    \interlinepenalty\@M
    \LARGE \bfseries #1\par\nobreak    \vskip 40\p@
  }}
\renewcommand{\@makeschapterhead}[1]{  \vspace*{20\p@}  {\parindent \z@ \raggedright
    \normalfont
    \interlinepenalty\@M
    \LARGE \bfseries  #1\par\nobreak
    \vskip 40\p@
  }}
\makeatother
\makeatletter
\def\@chapter[#1]#2{\ifnum \c@secnumdepth >\m@ne
                       \if@mainmatter
                         \refstepcounter{chapter}                         \typeout{\@chapapp\space\thechapter.}                         \addcontentsline{toc}{chapter}                                   {\protect \chaptername\ \thechapter.\ \ #1}                       \else
                         \addcontentsline{toc}{chapter}{#1}                       \fi
                    \else
                      \addcontentsline{toc}{chapter}{#1}                    \fi
                    \chaptermark{#1}                    \addtocontents{lof}{\protect\addvspace{10\p@}}                    \addtocontents{lot}{\protect\addvspace{10\p@}}                    \if@twocolumn
                      \@topnewpage[\@makechapterhead{#2}]                    \else
                      \@makechapterhead{#2}                      \@afterheading
                    \fi}
\makeatother
\makeatletter
\def\@sect#1#2#3#4#5#6[#7]#8{  \ifnum #2>\c@secnumdepth
    \let\@svsec\@empty
  \else
    \refstepcounter{#1}    \protected@edef\@svsec{\@seccntformat{#1}\relax}  \fi
  \@tempskipa #5\relax
  \ifdim \@tempskipa>\z@
    \begingroup
      #6{        \@hangfrom{\hskip #3\relax\@svsec}          \interlinepenalty \@M #8\@@par}    \endgroup
    \csname #1mark\endcsname{#7}    \addcontentsline{toc}{#1}{      \ifnum #2>\c@secnumdepth \else
        \protect\numberline{\S \csname the#1\endcsname}      \fi
      #7}  \else
    \def\@svsechd{      #6{\hskip #3\relax
      \@svsec #8}      \csname #1mark\endcsname{#7}      \addcontentsline{toc}{#1}{        \ifnum #2>\c@secnumdepth \else
          \protect\numberline{\csname the#1\endcsname}        \fi
        #7}}  \fi
  \@xsect{#5}}
\makeatother
\newtheorem{Pa}{Paper}[section]

\newtheorem{corollary}[Pa]{Corollary}

\newtheorem{definition}[Pa]{Definition}

\newtheorem{lemma}[Pa]{Lemma}

\begin{document}

\setcounter{endnote}{8}

\renewcommand{\contentsname}{Table of contents} \renewcommand{%
\bibname}{References}

\thispagestyle{empty}

\newpage

\begin{center}
{\large BELARUSIAN STATE UNIVERSITY}
\quad \\
\quad \\
\quad \\
\bigskip {\bf \LARGE Alexander Kiselev }\medskip\\
 {\large \emph{e-mail}\ \ aakiselev@yahoo.com}
\quad \\
\quad \\
\quad \\
{\bf \Huge Inaccessibility}
\quad \\
\quad \\
{\bf \Huge and}
\quad \\
\quad \\
{\bf \Huge Subinaccessibility}\\
\quad \\
\quad \\
{\Large In two parts} \\
\smallskip
{\Large Part II} \\
\quad \\
\quad \\
{\em \large Second edition} \\
\smallskip
{\em \large enriched and improved} \\
\quad \\
\quad \\
\quad \\
\quad \\
\quad \\
\quad \\
\quad \\
{Minsk} \\
\smallskip
{``Publishing center of BSU''} \\
\smallskip
{2010}

\end{center}

\newpage

\thispagestyle{empty}

\noindent UDK 510.227
\\

{\small \textbf{Kiselev Alexander.} Inaccessibility and
Subinaccessibility. In 2 pt. Pt. II / Alexander Kiselev. -- 2nd
ed., enrich. and improv. -- Minsk: Publ. center of BSU, 2010. --
\pageref{end} p. -- ISBN 978-985-476-596-9. }
\\

{\footnotesize The work presents the second part of the second
edition of its previous one published in 2000 under the same
title, containing the proof (in $ZF$) of the inaccessible
cardinals nonexistence, which is enriched and improved now.

This part contains applications of the subinaccessible cardinals
apparatus and its basic tools —-- theories of reduced formula
spectra and matrices, disseminators and others, which are used
here in this proof and are set forth now in their more transparent
and simplified form. Much attention is devoted to the explicit and
substantial development and cultivation of basic ideas, serving as
grounds for all main constructions and reasonings. The proof of
the theorem about inaccessible cardinals nonexistence is presented
in its detailed exposition. Several easy consequences of this
theorem and some well-known results are presented.

Ref. 38. }
\quad \\
\begin{center}
{R\;e\;f\;e\;r\;e\;e\;s:}
\\
Prof. {\em Petr P. Zabreiko;}
\\
Prof. {\em Andrei V. Lebedev}
\quad \\
\quad \\
\quad \\
\noindent Mathematics Subject Classification (1991): \\
03E05, 03E15, 03E35, 03E55, 03E60
\quad \\
\quad \\
\end{center}

\quad \\
\quad \\
\quad \\
\quad \\
\noindent {\scriptsize {\bf ISBN 978-985-476-596-9 (pt. 1)} \hfill
$\copyright$ Kiselev Alexander A., 2000}

\noindent {\scriptsize {\bf ISBN 978-985-476-597-6} \hfill
$\copyright$ Kiselev Alexander A., 2008, with modifications}

\hfill {\scriptsize $\copyright$ Kiselev Alexander A., 2010, with
modifications}

\newpage

\thispagestyle{empty}

\hfill \textit{To my mother Ann}

\newpage

\thispagestyle{empty}

\chapter*{Acknowledgements}
\hspace*{1em} The author sends his words of deep gratitude to
Hanna Ca\-lien\-do for understanding of the significance of the
theme, for hearty encouraging help and for the help in promotion
of the work.

The very special appreciation is expressed to Prof. Sergeiy
Kogalovskiy, who had taught the author the Hierarchy Theory, and
to Prof. Akihiro Kanamori for their invaluable and strong
encouragement that gave the necessary impetus for the completing
the work.

The author would like to express his very special gratitude to
Prof. Andrei Lebedev and Prof. Petr Zabreiko who have supported
the advancing of the work for many years; the very strong
intellectual and moral debt is expressed to both of them for the
spiritual and practical help.

Many deep thanks to Prof. Alexander Lapcevich, the Headmaster of
Minsk State Higher Aviation College, to Prof. Sergey Sizikov, the
Deputy Headmaster, to Prof. Anton Ripinskiy, the Dean of the Civil
Aviation Faculty, to Prof. Aleksey Kirilenko, the Chief of the
Sciences Department, for giving favorable conditions for the
fruitful work at the issue.

The deep gratitude is expressed to Prof.~Alexander Tuzikov and
Dr.~Yuri Prokopchuk who gave the author extensive and expert help
in typing of the text.

Prof. Vasiliy Romanchack gave the author the financial and moral
help during the most complicated period of investigations of the
theme and the author sends him many words of deep thankfulness.

Many thanks also goes to Ludmila Laptyonok who went through many
iterations of the difficult typing of previous author's works and
preparing this work.

But the most deep hearty thankfulness goes to Nadezhda Zabrodina
for the support and encouragement for a long time, without which
the work would be considerably hampered.

The host of other people who over the years provided promotion or
encouragement of the work, is too long to enumerate and thanks are
expressed to all of them.

\newpage
{} \thispagestyle{empty}

\tableofcontents

\newpage {} \thispagestyle{empty}

\newpage

\chapter*{Introduction}

\addcontentsline{toc}{chapter}{Introduction} %
\markboth{Introduction}{Introduction}

\setcounter{equation}{0}

\hspace*{1em} This work represents the direct continuation of its
previous Part~I \cite{Kiselev9} and constitutes the uniform text
with it. The author has considered it natural to organize this
work in such a manner, that it continues all enumerations
of~\cite{Kiselev9} and among them enumerations of chapters,
paragraphs, definitions, statements and formulas and even
enumerations of comments and references as well.
\\
Also formulas, notions or symbols used in this work without
explanations have been already introduced in
Part~I~\cite{Kiselev9} or generally accepted or used in remarkable
text of Jech~\cite{Jech} ``Lectures in Set Theory with Particular
Emphasis on the Method of Forcing'', providing many basic notions
and information and much more; therefore here they are assumed to
be known and will be used frequently without comments.
\\
So, it would be more convenient for the reader to familiarize
preliminarily with~\cite{Kiselev9} and with the main notions and
denotations of this work at least in outline.
\\
Anyway, it would be useful to get acquainted with the plan of all
the work beforehand and with the outline of developing basic ideas
as they are presented in~\cite{Kiselev9} on pp.~9--18.
\\
With this agreement in hand, the reader should remember that all
backward references to paragraphs with numbers less than 7 concern
Part~I~\cite{Kiselev9}, analogously for enumerated statements and
so on.

As to the content of this work and to the technical side of the
matter it should be noticed that it follows its previous
edition~\cite{Kiselev8} of 2000, but in the more systematic way.

Also it should be pointed out, that in this
edition~\cite{Kiselev8} and foregoing works the author tendered to
avoid the usage of new terminology, especially private notions and
symbols of his own, over certainly necessary, because he was
uneasy about difficulties and unacceptance it may cause for the
reader.

Nevertheless, the perception of the previous
edition~\cite{Kiselev8} by the readers showed, that apprehensions
of this kind are inappropriate and, so, this usage is unavoidable
all the same.

Therefore in the present work the author has taken another stand
and has considered more natural to involve all the system of his
own private concepts and definitions in the whole, which he has
developed since 1976, because it possesses the technical and
conceptual expressiveness and leads to the very essence of the
matter and, so, it would be too artificial to avoid its usage (see
comment 3~\cite{Kiselev9} as an example). Some statements have
received their strengthening; some details considered obvious in
the author's previous works, here have received their refinement;
some arguments have changed their places for the more suitable;
also some suitable redesignations are used.

But as for main constructions, one should note, that the present
work follows the edition~\cite{Kiselev8} of 2000 but in the more
clarified way; the main result of this work: the system
\[
    ZF+\exists k \hspace{2mm} (k \; \mbox{\it is weakly
    inaccessible cardinal})
\]
is inconsistent; all the reasonings  are carried out in this
theory. All weakly inaccessible cardinals become strongly
inaccessible in the constructive class \ $L$ \ and therefore the
reasonings are carried over to the standard countable basic model
\[
    \mathfrak{M}=(L_{\chi ^{0}},\; \in , \; =)
\]
of the theory
\[
    ZF+V=L+\exists k \hspace{2mm} (k \; \mbox{\it is weakly
    inaccessible cardinal}),
\]
and further \ $k$ \ is the smallest inaccessible cardinal in \
$\mathfrak{M}$. \ Actually only the formulas of the limited length
from this theory are used; moreover, the countability of this
model is required only for some technical convenience and it is
possible to get along without it (see
``Preliminaries''~\cite{Kiselev9}). In this model \ $\mathfrak{M}$
\ the so called matrix functions are constructed, possessing
simultaneously the two properties of monotonicity and
nonmonotonicity, that provides
\\

\noindent {\bf Main theorem } ($ZF$)
\\
\hspace*{1cm} {\it There are no weakly inaccessible cardinals.}
\\

\noindent It implies the nonexistence of strongly inaccessible
cardinals and therefore nonexistence of all other large cardinals.
These matrix functions are constructed and treated by means of the
elementary language from the formula classes (see definition
2.1~\cite{Kiselev9}) of some fixed level $>3$ over the standard
model
\[
    (L_k, \in, =)
\]
and further all constructions are carried out by means of this
language (if the opposite is not pointed out by the context).
\\

\noindent In addition in \S 12 some easy consequences of the Main
theorem and some well-known results are presented.

\setcounter{chapter}{1}

\chapter{Special Theory: Matrix Functions}

\markboth{\chaptername\ \thechapter. \ Special Theory}{\chaptername\
\thechapter. \ Special Theory} %

\setcounter{section}{6}

\section{Matrix \ $\protect\delta$\,-Functions}

\setcounter{equation}{0}

\hspace*{1em} Here we are going to start the further development
of the idea of the Main theorem proof and to modify the simplest
matrix functions \ $S^{<\alpha_1}_{\chi f} $ \ (see definition
5.14~\cite{Kiselev9}) in such a way that their new specialized
variants -- the so called \ $\alpha $-functions -- will provide
the required contradiction: they will possess the property of \
$\underline{\lessdot}$-monotonicity and at the same time will be
deprived of it.

Let us remind that the simplest matrix functions, which were
considered in \S~5~\cite{Kiselev9}, possess the property of
monotonicity, but it came out that the direct proof of the
required contradiction -- the proof of their nonmonotonicity -- is
hampered by the following obstacle: some essential properties of
lower levels of universe do not extend up to jump cardinals of
matrices on their carriers, which are values of the matrix
functions.
\\
In order to destroy this obstacle we shall equip such matrices
with their corresponding disseminators and as the result the
simplest matrix functions will be transformed to their more
complicated forms, \ $\alpha$-function.

However, the direct formation of these functions seems to be
considerably complicated and some their important singularities
unmotivated.
\\
Therefore in order to represent their introduction in the more
transparent way we shall beforehand undertake the second approach
to the idea of the Main theorem proof and turn attention to their
more simple forms, that is to the \ $\delta $-functions.
\\
To this end we shall apply results of \S 6~\cite{Kiselev9} for \
$m=n+1$ \ and the fixed level \ $n>3$, \ but the notion of
disseminator matrix should be sharpened; all disseminators in what
follows will be of of the level \ $n+1$ (see definition
6.9~\cite{Kiselev9}).

\begin{definition}
\label{7.1.} \hfill {} \newline
\hspace*{1em} Let
\begin{equation*}
\gamma < \alpha <\alpha _{1}\leq k.
\end{equation*}
\noindent \emph{1)}\quad We denote through \
$\mathbf{K}_{n}^{\forall <\alpha _{1}}(\gamma ,\alpha )$ \ the
formula:

\vspace{6pt}
\begin{equation*}
SIN_{n-1}^{<\alpha _{1}}(\gamma )\wedge \forall \gamma ^{\prime }\leq \gamma
\ (SIN_{n}^{<\alpha _{1}}(\gamma ^{\prime })\longrightarrow SIN_{n}^{<\alpha
}(\gamma ^{\prime }))~.
\end{equation*}
\vspace{0pt}

\noindent If this formula is fulfilled by the constants \ $\gamma
$, \ $\alpha $, \ $\alpha _{1}$, \ then we say that \ $\alpha $ \
conserves \ $SIN_{n}^{< \alpha _{1}}$-cardinals \ $\leq \gamma $ \
below \ $\alpha _{1}$.
\\
If \ $S$ \ is a matrix on a carrier \ $\alpha $ \ and its prejump
cardinal \ $\alpha _{\chi }^{\Downarrow }$ \ after \ $\chi$ \
conserves these cardinals, then we also say that \ $S$ \ on \
$\alpha $ \ conserves these cardinals below \ $\alpha _{1}$.
\newline

\noindent \emph{2)}\quad We denote through \
$\mathbf{K}_{n+1}^{\exists }(\chi ,\delta ,\gamma ,\alpha ,\rho
,S)$ \ the \ $\Pi _{n-2}$-formula:

\begin{equation*}
\sigma (\chi ,\alpha ,S)\wedge Lj^{<\alpha }(\chi )\wedge \chi
<\delta <\gamma <\alpha \wedge S\vartriangleleft \rho \leq \chi
^{+}\wedge \rho = \widehat{\rho }\wedge
\end{equation*}
\begin{equation*}
\wedge SIN_{n}^{<\alpha _{\chi }^{\Downarrow }}(\delta )\wedge
SIN_{n+1}^{<\alpha _{\chi }^{\Downarrow }}\left[ <\rho \right] (\delta ).
\end{equation*}
\vspace{0pt}

\noindent Here, remind, the \ $\Pi_{n-2}$-formula \ $\sigma(\chi,
\alpha, S)$ \ means that \ $S$ \ is the singular matrix on its
carrier \ $\alpha$ \ reduced to the cardinal \ $\chi$ \ (see
definition 5.7~\cite{Kiselev9}); \ $\delta$ \ is the disseminator
for \ $S$ \ on \ $\alpha$ \ with the base \ $\rho$ \ of the level
\ $n+1$ \ (definition 6.9~\cite{Kiselev9}); the upper indices \ $<
\alpha_{\chi}^{\Downarrow}$ \ mean the bounding of the formula
quantors under consideration by the prejump cardinal \
$\alpha_{\chi}^{\Downarrow}$ \ (see also definitions 2.3,
5.9~\cite{Kiselev9}); \ $\widehat{\rho}$ \ is the closure of \
$\rho$ \ under the pair function; and \ $Lj^{<\alpha}(\chi)$ \ is
the \ $\Delta_1$-property of the cardinal \ $\chi$ \ saturation
below \ $\alpha$ \ (see definition~6.9~4) \cite{Kiselev9}):
\[
    \chi < \alpha \wedge SIN_{n-1}^{<\alpha}(\chi) \wedge
    \Sigma rng\big(\widetilde{\mathbf{S}}_n^{sin \vartriangleleft \chi}\big)
    \in B_{\chi} \wedge \sup dom \big(\widetilde{\mathbf{S}}_n^{sin \vartriangleleft
    \chi}\big) = \chi.
\]
We denote through \ $\mathbf{K}^{<\alpha _{1}}(\chi ,\delta
,\gamma ,\alpha ,\rho ,S)$ \ the formula:

\vspace{6pt}
\begin{equation*}
\mathbf{K}_{n}^{\forall <\alpha _{1}}(\gamma ,\alpha _{\chi }^{\Downarrow
})\wedge \mathbf{K}_{n+1}^{\exists \vartriangleleft \alpha _{1}}(\chi
,\delta ,\gamma ,\alpha ,\rho ,S)\wedge \alpha <\alpha _{1}~.
\end{equation*}
\vspace{0pt}

\noindent \emph{3)}\quad If this formula is fulfilled by the
constants \ $\chi$, $\delta $, $\gamma $, $\alpha $, $\rho$, $S$,
$\alpha _{1}$, \ then we say that \ $\chi $, $\delta $, $\alpha $,
$\rho $, $S$ \ are strongly admissible for \ $\gamma $ \ below \
$\alpha _{1}$.
\newline If some of them are fixed or meant by the context, then
we say that the others are also strongly admissible for them (and
for \ $\gamma $) \ below \ $\alpha _{1}$. \newline

\noindent \emph{4)}\quad The matrix \ $S$ \ is called strongly
disseminator matrix or, briefly, \ $\delta $-matrix strongly
admissible on the carrier \ $\alpha $ \ for \ \mbox{$\gamma
=\gamma _{\tau }^{<\alpha _{1}}$} \ below \ $\alpha _{1}$, \ iff
it possesses some disseminator \ $\delta <\gamma $ \ with a base \
$\rho $ \ strongly admissible for them (also below \ $\alpha
_{1}$).
\newline In every case of this kind \ $\delta $-matrix
is denoted by the common symbol \ $\delta S$ \ or \ $S$.
\\
If \ $\alpha_1=k$, \ or \ $\alpha_1$ \ is pointed out by the
context, then the upper indices \ $<\alpha_1$, $\vartriangleleft
\alpha_1$ \ here and other mentionings about \ $\alpha_1$ \ are
dropped.
\\
\hspace*{\fill} $\dashv$
\end{definition}

\noindent Further up to the end of $\S$~7 the notions of
admissibility and of \ $\delta$-matrices will be considered to be
\textit{strongly} notions, so the term ``strongly'' will be
omitted. All matrices will be considered to be \
$\delta$-matrices; as the reducing cardinal \ $\chi $ \ in what
follows will be used the complete cardinal \ $\chi ^{\ast }$ \
(see definition~5.4~\cite{Kiselev9}) -- if the context will not
indicate some other case.
\\
Here one should pay attention also to the notion of the cardinal \
$\chi$ \ saturation below \ $\alpha$, \ that is to the \
$\Delta_1$-property \ $Lj^{<\alpha}(\chi)$; \ from lemma
5.5~\cite{Kiselev9} it follows, that \ $\chi^{\ast}$ \ is the
cardinal saturated below any \ $\alpha > \chi^{\ast}$, $\alpha \in
SIN_{n-2}$.
\\
The symbol \ $\chi ^{\ast }$ \ in notations and formula writings
will be often omitted for some shortening.
\\
Further every bounding cardinal \ $\alpha _{1}$ \ will belong to \
$SIN_{n-2}$ \ and hold the condition

\begin{equation}\label{e7.1}
\chi ^{\ast }<\alpha _{1} \leq k \wedge A_{n}^{\vartriangleleft
\alpha_{1}} ( \chi ^{\ast } ) =\left\| u_{n}^{\vartriangleleft
\alpha_{1}} ( \underline{l} ) \right\|,
\end{equation}

\noindent or \ $\alpha_1 = k$ \ (unless otherwise is specified by
the context).

The cardinal \ $\alpha_1 \le k$ \ here with this property will be
called {\em equinformative} (equally informative) with the
cardinal \ $\chi^{\ast}$.
\\
This term is introduced here because of the phenomenon: no \
$\Sigma_n$-proposition \ $\varphi(\underline{l})$ \ has
\textit{jump ordinals} after \ $\chi^{\ast}$ \ below \ $\alpha_1$
\ (see definition 2.4~\cite{Kiselev9}). \ It is not hard to see,
that it is equivalent to the following: for every generic
extension \ $\mathfrak{M}[l]$ \ every \ $\Pi_n$-proposition \
$\varphi(l)$ \ which holds in \ $\mathfrak{M}[l]$ \ below \
$\chi^{\ast}$, \ thereafter holds in this extension below \
$\alpha_1$ \ due to (\ref{e7.1}) and \ $\alpha_1 \in \Pi_{n-2}$; \
thus every \ $\Pi_n$-proposition \ $\varphi(l)$ \ holds or not in
both cases simultaneously for every generic extension \
$\mathfrak{M}[l]$ \ (see also comment 7~\cite{Kiselev9} to
illustrate the importance of this notion).
\\
One should pay attention to the important example of such
cardinal: the prejump cardinal \
$\alpha_{\chi^{\ast}}^{\Downarrow}$ \ after \ $\chi^{\ast}$ \ of
every matrix carrier \ $\alpha>\chi^{\ast}$. \ Besides, it will be
always assumed for \ $\chi^{\ast}$ \ and \ $\alpha_1$ \ that
\[
    \forall \gamma < \alpha_1 \exists \gamma^{\prime} \in \left[\gamma,
    \alpha_1\right[ \ SIN_n^{<\alpha_1}(\gamma^{\prime}) \wedge
    cf(\alpha_1) \ge \chi^{\ast +}
\]
for convenience of some formula transformations.
\\
The boundaries \ $ <\alpha _{1}$, \ $\vartriangleleft \alpha _{1}$
\ will be omitted, as usual, if \ $\alpha _{1}=k$, \ or \
$\alpha_1$ \ is meant by the context.

\begin{definition}
\label{7.2.} \hfill {} \newline
\hspace*{1em} Let \ $\chi ^{\ast }<\alpha _{1}$. \newline
\quad \newline
\emph{1)}\quad We call as the matrix \ $\delta $-function of the level \ $n$
\ below \ $\alpha _{1}$ \ reduced to \ $\chi ^{\ast }$ \ the \ function
\begin{equation*}
\delta S_{f}^{<\alpha _{1}}=(\delta S_{\tau }^{<\alpha _{1}})_{\tau }~
\end{equation*}
\noindent taking the value for \ $\tau$:
\begin{equation*}
\delta S_{\tau }^{<\alpha _{1}}=\min_{\underline{\lessdot }} \bigl
\{ S \vartriangleleft \chi^{\ast +} : \exists \delta ,\alpha ,\rho
< \gamma_{\tau+1}^{<\alpha_1} ~ \mathbf{K}^{<\alpha _{1}}(\delta
,\gamma _{\tau }^{<\alpha _{1}},\alpha ,\rho ,S) \bigr \};
\end{equation*}

\noindent \emph{2)}\quad the following accompanying ordinal
functions are defined below \ $\alpha _{1}$:
\begin{equation*}
\check{\delta}_{f}^{<\alpha _{1}}=(\check{\delta}_{\tau }^{<\alpha
_{1}})_{\tau };\quad \rho _{f}^{<\alpha _{1}}=(\rho _{\tau
}^{<\alpha _{1}})_{\tau };\quad \alpha _{f}^{<\alpha _{1}}=(\alpha
_{\tau }^{<\alpha _{1}})_{\tau }
\end{equation*}

\noindent taking the values:
\begin{equation*}
\check{\delta}_{\tau }^{<\alpha _{1}} = \min \{\delta <
\gamma_{\tau}^{<\alpha_1}: \exists \alpha ,\rho <
\gamma_{\tau+1}^{<\alpha_1} ~\mathbf{K}^{<\alpha _{1}}(\delta
,\gamma _{\tau }^{<\alpha _{1}},\alpha ,\rho ,\delta S_{\tau
}^{<\alpha _{1}})\};\quad
\end{equation*}
\vspace{-6pt}
\begin{equation*}
\rho _{\tau }^{<\alpha _{1}} = \min \{\rho < \chi^{\ast +} :
\exists \alpha < \gamma_{\tau+1}^{<\alpha_1} ~ \mathbf{K}^{<\alpha
_{1}}(\check{\delta}_{\tau }^{<\alpha _{1}},\gamma _{\tau
}^{<\alpha _{1}},\alpha ,\rho ,\delta S_{\tau }^{<\alpha _{1}})\};
\end{equation*}
\vspace{-6pt}
\begin{equation*}
\alpha _{\tau }^{<\alpha _{1}} = \min \{\alpha <
\gamma_{\tau+1}^{<\alpha_1} : \mathbf{K} ^{<\alpha
_{1}}(\check{\delta}_{\tau }^{<\alpha _{1}},\gamma _{\tau
}^{<\alpha _{1}},\alpha ,\rho _{\tau }^{<\alpha _{1}},\delta
S_{\tau }^{<\alpha _{1}})\}. \qquad
\end{equation*}

\noindent For each matrix \ $\delta S_{\tau }^{<\alpha _{1}}$ \
these functions define its generating disseminator \
$\check{\delta}_{\tau }^{<\alpha _{1}}<\gamma _{\tau }^{<\alpha
_{1}}$ \ along with its base \ $ \rho_{\tau }^{<\alpha _{1}}$ \
and its carrier \ \mbox{$\alpha _{\tau }^{<\alpha _{1}}$}.
\hspace*{\fill} $\dashv$
\end{definition}

\noindent Using lemma 6.8~\cite{Kiselev9} it is easy to see, that
here \ $\check{\delta}_{\tau}^{<\alpha _{1}}$ \ is the minimal
with the base
\[
    \rho_{\tau}^{<\alpha _{1}} = \widehat{\rho_1}, \
    \rho_1 = Od(\delta S_{\tau}^{<\alpha_{1}}),
\]
that is the closure of the ordinal \ $Od(\delta S_{\tau}^{<\alpha
_{1}})$ \ under the pair function; thereafter such disseminator is
called the {\em generating eigendisseminator of \ $\delta
S_{\tau}^{<\alpha _{1}}$} \ on \ $\alpha_{\tau}^{<\alpha _{1}}$ \
below \ $\alpha_1$ \ and is denoted through \
$\check{\delta}_{\tau}^{S < \alpha_1}$ \ (see also definition
6.9~2)~\cite{Kiselev9}), and its base \ $\rho_{\tau}^{<\alpha_1}$
\ is denoted through \ $\rho_{\tau}^{S <\alpha_1}$.

It is easy to obtain the following lemmas from these definitions
and lemmas 5.15, 5.16~\cite{Kiselev9}:

\begin{lemma}
\label{7.3.} \hfill {} \newline \hspace*{1em} For \ $\alpha
_{1}<k$ \ the formulas \ $\mathbf{K}_n^{\forall < \alpha_1}$, \
$\mathbf{K}^{<\alpha _{1}}$ \ belong to \ $\Delta _{1}$ \ and
therefore all functions
\begin{equation*}
\delta S_{f}^{<\alpha _{1}},\quad \check{\delta}_{f}^{<\alpha
_{1}},\quad \rho _{f}^{<\alpha _{1}},\quad \alpha _{f}^{<\alpha
_{1}}
\end{equation*}
are \ $\Delta _{1}$-definable through \ $\chi ^{\ast },\alpha
_{1}$.
\\
For \ $\alpha _{1}=k$ \ the formulas \ $\mathbf{K}_n^{\forall}$, \
$\mathbf{K}$ \ belong to \ $\Sigma _{n}$ \ and these functions are
\ $\Delta _{n+1}$-definable. \hspace*{\fill} $\dashv$
\end{lemma}

\begin{lemma}
\label{7.4.} \emph{(About \ $\delta $-function absoluteness)} \hfill {}
\newline
\hspace*{1em} Let \ $\chi ^{\ast }<\gamma _{\tau +1}^{<\alpha _{1}}<\alpha
_{2}<\alpha _{1}\leq k$, \quad $\alpha _{2}\in SIN_{n-2}^{<\alpha _{1}}$ \
and
\begin{equation*}
(\gamma _{\tau }^{<\alpha _{1}}+1)\cap SIN_{n}^{<\alpha _{2}}=(\gamma _{\tau
}^{<\alpha _{1}}+1)\cap SIN_{n}^{<\alpha _{1}},
\end{equation*}
then on the set
\begin{equation*}
\{\tau ^{\prime }:\ \ \ \chi ^{\ast }\leq \gamma _{\tau ^{\prime }}^{<\alpha
_{2}}\leq \gamma _{\tau }^{<\alpha _{1}}\}
\end{equation*}
the functions
\begin{equation*}
\delta S_{f}^{<\alpha _{2}},\quad \check{\delta}_{f}^{<\alpha
_{2}},\quad \rho _{f}^{<\alpha _{2}},\quad \alpha _{f}^{<\alpha
_{2}}
\end{equation*}
coincide respectively with the functions
\begin{equation*}
\delta S_{f}^{<\alpha _{1}},\quad \check{\delta}_{f}^{<\alpha
_{1}},\quad \rho _{f}^{<\alpha _{1}},\quad \alpha _{f}^{<\alpha
_{1}}.
\end{equation*}
\hspace*{\fill} $\dashv$
\end{lemma}

The following lemma and the reasoning proving it present the idea
which will be applied further in various significant typical
situations:

\begin{lemma}
\label{7.5.} \emph{(About disseminator)} \newline
\emph{1)}\quad Let
\begin{itemize}
\item[(i)]  \ $]\tau _{1},\tau _{2}[ \; \subseteq dom\bigl( \delta
S_{f}^{<\alpha _{1}}\bigr) ~,\quad \gamma _{\tau _{2}}\in
SIN_{n}^{<\alpha _{1}}$; \medskip

\item[(ii)]  \ $\tau _{3}\in dom\bigl( \delta S_{f}^{<\alpha _{1}}\bigr)
~,~~\tau _{2} \leq \tau _{3}$; \medskip

\item[(iii)]  \ $\check{\delta}_{\tau _{3}}^{<\alpha _{1}}<\gamma _{\tau
_{2}}^{<\alpha _{1}}$~. \medskip
\end{itemize}
Then
\begin{equation*}
\check{\delta}_{\tau _{3}}^{<\alpha _{1}}\leq \gamma _{\tau
_{1}}^{<\alpha _{1}}~.
\end{equation*}

\noindent \emph{2)}\quad Let \ $\delta$-matrix \ $S$ \ on its
carrier \ $\alpha$ \ be admissible for \
$\gamma_{\tau}^{<\alpha_1}$ \ along with its disseminator \
$\delta$ \ and base \ $\rho$ \ below \ $\alpha_1$, \ then:

\begin{itemize}
\item[(i)]  \ \ $\{\tau ^{\prime }: \ \delta
< \gamma _{\tau ^{\prime }}^{<\alpha _{1}}\leq \gamma _{\tau
}^{<\alpha _{1}}\}\subseteq dom\bigl( \delta S_{f}^{<\alpha
_{1}}\bigr);$ \medskip

\item[(ii)] this matrix \ $S$ \ along with the same \
$\delta$, \ $\rho$ \ possesses the minimal admissible carrier \
$\alpha^{\prime} \in \; ]\gamma _{\tau }^{<\alpha _{1}},\gamma
_{\tau +1}^{<\alpha _{1}}[$~.
\end{itemize}
\end{lemma}

\noindent \textit{Proof.} \ 1) The upper indices \ $< \alpha_{1}$,
$\vartriangleleft \alpha_{1}$ \ will be dropped. Let us consider
the matrix \ $S^3 = \delta S_{\tau_3}$ \ and \ $\check{\delta}^3 =
\check{\delta}_{\tau_3}$, \ $\rho^3 = \rho_{\tau_3}$. \ Suppose 1)
fails, then by $(iii)$
\begin{equation*}
\gamma_{\tau_{1}} < \check{\delta}^3< \gamma_{\tau_{2}} \quad
\mathrm{and} \quad \check{\delta}^3 =  \gamma_{\tau_{4}}
\end{equation*}
for some \ $\tau_{4} \in \; ] \tau_{1},\tau_{2} [ $. \ Let us
observe the situation below, standing on \ $\alpha^{3} =
\alpha_{\tau_{3}}^{\Downarrow}$\;. \ From $(i)$ and lemma~7.4 it
comes that
\begin{equation*}
\delta S_{f}^{<\alpha^{3}} \equiv \delta S_{f} \mathrm{\ \ on\ \ }
]\tau_1, \tau_2[
\end{equation*}
and the matrix \ $S^4 = \delta S_{\tau_{4}}^{<\alpha^{3}}=\delta
S_{\tau_{4}} $ \ on the carrier \
$\alpha_{\tau_{4}}^{<\alpha^{3}}=\alpha_{\tau_{4}} $ \ has the
disseminator
\begin{equation*}
\check{\delta}^{4} = \check{\delta}_{\tau_{4}}^{<\alpha^{3}} =
\check{\delta} _{\tau_{4}} < \gamma_{\tau_{4}} = \check{\delta}^3
\mathrm{~with~the~base~} \ \rho^4 = \rho_{\tau_4}^{<\alpha^3}.
\end{equation*}

\noindent Now the argument from the proof of lemma
6.6~\cite{Kiselev9} should be repeated. From \ $
\check{\delta}^{4} < \check{\delta}^{3} $ \ it comes that
\begin{equation*}
\rho^4<\rho^3 \mathrm{~and~that~is~why~} \check{\delta}^{4} \notin
SIN_{n+1}^{< \alpha^{3}} [<\rho_{\tau_{3}}]
\end{equation*}
and by lemma 6.6~\cite{Kiselev9} (for \ $m=n+1$) \ there exists
some \ \mbox{$\Sigma_{n}$-proposition} \ $\varphi(\alpha,
\overrightarrow{a}) $ \ with the train \ $\overrightarrow{a} $ \
of constants \ $< \rho_{\tau_{3}} $ \ and some ordinal \
$\alpha_{0} \in [ \check{\delta}^4, \alpha^3 [ $ \ such that
\begin{equation*}
\forall \alpha < \alpha_{0} \ \ \varphi^{\triangleleft \alpha^{3}}(\alpha,
\overrightarrow{a})\wedge \neg  \varphi^{\triangleleft
\alpha^{3}}(\alpha_{0}, \overrightarrow{a})~.
\end{equation*}
The disseminator \ $\check{\delta}^{3} $ \ restricts the
proposition \ $ \exists \alpha \neg \varphi(\alpha,
\overrightarrow{a}) $ \ below \ $\alpha^{3} $, \ so \ $\alpha_{0}
\in \; ] \check{\delta}^{4}, \check{\delta} ^{3} [ $. \ The \
$\Pi_{n+1}$-proposition
\begin{multline*}
    \forall \alpha, \gamma \Bigl( \neg \varphi(\alpha,
    \overrightarrow{a}) \longrightarrow \exists \gamma_1 \bigl (
    \gamma < \gamma_1 \wedge SIN_{n-1}(\gamma_1) \wedge  \
\\
    \wedge \exists \delta < \alpha \ \ \exists \alpha^{\prime}, \
    \mathbf{K} (\delta, \gamma_{1}, \alpha^{\prime}, \rho^4, S^4) \bigr)
    \Bigr)
\end{multline*}

\noindent is fulfilled below \ $\check{\delta}^{3} $ \ and hence \
$\check{ \delta}^{3} $ \ extends it up to \ $\alpha^{3} $, \
because
\[
    S^4 \vartriangleleft \rho_4 < \rho_3.
\]
Hence, for every \ $\gamma_{\tau}^{< \alpha^3} > \check{\delta}^3$
\ there appears \ $\delta$-matrix \ $S^4$ \ admissible on some
carrier
\[
    \alpha \in [ \gamma_{\tau}^{< \alpha^3}, \alpha^3 [
    \mathrm{\ \ for \ \ } \gamma_{\tau}^{< \alpha^3}
\]
along with its disseminator \ $\check{\delta}^4 < \check{\delta}
^3 $ \ and the base \ $\rho^4$.
\\
From here it follows that below \ $\alpha^3$ \ there are definable
the \textit{minimal} cardinal \ $\check{\delta}^m$ \ and the
\textit{minimal} base \ $ \rho^m$ with this property, that is
fulfilling the following statement \textit{below} \ $\alpha^3$:
\[
    \exists \gamma^m \forall \gamma > \gamma^m \bigl(
    SIN_{n-1}(\gamma) \rightarrow \exists \alpha^{\prime},
    S ~\mathbf{K}(\check{\delta}^m,
    \gamma, \alpha^{\prime}, \rho^m, S) \bigr),
\]
\begin{sloppypar}
\noindent that is there exists \ $\gamma^m < \alpha^3$ \ such that
for every \ \mbox{$\gamma_{\tau}^{< \alpha^3} \in \; ] \gamma^m,
\alpha^3 [$} \ there exists some \ $\delta$-matrix \ $ S$ \
admissible on some carrier \ \mbox{$\alpha \in [ \gamma_{\tau}^{<
\alpha^3}, \alpha^3 [$} \ for \ $\gamma_{\tau}^{< \alpha^3}$ \
below \ $\alpha^3$ \ along with its generating disseminator \
$\check{\delta}^m < \gamma^m$ \ with the base \ $\rho^m$.
\end{sloppypar}

\noindent Obviously, \ $\check{\delta}^m < \check{\delta}^3$. \
Since the minimal value \ $ \rho^m$ \ is definable below \
$\alpha^3$, \ by lemma 4.6~\cite{Kiselev9} about spectrum type, it
follows
\begin{equation*}
\rho^m < OT(\delta S_{\tau_3}) \leq Od(\delta S_{\tau_3}).
\end{equation*}
But then it implies the contradiction: there exist \
$\delta$-matrix \ $ S^m$ \ on some carrier \ $\alpha^m \in \; ]
\gamma_{\tau_3}, \alpha^3 [ $ \ admissible for \ $\gamma_{\tau_3}$
\ along with the disseminator \ $\check{\delta}^m <
\gamma_{\tau_3}$ \ and the base \ $\rho^m$ \ and by condition \
$\mathbf{K}^{\exists}_{n+1}$

\begin{equation*}
S^m \vartriangleleft \rho^m < OT(\delta S_{\tau_3})  \leq
Od(\delta S_{\tau_3}),
\end{equation*}
\vspace{0pt}

\noindent though \ $\delta S_{\tau_3}$ \ is \
$\underline{\lessdot}$-minimal by definition~7.2.
\\
\quad \\
Statement 2)~$(i)$ repeats lemma 5.17~2)~$(i)$~\cite{Kiselev9} and
follows from definition~7.2 immediately; while statement 2)~$(ii)$
one can establish easily by means of the argument of lemma
5.17~2)~$(ii)$~\cite{Kiselev9} \ proof for the matrix \ $S$ \
instead of \ $S_{\chi \tau}^{< \alpha_{1}} $ \ and for the formula
\ $\mathbf{K}$ \ instead of \ $\sigma$; \ we shall return to this
argument in $\S$~8 in the more important case. \hspace*{\fill}
$\dashv$
\\
\quad \\

\noindent The unrelativized function \ $\delta S_{f} $ \ really
does exist on the final subinterval of the inaccessible cardinal \
$k$ \ as it shows

\begin{lemma}
\label{7.6.} \emph{(About \ $\delta $-function definiteness)} \newline
\hspace*{1em} There exists an ordinal \ $\delta <k$ \ such that \ $\delta
S_{f}$ \ is defined on the set
\begin{equation*}
T = \{\tau :\delta <\gamma _{\tau } < k \}.
\end{equation*}
The minimal of such ordinals \ $\delta $ \ is denoted by \ $\delta
^{\ast }$, \ its successor in \ $SIN_{n}$ \ by \ $\delta ^{\ast
1}$ \ and the following corresponding ordinals are introduced:
\[
    \tau_1^{\ast }=\tau (\delta ^{\ast }),\quad \tau ^{\ast 1}=\tau
    (\delta ^{\ast 1}),
\]
\[
    \quad so~that~~
    \delta^{\ast}=\gamma_{\tau_1^{\ast}}, \quad \delta^{\ast
    1}=\gamma_{\tau^{\ast 1}},
\]
\[
    \quad and~~  \alpha ^{\ast
    1}=\alpha _{\tau ^{\ast 1}}^{\Downarrow },\quad \rho ^{\ast
    1}=\rho _{\tau _{{}}^{\ast 1}}.
\]
\end{lemma}

\noindent \textit{Proof} \ consists in the immediate application
of lemma~6.14~\cite{Kiselev9} for $\alpha_{1} = k, \ m=n+1, \ \chi
= \chi^{\ast} $. \hspace*{\fill} $\dashv$

\begin{lemma}
\label{7.7.}
\begin{equation*}
\delta ^{\ast }\in SIN_{n}\cap SIN_{n+1}^{<\alpha ^{\ast 1}}\left[ <\rho
^{\ast 1}\right] ^{`}.
\end{equation*}
\end{lemma}

\noindent \textit{Proof.} \ Let us consider the disseminator \
$\check{\delta}_{\tau^{\ast 1 }} $ \ with the base \ $\rho^{\ast
1} $ \ of the matrix \ $\delta S_{\tau^{\ast 1 }} $ \ on the
carrier \ $\alpha_{\tau^{\ast 1 }} $. \ Since
\begin{equation*}
\delta^{\ast 1} \in SIN_{n}, \quad \check{\delta}_{\tau^{\ast 1 }}
< \delta^{\ast 1}
\end{equation*}
and
\begin{equation*}
\check{\delta}_{\tau^{\ast 1 }} \in SIN_{n}^{ <
\alpha^{\ast1}}\cap SIN_{n+1}^{ < \alpha^{\ast 1}}[< \rho^{\ast
1}],
\end{equation*}
\vspace{0pt}

\noindent lemma 3.8~\cite{Kiselev9} implies \
$\check{\delta}_{\tau^{\ast 1 }} \in SIN_{n}$ \ and by
lemmas~7.5~2), 7.6 \ $\check{\delta}_{\tau^{\ast 1 }} =
\delta^{\ast} $. \hspace*{\fill}$\dashv$
\\

\begin{definition}
\label{7.8.} \hfill {}

1. The function \ $\delta S_{\tau}^{<\alpha_1}$ \ is called
monotone on an interval \ $[\tau_1, \tau_2[$ \ or on corresponding
interval \ $[\gamma_{\tau_1}^{<\alpha_1},
\gamma_{\tau_2}^{<\alpha_1}[$ \ below \ $\alpha_1$, \ iff \
$\tau_1+1<\tau_2$, \ $]\tau_1, \tau_2[\; \subseteq dom(\alpha
S_f^{<\alpha_1})$ \ and
\[
    \forall \tau^{\prime}, \tau^{\prime\prime}(\tau_1 < \tau^{\prime}
    < \tau^{\prime\prime} < \tau_2 \longrightarrow \delta
    S_{\tau^{\prime}}^{<\alpha_1} \; \underline{\lessdot } \; \delta
    S_{\tau^{\prime\prime}}^{<\alpha_1}).
\]

2. Thereafter the function \ $\delta S_f$ \ is called (totally)
monotone iff for \ $\tau_1^{\ast} = \tau(\delta^{\ast})$:
\[
    \forall \tau^{\prime}, \tau^{\prime\prime}(\tau_1^{\ast} < \tau^{\prime}
    < \tau^{\prime\prime} < k \longrightarrow \delta
    S_{\tau^{\prime}} \; \underline{\lessdot } \; \delta S_{\tau^{\prime\prime}} ).
\]
\hspace*{\fill}$\dashv$
\end{definition}

Some easy fragments of the matrix function \ $\delta S_f$ \
monotonicity comes from definition~7.2 and lemma 7.5~2)~$(ii)$ at
once:

\begin{lemma}
\label{7.9.} \emph{(About \ $\delta $-function monotonicity)}
\newline \hspace*{1em} Let
\[
    \tau_1 < \tau_2 \mathrm{~and~}
    \check{\delta}_{\tau_2}^{<\alpha_1} <
    \gamma_{\tau_1}^{<\alpha_1}.
\]

Then
\[
    \delta S_{\tau_1}^{<\alpha_1} \; \underline{\lessdot } \; \delta
    S_{\tau_2}^{<\alpha_1}.
\]
\hspace*{\fill}$\dashv$
\end{lemma}

Let us discuss the situation which arises.
\\
We have revealed above, that the simplest matrix function \
$S_{f}$ \ is \ $ \underline{ \lessdot }$-monotone, but for every \
$\tau>\tau ^{\ast }$ \ the prejump cardinal \ $\alpha ^{\Downarrow
}$ \ of \ $S_{\tau} $ \ on its corresponding carrier \ $\alpha \in
\; ]\gamma_{\tau}, \gamma_{\tau+1} [ $ \ do not conserve the
subinaccessibility of levels \ $\geq n$ \ of cardinals \ $\leq
\gamma _{\tau }$, \ and some other important properties of the
lower levels of the universe are destroyed when relativizing to \
$ \alpha ^{\Downarrow }$ \ (see lemmas 5.17, 5.18 and their
discussion in the end of $\S$5~\cite{Kiselev9}).

In order to overcome this obstruction we have supplied the values
of this function, matrices \ $S_{\tau }$, \ by disseminators of
the level \ $n+1$ \ and required the conservation of the
subinaccessibility of the level \ $n$ \ for cardinals \ $\leq
\gamma _{\tau }$, \ that is we passed to the \ $\delta $-function
\ $\delta S_{f}$.
\\

\noindent But now it involves the new complication: now with the
help of lemmas \ref{7.3.}-\ref{7.7.} above one can see, that after
this modification the \ $\delta $-function is deprived of the
property of total \ $\underline{\lessdot }$-monotonicity on
$[\tau_1^{\ast}, k[\;$, \ and just due to the fact that in many
cases the prejump cardinals \ $\alpha^{\Downarrow } $ \ of \
$\delta$-matrices carriers \ $\alpha$, \ vice versa, \textit{give
rise to the subinaccessibility of the level} \ $n$ \ of some
cardinals \ $\leq \gamma _{\tau }$ \ that become subinaccessible
(relatively to \ $\alpha^{\Downarrow} $), \ not being those in the
universe (Kiselev~\cite{Kiselev4}).
\newline

The way out of this new situation is pointed out by the following
discovery that affords the solution of the problem:

One can see that the matrix \ $\delta S_{\tau _{0}}$ \ breaking
the \ $ \underline{\lessdot}$-monotonicity on \ $[\tau_1^{\ast},
k[$ \ at the first time, that is for
\begin{equation*}
\tau _{0}=\sup \{ \tau :\delta S_{f}\quad \mathrm{{is}\quad
\underline{\lessdot}\mathrm{-monotone\ on} \quad \left]  \tau
_{1}^{\ast},\tau \right[ \; \} ~, }
\end{equation*}
is placed on some carrier \ $ \alpha_{\tau_0} \in \; ] \delta
^{\ast },\delta ^{\ast 1} [$ \ and also \ $\delta S_{\tau
_{0}}\vartriangleleft \rho ^{\ast 1}$ \ by lemma
3.2~\cite{Kiselev9}.
\\
Therefore from lemmas 7.7,~6.3~\cite{Kiselev9} (for \ $m=n+1$, \
$\alpha _{1}=\alpha ^{\ast 1}$) it follows that the disseminator \
$\check{\delta}_{\tau ^{\ast 1}}$ \ carries over precisely the
same situation, but below \ $\alpha ^{0}=\alpha _{\tau
_{0}}^{\Downarrow }$, \ that is:
\\
\quad \\
\textit{the class \ $SIN_{n}^{<\alpha ^{0}}$ \ contains
some cardinals \ $ \gamma _{\tau _{1}}^{<\alpha ^{0}}<\gamma
_{\tau _{2}}^{<\alpha ^{0}}$ such that
\begin{equation*}
\left] \tau_{1}, \tau _{2}\right[ \subseteq  dom ( \delta S_{f}^{<\alpha
^{0}} )
\end{equation*}
and again just the same matrix
\begin{equation*}
\delta S_{\tau ^{0}}=\delta  S_{\tau _{0}^{\prime }}^{<\alpha ^{0}}
\end{equation*}
is breaking the monotonicity of \ $\delta S_{f}^{<\alpha ^{0}}$ \
on \ $\left] \tau _{1},\tau _{2}\right[ $ \ for the first time for
some ordinal \ $\tau_{0}^{\prime} \in ]\tau_{1},\tau_{2} [\;$,} \
{\sl but below \ $\alpha^0$.}
\\

So, here we come to the third and final approach to the main idea:
\newline \quad \newline \textit{The following requirements should
be imposed on \ $\delta$-matrices:
\newline
1) they must possess the property of ``autoexorcizivity'', that is
of self-exclusion in such situations of monotonicity violation;
the matrices with this property (of {\sl ``unit characteristic''})
will have the priority over other matrices (of {\sl ``zero
characteristic''} respectively) during defining of the matrix
function;
\newline 2) one more requirement should be imposed on the matrices of {\sl
zero} characteristic, hampering their forming: their disseminator
data bases must increase substantially, when the proceeding part
of matrix function, that is have already been defined, contains
monotonicity violation, in order to correct this fault -- the
using of matrices of {\sl zero} characteristic;
\\
on this grounds the \ $\delta$-matrix function should receive
inconsistent properties of monotonicity and nonmonotonicity
simultaneously. }
\\

Obviously, all these reasons require the recursive definition of
the matrix function, setting its values depending on the
properties of its preceding values.
\\
We start to realize this idea from the following section.

\newpage

\section{Matrix \ $\protect\alpha$\,-Functions}
\setcounter{equation}{0}

\hspace*{1em} For the forthcoming recursive definition it is
necessary to complicate the previous formula \
$\mathbf{K}_{n+1}^{\exists }$ \ (definition \ref{7.1.}). But
beforehand certain subformulas are to be introduced in view to the
more clearness of the construction of this formula, where the
variable \ $X_1$ \ plays the role of the matrix function \ $\alpha
S_{f}^{< \alpha}$ \ and the variable \ $X_2$ \ plays the role of
the characteristic function \ $a_f^{< \alpha}$ \ forthcoming to be
defined below \ $\alpha$ \ both; the latter function assigns
corresponding characteristics (unit or zero) to reduced matrices
serving as values of \ $\alpha S_f^{< \alpha}$; \ these
characteristics of matrices on their carriers will take values
unit \ $a=1$ \ or zero \ $a=0$ \ according to the principle
sketched above.
\\
During introducing, these formulas will be accompanied by comments
on their sense, and after resulting definition~8.2 we shall
describe in outline how
it works as a whole.
\newline

All these formulas were used in the author's previous works
\cite{Kiselev1,Kiselev2,Kiselev3,Kiselev4,Kiselev5,Kiselev6,Kiselev7,Kiselev8},
but some of them were scattered over the text in their certain
forms (sometimes nonformalized, some others in semantic manner),
and here they are gathered together; also some suitable
redesignations are used.

In these formulas various cardinals from the classes \ $SIN_{n}$,
$SIN_{n-1}$, $SIN_{n-2}$, \ of subinaccessibility are used. It is
necessary to take in view that after \ $<$- \ or \
$\vartriangleleft$-bounding of these formulas by some cardinal \
$\alpha$ \ (see definition~2.3~\cite{Kiselev9}) there arise the
subinaccessibility classes of the same level, but bounded by this
\ $\alpha$; \ for example the \ $SIN_{n}$-subinaccessibility \
turns into the \ $SIN_{n}^{<\alpha}$-subinaccessibility, \ but
below \ $\alpha$; \ thus all formulas after that narrate about
corresponding situation below \ $\alpha$.
\\
Such transformations lean on definitions and on lemmas
3.3-3.8~\cite{Kiselev9}.
\\

\begin{definition}
\label{8.1.} \hfill {} \newline \hspace*{1em} The following
auxiliary formulas are introduced:
\\ \quad
\\
\emph{I.} Intervals of matrix function definiteness
\\ \quad
\\
\emph{1.0}\quad $A_0(\chi ,\tau _{1},\tau _{2},X_1)$:
\vspace{-6pt}
\begin{multline*}
    \tau_{1}+1<\tau _{2} ~\wedge~ \big(X_1\mbox{ is a function on
    }\left] \tau _{1},\tau _{2}\right[ \; \big)\wedge
\\
    \wedge \tau _{1} =\min \big\{\tau :\left] \tau ,\tau _{2}\right[
    \subseteq dom(X_1)\big\}\wedge
\\
    \wedge \exists \gamma^{1} \big( \chi \leq \gamma ^{1}=\gamma
    _{\tau _{1}} \wedge SIN_{n}(\gamma ^{1}) \big);
\end{multline*}

\noindent so, this formula means, that the interval \ $[\tau
_{1},\tau _{2}[$ \ takes up the special place in relation to the
matrix function \ $X_1$ \ domain: this function is defined on \
$]\tau _{1},\tau _{2}[$ \ and \ $\tau_{1}$ \ is the minimal
ordinal with this property; besides that the cardinal \ $ \gamma
_{\tau _{1}}$ \ belongs to \ $SIN_{n}$;
\\
due to this minimality \ $X_1$ \ always is not defined for such
ordinal \ $\tau_1$.
\quad \medskip %

\noindent \emph{1.1}\quad $A_{1}(\chi ,\tau _{1},\tau _{2},X_1)$:
\[
    A_0(\chi ,\tau _{1},\tau _{2},X_1)\wedge \exists \gamma^2
    \big( \gamma^2 = \gamma_{\tau_2} \wedge SIN_n (\gamma^2)\big);
\]

\noindent such interval \ $[\tau_1, \tau_2[$ \ and the
corresponding interval \ $[\gamma_{\tau_1}, \gamma_{\tau_2}[$ \
will be called the intervals of the function \ $X_1$ \
definiteness {\sl maximal to the left} (in \ $dom(X_1)$), \
maximal in the sense that there is no interval \ $]\tau^{\prime},
\tau_2[$ \ in \ $dom(X_1)$ \ with the lesser left end \
$\tau^{\prime} < \tau_1$; \ in addition it is still demanded that
\ $\gamma_{\tau_1} \in SIN_n$.
\\

\noindent \emph{1.2}\quad $A_{1.1}^M(\chi ,\tau _{1},\tau
_{2},X_1)$: \vspace{-6pt}
\[
    \quad A_{1}(\chi, \tau_{1}, \tau_{2}, X_1)\wedge \tau _{2}=
    \sup \big\{ \tau : A_{1}(\chi, \tau_{1}, \tau_{2}, X_1) \big\};
\]

\begin{sloppypar}
\noindent here the interval \ $]\tau _{1},\tau _{2}[$ \ is the
maximal (included in \ $dom(X_1)$), \ maximal in the sense that it
is not included in any other interval \ \mbox{$]\tau _{1}^{\prime
},\tau _{2}^{\prime }[\;\subseteq dom(X_1)$} \ such that \ $
\gamma _{\tau_{2}^{\prime}}\in SIN_{n}$; \ beyond this condition
it is still demanded that \ $\gamma _{\tau _{1}} \in SIN_{n}$; \
thus such interval \ $[\tau_1, \tau_2[$, \ and the corresponding
interval \ $[\gamma_{\tau_1}, \gamma_{\tau_2}[$ \ will be called
the {\sl maximal} intervals of the function \ $X_1$ \
definiteness.
\\
\quad \medskip %
\end{sloppypar}

\noindent \emph{1.3}\quad $A_{1.2}(\tau _{1},\tau _{2},\eta )$:
\vspace{-6pt}
\begin{multline*}
    \exists \gamma ^{1},\gamma ^{2} \Big(\gamma ^{1}=
    \gamma _{\tau_{1}}\wedge \gamma ^{2}=\gamma _{\tau _{2}}
    \wedge
\\
    \wedge \eta =OT\big(\big\{\gamma :\gamma ^{1}<
    \gamma <\gamma^{2}\wedge SIN_{n}(\gamma )\big\}\big)\Big);
\end{multline*}

\noindent here, remind, \ $OT$ \ denote the order type of the
specified set, therefore we shall call such ordinal \ $\eta $ \
the type of the interval \ $[\tau _{1},\tau _{2}[$ \ and also of
the corresponding interval \ $[\gamma _{\tau _{1}},\gamma _{\tau
_{2}}[$~.
\newline
\quad \medskip %

\noindent \emph{1.4}\quad $A_{2}(\chi ,\tau _{1},\tau _{2},\tau
_{3},X_1)$: \vspace{-6pt}
\begin{multline*}
    \ A_{1}(\chi ,\tau _{1},\tau _{3},X_1)\wedge
    \tau _{1}+1<\tau_{2}<\tau _{3} \wedge \tau_{2} =
\\
    \qquad = \sup \Big\{\tau <\tau _{3}:\forall \tau ^{\prime },
    \tau ^{\prime\prime }\big(\tau _{1}<\tau ^{\prime }
    <\tau ^{\prime \prime }<\tau \longrightarrow
    X(\tau ^{\prime })\underline{\lessdot } X
    (\tau^{\prime \prime })\big)\Big\};
\end{multline*}

\noindent so, here \ $\tau _{2}$ \ is the minimal index at which
the \ $ \underline{\lessdot }$-monotonicity of the matrix function
\ $X_1$ \ on \ $\left] \tau _{1},\tau _{3}\right[ $ \ fails.
\newline
\quad \medskip %

\noindent \emph{1.5}\quad $A_{3}(\chi ,\tau _{1},\tau _{1}^{\prime
},\tau _{2},\tau _{3},X_1,X_2)$: \newline

\noindent $\ A_{2}(\chi ,\tau _{1},\tau _{2},\tau _{3},X_1 )\wedge
\tau _{1}<\tau _{1}^{\prime }<\tau _{2}\wedge \big(X_2\mbox{\rm\
is a function on } ]\tau _{1},\tau _{3}[\ \big) \wedge $ \newline

\noindent \hfill $\ \wedge \tau _{1}^{\prime }=\min \big\{\tau \in
\;]\tau _{1},\tau _{2}[\;:X_1(\tau )\gtrdot X_1(\tau _{2})\wedge
X_2(\tau )=1\big\}$;
\newline

\noindent thus, here is indicated that the \ $\underline{\lessdot
}$ -monotonicity of the matrix function \ $X_1$ \ on \ $\left]
\tau _{1},\tau _{3}\right[ $ \ is broken first at the index \
$\tau _{2}$ \ and just because of the matrix \ $X_1(\tau
_{1}^{\prime })\gtrdot X_1(\tau _{2})$ \ for \ $\tau _{1}^{\prime
}\in \left] \tau _{1},\tau _{2}\right[ $ \ of unit characteristic.
\newline
\quad \medskip %

\noindent \emph{1.6.a}\quad $A_{4}^b(\chi ,\tau _{1},\tau
_{1}^{\prime },\tau _{2},\tau _{3},\eta ,X_1,X_2)$:
\newline

\qquad \qquad $A_{3}(\chi ,\tau _{1},\tau _{1}^{\prime },\tau
_{2},\tau _{3},X_1,X_2) \wedge A_{1.2}(\tau _{1},\tau _{3},\eta
);$
\\
\quad \\
\noindent \emph{1.6.a(i)}\quad $A_{4}^b(\chi ,\tau _{1}, \tau
_{2}, \eta, X_1, X_2)$:
\[
    \qquad \exists \tau_1^{\prime}, \tau_2^{\prime} \le \tau_2 \
    A_{4}^b(\chi, \tau_{1}, \tau_{1}^{\prime}, \tau_{2}^{\prime},
    \tau_{2}, \eta, X_1, X_2);
\]

\noindent \emph{1.6.b}\quad $A_{4}^{M b}(\chi ,\tau _{1},\tau
_{1}^{\prime },\tau _{2},\tau _{3},\eta ,X_1,X_2)$:
\[
    A_{4}^b(\chi, \tau_{1},\tau _{1}^{\prime },\tau _{2},\tau
    _{3},\eta,X_1,X_2) \wedge
    A_{1.1}^M(\chi ,\tau _{1},\tau _{3},X_1);
\]

\noindent \emph{1.6.b(i)}\quad $A_{4}^{M b}(\chi ,\tau _{1}, \tau
_{2}, \eta, X_1, X_2)$:
\[
    \qquad \exists \tau_1^{\prime}, \tau_2^{\prime} \le \tau_2 \
    A_{4}^{M b}(\chi, \tau_{1}, \tau_{1}^{\prime}, \tau_{2}^{\prime},
    \tau_{2}, \eta, X_1, X_2);
\]
in what follows every interval \ $[\tau _{1},\tau _{3}[$ \
possessing this property \ $A_{4}^b$ \ for some \ $\tau
_{1}^{\prime }$, \ $\tau _{2}$, \ $\eta $ \ and the corresponding
interval \ $[\gamma _{\tau _{1}},\gamma _{\tau _{3}}[$ \ will be
called the blocks of the type \ $\eta$, \ and if there in addition
holds \ $A_{1.1}^M(\chi ,\tau _{1},\tau _{3},X_1)$ \ -- \ then the
{\sl maximal blocks}.
\\
Such blocks are considered further as objectionable because of
there fatal defect: the violating of monotonicity. Because of that
we shall impose on such blocks certain hard conditions in order to
avoid their formation in the course of matrix function recursive
defining (see the condition \ $\mathbf{K} ^{0}$ \ below).
\end{definition}

Let us hold up for a little while this definition \ref{8.1.} to
explain the sense and direction of its subsequent part.

All formulas and notions introduced above and also forthcoming
will be used in the resulting recursive definition~\ref{8.2.} in
bounded forms, that is their variables and constants will be \
$<$- \ or \ $\vartriangleleft$-bounded by some corresponding
cardinal \ $\alpha_1$. \ In such cases their present formulations
are used, but with the added remark ``below \ $\alpha_1$''; \
respectively their designations are supplied by the upper index \
$<\alpha_1$ \ or \ $\vartriangleleft \alpha_1$.
\\
Thereby
\[
    A_1^{\vartriangleleft\alpha_1}(\chi, \tau_1, \tau_2, X_1)
\]
is the formula:

\vspace{-6pt}
\begin{eqnarray*}
    \tau_{1}+1<\tau _{2} ~\wedge~ \big(X_1\mbox{ is a function on
    }\left] \tau _{1},\tau _{2}\right[ \; \big)\wedge \qquad\qquad
\\
    \wedge \tau _{1} =\min \big\{\tau :\left] \tau ,\tau _{2}\right[
    \subseteq dom(X_1)\big\}\wedge \qquad\qquad\qquad\qquad
\\
    \quad \wedge \exists \gamma^{1}, \gamma^{2}
    \big( \chi \leq \gamma ^{1} \wedge
    \gamma ^{1}=\gamma _{\tau_{1}}^{<\alpha_1} \wedge
    \gamma^{2}=\gamma _{\tau _{2}}^{<\alpha_1} \wedge
\\
    \wedge SIN_{n}^{<\alpha_1}(\gamma ^{1})
    \wedge SIN_{n}^{<\alpha_1}(\gamma ^{2})     \big),
\end{eqnarray*}
\vspace{0pt}

\noindent which means that \ $\left]\tau_1, \tau_2 \right[$ \ is
the interval from domain of the function \ $X_1$ \ with the
minimal left end \ $\tau_1$, \ and in addition the corresponding
cardinals \ $\gamma_{\tau_1}^{<\alpha_1}$,
$\gamma_{\tau_2}^{<\alpha_1}$ \ are \
$SIN_n^{<\alpha_1}$-cardinals -- all it below \ $\alpha_1$.

\noindent Respectively,
\[
    A_4^{b \vartriangleleft\alpha_1}(\chi, \tau_1, \tau_1^{\prime},
    \tau_2, \tau_3, \eta, X_1, X_2)
\]
is the formula:
\[
    A_{3}^{\vartriangleleft\alpha_1}(\chi ,\tau _{1},\tau_{1}^{\prime },
    \tau _{2},\tau _{3},X_1,X_2) \wedge
    A_{1.2}^{\vartriangleleft\alpha_1}(\tau_{1},\tau_{3},\eta )
\]
which means that \ $[ \tau_1, \tau_3 [$ \ and \ $[
\gamma_{\tau_1}^{<\alpha_1}, \gamma_{\tau_3}^{<\alpha_1} [$ \ are
the \textit{blocks} below \ $\alpha_1$ of type \ $\eta$, \ that is
the interval \ $]\tau_1, \tau_3[$ \ is maximal to the left
included in \ $dom(X_1)$ \ and the cardinals \
$\gamma_{\tau_1}^{<\alpha_1}$, \ $\gamma_{\tau_3}^{<\alpha_1}$ \
are in \ $SIN_n^{<\alpha_1}$ \ both, and \
$\underline{\lessdot}$-monotonicity of \ $X_1$ \ on \ $]\tau_1,
\tau_3[$ \ is broken first at the index \ $\tau_2 \in \; ]\tau_1,
\tau_3[$ \ and just because of the matrix \ $X_1(\tau_1^{\prime})
\gtrdot X_1(\tau_2)$ \ of \textit{unit} characteristic for some \
$\tau_1^{\prime} \in \; ]\tau_1, \tau_2[$ \ -- \ and all it below
\ $\alpha_1$.

\noindent Is is not hard to see that all these and forthcoming
bounded formulas under consideration belong to the class \
$\Delta_1^1$ \ for any \ $\alpha_1 > \chi$, \ $\alpha_1 < k$, \
$\alpha_1 \in SIN_{n-2}$.

To introduce the notions forthcoming clearly it is convenient to
clarify in outline the principle regulating the assignment of
characteristics to matrices on their carriers and intersection of
these characteristics with each other, because the characteristic
function will play the leading role in recursive
definition~\ref{8.2.} of matrix function.
\\
So, the matrix \ $S$ \ on its carrier \ $\alpha$ --- and this
carrier \ $\alpha$ \ itself --- will receive \textit{zero}
characteristic \ $a=0$, \ if it participate in violation of matrix
function monotonicity in the following sense:
\\
there exist some interval of the matrix function definiteness
\[
    [ \gamma_{\tau_1}^{<\alpha_{\chi}^{\Downarrow}},
    \gamma_{\tau_3}^{<\alpha_{\chi}^{\Downarrow}} [
\]
below the prejump cardinal \ $\alpha_{\chi}^{\Downarrow}$ \ after
\ $\chi$ \ of this carrier \ $\alpha$, \ where occurs the
\textit{same matrix} \ $S$ \ as the value of the matrix function \
$X_1$, \ but already below \ $\alpha_{\chi}^{\Downarrow}$:
\[
    X_1(\tau_2) = S,
\]
for the index \ $\tau_2 \in ]\tau_1, \tau_3[$ \ which is the
minimal one violating the monotonicity of \ $X_1$ \ below \
$\alpha_{\chi}^{\Downarrow}$, \ that is when there holds
\[
    A_2^{\vartriangleleft \alpha_{\chi}^{\Downarrow}} ( \chi, \tau_1, \tau_2, \tau_3,
    X_1).
\]
And here comes the last refinement of this notion: in addition
there must be no admissible matrices for \ $\gamma_{\tau_1}$ \
below \ $\alpha_{\chi}^{\Downarrow}$ \ and all values of the
matrix function \ $X_1$ \ on the interval \ $]\tau_1, \tau_2]$ \
must be of \textit{unit} characteristic:
\\
\quad \\
\hspace*{8em} $ \forall \tau ( \tau_1 < \tau \le \tau_2
\rightarrow X_2(\tau)=1)$. \label{c9}
\endnote{
\ p. \pageref{c9}. \ This last refinement is not necessary and the
Main theorem proof can be conducted without it, but still it
should be accepted in order to shorten the reasoning forthcoming.
\\
\quad \\
} %
\\
\quad \\
Otherwise \ $S$ \ on \ $\alpha$ \ and \ $\alpha$ \ itself will
receive \textit{unit} characteristic \ $a=1$.
\\
And while the matrix function will receive its recursive
definition~\ref{8.2.} forthcoming, matrices of unit characteristic
will receive the priority over matrices of zero characteristic
--- to avoid the violation of monotonicity of this function.
\\
It is natural to realize the notion of ``priority'' in the sense:
when some value \ $X_1(\tau)$ \ of matrix function \ $X_1$ \ is on
definition and there are matrices \ $S^0$, \ $S^1$ \ of
characteristic zero and unit respectively that can be nominated as
such value, then just matrix \ $S^1$ \ should be assigned as \
$X_1(\tau)$.

But there will be certain rather specific cases, when zero
characteristic will be rejected by certain other conditions, when
\textit{zero matrix} \ $S$ \ on its carrier \ $\alpha$ \ will be
forbidden for nomination for a value of matrix functions; in every
such case we shall say, that \ $S$ \ on \ $\alpha$ \ {\em is
suppressed}.
\\
We use here the term ``suppression'', not ``nonpriority'', because
such suppression will happen not every time, but in special cases
depending on disposition of the carrier \ $\alpha$.

So, now we come to the description of cases, when such suppression
takes place.
\\
To organize these cases in the proper way and to formulate the
suppression condition one should notice, that formulas above in
this current definition~\ref{8.1.} must be used in the following
special way:
\\
Till now in all these formulas 1.0--1.6 b $(i)$ above the symbols \
$X_1$, $X_2$ \ were treated as functions defined on ordinals.
\\
But for recursive definition~\ref{8.2.} of matrix functions it is
necessary to use functions defined on {\em pairs} of ordinals.
Therefore let us introduce for such function \ $X$ \ another
function
\[
    X[\alpha] = \{(\tau,\eta) : ((\alpha,\tau),\eta) \in X \},
\]
so that
\[
    X(\alpha,\tau) = X[\alpha](\tau)
\]
for every pair \ $(\alpha,\tau) \in dom(X)$.
\\
Thereafter these formulas in definition~\ref{8.2.} and formulas
forthcoming will be used often for \ $X_1$, $X_2$ \ as functions:
\[
    X_1[\alpha^0], ~
    X_2[\alpha^0],
\]
where \ $\alpha^0$ \ are some ordinals.
\\
Now let us return to definition \ref{8.1.} with the aim to form
the so called suppression condition; it arises in connection with
coverings of cardinals by blocks of special kind and for this aim
the following band of conditions is needed:
\\
\quad \\

\em

\noindent \emph{II.} Suppression conditions
\\
\quad \\
\emph{2.1a.}\quad $A_{5.1}^{sc}(\chi, \gamma^m, \gamma, X_1,
X_2)$:
\[
    \gamma^m < \gamma \wedge \forall \gamma^{\prime} \in
    [\gamma^m, \gamma[ \ \exists \tau_1, \tau_2, \eta \ \big (
    \gamma^m \le \gamma_{\tau_1} \le \gamma^{\prime} <
    \gamma_{\tau_2} \le \gamma \wedge
\]
\[
    \qquad \qquad \qquad \qquad \qquad \qquad \qquad
    \wedge A_4^{Mb}(\chi, \tau_1, \tau_2, \eta, X_1, X_2) \big)
    \wedge
\]
\[
    \wedge \forall \gamma^{m \prime} < \gamma^m \neg \forall \gamma^{\prime\prime} \in
    [\gamma^{m \prime}, \gamma[ \ \exists \tau_1^{\prime}, \tau_2^{\prime},
    \eta^{\prime} \big ( \gamma^{m \prime} \le \gamma_{\tau_1^{\prime}} \le \gamma^{\prime\prime}
    < \gamma_{\tau_2^{\prime}} \le \gamma \wedge
\]
\[
    \qquad \qquad \qquad \qquad \qquad \qquad \qquad
    \wedge A_4^{Mb}(\chi, \tau_1^{\prime}, \tau_2^{\prime}, \eta^{\prime}, X_1, X_2)
    \big);
\]
here is indicated, that the interval \ $[\gamma^m, \gamma[$ \ with
the right end \ $\gamma$ \ is the union of maximal blocks and that
its left end \ $\gamma^m$ \ is the minimal one with this property;
such collection of intervals will be called the {\sl covering} of
the cardinal \ $\gamma$; \ it is easy to see, that under this
condition \ $\gamma^m$, $\gamma$ \ are \ $SIN_n$-cardinals;
\\
if one withdraw here the right end \ $\gamma$ \ it cause the
following condition:
\\
\quad \\
\emph{2.1b.}\quad $A_{5.1}^{sc}(\chi, \gamma^m, X_1, X_2)$:
\[
    \forall \gamma^{\prime} \ge \gamma^m \ \exists \tau_1, \tau_2, \eta \ \big (
    \gamma^m \le \gamma_{\tau_1} \le \gamma^{\prime} <
    \gamma_{\tau_2} \wedge \qquad \qquad \qquad
\]
\[
    \qquad \qquad \qquad \qquad \qquad \qquad \qquad
    \wedge A_4^{Mb}(\chi, \tau_1, \tau_2, \eta, X_1, X_2) \big)
    \wedge
\]
\[
    \wedge \forall \gamma^{m \prime} < \gamma^m \neg \forall \gamma^{\prime\prime}
    \ge \gamma^{m \prime} \ \exists \tau_1^{\prime}, \tau_2^{\prime},
    \eta^{\prime} \ \big ( \gamma^{m \prime} \le \gamma_{\tau_1^{\prime}} \le \gamma^{\prime\prime}
    < \gamma_{\tau_2^{\prime}} \wedge
\]
\[
    \qquad \qquad \qquad \qquad \qquad \qquad \qquad
    \wedge A_4^{Mb}(\chi, \tau_1^{\prime}, \tau_2^{\prime}, \eta^{\prime}, X_1, X_2)
    \big).
\]
In view to compose the suppression condition in a proper way the
following special conditions should be superimposed on such
coverings for the ordinals \ $\gamma^m < \gamma^{\ast} < \gamma$,
$\eta^{\ast}$:
\\
\quad \\
\emph{2.2.}\quad $A_{5.2}^{sc}(\chi, \gamma^m, \gamma^{\ast},
\eta^{\ast}, X_1, X_2)$:
\[
    A_{5.1}^{sc}(\chi, \gamma^m, \gamma^{\ast}, X_1, X_2) \wedge
    \qquad \qquad \qquad \qquad \qquad \qquad \qquad \qquad \qquad
\]
\[
    \wedge \forall \tau_1, \tau_2, \eta \ \big ( \gamma_{\tau_1} <
    \gamma_{\tau_2} \le \gamma^{\ast} \wedge
    A_4^{Mb}(\chi, \tau_1, \tau_2, \eta, X_1, X_2) \rightarrow
    \eta < \eta^{\ast} \big) \wedge
\]
\[
    \wedge \forall \eta < \eta^{\ast} \ \exists \gamma^{\prime} <
    \gamma^{\ast} \ \forall \tau_1^{\prime}, \tau_2^{\prime}, \eta^{\prime} \ \big (
    \gamma^{\prime} < \gamma_{\tau_2^{\prime}} \le \gamma^{\ast} \wedge
\]
\[
    \qquad \qquad \qquad \qquad \qquad \qquad
    \wedge A_4^{Mb}(\chi, \tau_1^{\prime}, \tau_2^{\prime}, \eta^{\prime}, X_1, X_2) \rightarrow
    \eta < \eta^{\prime} \big);
\]
in this case, when \ $A_{5.2}^{sc}$ \ holds, we shall say, that
covering types of the cardinal \ $\gamma^{\ast}$ \ are
nondecreasing up to \ $\eta^{\ast}$ \ {\sl substantially}; thereby
the ordinal \ $\eta^{\ast}$ \ is limit;
\\
\quad \\
\emph{2.3.}\quad $A_{5.3}^{sc}(\chi, \gamma^{\ast}, \gamma^1,
\gamma, \eta^{\ast}, X_1, X_2)$:
\[
    \exists \tau_1, \tau \ \Big (
    \gamma_{\tau_1} = \gamma^1 \wedge
    \gamma_{\tau} = \gamma \wedge A_4^b(\chi, \tau_1,
    \tau, \eta^{\ast}, X_1, X_2) \wedge
\]
\[
    \forall \tau_1^{\prime}, \tau_2^{\prime}, \eta^{\prime} \ \big (
    \gamma^{\ast} < \gamma_{\tau_2^{\prime}} \le \gamma^1 \wedge
    A_4^{Mb}(\chi, \tau_1^{\prime}, \tau_2^{\prime}, \eta^{\prime}, X_1, X_2) \rightarrow
\]
\[
    \qquad \qquad \qquad \qquad \qquad \qquad \qquad \qquad \qquad \qquad
    \rightarrow \eta^{\prime} = \eta^{\ast} \big) \Big);
\]
now these three conditions should be assembled in the following
\\
\quad \\
\noindent \emph{2.4.} \ Suppressing covering condition
\\
\quad \\
\hspace*{1.5em} $A_{5.4}^{sc}(\chi, \gamma, \eta^{\ast}, X_1,
X_2)$:
\[
    \exists \gamma^m, \gamma^{\ast}, \gamma^1 \ \Big (
    \gamma^m < \gamma^{\ast} < \gamma^1 < \gamma
    \wedge \eta^{\ast} < \chi^+ \wedge \quad
\]
\[
    \quad
    A_{5.1}^{sc}(\chi, \gamma^m, \gamma^1, X_1, X_2) \wedge
    A_{5.2}^{sc}(\chi, \gamma^m, \gamma^{\ast}, \eta^{\ast}, X_1, X_2)
    \wedge
\]
\[
    \qquad\qquad\qquad\qquad\qquad\qquad
    \wedge A_{5.3}^{sc}(\chi, \gamma^{\ast}, \gamma^1, \gamma,
    \eta^{\ast}, X_1, X_2) \Big);
\]
let us call the covering of \ $\gamma$ \ possessing this property
the {\sl suppressing covering} for \ $\gamma$ \ of the type \
$\eta^{\ast}$;
\\
so, these three conditions \ $A_{5.1}^{sc}-A_{5.3}^{sc}$ \ mean
together, that the covering of the cardinal \ $\gamma$ \ splits
into three parts: its types are nondecreasing up to the ordinal \
$\eta^{\ast} < \chi^+$ \ substantially to the left of \
$\gamma^{\ast}$, \ and then stabilizes from \ $\gamma^{\ast}$ \ up
to \ $\gamma^1$, \ that is the interval \ $[\gamma^{\ast},
\gamma^1[$ \ is covered by the maximal blocks of the constant type
\ $\eta^{\ast}$, \ also there is the block \ $[\gamma^1, \gamma[$
\ of the same type \ $\eta^{\ast} < \chi^+$; it is clear that
these conditions define the ordinals \ $\gamma^m < \gamma^{\ast} <
\gamma^1 < \gamma$, $\eta^{\ast}$ \ uniquely through \ $\gamma$ \
(if they exist);
\\
\quad \\
\emph{2.5.}\quad $A_{5.5}^{sc}(\chi, \gamma, \eta^{\ast}, \alpha,
X_1, X_2)$:
\[
    \forall \gamma^{\prime} \ \big( \gamma \le \gamma^{\prime} <
    \alpha \rightarrow \exists \tau_1^{\prime},
    \tau_2^{\prime}, \eta^{\prime} \big(
    \gamma_{\tau_1^{\prime}}^{<\alpha} \le \gamma^{\prime} <
    \gamma_{\tau_2^{\prime}}^{<\alpha} \wedge \qquad
\]
\[
    \qquad\qquad\qquad\qquad
    \wedge A_4^{Mb \vartriangleleft \alpha} (\chi, \tau_1^{\prime},
    \tau_2^{\prime}, \eta^{\prime}, X_1, X_2) \wedge
    \eta^{\prime} \ge \eta^{\ast}  \big ) \big );
\]
here is indicated, that the interval \ $[\gamma, \alpha[$ \ is
covered by maximal blocks below \ $\alpha$ \ of types \
$\eta^{\prime} \ge \eta^{\ast}$;
\\

\noindent Now all these conditions should be composed in the following
integrated
\\
\quad \\
\noindent \emph{2.6.} \ Resulting suppressing condition
\\
\quad \\
\hspace*{1.5em} $A_{5}^{S,0}(\chi, a, \gamma, \alpha, \rho, S,
X_1^0, X_2^0, X_1, X_2)$: \vspace{-6pt}
\begin{multline*}
    \qquad a=0 \wedge SIN_n(\gamma) \wedge \rho < \chi^+ \wedge
    \sigma(\chi, \alpha, S) \wedge
\\
    \wedge \exists \eta^{\ast}, \tau < \gamma \Big(
    \gamma = \gamma_{\tau} \wedge
    A_{5.4}^{sc}(\chi, \gamma, \eta^{\ast}, X_1^0, X_2^0) \wedge
\\
    \wedge \forall \tau^{\prime} \big( \tau < \tau^{\prime}
    \wedge SIN_n(\gamma_{\tau^{\prime}}) \rightarrow
\\
    \rightarrow \exists \alpha^{\prime}, S^{\prime} \big[
    \gamma_{\tau^{\prime}} < \alpha^{\prime} <
    \gamma_{\tau^{\prime}+1} \wedge SIN_n^{<\alpha_{\chi}^{\prime \Downarrow}}
    (\gamma_{\tau^{\prime}}) \wedge \sigma(\chi,
    \alpha^{\prime}, S^{\prime}) \wedge
\\
    \wedge A_{5.5}^{sc}(\chi, \gamma, \eta^{\ast},
    \alpha_{\chi}^{\prime \Downarrow},
    X_1[\alpha_{\chi}^{\prime \Downarrow}],
    X_2[\alpha_{\chi}^{\prime \Downarrow}]) \big] \big) \Big);
\end{multline*}
this last condition will superimpose on the matrix \ $S$ \ on its
carrier \ $\alpha$ \ the rather hard requirements and if it can be
realized, then only in the very special cases:
\\
the reduced matrix \ $S$ \ must be of \textit{zero} characteristic
on its carrier \ $\alpha$, the cardinal \ $\gamma$ \ must be \
$SIN_n$-{\sl subinaccessible}, the base \ $\rho$ \ must be {\sl
strictly less} than \ $\chi^+$, \ the cardinal \ $\gamma$ \ must
be {\sl covered by the suppressing covering} of the type \
$\eta^{\ast}$; \ moreover, for all \ $\gamma^{\prime} > \gamma$,
\mbox{$\gamma^{\prime} \in SIN_n$} \ there exist carriers \
$\alpha^{\prime} > \gamma^{\prime}$ \ with prejump cardinals \
$\alpha_{\chi}^{\prime \Downarrow}$ \ preserving \
$SIN_n$-cardinals \ $\le \gamma^{\prime}$ \ and with the interval
\ $[\gamma, \alpha_{\chi}^{\prime \Downarrow}[$ \ covered by the
maximal blocks of types \ $\eta^{\prime} \ge \eta^{\ast}$ \ below \ $\alpha_{\chi}^{\prime \Downarrow}$.

Further such occurences of the matrix \ $S$ \ of zero
characteristic on its corresponding carriers \ $\alpha$ \ will
fail while definition of the matrix function will go on, and, so, we
shall say that here zero matrix \ $S$ \ on \ $\alpha$ \ is {\sl
suppressed} for \ $\gamma$.
\\
Respectively, zero matrix \ $S$ \ on \ $\alpha$ \ with
disseminator \ $\delta$ \ and base \ $\rho$ \ is nonsuppressed for
\ $\gamma$ \ if this condition fails; \ thus any matrix \ $S$ \ on
\ $\alpha$ \ is {\sl nonsuppressed} if it is {\sl unit} or has the
base \ $\rho \ge \chi^{+}$ \ on \ $\alpha$, \ or \ $\gamma$ \ is
not \ $SIN_n$-cardinal; so one should always have in view the
corresponding cardinal \ $\gamma$.

\em

Let us suspend this definition \ref{8.1.} one more time for a
little bit to describe the direction of its final part operation.
\\
To construct the forthcoming definition \ref{8.2.} of matrix function
in the required way it should be managed by \ $\Pi_{n-2}$-formula
\[
    U_{n-2}(\mathfrak{n}, x, \chi, a, \delta, \gamma, \alpha, \rho, S)
\]
which is universal for the formula class \ $\Pi_{n-2}$ \ with
denoted free variables
\[
    x, \chi, a, \delta, \gamma, \alpha, \rho, S,
\]
and variable G\"{o}del number \ $\mathfrak{n}$ \ of such formulas
in basic model \ $\mathfrak{M}$ \ (see Tarski~\cite{Tarski}, also
Addison~\cite{Addison}).
\\
When this number \ $\mathfrak{n}$ \ and the variable \ $x$ \ will
take certain special value \ $\mathfrak{n}^{\alpha}$ \
simultaneously:
\[
    \mathfrak{n} = x = \mathfrak{n}^{\alpha}
\]
then this formula along with the \ $\Sigma_n$-formula \
$\mathbf{K}_n^{\forall}(\gamma, \alpha_{\chi}^{\Downarrow})$ \
will state, that \ $S$ \ is the \ $\alpha$-matrix reduced to \
$\chi$ \ on the carrier \ $\alpha$ \ of characteristic \ $a$ \ with
its disseminator \ $\delta$ \ and the base \ $\rho$, \ {\em
admissible} for \ $\gamma$ \ and obeying certain complex recursive
conditions; remind, the formula \
$\mathbf{K}_n^{\forall}(\gamma,\alpha)$ \
(definition~\ref{7.1.}\;) means, that the ordinal \ $\alpha$ \
preserves all \ $SIN_n$-cardinals \ $\le \gamma$.
\\
Nevertheless, until the value \ $\mathfrak{n}^{\alpha}$ \ will be
assigned to the variables \ $\mathfrak{n}$, \ $x$, \ this formula will
work in this definition~\ref{8.2.} with \ $\mathfrak{n} = x$:
\[
    U_{n-2}(x, x, \chi, a, \delta, \gamma, \alpha, \rho, S).
\]
Also further the following function restrictions are used:
\[
    X|\tau^0 = \big\{ (\tau, \eta) \in X: \tau < \tau^0 \big \};
\]
\[
    X|^1 \alpha^0 = \big\{ ((\alpha, \tau), \eta) \in X: \alpha < \alpha^0 \big
    \}.
\]
Now we return to definition \ref{8.1.} for the lase time. The
suppression condition \ $A_5^{S,0}$ \ will operate in the
following conjunction with the formula \ $U_{n-2}$, \ bearing all the
definition of matrix functions forthcoming:
\\

\em

\noindent \emph{III.} Bearing  conditions
\[
    3.1     \qquad U_{n-2}^{\ast}(\mathfrak{n}, x, \chi, a, \delta,
    \gamma, \alpha, \rho, S, X_1^0|\tau^{\prime}, X_2^0|\tau^{\prime},
    X_1|^1\alpha^0, X_2|^1\alpha^0): \qquad
\]
\[
    U_{n-2}(\mathfrak{n}, x, \chi, a, \delta, \gamma,
    \alpha, \rho, S) \wedge
    \qquad\qquad\qquad\qquad
\]
\[
    \quad \wedge \neg A_{5}^{S,0}(\chi, a, \gamma, \alpha, \rho, S,
    X_1^0|\tau^{\prime}, X_2^0|\tau^{\prime},
    X_1|^1\alpha^0, X_2|^1\alpha^0);
\]
this condition along with \ $\mathbf{K}_n^{\forall}(\gamma,
\alpha_{\chi}^{\Downarrow})$ \ after their \
$\vartriangleleft$-bounding by the cardinal \ $\alpha^0$ \ and for
the constants
\[
    x = \mathfrak{n}^{\alpha}, ~ \chi, ~ \delta, ~ \gamma, ~ \alpha, ~ \rho, ~ \tau^{\prime} < \alpha^0, ~
    S \vartriangleleft \rho
\]
will describe the following situation below \ $\alpha^0$: \ $S$ \
is the matrix reduced to \ $\chi$ \ on its carrier \ $\alpha$ \ of the
characteristic \ $a$ \ admissible for \ $\gamma$ \ along with its
disseminator \ $\delta$ \ and the base \ $\rho$, \ which is
nonsuppressed for \ $\gamma$ \ below \ $\alpha^0$  --- and one
should point out that this situation for any pair \ $(\alpha^0,
\tau^{\prime})$ \ will be determined by the functions
\[
    X_1^0|\tau^{\prime} = X_1[\alpha^0]|\tau^{\prime}, ~
    X_2^0|\tau^{\prime} = X_2[\alpha^0]|\tau^{\prime}, ~
    \mbox{ and } ~
    X_1[\alpha_{\chi}^{\prime \Downarrow}], ~
    X_2[\alpha_{\chi}^{\prime \Downarrow}]
\]
for various \ $\alpha_{\chi}^{\prime \Downarrow} < \alpha^0$; \
therefore the recursion mode provided by this condition will work
correctly.
\\
\quad \\
$3.2    \qquad A^0(x, \chi, \tau)$:
\\
\[
    \exists \gamma \Big( \gamma = \gamma_{\tau} \wedge \neg
    \exists a, \delta, \alpha, \rho, S \big( \mathbf{K}_n^{\forall}(\gamma,
    \alpha_{\chi}^{\Downarrow}) \wedge \qquad
\]
\[
    \qquad\qquad\qquad\qquad U_{n-2}(x, x, a, \delta, \gamma, \alpha, \rho, S) \big) \Big);
\]
this condition for \ $x=\mathfrak{n}^{\alpha}$ \ will mean, that
there is no \ $\alpha$-matrix \ $S$ \ on some carrier \ $\alpha$,
\ admissible for \ $\gamma_{\tau}$.
\\
$3.3    \qquad    A_2^0(x, \chi, \tau_1, \tau_2, \tau_3, X_1)$:
\\
\[
    A^0(x, \chi, \tau_1) \wedge A_2(\chi, \tau_1, \tau_2, \tau_3,
    X_1).
\]

\quad \\
\noindent \emph{IV.} Closing condition
\quad \\

This condition will help to close the diagonal reasoning providing the
final contradiction:

\noindent \emph{4.1} \qquad $ \Big( a=0 \longrightarrow \forall
\tau_1^{\prime}, \tau_1^{\prime\prime},
    \tau_2^{\prime}, \tau_3^{\prime}, \eta^{\prime} <
    \alpha_{\chi}^{\Downarrow} \big[
    \gamma_{\tau_1^{\prime}}^{<\alpha_{\chi}^{\Downarrow}} \le
    \delta <
    \gamma_{\tau_3^{\prime}}^{<\alpha_{\chi}^{\Downarrow}} \wedge$
\[
    \wedge A_4^{M b \vartriangleleft \alpha_{\chi}^{\Downarrow} } \big(
    \chi, \tau_1^{\prime}, \tau_1^{\prime\prime}, \tau_2^{\prime},
    \tau_3^{\prime}, \eta^{\prime}, X_1[\alpha_{\chi}^{\Downarrow}],
    X_2[\alpha_{\chi}^{\Downarrow}] \big) \rightarrow
    \eta^{\prime} < \rho \vee \rho = \chi^{+} \big] \Big);
\]
this formula has the following content for any {\textit zero}
matrix \ $S$ \ on the carrier \ $\alpha$ \ with the disseminator \
$\delta$ \ and base \ $\rho$:
\\
if this disseminator falls in any maximal block \ $ \big[
\gamma_{\tau_1^{\prime}}^{<\alpha_{\chi}^{\Downarrow}},
\gamma_{\tau_3^{\prime}}^{<\alpha_{\chi}^{\Downarrow}} \big[$ \
below the prejump cardinal \ $\alpha_{\chi}^{\Downarrow}$, \ that
is if
\[
    \gamma_{\tau_1^{\prime}}^{<\alpha_{\chi}^{\Downarrow}} \le
    \delta <
    \gamma_{\tau_3^{\prime}}^{<\alpha_{\chi}^{\Downarrow}},
\]
then this base \ $\rho$ \ has to arise substantially and to exceed
the type \ $\eta^{\prime}$ \ of this very block, or even to take its value
the greatest possible:
\[
    \eta^{\prime} < \rho \vee \rho = \chi^{+},
\]
for lack of anything better;
\\
therefore in similar cases the interval \ $[\tau_1^{\prime},
\tau_3^{\prime}[$ \ and the corresponding interval
\[
    \big[ \gamma_{\tau_1^{\prime}}^{<\alpha_{\chi}^{\Downarrow}},
    \gamma_{\tau_3^{\prime}}^{<\alpha_{\chi}^{\Downarrow}} \big[
\]
will be {\sl severe} for this {\sl zero} matrix \ $S$ \ on \
$\alpha$ \ and will {\sl hamper} the using of \ $S$ \ on \
$\alpha$ \ (with this disseminator \ $\delta$).
\\
\quad \\

\noindent \emph{V.} Equinformative condition
\\
\quad \\
\emph{5.1}\quad $A_{6}^{e}(\chi, \alpha^0)$:
\[
    \chi < \alpha^0 \wedge A_n^{\vartriangleleft \alpha^0}(\chi) =
    \| u_n^{\vartriangleleft \alpha^0}(\underline{l})\| \wedge
    SIN_{n-2}(\alpha^0) \wedge
\]
\[
    \qquad\qquad\qquad\qquad
    \wedge \forall \gamma < \alpha^0 \exists
    \gamma_1 \in [\gamma,\alpha^0[ \quad
    SIN_n^{< \alpha^0}(\gamma_1);
\]
the cardinal \ $\alpha^0$ \ here with this property is called,
remind, {\sl equinformative} with \ $\chi$.
\\
\hspace*{\fill} $\dashv$
\\

\em

The latter notion was used above several times (see also
\cite{Kiselev2,Kiselev3,Kiselev4,Kiselev5,Kiselev6,Kiselev7,Kiselev8},
\cite{Kiselev9}) and here it is
emphasized because of it special importance: every \
$\Pi_n$-proposition \ $\varphi(l)$ \ holds or not in any generic
extension \ $\mathfrak{M}[l]$ \ below \ $\chi$ \ and also in this extension below \
$\alpha^0$ \ simultaneously (see comment after (\ref{e7.1})); the
best example of such \ $\alpha^0$ -- any prejump cardinal \
$\alpha_{\chi}^{\Downarrow}$ \ after \ $\chi$ \ of any matrix
carrier \ $\alpha > \chi$ \ (if this cardinal is limit for the class \
$SIN_n^{< \alpha_{\chi}^{\Downarrow}}$).
\\

Now everything is ready to assemble all the fragments introduced above
in the following \textit{integrated} definition~\ref{8.2.} where
the variable matrix \ $S$ \ on its carrier \ $\alpha$ \ is under
consideration.
\\
Requirements which are superimposed there on \ $S$ \ on \ $\alpha$
\ and on its disseminator \ $\delta$ \ with the data base \ $\rho$ \
depend on the functions \ $X_i,~i=\overline{1,5}$, \ that already
are recursively defined below the prejump cardinal \
$\alpha_{\chi}^{\Downarrow}$; \ they are defined on the certain
subset of the domain
\begin{equation} \label{e8.1}
    \qquad\qquad
    \mathcal{A}_{\chi}^{\alpha_{\chi}^{\Downarrow}} = \Big\{
    (\alpha^0, \tau): \exists \gamma < \alpha^0 \big( \chi <
    \gamma=\gamma_{\tau}^{< \alpha^0} \wedge \qquad\qquad\qquad\qquad
\end{equation}
\[
    \qquad\qquad\qquad\qquad\qquad
    \wedge \alpha^0 \le \alpha_{\chi}^{\Downarrow} \wedge
    A_6^e(\chi,\alpha^0) \big) \Big\}
\]
and therefore the functions
\[
    X_i^0 = X_i[\alpha^0], ~ i=\overline{1,5}
\]
are defined on the corresponding subset of
\[
    \big\{ \tau: \gamma_{\tau} \in SIN_{n-1}^{<\alpha^0} \big\}
\]
for every cardinal \ $\alpha^0 \le \alpha_{\chi}^{\Downarrow}$ \
equinformative with \ $\chi$. \ This set \
$\mathcal{A}_{\chi}^{\alpha_{\chi}^{\Downarrow}}$ \ is considered
to be canonically ordered (with \ $\alpha^0$ \ as the first
component in this order and \ $\tau$ \ as the second).
\\

So, the variable \ $X_{2}^0$ \ will play here the role of
characteristic function \ $a_{f}^{<\alpha^0}$ \ defined below the
cardinal \ $\alpha^0$; \ $X_{1}^0$ \ will play the role of matrix
function \ $\alpha S_{f}^{<\alpha^0} $; \ $X_{3}^0$ \ -- the role
of disseminator function \ $\widetilde{\delta}_{f}^{<\alpha^0}$; \
$X_{4}^0$ \ -- the role of its data base function \ $\rho
_{f}^{<\alpha^0}$; \ $X_{5}^0$ \ -- the role of the carrier
function \ $\alpha _{f}^{<\alpha^0}$; \ all of them will be
defined below \ $\alpha^0$.
\\

After all these functions will be defined for all such cardinals
\[
    \alpha^0 \le \alpha_{\chi}^{\Downarrow}
\]
then in conclusion the resulting requirement will be superimposed
on the matrix \ $S$ \ on its carrier \ $\alpha$ \ under
consideration along with its disseminator \ $\delta$ \ and the data
base \ $\rho$ \ depending on location of this \ $\delta$, \ more
precisely -- depending on the maximal block
\[
    \big[ \gamma_{\tau_1^{\prime}}^{<\alpha_{\chi}^{\Downarrow}},
    ~ \gamma_{\tau_3^{\prime}}^{<\alpha_{\chi}^{\Downarrow}} \big[
\]
containing this disseminator \ $\delta$, \ that is already been
defined below \ $\alpha_{\chi}^{\Downarrow}$.
\\
And here, remind, is the closing requirement superimposed on the
matrix \ $S$ \ on its carrier \ $\alpha$, \ mentioned above:
\\
if \ $S$ \ on \ $\alpha$ \ has {\textit zero} characteristic and
its admissible disseminator \ $\delta$ \ falls in the maximal
block of type \ $\eta^{\prime}$ \ below \ $\alpha_{\chi}^{\Downarrow}$, \
then $\eta^{\prime} < \rho \wedge \rho = \chi^{+}$; \ thus, in such case
data base \ $\rho$ \ has to increase considerably and we shall see
that it is possible, but every time leads to some contradiction.
\\

In addition here is needed the formulation \ $L j^{<\alpha}(\chi)$
\ of the saturated cardinal \ $\chi$ \ notion (see argument
before~(\ref{e7.1}) or definition~6.9~4)~\cite{Kiselev9}\;);
remind also that \ $\widehat{\rho}$ \ denote the closure of \
$\rho$ \ under the pair function.

So, the recursive definition, based on the set \
$\mathcal{A}_{\chi}^{\alpha_{\chi}^{\Downarrow}}$, \ starts:
\label{c10}
\endnote{
\ p. \pageref{c10}. \ This definition has been presented earlier
(Kiselev~\cite{Kiselev6,Kiselev7,Kiselev8}) by means of the
uniform text, but now here it is split into parts to clarify its
structure.
\\
\quad \\
} %

\begin{definition}
\label{8.2.} \ \newline

\emph{1)}\quad Let
\[
    U_{n-2}(\mathfrak{n},x,\chi,a,\delta,\gamma,\alpha,\rho ,S)
\]
be the \ $\Pi _{n-2}$-formula universal for the class \ $\Pi
_{n-2} $ \ where \ $\mathfrak{n}$ \ is variable G\"{o}del number
of \ $\Pi _{n-2}$-formulas with the free variables \ $x$, $\chi $,
$a$, $\delta $, $\gamma $, $ \alpha $, $\rho $, $S$, \ and let
\[
    U_{n-2}^{\ast}(x,\chi,a,\delta, \gamma,
    \alpha ,\rho ,S,
    X_1^0|\tau^{\prime}, X_2^0|\tau^{\prime},
    X_1|^1\alpha^0, X_2|^1\alpha^0).
\]
be the formula
\begin{multline*}
    U_{n-2}(x,x,\chi,a,\delta,\gamma,\alpha,\rho,S) \wedge
\\
    \wedge \neg A_{5}^{S,0}(\chi, a, \gamma, \alpha, \rho, S,
    X_1^0|\tau^{\prime}, X_2^0|\tau^{\prime},
    X_1|^1\alpha^0, X_2|^1\alpha^0).
\end{multline*}

\quad \\
\emph{2)}\quad Let
\[
    A_7^{RC}(x, \chi, X_1, X_2, X_3, X_4, X_5, \alpha_{\chi}^{\Downarrow})
\]
be the following \ $\Delta_1$-formula, providing the required {\sl
recursive condition}:

\begin{equation*}
\begin{array}{l}
\bigwedge _{1\leq i\leq 5} \Big( ( X_{i}\mbox{\ \ is a
function})\wedge X_{i}\vartriangleleft \alpha_{\chi}^{\Downarrow
+} \wedge
\quad  \\
\quad  \\
\hspace{0.1cm}\wedge dom(X_{i}) \subseteq \Big \{ (\alpha^0,
\tau): \exists \gamma < \alpha^0 \big( \chi < \gamma =
\gamma_{\tau}^{< \alpha^0} \wedge \qquad
\quad  \\
\quad  \\
\hspace{4.5cm} \wedge \alpha^0 \le \alpha_{\chi}^{\Downarrow}
\wedge A_6^{e}(\chi, \alpha^0) \big) \Big\} \Big) \wedge
\end{array}
\end{equation*}

\vspace{-6pt}

\begin{equation*}
\begin{array}{l}
\hspace{0.1cm} \wedge \forall \alpha^0 \Big( \big( \alpha^0 \le
\alpha_{\chi}^{\Downarrow} \wedge A_6^{e}(\chi,\alpha^0) \big) \longrightarrow
\quad  \\
\quad  \\
\hspace{0.5cm} \longrightarrow \exists X_1^0, X_2^0, X_3^0, X_4^0,
X_5^0, X_1^{1,0}, X_2^{1,0} ~ \Big[\bigwedge_{1 \le i \le 5} X_i^0
= X_i[\alpha^0] \wedge
\quad  \\
\quad  \\
\hspace{4.5cm} \wedge X_1^{1,0} = X_1 |^1 \alpha^0 \wedge X_2^{1,0}
= X_2 |^1 \alpha^0 \wedge
\quad  \\
\quad  \\
\hspace{0.1cm} \wedge \forall \tau ^{\prime },\gamma ^{\prime },
\gamma^{\prime\prime} < \alpha^0 \Big( \chi <\gamma ^{\prime
}\wedge \gamma ^{\prime }=\gamma _{\tau ^{\prime }}^{<\alpha^0}
\wedge \gamma^{\prime\prime} =
\gamma_{\tau^{\prime}+1}^{<\alpha^0} \longrightarrow
\end{array}
\end{equation*}

\vspace{-6pt}

\begin{equation*}
\begin{array}{l}
\forall a^{\prime} \Big( X_{2}^0(\tau ^{\prime })=a^{\prime
}\leftrightarrow
\quad  \\
\quad  \\
\hspace{0.1cm} \leftrightarrow a^{\prime }=\max_{\leq } \big\{
a^{\prime \prime }:\exists \delta ^{\prime \prime },\alpha
^{\prime \prime },\rho ^{\prime \prime }< \gamma^{\prime \prime }
\exists S^{\prime \prime } \vartriangleleft \chi^{+} \big(
\mathbf{ K}_{n}^{\forall <\alpha^0}( \gamma ^{\prime},\alpha
_{\chi }^{\prime \prime \Downarrow }) \wedge
\quad  \\
\quad  \\
\wedge U_{n-2}^{\ast \vartriangleleft \alpha^0}(x,\chi ,a^{\prime \prime
}, \delta ^{\prime \prime },\gamma^{\prime}, \alpha ^{\prime
\prime }, \rho^{\prime \prime}, S^{\prime\prime}, X_1^0|\tau^{\prime}, X_2^0|\tau^{\prime},
X_1^{1,0}, X_2^{1,0}) \big) \big\} \Big) \wedge  \\
\end{array}
\end{equation*}

\vspace{-6pt}

\begin{equation*}
\begin{array}{l}
\wedge \forall S^{\prime } \Big( X_{1}^0(\tau ^{\prime
})=S^{\prime }\longleftrightarrow \exists a^{\prime }~ \Big(
a^{\prime} = X_{2}^0(\tau ^{\prime })\wedge  \\
\quad  \\
\hspace{0.5cm}\wedge S^{\prime } =\min_{\underline{\lessdot }}
\big\{ S^{\prime \prime } \vartriangleleft \chi^{+} :\exists
\delta ^{\prime \prime },\alpha ^{\prime \prime },\rho ^{\prime
\prime }< \gamma^{\prime\prime} \big( \mathbf{K}_{n}^{\forall
<\alpha^0} (\gamma ^{\prime },\alpha _{\chi }^{\prime \prime
\Downarrow })\wedge
\quad  \\
\quad  \\
\wedge U_{n-2}^{\ast \vartriangleleft \alpha^0}(x,\chi ,a^{\prime
},\delta ^{\prime \prime}, \gamma^{\prime}, \alpha ^{\prime
\prime}, \rho ^{\prime \prime },S^{\prime \prime },
X_1^0|\tau^{\prime}, X_2^0|\tau^{\prime},
X_1^{1,0}, X_2^{1,0}) \big) \big\} \Big) \Big) \wedge  \\
\end{array}
\end{equation*}

\vspace{-6pt}

\begin{equation*}
\begin{array}{l}
\wedge \forall \delta ^{\prime } \Big( X_{3}^0(\tau ^{\prime
})=\delta ^{\prime }\longleftrightarrow \exists a^{\prime
},S^{\prime } \Big( a^{\prime }=X_1^0(\tau ^{\prime })\wedge
S^{\prime }=X_2^0(\tau ^{\prime })\wedge
\quad  \\
\quad  \\
\hspace{0.5cm}\wedge \delta ^{\prime }=\min_{\leq } \big\{
\delta^{\prime \prime} < \gamma^{\prime} : \exists \alpha ^{\prime
\prime },\rho ^{\prime \prime }< \gamma^{\prime\prime} \big(
\mathbf{K}_{n}^{\forall <\alpha^0}(\gamma ^{\prime },\alpha _{\chi
}^{\prime \prime \Downarrow })\wedge
\quad  \\
\quad  \\
\wedge U_{n-2}^{\ast \vartriangleleft \alpha^0}(x,\chi ,a^{\prime
},\delta ^{\prime \prime }, \gamma^{\prime}, \alpha ^{\prime
\prime },\rho ^{\prime \prime },S^{\prime }, X_1^0|\tau^{\prime},
X_2^0|\tau^{\prime}, X_1^{1,0}, X_2^{1,0}) \big) \big\} \Big) \Big) \wedge \\
\end{array}
\end{equation*}

\vspace{-6pt}

\begin{equation*}
\begin{array}{l}
\wedge \forall \rho ^{\prime }~ \Big( X_{4}^0(\tau ^{\prime
})=\rho ^{\prime }\longleftrightarrow \exists a^{\prime
},S^{\prime },\delta ^{\prime } \Big( a^{\prime }=X_1^0(\tau
^{\prime}) \wedge S^{\prime }=X_2^0(\tau ^{\prime})\wedge
\quad  \\
\quad  \\
\hspace{0.1cm}\wedge \delta ^{\prime }=X_{3}^0(\tau ^{\prime })
\wedge \rho ^{\prime }=\min_{\leq } \big \{ \rho ^{\prime \prime }
< \chi^{+} :\exists \alpha^{\prime\prime} < \gamma^{\prime\prime}
\big( \mathbf{K}_{n}^{\forall <\alpha^0}(\gamma ^{\prime },\alpha
_{\chi }^{\prime \prime \Downarrow })\wedge
\quad  \\
\quad  \\
\hspace{0.2cm}\wedge U_{n-2}^{\ast \vartriangleleft
\alpha^0}(x,\chi ,a^{\prime },\delta ^{\prime }, \gamma^{\prime},
\alpha ^{\prime\prime}, \rho ^{\prime \prime }, S^{\prime },
X_1^0|\tau^{\prime}, X_2^0|\tau^{\prime}, X_1^{1,0}, X_2^{1,0})
\big) \big\} \Big) \Big) \wedge
\end{array}
\end{equation*}

\vspace{-6pt}

\[
\wedge \forall \alpha ^{\prime } \Big( X_{5}^0(\tau ^{\prime
})=\alpha ^{\prime }\longleftrightarrow \exists a^{\prime
},S^{\prime },\delta ^{\prime },\rho ^{\prime } \Big(
a^{\prime}=X_1^0(\tau ^{\prime })\wedge S^{\prime }=X_2^0(\tau
^{\prime })\wedge \qquad \qquad \qquad
\]
\[
\hspace{0.1cm}\delta ^{\prime }=X_3^0(\tau ^{\prime })\wedge
\rho^{\prime }=X_4^0(\tau ^{\prime })\wedge \alpha ^{\prime
}=\min_{\leq} \big \{ \alpha^{\prime\prime}<
\gamma^{\prime\prime}: \mathbf{K}_{n}^{\forall <\alpha^0}(\gamma
^{\prime },\alpha _{\chi}^{\prime \prime \Downarrow })\wedge
\]
\[
\hspace{0.2cm}\wedge U_{n-2}^{\ast \vartriangleleft
\alpha^0}(x,\chi , a^{\prime }, \delta ^{\prime },
\gamma^{\prime}, \alpha^{\prime \prime}, \rho ^{\prime },S^{\prime
}, X_1^0|\tau^{\prime}, X_2^0|\tau^{\prime}, X_1^{1,0}, X_2^{1,0}
)  \big\} \Big) \Big) \Big) \Big]
\Big). \quad  \\
\]

\noindent \emph{3)}\quad We denote through

\begin{equation*}
\alpha \mathbf{K}_{n+1}^{\exists }(x,\chi ,a,\delta ,\gamma ,\alpha ,\rho ,S)
\end{equation*}
\vspace{0pt}

\noindent the \ $\Pi _{n-2}$-formula which is equivalent to the
following formula:

\vspace{6pt}
\begin{equation*}
\begin{array}{l}
(a=0\vee a=1)\wedge \sigma (\chi ,\alpha ,S)\wedge Lj^{<\alpha
}(\chi )\wedge \chi <\delta <\gamma <\alpha \wedge
\quad  \\
\quad  \\
\hspace{0.5cm}\wedge S\vartriangleleft \rho \leq \chi ^{+}\wedge
\rho = \widehat{\rho }\wedge SIN_{n}^{<\alpha _{\chi }^{\Downarrow
}}(\delta )\wedge SIN_{n+1}^{<\alpha _{\chi }^{\Downarrow }}\left[
<\rho \right] (\delta )\wedge
\end{array}
\end{equation*}

\vspace{-6pt}

\begin{equation*}
\begin{array}{l}
\hspace{1.0cm}\wedge \forall \gamma <\alpha _{\chi }^{\Downarrow
}~\exists \gamma ^{\prime }\in [ \gamma ,\alpha _{\chi
}^{\Downarrow }[ ~SIN_{n}^{<\alpha _{\chi }^{\Downarrow }}(\gamma
^{\prime })\wedge cf(\alpha _{\chi }^{\Downarrow })\geq \chi
^{+}\wedge
\quad  \\
\quad  \\
\wedge \exists X_{1},X_{2},X_{3},X_{4},X_{5}\Bigl\{ A_7^{RC}(x,
\chi, X_1, X_2, X_3, X_4, X_5, \alpha_{\chi}^{\Downarrow}) \wedge
\quad  \\
\quad  \\
\Big( a=0\longleftrightarrow \exists \tau_1^{\prime},
\tau_2^{\prime}, \tau_3^{\prime}<\alpha_{\chi}^{\Downarrow} ~
\big( A_{2}^{0\vartriangleleft \alpha _{\chi }^{\Downarrow }}
(\chi, \tau_1^{\prime}, \tau_2^{\prime}, \tau_3^{\prime},
X_1[\alpha_{\chi}^{\Downarrow}] ) \wedge
\quad  \\
\quad  \\
\hspace{0.2cm} \wedge \forall \tau^{\prime\prime} \big(
\tau_1^{\prime} < \tau^{\prime\prime} \le \tau_2^{\prime}
\rightarrow X_2[\alpha_{\chi}^{\Downarrow}](\tau^{\prime\prime})=1
\big) \wedge X_1[\alpha_{\chi}^{\Downarrow}] (\tau_{2}^{\prime})
=S \big) \Big)\wedge
\end{array}
\end{equation*}

\vspace{-6pt}

\begin{equation*}
\begin{array}{l}
\wedge \Big( a=0\longrightarrow \forall \tau_1^{\prime},
\tau_1^{\prime\prime}, \tau_2^{\prime}, \tau_3^{\prime},
\eta^{\prime} <\alpha _{\chi }^{\Downarrow }~\big[
\gamma_{\tau_{1}^{\prime}}^{<\alpha _{\chi }^{\Downarrow }} \leq
\delta < \gamma_{\tau_{3}^{\prime}}^{<\alpha _{\chi
}^{\Downarrow }}\wedge  \\
\quad  \\
\wedge A_{4}^{M b \vartriangleleft \alpha _{\chi }^{\Downarrow }}
(\chi ,\tau_1^{\prime}, \tau_1^{\prime\prime}, \tau_2^{\prime},
\tau_3^{\prime}, \eta^{\prime}, X_1[\alpha_{\chi}^{\Downarrow}],
X_2[\alpha_{\chi}^{\Downarrow}]) \rightarrow \\
\quad  \\
\qquad\qquad\qquad\qquad\qquad\qquad\qquad\qquad \rightarrow
\eta^{\prime} < \rho \vee \rho =\chi ^{+}\big] \Big) \Bigr\}.
\end{array}
\end{equation*}
\vspace{6pt}

Let us denote by \ $\mathbf{K}^{0}(\chi ,a,\delta ,\alpha ,\rho)$
\ the last conjunctive constituent in big curly brackets $\{\;,\}$
in the latter formula, that is the closing condition:
\[
    \Bigl(a=0\longrightarrow \forall \tau_1^{\prime},
    \tau_1^{\prime\prime}, \tau_2^{\prime}, \tau_3^{\prime}, \eta^{\prime}
    <\alpha _{\chi }^{\Downarrow } ~ \bigl[
    \gamma_{\tau_{1}^{\prime}}^{<\alpha _{\chi }^{\Downarrow }} \leq
    \delta < \gamma_{\tau_{3}^{\prime}}^{<\alpha _{\chi
    }^{\Downarrow }} \wedge  \qquad \qquad \\
    \quad  \\
\]
\[
    \wedge A_{4}^{M b\vartriangleleft \alpha _{\chi }^{\Downarrow
    }}(\chi ,\tau_1^{\prime}, \tau_1^{\prime\prime}, \tau_2^{\prime},
    \tau_3^{\prime}, \eta^{\prime}, X_1[\alpha_{\chi}^{\Downarrow}],
    X_2[\alpha_{\chi}^{\Downarrow}]) \rightarrow \qquad \qquad
\]
\hspace*{17em} $\rightarrow \eta^{\prime} < \rho \vee \rho =
\chi^{+}\bigr]\Bigr).$\label{c11}
\endnote{
\ p. \pageref{c11}. \ This closure condition \ $\mathbf{K}^0$ \
operates as the condition \ $\mathbf{K}^0$ \ used earlier
(Kiselev~\cite{Kiselev6,Kiselev7,Kiselev8}), but in the more
managing way, because now it manages quite well without the
subformula \
$(\gamma_{\tau_3^{\prime}}^{<\alpha_{\chi}^{\Downarrow}} = \gamma
\rightarrow \lim(\gamma))$, \ which caused the significant proof
complication.
\\
\quad \\
} %
\quad \\

\noindent The functions \ $X_1[\alpha_{\chi}^{\Downarrow}],
X_2[\alpha_{\chi}^{\Downarrow}]$ \ are not mentioned here in the
denotation of \ $\mathbf{K}^{0}$ \ for brevity taking into
account, that they are defined uniquely in the preceding part of
this formula \ $\alpha \mathbf{K}_{n+1}^{\exists}$.

\emph{4)}\quad The formula \ $\alpha \mathbf{K}_{n+1}^{\exists }$
\ is the \ $\Pi _{n-2}$-formula and thereby it receives its
G\"{o}del number \ $\mathfrak{n} ^{\alpha }$, \ that is:
\begin{equation*}
\alpha \mathbf{K}_{n+1}^{\exists }(x,\chi ,a,\delta ,\gamma ,\alpha ,\rho
,S)\longleftrightarrow U_{n-2}(\mathfrak{n}^{\alpha },x,\chi ,a,\delta
,\gamma ,\alpha ,\rho ,S).
\end{equation*}
Let us assign the value \ $\mathfrak{n}^{\alpha }$ \ to the
variable \ $x$ \ in this equivalence and everywhere further,
thereafter these \ $\mathfrak{n}^{\alpha}$, $x$ \ will be dropped
in the notations.
\\
We denote through \ $\alpha \mathbf{K}^{<\alpha _{1}}(\chi
,a,\delta ,\gamma ,\alpha ,\rho ,S)$ \ the \ $\Delta_1$-formula:
\[
    \mathbf{K}_{n}^{\forall <\alpha _{1}}(\gamma ,\alpha _{\chi }^{\Downarrow
    })\wedge \alpha \mathbf{K}_{n+1}^{\exists \vartriangleleft \alpha _{1}}
    (\chi ,a,\delta ,\gamma ,\alpha ,\rho
    ,S)\wedge \alpha <\alpha _{1},
\]
and, respectively, through \ $\alpha \mathbf{K}^{\ast <\alpha
_{1}}(\chi,a,\delta,\gamma, \alpha, \rho, S)$ -- the formula,
which is obtained from the formula \ $\alpha
\mathbf{K}^{<\alpha_1}$ \ through joining to it the conjunctive
condition of the matrix \ $S$ \ {\sl nonsuppression} on \ $\alpha$
\ for \ $\gamma$ \ (see definition \ref{8.1.}~2.6\;), but below \
$\alpha_1 < k$ \ (as it was done above in points 2), 3) for \
$\alpha_{\chi}^{\Downarrow}$, \ but now for \ $\alpha_1 < k$ \
instead of \ $\alpha_{\chi}^{\Downarrow}$) -- in the following
way:
\[
    \mathbf{K}_{n}^{\forall <\alpha _{1}}(\gamma ,\alpha _{\chi }^{\Downarrow
    })\wedge \alpha \mathbf{K}_{n+1}^{\exists \vartriangleleft  \alpha _{1}}
    (\chi,a,\delta,\gamma, \alpha, \rho, S) \wedge \alpha <\alpha_{1}
    \wedge
\]
\[
    \wedge \neg \Big( a=0 \wedge SIN_n^{<\alpha_1}(\gamma) \wedge
    \rho < \chi^+ \wedge \sigma(\chi, \alpha, S) \wedge \qquad\qquad\qquad
\]
\[
    \wedge \exists X_1, X_2, X_3, X_4, X_5 \Big(
    A_7^{RC}(\mathfrak{n}^{\alpha}, \chi, X_1, X_2, X_3, X_4, X_5,
    \alpha_1) \wedge
\]
\[
    \wedge \exists \eta^{\ast}, \tau < \gamma \Big( \gamma =
    \gamma_{\tau}^{<\alpha_1} \wedge A_{5.4}^{sc \vartriangleleft
    \alpha_1} (\chi, \gamma, \eta^{\ast}, X_1[\alpha_1]|\tau,
    X_2[\alpha_1]|\tau) \wedge
\]
\[
    \wedge \forall \tau^{\prime} \big( \tau < \tau^{\prime} \wedge
    SIN_n^{<\alpha_1}(\gamma_{\tau^{\prime}}^{<\alpha_1})
    \rightarrow
    \qquad\qquad\qquad \qquad\qquad\qquad
\]
\[
    \exists \alpha^{\prime}, S^{\prime} \big[
    \gamma_{\tau^{\prime}}^{<\alpha_1} < \alpha^{\prime} <
    \gamma_{\tau^{\prime}+1}^{<\alpha_1} \wedge
    SIN_n^{<\alpha_{\chi}^{\prime\Downarrow}}(\gamma_{\tau^{\prime}}^{<\alpha_1})
    \wedge \sigma(\chi, \alpha^{\prime}, S^{\prime}) \wedge
\]
\[
    \qquad \qquad \qquad
    \wedge A_{5.5}^{sc \vartriangleleft \alpha_1}(\chi, \gamma, \eta^{\ast},
    \alpha_{\chi}^{\prime \Downarrow},
    X_1[\alpha_{\chi}^{\prime \Downarrow}],
    X_2[\alpha_{\chi}^{\prime \Downarrow}]) \big] \big) \Big)
    \Big) \Big);
\]
here is stated the admissibility of \ $S$ \ on \ $\alpha$, \ and
in addition --- its nonsuppressibility for \ $\gamma$ \ below \
$\alpha_1$. \ So, if \ $\alpha \mathbf{K}^{<\alpha_1}$ \ holds,
but \ $\alpha \mathbf{K}^{\ast <\alpha_1}$ \ fails, then \ $S$ \
on \ $\alpha$ \ is admissible but suppressed (all it below \
$\alpha_1$).\label{c12}
\endnote{
\ p. \pageref{c12}. \ This notion could be introduced in the
nonrelativized form for \ $\alpha_1 = k$ \ as well, but it is not
used in what follows; besides that, in this form it requires the
more complicated non-elementary language over \ $L_k$.
\\
\quad \\
} %

\emph{5)}\quad If the formula \ $\alpha \mathbf{K}^{<\alpha
_{1}}(\chi ,a,\delta ,\gamma ,\alpha ,\rho ,S)$ \ is fulfilled by
the constants \ $\chi $, $a$, $\delta $, $\alpha $, $\gamma $,
$\rho $, $S$, $\alpha _{1}$, \ then we say that \ $\chi $, $a$,
$\delta $, $\alpha $, $\rho $, $S$ \ are admissible {\sl very
strongly} for \ $\gamma $ \ below \ $\alpha _{1}$.
\\
If some of them are fixed or meant by the context, then we say
that others are admissible {\sl very strongly} for them (and for \
$\gamma $) below \ $\alpha _{1}$. \ So, by
\[
    \alpha \mathbf{K}^{<\alpha_1}(\chi, \gamma, \alpha, S)
\]
will be denoted the formula
\[
    \exists a, \delta, \rho < \gamma ~ \alpha
    \mathbf{K}^{<\alpha_1}(\chi, a, \delta, \gamma, \alpha, \rho, S)
\]
meaning that \ $S$ \ on \ $\alpha$ \ is admissible {\sl very
strongly} for \ $\gamma$ \ below \ $\alpha_1$.

\emph{6)} \ The matrix \ $S$ \ is called autoexorcizive or,
briefly, \ \mbox{$\alpha $-matrix} admissible {\sl very strongly}
on the carrier \ $\alpha $ \ of the characteristic \ $a$ \ for \
$\gamma = \gamma_{\tau}^{<\alpha_1}$ \ below \ $\alpha _{1}$, \
iff it possesses on \ $\alpha$ \ some disseminator \
$\delta<\gamma$ \ with a base \ $\rho$ \ admissible very strongly
for them also below \ $ \alpha _{1}$.

In every case of this kind \ $\alpha $-matrix is denoted by the
general symbol \ $\alpha S$ \ or \ $S$.

If \ $a_1=k$, \ or \ $\alpha_1$ \ is pointed out by the context,
then the upper indices \ $< \alpha_1, \vartriangleleft \alpha_1$ \ and
other mentionings about \ $\alpha_1$ \ are dropped.

Further all notions of admissibility will be always considered to
be very strongly, so the term ``very strongly'' will be omitted in
what follows. \hspace*{\fill} $\dashv$
\end{definition}

\noindent Thus here the bounded formula
\[
    \alpha \mathbf{K}^{\ast < \alpha _{1}}
    (\chi,a,\delta,\gamma, \alpha, \rho, S)
\]
arises from \ $\alpha \mathbf{K}^{<\alpha_1}$ \ by adding the
condition of nonsuppression of the matrix \ $S$ \ on \ $\alpha$ \
for \ $\gamma$ \ below \ $\alpha_1$ \ which is obtained from the
condition \ $\neg A_5^{S,0}$ \ in the way indicated above through
its \ $\vartriangleleft$-bounding by the cardinal \ $\alpha_1$ \
(that is through \ $\vartriangleleft$-bounding its individual
variables by \ $\alpha_1$ \ and through replacing its constituents
\ $SIN_n(\gamma)$ \ with \ $SIN_n^{<\alpha_1}(\gamma)$).
\\
So, any matrix \ $S$ \ on \ $\alpha$ \ if \textit{suppressed} for
the cardinal \ $\gamma$ \ below \ $\alpha_1$, \ then when \
$\gamma$ \ is \ $SIN_n^{< \alpha_1}$-cardinal, and when this \ $S$
\ has \textit{zero} characteristic on \ $\alpha$ \ and the base \
$\rho < \chi^+$ \ below \ $\alpha_1$.

Everywhere further \ $\chi =\chi ^{\ast }<\alpha _{1}$; \ we shall
often omit the notations of the functions \ $X_1=\alpha S_{f}^{<\alpha
_{1}} $, \ $X_2=a_{f}^{<\alpha _{1}}$ \ and the symbols \ $\chi
^{\ast }$, \ $\mathfrak{n} ^{\alpha }$ \ in writings of all
formulas
\begin{equation*}
A_{0}-A_5^{S,0}, \quad A^0, \quad A_2^0, \quad \alpha
\mathbf{K}_{n+1}^{\exists }, \quad \mathbf{K }^0, \quad \alpha
\mathbf{K}^{<\alpha _{1}}, \quad \alpha \mathbf{K}^{\ast <\alpha
_{1}}
\end{equation*}
from definitions \ref{8.1.}, \ref{8.2.} and of other denotations
for some shortening (if it will not perform a misunderstanding);
for example any prejump cardinals \
$\alpha_{\chi^{\ast}}^{\Downarrow}$ \ will be denoted through \
$\alpha^{\Downarrow}$ \ and so on.

\noindent Concerning these formulas it should be pointed out, that
definition \ref{8.2.} has been constructed with the aim to receive
the key formula \ $\alpha \mathbf{K}_{n+1}^{\exists}$ \ of the
class \ $\Pi_{n-2}$. \ To this end all constituent formulas were \
$\vartriangleleft$-bounded by ordinals $\alpha^0$ \ or $\alpha^{\Downarrow}$.
\\
But in what follows these boundaries will be often dropped without
loss of their content, because their individual variables and
constants are in fact will be bounded by ordinals pointed out in
the context during their applications.
\\

Clearly, variables \ $X_i,~ i=\overline{1,5}$ \ are defined in
definition \ref{8.2.} uniquely through all their parameters, thus
similar functions can be defined recursively following this
construction by the same recursion on the similar set (remind the set
(\ref{e8.1}))
\[
    \mathcal{A} = \bigl\{ (\alpha_1,\tau): \exists \gamma<\alpha^1
    (\chi^{\ast}<\gamma=\gamma_{\tau}^{<\alpha_1} \wedge
    A_6^e(\chi,\alpha_1) ) \bigr\}
\]
of pairs \ $(\alpha_1,\tau)$ \ (ordered canonically as in the
former case, with \ $\alpha_1$ \ as the first component in this
order and \ $\tau$ \ as the second).

\begin{definition}
\label{8.3.} \ \\
Let \ $\chi^{\ast}<\alpha_1$.
\\
\emph{1)}\quad We call as the characteristic function of the level
\ $n$ \ below \ $\alpha _{1}$ \ reduced to \ $\chi ^{\ast }$ \ the
function
\[
    a_{f}^{<\alpha _{1}}=(a_{\tau }^{<\alpha _{1}})_{\tau }
\]
taking the values: \vspace{-6pt}
\begin{multline*}
    a_{\tau }^{<\alpha _{1}} =
\\
    = \max_{\leq } \bigl \{ a: \exists \delta ,\alpha,
    \rho < \gamma_{\tau+1}^{<\alpha_1} ~ \exists S \vartriangleleft \chi^{\ast +}
    ~ \alpha \mathbf{K}^{\ast <\alpha _{1}}(a,\delta,
    \gamma_{\tau }^{<\alpha _{1}}, \alpha ,\rho ,S)
    \bigr\};
\end{multline*}

\noindent \emph{2)}\quad we call as the matrix autoexorcizive (in
monotonicity violation) function or, briefly, \ $\alpha $-function
of the level \ $n$ \ below $\alpha _{1}$ \ reduced to \ $\chi
^{\ast }$ \ the function
\[
    \alpha S_{f}^{<\alpha _{1}}=(\alpha S_{\tau }^{<\alpha _{1}})_{\tau }
\]
taking the values \vspace{-6pt}
\begin{multline*}
    \alpha S_{\tau }^{<\alpha _{1}}=
\\
    = \min_{\underline{\lessdot }}
    \bigl\{ S \vartriangleleft \chi^{\ast +}:
    \exists \delta ,\alpha ,\rho < \gamma_{\tau+1}^{<\alpha_1} ~  \alpha
    \mathbf{K}^{\ast <\alpha_{1}} (a_{\tau}^{<\alpha _{1}},\delta ,
    \gamma _{\tau }^{<\alpha _{1}}, \alpha, \rho,
    S) \bigr\};
\end{multline*}

\noindent \emph{3)} the following accompanying ordinal functions
are defined below~$\alpha _{1}$

\hspace*{1.5em} the floating disseminator function \hspace{\stretch{0.5000}}
$\widetilde{\delta }_{f}^{<\alpha _{1}}=(\widetilde{\delta }_{\tau
}^{<\alpha _{1}})_{\tau }$, \hspace*{1em}

\hspace*{1.5em} its data base function \hspace{\stretch{0.5000}} $\rho
_{f}^{<\alpha _{1}}=(\rho _{\tau }^{<\alpha _{1}})_{\tau }$, \hspace*{1em}

\hspace*{1.5em} the carrier function \hspace{\stretch{0.5000}} $\alpha
_{f}^{<\alpha _{1}}=(\alpha _{\tau }^{<\alpha _{1}})_{\tau }$, \hspace*{1em}

\hspace*{1.5em} the generating disseminator function
\hspace{\stretch{0.5} } $\check{\delta}_{f}^{<\alpha
_{1}}=(\check{\delta}_{\tau }^{<\alpha _{1}})_{\tau }$,
\hspace*{1em} \newline \quad \newline taking the values for \
$a_{\tau} = a_{\tau }^{<\alpha _{1}}$, $S_{\tau} = \alpha S_{\tau
}^{<\alpha _{1}}$:
\begin{multline*}
    \widetilde{\delta }_{\tau }^{<\alpha _{1}} = \vspace{-12pt}
\\
    = \min_{\leq }
    \bigl \{ \delta < \gamma_{\tau}^{< \alpha_1}:
    \exists \alpha ,\rho <
    \gamma_{\tau+1}^{< \alpha_1} ~ \alpha
    \mathbf{K}^{\ast <\alpha_{1}}(a_{\tau},
    \delta, \gamma_{\tau }^{<\alpha_{1}},
    \alpha, \rho, S_{\tau}) \bigr \};
\end{multline*}

\vspace{-24pt}
\begin{multline*}
    \rho _{\tau }^{<\alpha _{1}} = \vspace{-12pt}
\\
    = \min_{\leq }
    \bigl \{ \rho \le \chi^{\ast +} :
    \exists \alpha < \gamma_{\tau+1}^{< \alpha_1}
    ~ \alpha \mathbf{K}^{\ast <\alpha_{1}}(a_{\tau},
    \widetilde{\delta }_{\tau}^{<\alpha _{1}},
    \gamma_{\tau}^{<\alpha_{1}},
    \alpha, \rho, S_{\tau})  \bigr \};
\end{multline*}

\vspace{-24pt}
\begin{multline*}
    \alpha _{\tau }^{<\alpha _{1}} = \vspace{-12pt}
\\
    = \min_{\leq }
    \bigl \{ \alpha < \gamma_{\tau+1}^{< \alpha_1} :
    \alpha \mathbf{K}^{\ast <\alpha _{1}}(a_{\tau},
    \widetilde{\delta }_{\tau}^{<\alpha _{1}},
    \gamma_{\tau}^{<\alpha_{1}},
    \alpha ,\rho_{\tau }^{<\alpha _{1}},
    S_{\tau})  \bigr \};
\end{multline*}

\noindent and for \ $\alpha^1 =
\alpha_{\tau}^{<\alpha_1^{\Downarrow}}$:

\vspace{-12pt}
\begin{multline*}
    \check{\delta}_{\tau }^{<\alpha _{1}}= \min_{\leq }
    \bigl \{ \delta < \gamma_{\tau}^{< \alpha_1} :
    SIN_{n}^{<\alpha ^{1}}(\delta )\wedge SIN_{n+1}^{<\alpha ^{1}}
    \left[ <\rho_{\tau }^{<\alpha _{1}}\right] (\delta )) \bigr \};
\end{multline*}
\vspace{0pt}

\noindent The value \ $a_{\tau }^{<\alpha _{1}}$ \ is called,
remind, the characteristic of the matrix \ $\alpha S_{\tau
}^{<\alpha _{1}}$ \ on the carrier \ $\alpha_{\tau}^{<\alpha_1}$,
\ and of this carrier itself.

\noindent All the functions
\[
    a_{f}^{<\alpha _{1}}, \quad
    \widetilde{\delta }_{f}^{<\alpha _{1}},\quad
    \check{\delta}_{f}^{<\alpha _{1}},\quad
    \rho_{f}^{<\alpha _{1}}
\]
are called, for some brevity, the accessories of the functions
\[
    \alpha_f^{<\alpha_1}, \quad \alpha S_f^{<\alpha_1},
\]
and their values for the index \ $\tau$ \ are called also the
accessories of the values
\[
    \alpha_{\tau}^{<\alpha_1}, \quad \alpha S_{\tau}^{<\alpha_1};
\]
similarly the function \ $\alpha_f^{<\alpha_1}$ \ is called the
accessory of \ $\alpha S_f^{<\alpha_1}$, \ and its value \
$\alpha_{\tau}^{<\alpha_1}$ \ -- the accessory of the matrix \
$\alpha S_{\tau}^{<\alpha_1}$ \ below \ $\alpha_1$, \ and so on.
\\
\hspace*{\fill} $\dashv$

\end{definition}

\noindent The notion of characteristic is introduced in the
general case:

\begin{definition}
\label{8.4.} \hfill {} \newline
\hspace*{1em} We call as a characteristic of a matrix \ $S$ \ on a carrier \
$\alpha >\chi ^{\ast }$ \ the number \ $a(S,\alpha )=a$ \ defined in the
following way:
\[
    \hspace{0.4cm} \bigl(a=1\vee a=0 \bigr)\wedge \Big( a=0\longleftrightarrow
    \qquad\qquad\qquad\qquad\qquad\qquad\qquad\qquad
\]
\[
    \hspace{0.8cm} \longleftrightarrow \exists \tau_{1}^{\prime},
    \tau_{2}^{\prime},\tau _{3}^{\prime} < \alpha^{\Downarrow} \
    \big( A_{2}^{0 \vartriangleleft \alpha^{\Downarrow}}(\tau_{1}^{\prime},
    \tau_{2}^{\prime},\tau _{3}^{\prime},
    \alpha S_{f}^{<\alpha ^{\Downarrow }}) \wedge
\]
\[
    \qquad \qquad \qquad
    \wedge \forall \tau^{\prime\prime} \big( \tau_1^{\prime} <
    \tau^{\prime\prime} \le \tau_2^{\prime} \rightarrow
    a_{\tau^{\prime\prime}}^{<\alpha^{\Downarrow}}=1 \big) \wedge
    \alpha S_{\tau_2^{\prime}}^{<\alpha ^{\Downarrow }}=S \big) \Big).
\]

\noindent The matrix \ $S$ \ on its carrier \ $\alpha$ \ is called
the unit matrix on \ $\alpha$ \ iff it has the unit characteristic on
\ $\alpha$; \ otherwise it is called zero matrix on \ $\alpha$.
\hspace*{\fill} $\dashv$
\end{definition}

\noindent Thereafter when the \ $\alpha $-function \ $\alpha
S_{f}^{<\alpha _{1}}$ \ is defined, the priority belongs to \
$\alpha$-matrices possessing the greater characteristic.
\\
This circumstance, although making possible the solution of the
inaccessibility problem, complicates considerably the matrix
function theory as a whole because the restriction reasoning does
not work now freely: a situation concerning zero characteristic
may not be carried over to the part of the universe below which is
determined by unit characteristic, or by other reasons connected
with suppressibility.
\\
\quad \\

\begin{sloppypar}
Definition~\ref{8.3.} of \ $\alpha $-function and accompanying
ordinal functions follows the recursive definition~\ref{8.2.} and
since the functions \ \mbox{$X_i, \; i=\overline{1,5}$} \ are
defined in the formula \ $\alpha \mathbf{K}_{n+1}^{\exists}$ \ in
its subformula \ $\ A_7^{RC}$\ uniquely through its parameters by
this recursion, it is easy to see that functions \
$X_i[\alpha^0]$, $i=\overline{1,5}$, \ in definition \ref{8.2.}
coincide with corresponding functions
\end{sloppypar}
\begin{equation} \label{e8.2}
    a_{f }^{<\alpha^0},\quad
    \alpha S_{f }^{<\alpha^0},\quad
    \widetilde{\delta }_{f }^{<\alpha^0},\quad
    \rho_{f }^{<\alpha^0},\quad
    \alpha_{f }^{<\alpha^0},\quad
\end{equation}
for every cardinal \ $\alpha^0$, \ equinformative
with \ $\chi^{\ast}$.
\\
Cause of that we shall use their notations (\ref{e8.2}) instead of
corresponding notations of these functions \ $X_i[\alpha^0], ~
i=\overline{1,5}$ \ in formulas from definition~\ref{8.1.} that is
using these formulas but for the functions \ \mbox{$X_i[\alpha^0],
\; i=\overline{1,5}$}, \ replaced with corresponding functions
(\ref{e8.2}) for \ $\alpha^0=\alpha_1$; \ we shall even omit them
often for some brevity, when it will not cause misunderstanding
and when the context will point out them clearly.
\\
For instance, the formula \ $ A_0^{\vartriangleleft
\alpha_1}(\tau_1, \tau_2, \alpha S_f^{<\alpha_1})$ \ means, that
below \ $\alpha_1$ \ there holds
\[
    \tau_1+1 < \tau_2 \wedge (\alpha S_f^{<\alpha_1} \mbox{ is the
    function on } \left] \tau_1, \tau_2 \right[ ) \wedge
\]
\[
    \wedge \tau_1 = \min \bigl \{ \tau: \left] \tau, \tau_2 \right[ \subseteq
    dom(\alpha S_{f }^{<\alpha_1} )\wedge
\]
\[
    \wedge \chi^{\ast} \le  \gamma_{\tau_1}^{<\alpha_1} \wedge \gamma_{\tau_1}^{<\alpha_1}
    \in SIN_n^{<\alpha_1};
\]
the formula \ $ A_1^{\vartriangleleft \alpha_1}(\tau_1,
\tau_2, \alpha S_f^{<\alpha_1})$ \ means, that below \ $\alpha_1$ \ there holds
\[
    A_0^{\vartriangleleft \alpha_1}(\tau_1, \tau_2, \alpha S_f^{< \alpha_1}) \wedge
    \gamma_{\tau_2}^{<\alpha_1} \in SIN_n^{<\alpha_1};
\]
the formula
\[
    A_2^{0 \vartriangleleft \alpha^{\Downarrow}}(\tau_1, \tau_2, \tau_3,
    \alpha S_f^{<\alpha^{\Downarrow}}) \wedge
    \forall \tau \in \; ]\tau_1, \tau_2] ~
    a_{\tau}^{<\alpha^{\Downarrow}}=1
    \wedge \alpha S_{\tau_2}^{<\alpha^{\Downarrow}} = S
\]
means that here \ $\alpha^{\Downarrow}$ \ is the prejump cardinal
of \ $\alpha$ \ after \ $\chi^{\ast}$, \ and there is no \
$\alpha$-matrices admissible for \ $\gamma_{\tau_1}^{<\alpha_1}$ \
below \ $\alpha^{\Downarrow}$, \ and below \ $\alpha^{\Downarrow}$
\ there holds \
\[
    A_1^{\vartriangleleft \alpha^{\Downarrow}} (\tau_1, \tau_3, \alpha
    S_f^{<\alpha^{\Downarrow}}),
\]
where \ $\tau_2 \in \left] \tau_1, \tau_3 \right[$ \ is the first
ordinal at which monotonicity on \ $]\tau_1, \tau_3[$ \ of the
matrix function \ $\alpha S_f^{<\alpha^{\Downarrow}}$ \ fails, but
already below \ $\alpha^{\Downarrow}$, \ and, moreover, \ $\alpha
S_{\tau_2}^{<\alpha^{\Downarrow}} = S$ \ and all matrices \
$\alpha S_{\tau}^{<\alpha^{\Downarrow}}$ \ are of unit
characteristic on \ $]\tau_1, \tau_2]$ \  -- \ and so on.

\noindent Next, two easy remarks should be done:

1. All intervals \ $[ \gamma_{\tau_1}^{<\alpha_1},
\gamma_{\tau_2}^{<\alpha_1} [$ \ of definiteness below \
$\alpha_1$, \ considered in definition~\ref{8.1.} \ for the
functions
\[
    X_1 = \alpha S_f^{< \alpha_1}, \quad X_2 = a_f^{< \alpha_1},
\]
were of different types and were defined by different conditions,
but all of them include the condition of the interval \
$[\gamma_{\tau_1}^{<\alpha_1}, \gamma_{\tau_2}^{<\alpha_1}[$ \
maximality to the left:
\[
    A_0^{\vartriangleleft \alpha_1}(\tau_1, \tau_2, \alpha
    S_f^{<\alpha_1})
\]
which states, among other things, that the matrix function \
$\alpha S_f^{<\alpha_1}$ \ below \ $\alpha_1$ \ is defined on the
interval \ $\left] \tau_1, \tau_2 \right[$ \ and the ordinal \
$\tau_1$ \ is the \textit{minimal} one with this property and,
moreover, \ $\gamma_{\tau_1}^{<\alpha_1}$ \ is the \
$SIN_n^{<\alpha_1}$-cardinal. Due to this minimality it is not
hard to see, that \ $\alpha S_f^{<\alpha_1}$ \ is not defined for
this ordinal \ $\tau_1$ \ itself!

2. The notions of admissibility, priority and nonsuppression
should be distinguished. One can imagine two matrices \
$S^{\prime}, S^{\prime\prime}$ \ on their carriers \
$\alpha^{\prime}, \alpha^{\prime\prime}$ \ respectively along with
their corresponding accessories, both admissible for one cardinal
\ $\gamma_{\tau}^{<\alpha_1}$; when \ $S^{\prime}$ \ is of unit
characteristic on \ $\alpha^{\prime}$ \ it is always nonsuppressed
and has the priority over \ $S^{\prime\prime}$ \ of zero
characteristic on \ $\alpha^{\prime\prime}$. \ But even when there
is no such matrix \ $S^{\prime}$, \ still the matrix \
$S^{\prime\prime}$ \ on \ $\alpha^{\prime\prime}$ \ can be
suppressed, if there holds the suppression condition \ $A_5^{S,0}$
\ below \ $\alpha_1$; \ in any case every matrix, being
suppressed, can not be the value of the matrix function \ $\alpha
S_f^{<\alpha_1}$.
\\
So, for the interval $[\gamma_{\tau_1}^{<\alpha_1},
\gamma_{\tau_2}^{<\alpha_1}[$ \ maximal to the left below \
$\alpha_1$ \ there can be no value \ $\alpha S_{\tau}^{<\alpha_1}$
\ for \ $\tau=\tau_1$, \ but still it does not exclude the
existence of some matrix \textit{only admissible} (but suppressed)
for \ $\gamma_{\tau_1}^{<\alpha_1}$ \ below \ $\alpha_1$.
\\
\quad \\

And now, with all these comments in hand, let us see how
definition~\ref{8.2.} -- and, hence, definition~\ref{8.3.} --
works below \ $\alpha_1$ (we consider, remind, the most important
case when \ $\chi=\chi^{\ast}$, \ $\mathfrak{n} =
\mathfrak{n}^{\alpha}$).
\\

I. So, in the third part in the beginning of the formula
\[
    \alpha \mathbf{K}_{n+1}^{\exists < \alpha_1}
    (a, \delta, \gamma, \alpha, \rho, S)
\]
it is stated, that \ $S$ \ is the \ $\delta$-matrix on its carrier
\ $\alpha > \chi^{\ast}$, $\alpha < \alpha_1$ \ reduced to \
$\chi^{\ast}$ \ with the disseminator \ $\delta<\gamma$ \ and base
\ $\rho$:
\[
    S \vartriangleleft \rho = \widehat{\rho} \le \chi^{\ast +};
\]
the prejump cardinal \ $\alpha^{\Downarrow} =
\alpha_{\chi^{\ast}}^{\Downarrow}$ \ is limit for \
$SIN_n^{<\alpha^{\Downarrow}}$ \ and has the cofinality \ $\ge
\chi^{\ast +}$; \ the disseminator \ $\delta$ \ has
subinaccessibility below \ $\alpha^{\Downarrow}$ \ of the level \
$n$ \ and even of the level \ $n+1$ \ with the base \ $\rho$, \ that
is
\[
    \delta \in SIN_n^{<\alpha^{\Downarrow}} \cap
    SIN_{n+1}^{<\alpha^{\Downarrow}}[<\rho].
\]

II. Then below \ $\alpha^{\Downarrow}$ \ there are defined the
functions \ $X_i, \ i = \overline{1,5}$ \ on pairs \ $(\alpha^0,
\tau^{\prime}) \in
\mathcal{A}_{\chi^{\ast}}^{\alpha^{\Downarrow}}$, \ where
cardinals \ $\alpha^0 \in \left] \chi^{\ast}, \alpha^{\Downarrow}
\right]$ \ are equinformative with \ $\chi^{\ast}$ \ and there
exist cardinals \ $\gamma_{\tau^{\prime}}^{<\alpha^0}$.
\\
All these functions  are recursively defined through definition of the
functions \ $X_i^0, \ i = \overline{1,5}$, \ by means of the
recursive condition \ $A_7^{RC}$:
\[
    X_1[\alpha^0] = \alpha S_f^{<\alpha^0}, \quad
    X_2[\alpha^0] = a_f^{<\alpha^0}, \quad
    X_3[\alpha^0] = \widetilde{\delta}_f^{<\alpha^0},
\]
\[
    X_4[\alpha^0] = \rho_f^{<\alpha^0}, \quad
    X_5[\alpha^0] = \alpha_f^{<\alpha^0}.
\]
The aim of this definition -- to receive the resulting matrix
function \ $\alpha S_f^{<\alpha^0}$, \ but the first it is
introduced just the characteristic function
\[
    X_2^0 = X_2[\alpha^0] = a_f^{<\alpha^0}.
\]
This function accepts the \textit{maximal} possible values, unit
or zero, that are characteristics of admissible matrices below \
$\alpha^0$, \ \textit{but only not zero characteristic of
suppressed zero matrices} \ $S^{\prime\prime}$ \ on their carriers
\ $\alpha^{\prime\prime}$, \ that satisfies the suppression
condition below \ $\alpha^0$:
\[
    A_5^{S,0 \vartriangleleft \alpha^0}(0, \gamma_{\tau^{\prime}}^{<\alpha^0}, \alpha^{\prime\prime},
    \rho^{\prime\prime}, S^{\prime\prime}, X_1^0|\tau^{\prime},
    X_2^0|\tau^{\prime}, X_1|^1\alpha^0, X_2|^1\alpha^0),
\]
where the functions here
\[
    X_1^0|\tau^{\prime}, \quad X_2^0|\tau^{\prime}, \quad
    X_1|^1\alpha^0, \quad X_2|^1\alpha^0
\]
are already defined. And everywhere further such suppressed zero
matrices are systematically rejected.
\\
After the characteristic function \ $X_2[\alpha^0] = a_f^{<\alpha_0}$ \ is defined,
all remained functions
\[
    X_1[\alpha^0], \quad X_i[\alpha^0], \quad i=\overline{3,5}
\]
are defined one by one successively through the minimization of
their admissible and nonsuppressed values.
\\
So, the next it is defined the matrix function \ $X_1[\alpha^0] =
\alpha S_f^{<\alpha^0}$, \ after that the corresponding
disseminator function \ $X_3[\alpha^0] =
\widetilde{\delta}_f^{<\alpha^0}$, \ then the data base function \
$X_4[\alpha^0] = \rho_f^{<\alpha^0}$, \ and, in last turn, the
carrier function \ $X_5[\alpha^0] = \alpha_f^{<\alpha^0}$ \ is
defined.
\\
The values of every subsequent of these functions depend
essentially on the values of the previous ones.

III. After all these functions are constructed for every
\[
    \alpha^0 \in ] \chi^{\ast}, \alpha^{\Downarrow}[\;,
\]
the definition passes to the cardinal
\[
    \alpha^0 = \alpha^{\Downarrow}
\]
and after that defines the characteristic of the matrix \ $S$ \ on its carrier \
$\alpha$ \ itself:
\\
$S$ \ on \ $\alpha$ \ receives \textit{zero characteristic}, if it
participate in the following violation of the matrix function
\[
    X_1[\alpha^{\Downarrow}] = \alpha S_f^{<\alpha^{\Downarrow}}
\]
monotonicity below \ $\alpha^{\Downarrow}$, \ when below \
$\alpha^{\Downarrow}$ \ there holds the condition
\begin{multline*}
    \exists \tau_1^{\prime}, \tau_2^{\prime}, \tau_3^{\prime} <
    \alpha^{\Downarrow} \big(A_2^{0 \vartriangleleft \alpha^{\Downarrow}}
    (\tau_1^{\prime}, \tau_2^{\prime}, \tau_3^{\prime},
    \alpha S_f^{<\alpha^{\Downarrow}} ) \wedge
\\
    \wedge \forall \tau^{\prime\prime} \in \; ]\tau_1^{\prime}, \tau_2^{\prime}] ~
    a_{\tau^{\prime\prime}}^{<\alpha^{\Downarrow}}=1 \wedge
    \alpha S_{\tau_2^{\prime}}^{<\alpha^{\Downarrow}} = S \big);
\end{multline*}
otherwise \ $S$ \ on \ $\alpha$ \ receives
\textit{unit characteristic}.

IV. And in the last turn this definition forms the \textit{closing
condition} for \ $S$ \ on \ $\alpha$:
\\
If \ $S$ \ is \textit{zero} matrix on \ $\alpha$ \ and its
admissible disseminator \ $\delta$ \ falls in some maximal block
of type \ $\eta^{\prime}$ \ below \ $\alpha^{\Downarrow}$
\[
    [ \gamma_{\tau_1^{\prime}}^{<\alpha^{\Downarrow}},
    \gamma_{\tau_3^{\prime}}^{<\alpha^{\Downarrow}} [
\]
{\sl severe} for \ $S$ \ on \ $\alpha$, \ that is if there holds
\[
    \gamma_{\tau_1^{\prime}}^{<\alpha^{\Downarrow}} \le \delta <
    \gamma_{\tau_3^{\prime}}^{<\alpha^{\Downarrow}} \wedge
    A_4^{M b \vartriangleleft \alpha^{\Downarrow}}
    (\tau_1^{\prime}, \tau_1^{\prime\prime}, \tau_2^{\prime},
    \tau_3^{\prime}, \eta^{\prime}, \alpha S_f^{<\alpha^{\Downarrow}},
    a_f^{<\alpha^{\Downarrow}} )
\]
below \ $\alpha^{\Downarrow}$, \ then there is demanded the
admissible data base \ $\rho$ \ of the disseminator \ $\delta$ \
of \ $S$ \ on \ $\alpha$ \ but only such that
\[
    \eta^{\prime} < \rho \vee \rho = \chi^{\ast +}.
\]
So, this case \textit{hampers} the using of such \ $S$ \ on \
$\alpha$ \ considerably; besides that \ $S$ \ on \ $\alpha$ \ must
be nonsuppressed; in any other cases no requirements are inflicted
on \ $S$ \ on \ $\alpha$.
\\
But remind, that the base \ $\rho = \chi^{\ast +}$ \ and every
\textit{unit} matrix are always admissible and nonsuppressed;
every matrix is nonsuppressed for \ $\gamma \notin SIN_n$ \ in any
case.

After that this definition forms the conjunction \ $\alpha
\mathbf{K}^{<\alpha_1}$:
\[
    \mathbf{K}_n^{\forall < \alpha_1} (\gamma,
    \alpha^{\Downarrow}) \wedge \alpha
    \mathbf{K}_{n+1}^{\exists \vartriangleleft \alpha_1}
    (a, \delta, \gamma, \alpha, \rho, S)
    \wedge \alpha < \alpha_1
\]
where is required in addition, as usual, that \
$\alpha^{\Downarrow}$ \ preserves all \
$SIN_n^{<\alpha_1}$-cardinals \ $\le \gamma$ \ below \ $\alpha_1$;
\ and, at last, there arises the formula \ $\alpha
\mathbf{K}^{\ast <\alpha_1}$ \ from this formula \ $\alpha
\mathbf{K}^{<\alpha_1}$ \ under the requirement of the
nonsuppression of zero matrix \ $S$ \ on \ $\alpha$ \ below \
$\alpha_1$.
\\

Since definition \ref{8.3.} of the matrix \ $\alpha$-function and
of the accompanying functions follows definition \ref{8.2.}, there
holds the next obvious lemma, which actually repeats this
definition. Here is used the notion of generating
eigendisseminator \ $\check{\delta}^S$ \ for arbitrary matrix \
$S$ \ on a carrier \ $\alpha$, \ that is, remind, the minimal
disseminator for \ $S$ \ on \ $\alpha$ \ with the minimal possible
base \ $\rho^S = \widehat{\rho_1}$, $\rho_1 = Od(S)$ \ (see
\cite{Kiselev9}, \cite{Kiselev8}).

\begin{lemma}
\label{8.5.} \hfill {} \newline \hspace*{1em} Let \ $S$ \ be an
arbitrary \ $\alpha$-matrix reduced to \ $\chi ^{\ast }$ \ of
characteristic \ $a$ \ on a carrier \ $\alpha < \alpha_1$, {\sl
admissible} for \ $\gamma _{\tau }^{<\alpha _{1}}$ \ along with
its disseminator \ $\widetilde{ \delta }$, \ generating
disseminator \ $\check{\delta}$ \ with a base \ $\rho $, and
generating eigendisseminator \ $\check{\delta}^S$ \ below \
$\alpha _{1}$, \ then for the prejump cardinal \ $\alpha
^{\Downarrow }$ \ after \ $\chi^{\ast}$ \ there holds below \
$\alpha_1$:
\newline

\noindent \medskip \emph{1)}\quad $\forall \gamma \leq \gamma
_{\tau }^{<\alpha _{1}}(SIN_{n}^{<\alpha _{1}}(\gamma
)\longrightarrow SIN_{n}^{<\alpha ^{\Downarrow }}(\gamma ))$~;
\newline

\noindent \medskip \emph{2)}\quad $\chi ^{\ast }<\widetilde{\delta }<\gamma
_{\tau }^{<\alpha _{1}}<\alpha ^{\Downarrow }\wedge S\vartriangleleft \rho
\leq \chi ^{\ast +}\wedge \rho =\widehat{\rho }$~; \newline

\noindent \medskip \emph{3)}\quad $\widetilde{\delta }\in
SIN_{n}^{<\alpha ^{\Downarrow }}\cap SIN_{n+1}^{<\alpha
^{\Downarrow }}\left[ <\rho \right] $; \ analogously for \
$\check{\delta}$; \newline

\noindent \medskip \emph{4)}\quad $\sup SIN_{n}^{<\alpha
^{\Downarrow }}=\alpha ^{\Downarrow}\wedge cf(\alpha ^{\Downarrow
})\geq \chi ^{\ast +}$; \newline

\noindent \medskip \emph{5)}\quad $a=0\longleftrightarrow \exists
\tau_{1}^{\prime}, \tau_{2}^{\prime}, \tau_{3}^{\prime}
~\big(A_{2}^{0 \vartriangleleft
\alpha^{\Downarrow}}(\tau_{1}^{\prime}, \tau_{2}^{\prime},
\tau_{3}^{\prime}, \alpha S_f^{<\alpha^{\Downarrow}}) \wedge $
\newline

\noindent \medskip \hspace{10em} $ \wedge \forall
\tau^{\prime\prime} \in \; ]\tau_1^{\prime}, \tau_2^{\prime}] ~
a_{\tau^{\prime\prime}}^{<\alpha^{\Downarrow}}=1 \wedge \alpha
S_{\tau_{2}^{\prime}}^{<\alpha ^{\Downarrow }}=S\big)$;
\newline

\noindent \medskip \emph{6)}\quad $a=0\longrightarrow \forall
\tau_{1}^{\prime}, \tau_{1}^{\prime\prime}, \tau_{2}^{\prime},
\tau_{3}^{\prime}, \eta^{\prime} \bigl[ \gamma_{\tau_1^{\prime}}^{<\alpha
^{\Downarrow }}\leq \widetilde{\delta} <
\gamma_{\tau_3^{\prime}}^{<\alpha ^{\Downarrow }}\wedge $
\newline

\noindent \medskip \hspace{1em} $\wedge A_{4}^{M b
\vartriangleleft \alpha^{\Downarrow} } (\tau_{1}^{\prime},
\tau_{1}^{\prime\prime}, \tau_{2}^{\prime}, \tau_{3}^{\prime},
\eta^{\prime}, \alpha S_f^{<\alpha^{\Downarrow}},
a_f^{<\alpha^{\Downarrow}} ) \longrightarrow \eta^{\prime} <\rho \vee \rho
=\chi ^{\ast +} \bigr]$;
\newline

\noindent \medskip \emph{7)}\quad (i) \ $\check{\delta}^S \le
\check{\delta}\leq \widetilde{\delta }<\gamma _{\tau }^{<\alpha
_{1}}$; \newline

\noindent \medskip \qquad (ii) \ if \ $\widetilde{\delta }$ \ is
the {\sl minimal} floating disseminator of \ $S$ \ on \ $\alpha $ \ with
the minimal base \ $\rho$ \ admissible for \ $\gamma _{\tau
}^{<\alpha _{1}}$ \ along with \ $\rho $, \ then:
\[
    a=1\longrightarrow \widetilde{\delta } = \check{\delta}^S
    \wedge \rho = \rho^S = \widehat{\rho_1}, \wedge \rho_1 = Od(S),
\]
that is when \ $S$ \ is the unit matrix on \ $\alpha$, \ then \
$\widetilde{\delta}$ \ is the generating eigendisseminator \
$\check{\delta}^S$ \ of \ $S$ \ on \ $\alpha$ \ with the base \
$\rho^S$;
\\

\noindent \medskip \emph{8)}\quad there exist the minimal carrier
\ $\alpha^{\prime} < \gamma_{\tau+1}^{<\alpha_1}$ \ of \ $S$ \ of
the same characteristic \ $a$ \ admissible for \
$\gamma_{\tau}^{<\alpha_1}$ \ along with the same accessories \
$\widetilde{\delta}$, $\rho$ \ below \ $\alpha_1$:
\[
    \gamma _{\tau }^{<\alpha _{1}}<\alpha^{\prime} <
    \gamma _{\tau +1}^{<\alpha_{1}}~;
\]

\noindent analogously for nonsuppressibility  of \ $S$ \ for \
$\gamma_{\tau}$ \ along with its accessories.
\end{lemma}

\noindent \textit{Proof.} \ It remains to prove the last two
statements; the upper index \ $< \alpha_1$ \ will be dropped.
\\
So, let us consider the matrix \ $S$ \ of characteristic \ $a$ \
on its carrier \ $\alpha < \alpha_1$ \ admissible for \
$\gamma_{\tau}^{<\alpha_1}$ \ along with its disseminator \
$\widetilde{\delta}$ \ and base \ $\rho$. \ Statement 7)~$(i)$ is
obvious; as to 7)~$(ii)$ let us consider \ $a=1$, \ then the base
\[
    \rho = \rho^S = \widehat{\rho_1}, \ \rho_1 = Od(S)
\]
along with the minimal disseminator
\[
    \check{\delta}^{S} \in SIN_n^{<\alpha^{\Downarrow}} \cap
    SIN_{n+1}^{<\alpha^{\Downarrow}}[<\rho^S]
\]
evidently fulfill all requirements of the condition
\[
    \alpha \mathbf{K}(a, \check{\delta}^S, \gamma, \alpha, \rho^S, S)
\]
up to the last it conjunctive constituent \ $\mathbf{K}^0$.
\\
But the latter is fulfilled also, because for \ $a=1$ \ its
premise fails.
\\
Thus the whole \ $\alpha \mathbf{K}$ \ is fulfilled and \
$\widetilde{\delta} = \check{\delta}^S$, $\rho=\rho^S$.

\noindent Turning to the proof of 8) it is not hard to apply lemma
3.2~\cite{Kiselev9} (about restriction) just as it was done in the
lemma 5.17~2)~$(ii)$ proof. Nevertheless, this application
presents the typical reasoning, which will be used further in
various important cases, so one should accept it in details.
\\
First, it was assumed above that \ $\alpha_1$ \ is limit for the
class \ $SIN_{n-1}^{<\alpha_1}$ (remind the convention after
(\ref{e7.1})), therefore it always exist \
$\gamma_{\tau+1}^{<\alpha_1}$ \ for every \
$\gamma_{\tau}^{<\alpha_1}$.
\\
Next, suppose that the matrix \ $S$ \ with the disseminator \
$\delta$ \ and base \ $\rho$ \ on the carrier
\[
    \alpha \in \; ]\gamma_\tau^{<\alpha_1}, \alpha_1 [
\]
is admissible for \ $\gamma_{\tau}^{<\alpha_1}$ \ below \
$\alpha_1$, \ then it holds the following proposition \
$\varphi(\chi^{\ast}, \delta, \gamma_{\tau}^{<\alpha_1}, \rho,
S)$:
\[
    \exists \alpha^{\prime} ~ ( \gamma_{\tau}^{<\alpha_1} < \alpha^{\prime}
    \wedge \alpha \mathbf{K}(\delta, \gamma_{\tau}^{<\alpha_1},
    \alpha^{\prime}, \rho, S))
\]
below \ $\alpha_1$, \ that is after its \
$\vartriangleleft$-bounding by the cardinal \ $\alpha_1$. \ This
proposition \ $\varphi$ \ itself is from the class \ $\Sigma_n$, \
because it includes \ $\Sigma_n$-formula \
$\mathbf{K}_n^{\forall}$. \ But let us use the cardinal
\[
    \gamma_{\tau^n}=\sup\left\{\gamma \le \gamma_\tau^{<\alpha_1} :
    \gamma \in SIN_n^{<\alpha_1} \right\};
\]
by lemma~3.4~\cite{Kiselev9} \ $\gamma_{\tau^n}$ \ belongs to \
$SIN_n^{<\alpha_1}$ \ too. Now let us replace in the formula \
$\alpha \mathbf{K}$ \ its subformula \ $\mathbf{K}_n^{\forall}$ \
with the \ $\Delta_1$-formula
\[
    SIN_n^{<\alpha^{\Downarrow}}(\gamma_{\tau^n}),
\]
then the \ $\Sigma_n$-formula \ $\alpha \mathbf{K}$ \ turns into
some \ $\Pi_{n-2}$-formula, which we shall denote through \
$\alpha \mathbf{K}_{n-2}$. \ Consequently, the formula \ $\varphi$
\ turns into some \ $\Sigma_{n-1}$-formula \
$\varphi_{n-2}(\chi^{\ast}, \delta, \gamma_{\tau^n},
\gamma_{\tau}^{<\alpha_1}, \rho, S)$:
\[
    \exists \alpha^{\prime} ~ (\gamma_{\tau}^{<\alpha_1} < \alpha^{\prime} \wedge
    \alpha \mathbf{K}_{n-2}(\delta, \gamma_{\tau^n}, \gamma_{\tau}^{<\alpha_1},
    \alpha^{\prime}, \rho, S))
\]
precisely with the same content below \ $\alpha_1$, \ and there
holds
\[
    \varphi_{n-2}^{\vartriangleleft \alpha_1}(\chi^{\ast}, \delta,
    \gamma_{\tau^n}, \gamma_{\tau}^{<\alpha_1}, \rho, S).
\]
The last proposition contains individual constants
\[
    \chi^{\ast}, \delta, \gamma_{\tau^n}, \gamma_{\tau}^{<\alpha_1}, \rho, S
\]
less then the \ $SIN_{n-1}^{<\alpha_1}$-cardinal \
$\gamma_{\tau+1}^{<\alpha_1}$ \ and therefore this cardinal
restricts this proposition by lemma 3.2~\cite{Kiselev9} (where \
$n$ \ replaced with \ $n-1$), that is there holds the formula
\[
    \exists \alpha^{\prime} \in [\gamma_{\tau}^{<\alpha_1},
    \gamma_{\tau+1}^{<\alpha_1} [ \quad
    \alpha \mathbf{K}_{n-2}^{<\alpha_1}(\delta,
    \gamma_{\tau^n}, \gamma_{\tau}^{<\alpha_1},
    \alpha^{\prime}, \rho, S)
\]
and \ $S$ \ receives its carrier \ $\alpha^{\prime} \in
[\gamma_{\tau}^{<\alpha_1}, \gamma_{\tau+1}^{<\alpha_1}[$ \
admissible for \ $\gamma_{\tau}^{<\alpha_1}$ \ along with its
previous disseminator and data base.
\\
The part 8) for nonsuppressibility will not be used up to \S 11
and there we shall return to it once more.
\\
\hspace*{\fill} $\dashv$
\\

It is not hard to see that the functions introduced above in
definition \ref{8.3.} possess many simple properties of \ $\delta
$-functions and their accompanying functions, so the proofs of the
following three lemmas are analogous to those of lemmas
\ref{7.3.}, \ref{7.4.} (or lemmas~5.16, 5.15 \cite{Kiselev9}) and
lemma \ref{7.5.}:

\begin{lemma}
\label{8.6.} \hfill {} \newline \hspace*{1em} For \ $\alpha
_{1}<k$ \ the formulas \ $\alpha \mathbf{K}^{<\alpha _{1}}$, \
$\alpha \mathbf{K}^{\ast <\alpha _{1}}$ \ belong to \ $\Delta
_{1}$ \ and therefore all functions from definition~\ref{8.3.}:
\begin{equation*}
a_{f}^{<\alpha _{1}}, ~~\alpha S_{f}^{<\alpha _{1}},
~~\widetilde{\delta } _{f}^{<\alpha _{1}}, ~~\rho _{f}^{<\alpha
_{1}}, ~~\alpha _{f}^{<\alpha _{1}}, ~~\check{\delta}_{f}^{<\alpha
_{1}}
\end{equation*}
are $\Delta _{1}$-definable through \ $\chi ^{\ast },\alpha _{1}$.
\ For \ $\alpha _{1}=k$ \ the formula \ $\alpha \mathbf{K}$ \
belongs to \ $\Sigma _{n+1}$.
\\
\hspace*{\fill} $\dashv$
\end{lemma}

\begin{lemma}
\label{8.7.} \emph{(About \ $\alpha $-function absoluteness)} \newline
\hspace*{1em} Let \ $\chi ^{\ast }<\gamma _{\tau +1}^{<\alpha _{1}}<\alpha
_{2}<\alpha _{1}\leq k$, \quad $\alpha _{2}\in SIN_{n-2}^{<\alpha _{1}}$ \
and
\begin{equation*}
(\gamma _{\tau }^{<\alpha _{1}}+1)\cap SIN_{n}^{<\alpha _{2}}=(\gamma _{\tau
}^{<\alpha _{1}}+1)\cap SIN_{n}^{<\alpha _{1}}.
\end{equation*}

\noindent \medskip \emph{1)}\quad Then on the set
\[
    T = \{\tau^{\prime }:~\chi ^{\ast }\leq
    \gamma _{\tau ^{\prime }}^{<\alpha_{2}}\leq
    \gamma _{\tau }^{<\alpha _{1}}\}
\]
the admissibility below \ $\alpha_2$ \ coincides with the
admissibility below \ $\alpha_1$: \ for every \ $\tau^{\prime} \in
T$ \ and a matrix \ $S^{\prime}$ \ on its carrier \
\mbox{$\alpha^{\prime} \in \; ] \gamma_{\tau^{\prime}}^{<\alpha_2},
\gamma_{\tau+1}^{<\alpha_1}[$}
\[
    \alpha \mathbf{K}^{<\alpha_2}(
    \gamma_{\tau^{\prime}}^{<\alpha_2}, \alpha^{\prime},
    S^{\prime}) \leftrightarrow
    \alpha \mathbf{K}^{<\alpha_1}(
    \gamma_{\tau^{\prime}}^{<\alpha_2}, \alpha^{\prime},
    S^{\prime});
\]
\emph{2)}\quad on the set
\[
    \big \{ \tau^{\prime}: \chi^{\ast} \le \gamma_{\tau^{\prime}}^{<\alpha_2}
    \le \gamma_{\tau}^{<\alpha_1} \wedge ( a_{\tau^{\prime}}^{<\alpha_2} = 1 \vee
    \neg SIN_n^{<\alpha_2}(\gamma_{\tau^{\prime}}^{<\alpha_2}) ) \big \}
\]
the functions (\ref{e8.2}) below \ $\alpha^0 =
\alpha_{2}$ \ coincide respectively with the functions
(\ref{e8.2}) below \ $\alpha^0 = \alpha _{1}$.
\\
\hspace*{\fill} $\dashv$
\end{lemma}

\begin{lemma}
\label{8.8.} \emph{(About disseminator)} \newline
\emph{1)}\quad Let

\quad \medskip (i) \quad $] \tau_{1}, \tau_{2}[\; \subseteq dom
\bigl( \alpha S_{f}^{< \alpha_{1}} \bigr ) $, \quad
$\gamma_{\tau_{2}} \in SIN_{n}^{< \alpha_{1}}$;

\quad \medskip (ii) \quad $\tau_{3} \in dom \bigl ( \alpha
S_{f}^{< \alpha_{1}} \bigr ),\quad \tau_{2} \leq \tau_{3}$;

\quad \medskip (iii) \quad $\widetilde{\delta}_{\tau_{3}}^{< \alpha_{1}} <
\gamma_{\tau_{2}}^{< \alpha_{1}} $ \ and \ $a_{\tau_{3}}^{< \alpha_{1}}=0$.
\newline
Then
\begin{equation*}
\widetilde{\delta}_{\tau_{3}}^{ < \alpha_{1}} \leq
\gamma_{\tau_{1}}^{ < \alpha_{1}}.
\end{equation*}
Analogously for \ $\check{\delta}_{\tau_{3}}^{ < \alpha_{1}}$.
\\
\quad \\
\emph{2)}\quad Let \ $\alpha$-matrix \ $S$ \ of characteristic \
$a$ \ on a carrier \ $\alpha$ \ be admissible for \
$\gamma_{\tau}^{<\alpha_1}$ \ along with its disseminator \
$\widetilde{\delta}$ \ and base \ $\rho$ \ below \ $\alpha_1$, \
then

\[
    \{ \tau^\prime : \widetilde{\delta}  <
    \gamma_{\tau^\prime}^{ < \alpha_{1}} < \gamma_{\tau}^{ <
    \alpha_{1}} \} \subseteq dom \bigl ( \alpha S_{f}^{< \alpha_{1}}
    \bigr ).
\]
\end{lemma}

\noindent \textit{Proof.} \ 1) The reasoning forthcoming is
analogous to the proof of lemma 7.5~1), but now the special
properties of matrix disseminator of unit or zero characteristic
involves the singular situation. Therefore, here one should use
the following argument that will be applied further in various
significant typical situations; here it is presented in outline;
the upper indices \ $< \alpha_{1} $, $\vartriangleleft \alpha_{1}$
\ will be dropped for shortness.
\\
Suppose 1) fails; let us consider the matrix \ $S^3 = \alpha
S_{\tau_{3}} $ \ of characteristic \ $a^3 = a_{\tau_{3}}=0 $ \ on
the carrier \ $\alpha_{\tau_3}$ \ with the prejump cardinal \
$\alpha^3 = \alpha_{\tau_3}^{\Downarrow}$, \ possessing the
disseminators \ $\check{\delta}^3 = \check{\delta}_{\tau_{3}}$, \
$\widetilde{\delta}^3 = \widetilde{\delta}_{\tau_{3}}$ \ with the
base \ $\rho^3 = \rho_{\tau_{3}}$, \ and suppose that
\\
\begin{equation} \label{e8.3}
\gamma_{\tau_1} < \widetilde{\delta}^3 < \gamma_{\tau_2}.
\end{equation}
One should consider here the minimal ordinal \
$\tau_1$ \ with the property $(i)$.
\\
By definition \ref{8.3.} \ the proposition \ $\alpha \mathbf{K} $
\ holds and hence there holds the proposition \ $\mathbf{K}^{0}$:

\begin{equation*}
\begin{array}{l}
a^3 = 0 \longrightarrow \forall \tau_1^{\prime},
\tau_1^{\prime\prime}, \tau_2^{\prime}, \tau_3^{\prime},
\eta^{\prime} < \alpha^3 ~ \bigl[ \gamma
_{\tau_1^\prime}^{<\alpha^3} \leq
\widetilde{\delta}^3 < \gamma_{\tau_3^\prime}^{<\alpha^3} \wedge \\
\quad \\
~ \wedge A_{4}^{M b \vartriangleleft \alpha^3} ( \tau_1^{\prime},
\tau_1^{\prime\prime}, \tau_2^{\prime}, \tau_3^{\prime},
\eta^{\prime}, \alpha S_f^{<\alpha^3}, a_f^{<\alpha^3})
\longrightarrow \eta^{\prime} < \rho^3 \vee \rho^3 = \chi^{\ast +}
\bigr].
\end{array}
\end{equation*}

\noindent Suppose that there exist some ordinals \ $\tau_1^{\prime},
\tau_1^{\prime\prime}, \tau_2^{\prime}, \tau_3^{\prime},
\eta^{\prime} < \alpha^3$ \ that fulfill the premise of this
proposition: \vspace{6pt}
\begin{equation} \label{e8.4}
\gamma_{\tau_{1}^{\prime}}^{< \alpha^{3}} \leq
\widetilde{\delta}^3 < \gamma_{\tau_{3}^{\prime}}^{< \alpha^{3}}
\wedge A_{4}^{M b \vartriangleleft \alpha^3} ( \tau_1^{\prime},
\tau_1^{\prime\prime}, \tau_2^{\prime}, \tau_3^{\prime},
\eta^{\prime}, \alpha S_{f}^{< \alpha^{3}}, a_{f}^{< \alpha^{3}}
)~.
\end{equation}
\vspace{0pt}

\noindent It should be pointed out again that due to \ $A_{4}^{M b
\vartriangleleft \alpha^3} $ \ these ordinals are defined through
\ $\widetilde{\delta}^3, \alpha^{3} $ \ \textit{uniquely}. Since \
$\gamma_{\tau_{2}} \in SIN_{n} $ \ and \ $\gamma_{\tau_{1}}$ \ is
the minimal one can see, that due to supposition (\ref{e8.3})
there comes
\begin{equation} \label{e8.5}
    \gamma_{\tau_{1}^{\prime}} < \check{\delta}^3 =
    \widetilde{\delta}^3 < \gamma_{\tau_{2}}
\end{equation}
as a result of the minimizing of the disseminator \
$\widetilde{\delta}_{\tau_{3}} $ \ inside the interval \ $[
\gamma_{\tau_{1}^{\prime}}, \gamma_{\tau_{2}} [ $ \ according to
definition~8.3. Now we come to the situation from the proof of
lemma~7.5~1) and it remains to repeat its arguments and
to use the \ $\underline{\lessdot}$-minimal matrix \ $S^m \lessdot
S^3$ \ on some carrier \ $\alpha^m \in ] \gamma_{\tau_3}, \alpha^3
[$ \ of characteristic \ $a^m$, \ admissible and nonsuppressed for
\ $\gamma_{\tau_3}$ \ along with its minimal disseminator \
$\widetilde{\delta}^m < \gamma_{\tau_2}$ \ and base \ $\rho^m <
Od(S^3)$, \ because the suppression of \ $S^m$ \ for \
$\gamma_{\tau_3}$ \ implies the suppression of the matrix \ $S^3$
\ itself for \ $\gamma_{\tau_3}$, \ though it is nonsuppressed by
definition (below \ $\alpha_1$).
\\
Now there comes the contradiction: since \ $S^m \lessdot S^3$ \
and \ $a^3=0$ \ then by definition \ref{8.3.} the matrix \ $S^3$ \
cannot be the minimal value \ $\alpha S_{\tau_3}$.

If there is no such ordinals \ $\tau_1^{\prime},
\tau_1^{\prime\prime}, \tau_2^{\prime}, \tau_3^{\prime},
\eta^{\prime}$, \ the proposition \ $\mathbf{K}^{0} $ \ survives
evidently under minimizing the disseminator \
$\widetilde{\delta}_{\tau_{3}} $ \ inside \ $[ \gamma_{\tau_{1}},
\gamma_{\tau_{2}} [ $ \ and so \ $\widetilde{\delta}_{\tau_{3}}
\leq \gamma_{\tau_{1}} $, \ otherwise again there holds \
$\gamma_{\tau_{1}} < \check{\delta}_{\tau_{3}} =
\widetilde{\delta}_{\tau_{3}} < \gamma_{\tau_{2}} $ \ and the same
reasoning provides the same contradiction.

Turning to 2), one should simply notice, that this statement
repeats the previous lemmas~5.17~2) \cite{Kiselev9}, 7.5~2) in the
following form:
\\
the matrix \ $S$ \ being admissible for \ $\gamma_{\tau}$ \ on its
carrier \ $\alpha_\tau$, \ by lemma \ref{8.5.}~8) and definition
\ref{8.2.} remains still admissible and nonsuppressed for every \
$\gamma_{\tau^\prime} < \gamma_{\tau}$, \ such that \
$\widetilde{\delta} < \gamma_{\tau^\prime}$, \ along with the same
accompanying ordinals \ $a$, $\widetilde{\delta}$, $\rho$,
$\alpha$. \ For the unit characteristic \ $a=1$ \ it is obvious;
but for \ $a=0$ \ this lemma will be used only in \S 11 and there
we shall turn to its proof in the exposition detailed more.
\hspace*{\fill} $\dashv$
\\
\quad \\

The following lemmas confirm the further extension of the \
\mbox{$\delta$-functions} theory on \ $\alpha$-functions and are
analogous to lemmas~7.6, 7.7 about \ $\delta$-function
definiteness on the final subinterval of the inaccessible cardinal
\ $k$.
\\
So, the next lemma shows, that there exists the cardinal \ $\delta
<k$ \ such that
\begin{equation*}
\{ \tau: \delta < \gamma _{\tau } < k \} \subseteq  dom ( \alpha
S_{f} );
\end{equation*}
more precisely:
\\
\begin{lemma}
\label{8.9.} \emph{(About \ $\alpha $-function definiteness)}
\newline \hspace*{1em} There exist cardinals \ $\delta < \gamma <
k$ \ such that for every \ \mbox{$SIN_n$-}car\-di\-nal \ $\alpha_1
> \gamma$, $\alpha_1 < k$ \ limit for \ $SIN_n \cap \alpha_1$ \ the function \ $\alpha
S_f^{<\alpha_1}$ \ is defined on the nonemty set
\[
    T^{\alpha_1} = \{\tau: \delta < \gamma_{\tau}^{<\alpha_1} < \alpha_1\}.
\]
The minimal of such cardinals \ $\delta$ \ is denoted by \ $\alpha
\delta ^{\ast }$, \ its successor in \ $SIN_n$ \ by \ $\alpha
\delta ^{\ast 1}$ \ and the following corresponding ordinals are
introduced:
\begin{align*}
    & \alpha \tau_1^{\ast} = \tau(\alpha \delta^{\ast}),
    \ \alpha \tau^{\ast 1} = \tau(\alpha \delta^{\ast 1}),
\\
    & \mathit{so\ that\ \ } \alpha \delta^{\ast} = \gamma_{\alpha
    \tau_1^{\ast}}, \ \alpha \delta^{\ast 1} = \gamma_{\alpha
    \tau^{\ast 1}},
\\
    & \mathit{and\ \ } \alpha^{\ast 1} = \alpha_{\alpha \tau^{\ast
    1}}^{< \alpha_1 \Downarrow}, \ \alpha \rho^{\ast 1} = \rho_{\alpha \tau^{\ast
    1}}^{< \alpha_1}.
\end{align*}
\end{lemma}

\noindent \textit{Proof} consists in the application of lemma~6.14
\cite{Kiselev9} just as it was done in the proof of lemma~7.6 but
for the \textit{greater} reducing cardinal
\begin{equation*}
\chi = ( \chi^{\ast})^{+ \omega_{0}} \mathrm{\quad and \quad}  \alpha_{1} =
k, \quad m=n+1.
\end{equation*}
The resulting function \ $\mathfrak{A}$, \ being defined on the
set
\begin{equation*}
T = \{ \tau : \gamma_{\tau_{0}} \leq \gamma_{\tau} < k \},
\end{equation*}
should be treated in the following way:
\\
Let us consider by lemma~6.14 \cite{Kiselev9} the matrix \
$S_{\tau}^{1} = \mathfrak{A} (\tau) $ \ reduced to the cardinal \
$\chi=(\chi^{\ast})^{+ \omega_0}$ \ on the carrier \
\mbox{$\alpha_{\tau}^{1}> \gamma_{\tau} $;} \ it has the
admissible generating eigendisseminator \ $
\check{\delta}_{\tau}^{1} < \gamma_{\tau} $ \ with the base \
$\rho_{\tau}^{1} \vartriangleright S_{\tau}^{1} $. \ One can see
that \ $\rho_{\tau}^{1} > \chi^{\ast+} $ \ and so \
$\check{\delta}_{\tau}^{1} $ \ can be considered as the
disseminator for \ $S_{\tau}^1$ \ on \ $\alpha_{\tau}^1$ \ with
the base \ $\chi^{*+}$.
\\
Now let us turn to the prejump cardinal
\[
    \alpha^{1} = \alpha_{\tau}^{1 \Downarrow};
\]
by the same lemma \ $cf(\alpha^{1}) \geq \chi^{\ast +} $ \ and it
is possible to introduce the matrix \ $S_{\tau}^{2} $ \
\textit{reduced to} \ $ \chi^{\ast} $ \ possessing the same
prejump cardinal \ $\alpha^{1} $ \ and hence the same disseminator
\ $\delta_{\tau}^{1} $ \ with the same base \ $ \chi^{\ast +} $ \
using lemma 5.12~\cite{Kiselev9} in the following way:
\\
If there holds the proposition
\begin{equation} \label{e8.6}
\exists \alpha \in [ \alpha^{1}, \alpha_{\tau}^{1}  [ \ \ \sigma
(\chi^{\ast}, \alpha)
\end{equation}
then let \ $S_{\tau}^{2} $ \ be the matrix reduced to \
$\chi^{\ast} $ \ on the minimal carrier \ $\alpha_{\tau}^{2}$ \
and produced by the cardinal \ $\alpha^{1}$, \ so that \ $\alpha^1
= \alpha_\tau ^{2 \Downarrow} $ \ (just as it was done in the
proof of lemma~6.12 \cite{Kiselev9} by the cardinal \ $\alpha_0$,
\ playing the role of \ $\alpha^1$ \ here).
\\
In the opposite case, when (\ref{e8.6}) fails, one can see that,
since the proposition of lemma~5.12 \cite{Kiselev9} is fulfilled
below \ $\alpha^{1} $, \ the matrix \ $S_{\tau}^{1} $ \ protects
the jump cardinal \ $\alpha_{\tau}^{\downarrow} $ \ (and, hence, \
$ \alpha^{1}$) \ which is preserves under the reducing the matrix
\ $S_{\tau}^{1} $ \ on the carrier \ $\alpha_{\tau}^{1} $ \ to \
$\chi^{\ast} $ ; \ so we can define the matrix (see definitions
4.1, 5.1, 5.5~\cite{Kiselev9})
\begin{equation*}
S_{\tau}^{2} \Rightarrow \widetilde{\mathbf{S}}_{n}^{\sin \triangleleft
\alpha_{\tau}^{1}} \overline{\overline{\lceil }}  \chi^{\ast} \mathrm{\quad
on\ the\ carrier \quad}  \alpha_{\tau}^{2} = \sup dom \bigl (
\widetilde{\mathbf{S}}_{n}^{\sin \triangleleft \alpha_{\tau}^{1}}
\overline{ \overline{\lceil }}  \chi^{\ast} \bigr ) .
\end{equation*}
This matrix \ $S_\tau^2$ \ is singular on the carrier \
$\alpha_\tau^2$: conditions 1), 3) of definition
5.7~\cite{Kiselev9} are obvious, while condition 2) one can verify
with the help of the \textit{splitting method}, repeating the
argument from the proof of lemma 5.12~\cite{Kiselev9} (where \
$\alpha_1$, \ $\chi$ \ are replaced with \ $\alpha_\tau^2$, \
$\chi^\ast$ \ respectively) literally.
\\
In any case \ $\alpha^{1} = \alpha_{\tau}^{2 \Downarrow} $ \ and \
$ S_{\tau}^{2} $ \ is found to be admissible on \
$\alpha_{\tau}^{2} $ \ for \ $\gamma_{\tau} $ \ along with the
same disseminator \ $\check{\delta}_{\tau}^{1} $ \ and its base \
$\chi^{\ast +} $, \ because all conditions of \ $\mathbf{K}^0$ \
from definition \ref{8.2.} trivially holds when \ $\rho=\chi^{\ast
+}$. \ Also such matrix \ $S_{\tau}^2$ \ on its carrier \
$\alpha_{\tau}^2$ \ is nonsuppressed due to its base \
$\rho=\chi^{\ast +}$. \ It can be unit or zero, but in any case
there exist some \ $\alpha$-matrix reduced to \ $\chi^{\ast}$ \
admissible and nonsuppressed for \ $\gamma_{\tau}$ \ along with
its accessories.
\\
Now one should take any cardinal \ $\gamma$ \ great enough and
such that for any \ $\gamma_{\tau} > \gamma$ \ there exist some
matrix \ $S_{\tau}^2$ \ with the base \ $\rho = \chi^{\ast +}$; \
it is admissible and nonsuppressed for \ $\gamma_{\tau}$ \ below \
$\alpha_1$ \ for any \ $\alpha_1 \in SIN_n$ \ by definition.
\\
So, after the minimizing such resulting matrices and their
accompanying ordinals according to definition~\ref{8.3.} there
appears the function \ $\alpha S_{f}^{<\alpha_1} $ \ and
accompanying ordinal functions defined on \ $T^{\alpha_1}$ \ for
any \ $\alpha_1 \in SIN_n$, $\alpha_1 > \gamma$.
\\
\hspace*{\fill} $ \dashv$
\\

\noindent In conclusion of this section repeating the proof mode
of lemma~\ref{7.7.} it is easy to draw out

\begin{lemma}
\label{8.10.}
\begin{equation*}
\alpha \delta ^{\ast } \in SIN_{n} \cap SIN_{n+1}^{<\alpha^{\ast
1}}[< \alpha \rho^{\ast 1}].
\end{equation*}
\end{lemma}
\noindent \textit{Proof.} \ Let us use the designations from the
previous lemma~\ref{8.9.}.
\\
First starts lemma~\ref{7.7.} proof mode, treating the
disseminator \ $\widetilde{\delta}_{\alpha \tau^{\ast 1}}$ \ with
the base \ $\alpha \rho^{\ast 1}$ \ of the matrix \ $\alpha
S_{\alpha \tau^{\ast 1}}$ \ on the carrier \ $\alpha_{\alpha
\tau^{\ast 1}}$ \ with the prejump cardinal \ $\alpha^{\ast 1} =
\alpha_{\alpha \tau^{\ast 1}}^{\Downarrow}$. \ Since
\[
    \alpha \delta^{\ast 1} \in SIN_n, \
    \widetilde{\delta}_{\alpha \tau^{\ast 1}} <
    \alpha \delta^{\ast 1}
\]
and
\[
    \widetilde{\delta}_{\alpha \tau^{\ast 1}} \in
    SIN_n^{<\alpha^{\ast 1}} \cap  SIN_{n+1}^{<\alpha^{\ast 1}}
    [< \alpha \rho^{\ast 1}],
\]
lemma~3.8 implies \ $\widetilde{\delta}_{\alpha \tau^{\ast 1}} \in
SIN_n$.

\noindent Now suppose, that  this lemma~\ref{8.10.} is wrong and
\[
    \alpha \delta^{\ast} \notin SIN_n,
\]
then
\[
    \widetilde{\delta}_{\alpha \tau^{\ast 1}} <
    \alpha \delta^{\ast} = \gamma_{\alpha \tau_1^{\ast}}.
\]
Due to lemma 3.2~\cite{Kiselev9} it makes possible to restrict the
\ \mbox{$\Sigma_{n-1}$-}pro\-po\-si\-tion about existence of
matrix \ $\alpha S_{\alpha \tau^{\ast 1}}$ \ admissible carrier
along with the same \ $\widetilde{\delta}_{\alpha \tau^{\ast 1}}$,
\ $\alpha \rho^{\ast 1}$ \ to the \ $SIN_{n-1}$-cardinal \
$\gamma_{\alpha \tau_1^{\ast}+1}$, \ just as it was done in the
proof of lemma~\ref{8.5.}~8).
\\
Then the matrix \ $\alpha S_{\alpha \tau^{\ast 1}}$ \ receives
again some its carrier
\[
    \alpha^{\prime} \in \; ] \gamma_{\alpha \tau^{\ast 1}},
    \gamma_{\alpha \tau^{\ast 1} + 1} [
\]
admissible for \ $\gamma_{\alpha \tau^{\ast 1}}$ \ along with its
previous disseminator and data base.
\\
But due to the minimality of \ $\alpha \tau_1^{\ast}$ \ there
holds
\[
    \alpha \tau_1^{\ast} \notin dom(\alpha S_f).
\]
It can be only when \ $\alpha S_{\tau^{\ast 1}}$ \ on \
$\alpha^{\prime}$ \ is admissible but suppressed for \
$\gamma_{\alpha \tau_1^{\ast}}$; \ in its turn it can be only when
\[
    \alpha \delta^{\ast} = \gamma_{\alpha \tau_1^{\ast}} \in SIN_n
\]
contrary to the supposition.

As to the rest part of this lemma:
\[
    \alpha \delta^{\ast} \in SIN_n^{<\alpha^{\ast 1}} [ < \alpha
    \rho^{\ast 1} ],
\]
it is not needed in what follows up to \S 11 and therefore we
shall return to it there.  \hspace*{\fill} $\dashv$

\newpage

\section{Analysis of \ $\protect\alpha $\,-Function Monotonicity}

\setcounter{equation}{0}

\hspace*{1em} Here the first component of the required
contradiction -- the monotonicity of \ $\alpha$-function -- is
investigated in various important cases.

As we shall see, this property is rather strong; in particular,
any interval \ $[ \tau_1, \tau_2 [$ \ of its monotonicity can not
be ``too long'', -- the corresponding interval \ $]
\gamma_{\tau_1}, \gamma_{\tau_2} [$ \ can not contain any \
\mbox{$SIN_n$-cardinals}, and if \ $\gamma_{\tau_2} \in SIN_n$, \
then this function receive some constant characteristics and
\textit{stabilizes} on such \ $[ \tau_1, \tau_2 [$.
\\
We start with the latter situation:

\begin{definition}
\label{9.1.} \hfill {}

The function \ $\alpha S_{f}^{<\alpha_{1}}$ \ is called monotone
on interval \ $\left[ \tau _{1},\tau _{2}\right[ $ \ and on the
corresponding interval \ $[ \gamma_{\tau_1}^{<\alpha_1},
\gamma_{\tau_2}^{<\alpha_1} [$ \ below \ $ \alpha _{1}$ \ iff \
$\tau _{1}+1<\tau _{2}$, \ $\left] \tau _{1},\tau _{2} \right[
\subseteq dom ( \alpha S_{f}^{<\alpha _{1}} ) $ \ and

\begin{equation*}
\forall \tau ^{\prime },\tau ^{\prime \prime }  ( \tau _{1} <
\tau^{\prime }<\tau ^{\prime \prime } <  \tau _{2}\longrightarrow
\alpha S_{\tau ^{\prime }}^{<\alpha _{1}} \underline{\lessdot
}\alpha S_{\tau ^{\prime \prime }}^{<\alpha _{1}} )~.
\end{equation*}
\vspace{-12pt} \hspace*{\fill} $\dashv$
\end{definition}

\noindent To operate with this notion it is suitable to use the
following \ \mbox{$\Delta_1$-formulas}, which play the main role
in this section:
\\

\noindent \quad $A_0^{1 \vartriangleleft \alpha_1}(\chi, \tau_1,
\tau_2, \alpha S_f^{<\alpha_1})$: \vspace{-6pt}
\begin{multline*}
    A_0^{\vartriangleleft \alpha_1}(\chi, \tau_1, \tau_2, \alpha S_f^{<\alpha_1}) \wedge
    \forall \tau^{\prime}, \tau^{\prime\prime} \big( \tau_1 < \tau^{\prime}
    < \tau^{\prime\prime} < \tau_2 \rightarrow
\\
    \rightarrow \alpha S_f^{<\alpha_1}(\tau^{\prime})
    \underline{\lessdot} \alpha
    S_f^{<\alpha_1}(\tau^{\prime\prime}) \big);
\end{multline*}
so, here is stated, that \ $\alpha S_f^{<\alpha_1}$ \ is defined
on the interval \ $] \tau_1,\tau_2 [$ \ with the property \ $A_0$
\ (remind definition \ref{8.1.}~1.0\;) and, moreover, it is
monotone on the interval \ $[ \tau_1,\tau_2[\;$; \ thus we shall
name it and the corresponding interval \ $[
\gamma_{\tau_1},\gamma_{\tau_2} [$ \ the intervals of the function
\ $\alpha S_f^{<\alpha_1}$ \ monotonicity;
\\

\noindent \quad $A_1^{1 \vartriangleleft \alpha_1}(\chi, \tau_1,
\tau_2, \alpha S_f^{<\alpha_1})$:
\\
\[
    \qquad \qquad A_0^{1 \vartriangleleft \alpha_1}(\chi, \tau_1, \tau_2, \alpha S_f^{<\alpha_1})
    \wedge \exists \gamma^2 \big(
    \gamma^2 = \gamma_{\tau_2} \wedge SIN_n^{<\alpha_1}(\gamma^2) \big);
\]
\quad %

\noindent further the function \ $\alpha S_f^{<\alpha_1}$ \ will
be omitted in such notations for brevity (if it will be pointed
out in the context).

Now it is not still all ready to prove the total monotonicity of
matrix function \ $\alpha S_f$, \ but some its fragments are clear
quite analogous to lemmas 5.17~1)~\cite{Kiselev9}, \ref{7.9.}. For
instance, from lemma 3.2~\cite{Kiselev9} it comes directly

\begin{lemma}
\label{9.2.} \emph{(About $\alpha $-function monotonicity)} \newline
\hspace*{1em} Let
\begin{equation*}
\tau _{1}<\tau _{2}, ~~ a_{\tau _{2}}^{<\alpha _{1}}= 1  ~~ and ~~
\widetilde{\delta }_{\tau _{2}}^{<\alpha _{1}} <  \gamma_{\tau _{1}},
\end{equation*}
\noindent then
\begin{equation*}
\alpha S_{\tau _{1}}^{<\alpha _{1}}\underline{\lessdot }  \ \ \alpha S_{\tau
_{2}}^{<\alpha_{1}}  \quad and \quad a_{\tau_{1}}^{< \alpha_{1}} = 1 .
\end{equation*}

\noindent Analogously for zero characteristic \ $a_{\tau
_{1}}^{<\alpha _{1}} = a_{\tau _{2}}^{<\alpha _{1}} = 0 $.
\\
\hspace*{\fill} $\dashv$
\end{lemma}

\begin{lemma}
\label{9.3.} \emph{(About \ $\alpha $-function stabilization)}
\newline \hspace*{1em} Let
\\
(i) \ $\alpha S_{f}^{<\alpha _{1}}$ \ be monotone on \ $[ \tau
_{1},\tau _{2} [$ \ below \ $\alpha_1$:
\[
    A_1^{1 \vartriangleleft \alpha_1}(\tau_1, \tau_2);
\]
(ii) \ $\gamma _{\tau _{2}}^{<\alpha_{1}}$ \ be a successor in \
$SIN_{n}^{<\alpha _{1}}$. \newline Then \ $\alpha S_{f}^{<\alpha
_{1}}$ \ stabilizes on \ $[ \tau _{1},\tau _{2} [ $, \ that is
there exist  \ $ S^{0}$ \ and \ $\tau _{0}\in \; ] \tau_{1},\tau
_{2} [ \;$ \ such that

\vspace{6pt}
\begin{equation*}
\forall \tau \in \left[ \tau _{0},\tau _{2}\right[ \quad \alpha
S_{\tau }^{<\alpha_{1}}=S^{0}.
\end{equation*}
\vspace{0pt}

\noindent The least of such ordinals \ $\tau _{0}$ \ is called the
{\sl stabilization ordinal} of \ $\alpha S_{f}^{<\alpha _{1}}$ \
for \ $\tau_2$ \ below \ $\alpha_1$ \ and denoted  through \
$\tau_2^{s<\alpha_1} $.
\end{lemma}

\noindent \textit{Proof} represents once more the typical
application of lemma 3.2~\cite{Kiselev9}; we shall omit the upper
indices \ $<\alpha _{1}$, $\vartriangleleft \alpha _{1}$. \
Suppose this lemma is wrong; let us consider the ordinal
\[
    \rho ^{0}= \sup \{ Od ( \alpha S_{\tau } ) :  \tau _{1}<\tau
    <\tau_{2}\}.
\]
Let us apply the mode of reasoning used above in the proof of
lemma \ref{8.5.}~8) and introduce the cardinals
\[
    \gamma_{\tau_2^n} = \sup \{ \gamma < \gamma_{\tau_2}: \gamma
    \in SIN_n \};
\]
\[
    \gamma_{\tau_{1,2}^n} = \max \{ \gamma_{\tau_1}, \gamma_{\tau_2^n}
    \}.
\]

\noindent Then one should repeat definition~\ref{8.3.} of the
matrix function \ $\alpha S_f$ \ and its accompanying ordinal
functions below \ $\alpha_1$ \ on the set

\begin{equation} \label{e9.1n}
    T_{\tau_{1,2}^n}^{\alpha_1} = \{ \tau: \gamma_{\tau_{1,2}^n}
    < \gamma_{\tau} < \alpha_1 \}
\end{equation}
but preserving \textit{only $SIN_n$-cardinals} \ $\le
\gamma_{\tau_2^n}$; \ it can be done in the following way:
\\
Definition~\ref{8.3.} is based on the formula
\begin{equation} \label{e9.2n}
\alpha \mathbf{K}^{\ast}(a, \delta, \gamma_{\tau}, \alpha, \rho,
S)
\end{equation}
below \ $\alpha_1$ \ (see definition \ref{8.2.}~4)\;), which
means, that \ $S$ \ is the \ $\alpha$-matrix on its carrier \
$\alpha$ \ of characteristic \ $a$ \ with the disseminator \
$\delta$ \ and base \ $\rho$ \ admissible for \ $\gamma_{\tau}$ \
and, moreover, nonsuppressed on this \ $\alpha$ \ for \
$\gamma_{\tau}$ \ below \ $\alpha_1$; \ but since for every \
$\tau \in T_{\tau_{1.2}^n}^{<\alpha_1}$ \ there holds \
$\gamma_{\tau} \notin SIN_n$, \ the nonsuppression condition \
$\neg A_5^{S,0}$ \ in \ $\alpha \mathbf{K}^{\ast}$ \ holds on and
it can be dropped, and then \ $\alpha \mathbf{K}^{\ast}$ \ turns
into the formula $\alpha \mathbf{K}$.
\\
This formula is from the class \ $\Sigma_n$, \ because it includes
the \ $\Sigma_n$-formula \ $\mathbf{K}_n^{\forall}$. \ But let us
use the cardinal \ $\gamma_{\tau_2^n}$ \ and replace in
formula~(\ref{e9.2n}) its subformula \ $\mathbf{K}_n^{\forall}$ \
with \ $\Delta_1$-formula
\[
    SIN_n^{<\alpha^{\Downarrow}}(\gamma_{\tau_2^n}),
\]
then the \ $\Sigma_n$-formula (\ref{e9.2n}) turns into some \
$\Pi_{n-2}$-formula, which we shall denote through
\[
    \alpha \mathbf{K}_{n-2}^{\ast 1}(a, \delta, \gamma_{\tau},
    \alpha, \rho, S).
\]
So, the matrix function defined on the set \
$T_{\tau_{1,2}^n}^{\alpha_1}$ \ (\ref{e9.1n}) as in
definition~\ref{8.3.}, but through the formula (\ref{e9.2n})
replaced with \ $\alpha\mathbf{K}_{n-2}^{\ast 1}$, \ evidently
coincides with the function \ $\alpha S_f$ \ on the interval \
$]\tau_{1,2}^n, \tau_2 [\;$; \ we shall denote it by \ $\alpha
S_f^1$.
\\
Now, since \ $\alpha S_f^1$ \ is monotone on \
$]\tau_{1,2}^n,~\tau_2[$ \ but is not stabilized on this interval,
the ordinal \ $\rho_0$ \ is \textit{limit} and there holds the
following proposition below \ $\gamma_{\tau_2}$:

\vspace{6pt}
\begin{equation*}
\forall \tau ~ ( \tau_{1,2}^n < \tau \longrightarrow \exists S~ (
S=\alpha S_{\tau }^{1}\wedge S\vartriangleleft  \rho ^{0} ) )~.
\end{equation*}
\vspace{0pt}

\noindent It can be formulated in the \ $\Pi_n$-form:

\begin{equation*}
\forall \tau, \gamma^{\prime}, \gamma^{\prime\prime} \Bigl[
\gamma_{\tau_{1,2}^n} < \gamma^{\prime} = \gamma_{\tau} <
\gamma^{\prime\prime} = \gamma_{\tau+1} \rightarrow
\qquad\qquad\qquad\qquad\qquad\qquad
\end{equation*}
\begin{equation} \label{e9.3n}
~ \rightarrow \Big( \exists \delta, \alpha, \rho <
\gamma^{\prime\prime} \exists S \vartriangleleft \rho ~ \big(
\alpha \mathbf{K}_{n-2}^{\ast 1} (1, \delta, \gamma^{\prime},
\alpha, \rho, S) \wedge S \vartriangleleft \rho^0 \big) \vee
\end{equation}
\begin{equation*}
~ \vee \big( \exists \delta, \alpha, \rho < \gamma^{\prime\prime}
\exists S \vartriangleleft \rho ~ \alpha \mathbf{K}_{n-2}^{\ast 1}
(0, \delta, \gamma^{\prime}, \alpha, \rho,  S) \wedge S
\vartriangleleft \rho^0 \wedge \quad \quad
\end{equation*}
\begin{equation*}
\quad \wedge \forall \delta^{\prime}, \alpha^{\prime},
\rho^{\prime} < \gamma^{\prime\prime} \forall
S^{\prime}\vartriangleleft \rho^{\prime} ~ \neg \alpha
\mathbf{K}_{n-2}^{\ast 1} (1, \delta^{\prime}, \gamma^{\prime},
\alpha^{\prime}, \rho^{\prime}, S^{\prime}) \big) \Big) \Bigr ] ~.
\end{equation*}
\vspace{0pt}

\noindent Now it comes the contradiction:
\\
On one hand, the \ $SIN_n$-cardinal \ $\gamma _{\tau _{2}}$ \
extends this proposition (\ref{e9.3n}) up to \ $\alpha _{1}$ \ and
as a result the matrix \ $ \alpha S_{\tau
_{2}}^{1}\vartriangleleft \rho ^{0}$ \ arises.
\\
But, on the other hand, \ $\rho ^{0} $ \ is the limit ordinal and there
exists
\[
    \tau _{1,3}^{n} \in ] \tau _{1,2}^{n}, \tau _{2} [
    \mathrm{~such~that~}
    \alpha S_{\tau _{1,3}^{n}}^{1} \gtrdot \alpha S_{\tau _{2}}^{1}.
\]
That is why below \ $\gamma _{\tau _{2}}$ \ the next proposition
holds:

\vspace{6pt}
\begin{equation*}
\forall \tau ~ \big( \tau _{1,3}^{n }<\tau \longrightarrow \forall
S ( S=\alpha S_{\tau }^{1}\longrightarrow S \gtrdot \alpha
S_{\tau_{2}}^{1} ) \big)~.
\end{equation*}
\vspace{0pt}

\noindent It also can be formulated in \ $\Pi_{n}$-form:

\begin{equation*}
\forall \tau, \gamma^{\prime}, \gamma^{\prime\prime} \Bigl[
\gamma_{\tau_{1,3}^n} < \gamma^{\prime} = \gamma_{\tau} <
\gamma^{\prime\prime} = \gamma_{\tau+1} \rightarrow \quad \quad
\quad \quad \quad \quad \quad \quad \quad \quad \quad \quad \quad
\end{equation*}
\begin{equation*}
\rightarrow \Big( \forall \delta, \alpha, \rho <
\gamma^{\prime\prime} ~ \forall S \vartriangleleft \rho \big(
\alpha \mathbf{K}_{n-2}^{\ast 1} (1, \delta, \gamma^{\prime},
\alpha, \rho, S) \rightarrow \alpha S_{\tau_2}^{1} \lessdot S
\big) \wedge
\end{equation*}
\begin{equation} \label{e9.4}
\quad \wedge \forall \delta, \alpha, \rho < \gamma^{\prime\prime}
~ \forall S \vartriangleleft \rho  \big( \alpha
\mathbf{K}_{n-2}^{\ast 1} (0, \delta, \gamma^{\prime}, \alpha,
\rho, S) \wedge
\end{equation}
\begin{equation*}
\wedge \forall \delta^{\prime}, \alpha^{\prime},  \rho^{\prime} <
\gamma^{\prime\prime} ~ \forall S^{\prime} \vartriangleleft
\rho^{\prime}  \neg \alpha \mathbf{K}_{n-2}^{\ast 1} (1,
\delta^{\prime}, \gamma^{\prime}, \alpha^{\prime}, \rho^{\prime},
S^{\prime}) \rightarrow
\end{equation*}
\begin{equation*}
\qquad \qquad \qquad \qquad \qquad \qquad \qquad \qquad \qquad
\qquad \rightarrow \alpha S_{\tau_2}^{1} \lessdot S \big) \Big)
\Bigr ] ,
\end{equation*}

\vspace{6pt}

\noindent which \ $\gamma _{\tau _{2}}$ \ extends up to \ $\alpha
_{1} $ \ and therefore
\[
    \alpha S_{\tau _{2}}^{1} \lessdot \alpha S^1_{\tau _{2}}.
\]
\hspace*{\fill} $\dashv$

\quad \\
Let us remind that the symbols \ $\mathfrak{n}^{\alpha}$,
\ $\chi^{\ast}$, \ $ \alpha S_{f}^{<\alpha _{1}}$, \
$a_{f}^{<\alpha _{1}}$ \ in writings of formulas can be often
omitted for shortening. Besides, the usual condition of
equinformativeness \ $A_6^e(\alpha_1)$
\begin{multline*}
    \chi^{\ast} < \alpha_1 \wedge A_n^{<\alpha_1} (\chi^\ast) =  \|
    u_n^{\vartriangleleft \alpha_1} ( \underline{l}) \| \wedge SIN_{n-2}
    (\alpha_1) \wedge
\\
    \wedge \forall \gamma < \alpha_1 \exists \gamma_1 \in [\gamma,
    \alpha_1[ ~ SIN_n^{<\alpha_1}(\gamma_1)
\end{multline*}
\vspace{0pt}

\noindent is always superimposed on the bounding cardinals \
$\alpha_1$.
\newline

The stabilization property is very important for what follows; moreover, it
comes out that the analogous attribute arises for the characteristic
function, which play the crucial role further.
\\
Complicating in a certain way the reasoning from the proof of
lemma~\ref{9.3.} it is possible to prove the similar
characteristic property:

\begin{lemma}
\label{9.4.} \emph{(About characteristic stabilization)} \newline
\hspace*{1em} Let \medskip

(i) \quad $A_{1}^{1 \vartriangleleft \alpha_1} ( \tau _{1},
\tau_{2} ) ;$
\medskip

(ii) \quad $\forall \tau <\tau _{2}\quad \exists \tau^{\prime} \in \left[
\tau ,\tau_{2}\right[ \quad a_{\tau ^{\prime }}^{<\alpha _{1}}=1$; \medskip
\newline
Then
\begin{equation*}
\forall \tau ^{\prime }\in \left] \tau _{1},  \tau_{2}\right[
\quad a_{\tau^{\prime} }^{<\alpha _{1}}=1 .
\end{equation*}
In this case we shall say that the unit characteristic stabilizes
on \ $[ \tau _{1},\tau _{2} [ $ \ below \ $\alpha _{1}$.
\\
Analogously for zero characteristic.
\end{lemma}

\noindent \textit{Proof} is carried out by the induction on the
pair \ $( \alpha_{1}, \tau_{2} )$; \ (remind, the set of such
pairs is considered to be canonically ordered as above, with \
$\alpha_1$ \ as the first component and \ $\tau_2$ \ as the
second). Suppose this pair is minimal violating the lemma. It is
not hard to see that \ $\gamma_{\tau_{2}}^{< \alpha_{1}} $ \ is
the successor in \ $SIN_{n}^{< \alpha_{1}} $; \ precisely this
case is used further. Recall, that for matrices of unit
characteristic on their carriers the suppression condition \
$A_5^{S,0}$ \ fails and it can be dropped for these matrices; so
for the unit characteristic \ $a=1$ \ the formula \ $\alpha
\mathbf{K}^{\ast <\alpha_1}$ \ is equivalent to $\alpha
\mathbf{K}^{<\alpha_1}$. \ The upper indices \ $< \alpha_1$,
$\vartriangleleft \alpha_1$ \ will be dropped for some brevity.
\\
By the previous lemma there exist the stabilization ordinal \
$\tau_{2}^s $ \ of \ $\alpha S_{f} $ \ on \ $[ \tau_{1}, \tau_{2}
[ $, \  and the matrix \ $S^0$ \ such that

\begin{equation*}
\forall \tau \in [ \tau_{2}^s, \tau_{2}[ \quad  \alpha S_{\tau} =
S^{0} .
\end{equation*}

\noindent According \ to the \ condition \ of \ this \ lemma \
there \ exists \ the minimal \ $\tau^{1} \in [ \tau_{2}^s,
\tau_{2}[ $ \ such that \ $ a_{\tau^{1}} = 1 $. \ The further
reasoning splits into two parts:
\\
1. \ First, let us prove that \ $\forall \tau \in [ \tau^{1},
\tau_{2}[ \ \ a_{\tau} = 1 $. \newline Suppose it is wrong, then
there exists the minimal \ $\tau^{0} \in ] \tau^{1}, \tau_{2}[ $ \
providing \ $a_{\tau^{0}} = 0 $; \ thus \ $a_{\tau} \equiv 1 $ \
on \ $[ \tau^{1}, \tau^{0}[ $ \ (so, take notice, the matrix \
$S^0$ \ on different admissible carriers can possess different
characteristic here). Let us consider the following subcases for
\begin{equation*}
S^0 = \alpha S_{\tau^0}, \quad \check{\delta}^{0} =
\check{\delta}_{\tau^{0}} , \quad \alpha^{0} =
\alpha_{\tau^{0}}^{\Downarrow} ~:
\end{equation*}
\newline
1a. \quad $\check{\delta}^{0} \notin SIN_{n} $. \ Since \
$\gamma_{\tau_{1}} \in SIN_{n} $, \ lemma~3.8~1) \cite{Kiselev9}
implies
\begin{equation*}
\gamma_{\tau_1} < \check{\delta}^{0}, \quad \check{\delta}^{0} \in
\left ( SIN_{n}^{< \alpha^{0}} - SIN_{n} \right )
\end{equation*}
and we can use the cardinal
\begin{equation*}
\gamma_{\tau^{2}} = min \left (SIN_{n}^{< \alpha^{0}} -  SIN_{n} \right ) .
\end{equation*}
Due to the same lemma 3.8 \cite{Kiselev9} it is not hard to see
that \ $\gamma_{\tau^{2}} $ \ is the \textit{successor} in \ $
SIN_{n}^{< \alpha^{0}} $ \ of some cardinal
\begin{equation*}
\gamma_{\tau^{3}} \geq \gamma_{\tau_{1}}, \quad  \gamma_{\tau^{3}}
\in SIN_{n} \mathrm{\ \ below\ \ } \alpha_1,
\end{equation*}
and the function \ $\alpha S_{f}^{< \alpha^{0}} $ \ is monotone on
the interval \ $[ \tau^3, \tau^2 [\;$. \ Since \ $a_{\tau} \equiv
1 $ \ on \ $[ \tau^{1}, \tau^{0} [ $ \ lemma 3.2~\cite{Kiselev9}
provides that the interval \ $ ]\gamma_{\tau_{1}},
\gamma_{\tau^{2}} [ $ \ contains admissible carriers of matrices
of unit characteristic disposed cofinally to \ $\gamma_{\tau^{2}}
$ \ because the \ $SIN_n^{<\alpha^0}$-cardinal \ $\gamma_{\tau^2}$
\ restricts the \ $\Sigma_n$-proposition about the existence of
such carriers. After that the cardinal \ $\gamma_{\tau^{2}} $ \
\textit{extends} unit characteristic up to \ $\alpha_0$, \ and,
so, \ $S^0$ \ on \ $\alpha_{\tau_0}$ \ becomes unit matrix
contrary to the supposition.

\noindent This argument mode consists in restrictions and
extensions applied in turns and therefore we shall call it the
\textit{restriction-and-extension method}. It will be used further
often enough in various forms and is typical in disseminator
theory, therefore one should consider it more in details:
\\
Let \ $\gamma < \gamma_{\tau_2}$ \ be an arbitrary cardinal; there
exist the unit matrix \ $S^0$ \ on some carrier \ $\alpha >
\gamma$ \ and it remains unit below \ $\alpha^0$ \ due to
lemma~\ref{8.7.} about absoluteness. Now the reasoning passes to
the situation below \ $\alpha^0$; \ there holds the proposition
below \ $\alpha^0$:
\[
    \exists \delta, \gamma_{\tau}, \alpha, \rho ~ \bigl( \gamma <
    \gamma_{\tau} \wedge \alpha \mathbf{K}(1, \delta, \gamma_{\tau}, \alpha,
    \rho, S^0) \bigr),
\]
it belongs to \ $\Sigma_n$ \ and contains only constants
\[
    \chi^{\ast}, ~ \gamma < \gamma_{\tau^2}, ~ S^0 \vartriangleleft \chi^{\ast
    +} < \gamma_{\tau^2}.
\]
Thus the \ $\Pi_n$-cardinal \ $\gamma_{\tau_2}$ \ below \
$\alpha_0$ \ restricts it by lemma 3.2~\cite{Kiselev9}, that is it
holds after its bounding by $\gamma_{\tau^2}$:
\[
    \exists \delta, \gamma_{\tau}^{<\gamma_{\tau^2}}, \alpha,
    \rho < \gamma_{\tau^2} ~ \bigl( \gamma <
    \gamma_{\tau}^{<\gamma_{\tau^2}} \wedge \alpha
    \mathbf{K}^{\vartriangleleft
    \gamma_{\tau^2}}(1, \delta, \gamma_{\tau}^{<\gamma_{\tau^2}},
    \alpha, \rho, S^0) \bigr).
\]
But here the upper indices \ $<\gamma_{\tau^2}$, \
$\vartriangleleft \gamma_{\tau^2}$ \ can be dropped due to the \
$\Pi_n^{<\alpha^0}$-subinaccessibility of \ $\gamma_{\tau^2}$ \
and as a result there appear admissible carriers
\[
    \alpha \in \; ] \gamma, \gamma_{\tau_2} [
\]
of the matrix \ $S^0$ \ of \textit{unit} characteristic on such \
$\alpha$ \ below \ $\alpha^0$ \ for arbitrary \ $\gamma <
\gamma_{\tau_2}$.

\noindent Then by the inductive hypothesis \ $a_{\tau} \equiv 1 $
\ on \ $] \tau_{1}, \tau^{2} [ $, \ and below \ $\gamma_{\tau^{2}}
$ \ there is fulfilled the proposition
\begin{equation*}
\forall \tau \ \ (\gamma_{\tau^{3}} <  \gamma_{\tau}\longrightarrow a_{\tau}
= 1 )
\end{equation*}
that can be formulated in the \ $\Pi_{n} $-form:
\begin{multline*}
\forall \gamma \Bigl( \gamma_{\tau^{3}} < \gamma  \wedge
SIN_{n-1}(\gamma) \longrightarrow
\\
\longrightarrow \exists \delta, \alpha, \rho, S \Bigl(
SIN_{n}^{<\alpha^{\Downarrow}} (\gamma_{\tau^{3}}) \wedge \alpha
\mathbf{K}_{n+1}^{\ast \exists} (1, \delta, \gamma, \alpha,  \rho,
S) \Bigr) \Bigr) \ .
\end{multline*}
\vspace{0pt}

\noindent The cardinal \ $\gamma_{\tau^{2}} \in SIN_{n}^{<
\alpha^{0}} $ \ extends this last proposition up to \ $\alpha^{0}
$, \ and below \ $ \alpha^{0} $ \ there appears the matrix of unit
characteristic on some carrier \ $\in \; ]\gamma_{\tau^{0}},
\alpha^{0} [ $ \ admissible together with its disseminator \
$<\gamma_{\tau_{0}} $ \ and its base for \ $ \gamma_{\tau_{0}} $.
\\
Thus, \ $a_{\tau_{0}}=1 $ \ contrary to the assumption, and we
turn to the next subcase:
\\

\noindent 1b. \quad $\check{\delta}^{0} \in SIN_{n} $. \ Since \
$a_{\tau} \equiv 1 $ \ on \ $[ \tau^{1}, \tau^{0} [ $, \ there
exist the matrix
\[
    S^0 = \alpha S_{\tau^{1,0}} \mathrm{\ on\ the\ carrier\ }
    \alpha_{\tau^{1,0}} \in [ \check{\delta}^0, \gamma_{\tau^0}[
\]
of \textit{unit characteristic} \ $a_{\tau^{1,0}}=1$ \ and one can
reveal the situation below \ $\alpha^{1,0} =
\alpha_{\tau^{1,0}}^{\Downarrow}$ \ in the following way.
\\
The reasoning forthcoming is applied further subsequently, thus it
is necessary to dwell upon it.
\\
We start with \ $S^0$ \ on \ $\alpha_{\tau_0}$. \ By
lemma~\ref{8.5.}~5) zero characteristic of \ $S^0$ \ on \
$\alpha_{\tau^0}$ \ means that there holds

\begin{equation}  \label{e9.5}
    \qquad
    \exists \tau_1^{\prime}, \tau_2^{\prime}, \tau_3^{\prime} <
    \alpha^0 ~ \big( A_2^{0 \vartriangleleft \alpha^0} (\tau_1^{\prime},
    \tau_2^{\prime}, \tau_3^{\prime}, \alpha S_f^{<\alpha^0} )
    \wedge \qquad \qquad \qquad
\end{equation}
\[
    \qquad \qquad \qquad
    \wedge \forall \tau^{\prime\prime} \in \; ]\tau_1^{\prime},
    \tau_2^{\prime}] ~ a_{\tau^{\prime\prime}}^{<\alpha^0} = 1
    \wedge \alpha S_{\tau_2^{\prime}}^{<\alpha^0} = S^0 \big).
\]
\vspace{0pt}

\noindent Thus there can be used some ordinals
\[
    \tau_1^{\prime} < \tau_2^{\prime} < \tau_3^{\prime} <
    \alpha^0
\]
such that there holds

\begin{equation}  \label{e9.6}
    A_2^0( \tau_{1}^{\prime}, \tau_{2}^{\prime}, \tau_{3}^{\prime},
    \alpha S_{f})
    \wedge \forall \tau^{\prime\prime} \in \; ]\tau_1^{\prime},
    \tau_2^{\prime}] ~ a_{\tau^{\prime\prime}}^{<\alpha^0} = 1
    \wedge \alpha S_{\tau_{2}^{\prime}} = S^{0}
\end{equation}

\noindent below\ $\alpha^0$, \ that is after \
$\vartriangleleft$-bounding by the cardinal \ $\alpha^0$.

\noindent The key role will be played here by the so called
\textit{mediator}: it is some \ $SIN_n^{<\alpha^0}$-cardinal \
$\gamma^0$ \ such that

\begin{equation}  \label{e9.7}
    \gamma_{\tau_1^{\prime}}^{<\alpha^0} <
    \gamma_{\tau_2^{\prime}}^{<\alpha^0} <
    \gamma_{\tau_3^{\prime}}^{<\alpha^0} <
    \gamma^0 < \alpha^0
\end{equation}

\noindent which exist due to lemma~\ref{8.5.}~4). By
lemma~\ref{8.7.} about absoluteness of the admissibility and of
unit values of the matrix \ $\alpha$-function these values and
their accessories below \ $\alpha^0$ \ and below \ $\gamma^0$ \
coincide on the set
\[
    \bigl \{ \tau: \gamma_{\tau}^{<\alpha^0} < \gamma^0
    \wedge a_{\tau}^{<\alpha^0} = 1
    \bigr \}
\]
and therefore (\ref{e9.7}), (\ref{e9.6}) imply the following \
$\Sigma_{n+1}$-proposition below \ $\alpha^0$:
\[
    \exists \gamma^0 \exists \tau_1^{\prime}, \tau_2^{\prime},
    \tau_3^{\prime} < \gamma^0 \big( SIN_n(\gamma^0) \wedge
    \gamma_{\tau_1^{\prime}}^{<\gamma^0} <
    \gamma_{\tau_2^{\prime}}^{<\gamma^0} <
    \gamma_{\tau_3^{\prime}}^{<\gamma^0} < \gamma^0 \wedge
    \qquad
\]
\begin{equation} \label{e9.8}
    \wedge  A_2^{0 \vartriangleleft \gamma^0} (\tau_1^{\prime},
    \tau_2^{\prime}, \tau_3^{\prime}, \alpha S_f^{<\gamma^0})
    \wedge \forall \tau^{\prime\prime} \in \; ]\tau_1^{\prime},
    \tau_2^{\prime}] ~ a_{\tau^{\prime\prime}}^{<\gamma^0} = 1
    \wedge
    \quad
\end{equation}
\[
    \qquad \qquad \qquad \qquad \qquad \qquad \qquad \qquad
    \wedge \alpha S_{\tau_2^{\prime}}^{<\gamma^0} = S^0 \big).
\]

\noindent Due to lemma~\ref{8.5.}~3)
\[
    \check{\delta}^0 \in SIN_{n+1}^{<\alpha^0}[< \rho_{\tau^0}]
\]
and then by lemma 3.2~\cite{Kiselev9} there exist some \
$\gamma^0$ \ with the property (\ref{e9.8}), but \textit{already
below} \ $\check{\delta}^0$. \ From here and lemma
3.8~\cite{Kiselev9} it follows that the \
$SIN_n$-subinaccessibility of \ $\check{\delta}^0$ \ draws the
existence of \ $SIN_n$-cardinal \ $\gamma^{0 \prime} <
\check{\delta}^0$ \ with the same property (\ref{e9.8}); note,
that \ $\gamma^{0 \prime}$ \ possesses the same \
$SIN_n$-subinaccessibility as \ $\check{\delta}^0$. \ So, for some
cardinals
\begin{equation} \label{e9.9}
    \gamma_{\tau_1^{\prime\prime}}^{<\gamma^{0 \prime}} <
    \gamma_{\tau_2^{\prime\prime}}^{<\gamma^{0 \prime}} <
    \gamma_{\tau_3^{\prime\prime}}^{<\gamma^{0 \prime}} <
    \gamma^{0 \prime}
\end{equation}
there holds
\begin{equation} \label{e9.10}
    \qquad \qquad
    A_2^{0 \vartriangleleft \gamma^{0 \prime}} (\tau_1^{\prime\prime},
    \tau_2^{\prime\prime}, \tau_3^{\prime\prime}, \alpha S_f^{<\gamma^{0 \prime}})
    \wedge
    \qquad\qquad\qquad\qquad\qquad\qquad\qquad
\end{equation}
\[
    \qquad\qquad
    \wedge \forall \tau^{\prime\prime\prime} \in \; ]\tau_1^{\prime\prime},
    \tau_2^{\prime\prime}] ~ a_{\tau^{\prime\prime\prime}}^{<\gamma^{0\prime}} = 1
    \wedge \alpha S_{\tau_2^{\prime\prime}}^{<\gamma^{0 \prime}} = S^0
\]

\noindent Since \ $\gamma^{0 \prime}$ \ is the \ $SIN_n$-cardinal
everywhere in (\ref{e9.9}), (\ref{e9.10}) \ $\vartriangleleft$-
and \ $<$-boundaries by \ $\gamma^{0 \prime}$ can be dropped by
the same lemmas 3.8, \ref{8.7.}.

From this place one have to repeat the reasoning above but in the
reserve direction, and not for \ $\alpha^0$, \ \textit{but for} \
$\alpha^{1,0}$. \ Then (\ref{e9.9}), (\ref{e9.10}) draw
(\ref{e9.5}) where \ $\alpha^0$ \ is replaced with \
$\alpha^{1,0}$ \ and thereby \ $S^0$ \ on \ $\alpha_{\tau^{1,0}}$
\ receives \textit{zero characteristic} contrary to the
assumptions.

\noindent 2. So, statement 1 is proved; it remains to examine the
ordinal
\[
    \tau^{1,2} = \min \bigl \{\tau \in [ \tau_1, \tau_2 [ : \forall
    \tau^{\prime} \in ] \tau, \tau_2 [ ~ a_{\tau^{\prime}} = 1 \bigr \}
\]
and to prove that it coincides with \ $\tau_1$.
\\
Suppose it is wrong and \ $\tau_1 < \tau^{1,2}$, \ then one should
examine two unit matrices
\[
    S^1 = \alpha S_{\tau^{1,2}}, ~ S^2 = \alpha S_{\tau^{1,2}+1}
\]
and treat the matrix \ $S^2$ \ on its carrier \
$\alpha_{\tau^{1,2}+1}$ \ with its generating disseminator \
$\check{\delta}^2 = \check{\delta}_{\tau^{1,2}+1} $. \ By
lemmas~\ref{8.5.}~7)~$(ii)$, \ \ref{8.8.}~2)
\[
    \gamma_{\tau_1} \le \check{\delta}^2 = \widetilde{\delta}_{\tau^{1,2}+1}
\]
and cause of that only three subcases arises:

\noindent 2a. \quad $\gamma_{\tau_1} = \check{\delta}^2$. \ Then
by definition~\ref{8.3.}
\[
    \forall \tau \in ] \tau_1, \tau^{1,2} [ ~ a_{\tau} = 1
\]
contrary to the supposition.

\noindent 2b. \quad $\gamma_{\tau_1} < \check{\delta}^2$, \
$\check{\delta}^2 \notin SIN_n$. \ Then the
restriction-and-extension reasoning works, literally as it was in
subcase 1a. of this proof above, but for
\[
    S^2, ~ \check{\delta}^2 ~ \mbox{ instead of }
    ~ \alpha S_{\tau_0}, ~ \check{\delta}^0
\]
and again it comes \ $a_{\tau} \equiv 1$ \ on \ $] \tau_1,
\tau^{1,2} ]$.

\noindent 2c. \quad $\gamma_{\tau_1} < \check{\delta}^2$, \
$\check{\delta}^2 \in SIN_n$. \ Here again the
restriction-and-extension method works, but in slightly another
manner. First by lemma 3.2~\cite{Kiselev9} matrix \ $S^1$ \
receives its admissible carriers of unit characteristic disposed
cofinally to \ $\check{\delta}^2$, \ so by the inductive
hypothesis
\[
    a_{\tau} \equiv 1 \mathrm{\ on\ the\ set\ }
    \{ \tau : \gamma_{\tau_{1}} < \gamma_{\tau} <
    \check{\delta}^2 \}.
\]
Then below \ $\check{\delta}^2$ \ the following \ $ \Pi_{n+1}
$-proposition holds

\begin{equation*}
\forall \gamma \big( \gamma_{\tau_1} < \gamma  \wedge
SIN_{n-1}(\gamma) ~ \rightarrow ~ \exists \delta,  \alpha, \rho, S
\quad \alpha \mathbf{K} (1, \delta,  \alpha, \gamma, \rho, S)
\big)
\end{equation*}

\noindent which is extended by this disseminator up to \ $
\alpha_{\tau^{1,2}+1}^{\Downarrow} $ \ according to lemma
6.6~\cite{Kiselev9} (for \ $ m=n+1 $ , $\delta =
\check{\delta}^2$, $\alpha_{1} =
\alpha_{\tau^{1,2}+1}^{\Downarrow}$) \ and again it comes \ $
a_{\tau} \equiv 1 $ \ on the same set \ $]\tau_1, \tau^{1,2}]$.
\\
In every case it implies \ $\tau_1 = \tau^{1,2}$. \hspace*{\fill}
$\dashv$
\\

The following important lemma will be proved again by means of the
\textit{restric\-tion-and-extension method} but in some
synthesized form.

However, beforehand the following rather suitable notion should be
introduced using the notions of reduced spectra and matrices
(remind definitions~4.1, 5.1~\cite{Kiselev9}).
\\
In what follows the main technical mode of reasonings will consist
in the examination of some matrix \ $S$ \ under consideration on
its \textit{different carriers in turns}. Such transition of the
reduced matrix \ $S$ \ from one its carrier \ $\alpha$ \ over to
another its carrier \ $\alpha^1$ \ will be called the
\textit{carrying over} of the matrix \ $S$ \ from \ $\alpha$ \ to
\ $\alpha^1$.
\\
This technique will be frequent enough to be used further and was
already used above in the proofs of lemmas \ref{7.5.}, \ref{8.8.},
\ref{9.4.}.
\\
During such carrying over of reduced matrix \ $S$ \ from \
$\alpha$ \ to \ $\alpha^1$ \ some properties of the universe
bounded by jump or prejump cardinals of \ $S$ \ on \ $\alpha$ \
can be preserved and thereby they will be called the \textit{inner
properties} of \ $S$; \ other properties of \ $S$ \ may be lost
and they will be called the \textit{outer properties}.

More precisely: a property or attribute of matrix \ $S$ \ reduced
to \ $\chi^{\ast}$ \ on its carrier \ $\alpha$ \ will be called
the \textit{inner} property or attribute of this \ $S$ \ (on \
$\alpha$) \ if it is definable below some jump or prejump cardinal
of the spectrum
\[
    dom  \bigl( \widetilde{\mathbf{S}}_{n}^{\sin \vartriangleleft
    \alpha}\overline{\overline{\lceil}}\chi^{\ast} \bigr )
\]
through its some other jump of prejump cardinals; analogously for
other objects from \ $L_{\alpha}$; \ in all other cases they will
be called the \textit{outer properties or attributes or objects}
of \ $S$.
\\
These notions are activated by lemma~5.11~\cite{Kiselev9} about
matrix informativeness which means that such \textit{inner}
properties are preserved while matrix \ $S$ \ is carried over from
one its carrier to any other one.

Here is very important example of the \textit{outer} property --
the property of \textit{characteristic}; it involves the
\textit{whole} matrix \ $S$ \ on its carrier \ $\alpha$, \ but not
only its some jump cardinals.
\\
Really, take any matrix \ $S$ \ on its carrier \ $\alpha$ \ of
\textit{zero} characteristic (if such exist), then by
lemma~\ref{8.5.}~5) there holds \vspace{-6pt}
\begin{multline*}
    \exists \tau_1^{\prime}, \tau_2^{\prime}, \tau_3^{\prime} <
    \alpha^{\Downarrow} \bigl(
    A_2^{0 \vartriangleleft \alpha^{\Downarrow}}
    ( \tau_1^{\prime}, \tau_2^{\prime}, \tau_3^{\prime},
    \alpha S_f^{< \alpha^{\Downarrow}}) \wedge
\\
    \wedge \forall \tau^{\prime\prime} \in \; ]\tau_1^{\prime},
    \tau_2^{\prime}] ~ a_{\tau^{\prime\prime}}^{<\alpha^{\Downarrow}} = 1
    \wedge \alpha S_{\tau_2^{\prime}}^{< \alpha^{\Downarrow}} = S \bigr),
\end{multline*}
where \ $S$ \ receives the \textit{lesser} carrier \
$\alpha_{\tau_2^{\prime}}^{< \alpha^{\Downarrow}}$, \ already of
the \textit{unit} characteristic due to the condition \ $\forall
\tau^{\prime\prime} \in \; ]\tau_1^{\prime}, \tau_2^{\prime}] ~
a_{\tau^{\prime\prime}}^{<\alpha^{\Downarrow}} = 1$.

But other matrix properties used in what follows are
\textit{inner}, and one of them realizes the
restriction-and-extension reasoning in the following lemma.

\noindent This lemma uses the suitable function, which was already
used in the proof of lemma~9.3:
\[
    Od \alpha S_{f}^{<\alpha _{1}} ( \tau _{1},\tau _{2} )  = \sup \{
    Od ( \alpha S_{\tau }^{<\alpha _{1}} ) :  \tau_{1} < \tau <\tau
    _{2}\};
\]
it will be applied to forming the so called \textit{stairways} ---
collections of intervals, which will be the main technical tools
in the Main theorem proof. To this end the following formulas
below \ $\alpha_1$ \ are needed:

\noindent 1.\quad $A_{1.1}^{m \vartriangleleft \alpha_1}(\tau_1,
\tau_2, \alpha S_f^{<\alpha_1})$:
\\
\[
    A_1^{1 \vartriangleleft \alpha_1}
    (\tau_1, \tau_2, \alpha S_f^{<\alpha_1})
    \wedge \tau_2 = \sup \big \{\tau: A_1^{1 \vartriangleleft \alpha_1}
    (\tau_1, \tau, \alpha S_f^{<\alpha_1}) \big \};
\]
\quad \\
here the interval \ $\left[ \tau_1,\tau_2 \right[$ \ is the
maximal of monotonicity intervals with the left \
$SIN_n^{<\alpha_1}$-end \ $\gamma_{\tau_1}^{<\alpha_1}$ \ and with
right \ $SIN_n^{<\alpha_1}$-ends, thus we shall call it and the
corresponding interval \ $[ \gamma_{\tau_1}^{<\alpha_1},
\gamma_{\tau_2}^{<\alpha_1} [$ \ the \textit{maximal intervals} of
the function \ $\alpha S_f^{<\alpha_1}$ \ \textit{monotonicity}
below \ $\alpha_1$.
\\
\quad %

\noindent 2.\quad $A_{1.1}^{m 1 \vartriangleleft \alpha_1}(\tau_1,
\tau_2, \alpha S_f^{<\alpha_1}, a_f^{<\alpha_1})$:
\\
\[
    \
    A^{0 \vartriangleleft \alpha_1}(\tau_1) \wedge
    A_{1.1}^{m \vartriangleleft \alpha_1}(\tau_1, \tau_2, \alpha S_f^{<\alpha_1})
    \wedge \forall \tau \big( \tau_1 < \tau < \tau_2
    \rightarrow a_{\tau}^{<\alpha_1} = 1 \big);
\]
\quad \\
in addition to \ $A_{1.1}^{m \vartriangleleft \alpha_1}$ \ here is
stated, that there is no \ $\alpha$-matrices admissible for \
$\gamma_{\tau_1}^{<\alpha_1}$ \ below \ $\alpha_1$ \ and the
function \ $\alpha S_f^{<\alpha_1}$ \ has on \ $]\tau_1,\tau_2[$ \
the values \ $\alpha S_{\tau}^{<\alpha_1}$ \ only of unit
characteristic \ $a_{\tau}^{<\alpha_1}=1$; \ in such cases the
unit characteristic \ $a=1$ \ stabilizes on the interval \
$[\tau_1,\tau_2[$ \ and on the corresponding interval \ $[
\gamma_{\tau_1}^{<\alpha_1}, \gamma_{\tau_2}^{<\alpha_1} [$ \
below \ $\alpha_1$ \ by lemma~\ref{9.4.}
\\
\quad \medskip %

\noindent 3.\quad $A_{1.1}^{st \vartriangleleft \alpha_1}(\tau_1,
\tau_{\ast}, \tau_2,\alpha
S_f^{<\alpha_1},a_f^{<\alpha_1})$:\newline
\[
    A_{1.1}^{m 1 \vartriangleleft \alpha_1}(\tau_1, \tau_{\ast},
    \alpha S_f^{<\alpha_1}, a_f^{<\alpha_1})
    \wedge \tau_1 < \tau_{\ast} \le \tau_2 \wedge
    A_{1}^{\vartriangleleft \alpha_1}(\tau_1, \tau_2, \alpha S_f^{<\alpha_1});
\]
\quad \\
here is indicated, that the function \ $\alpha S_f^{<\alpha_1}$ \
is defined on the interval \ $]\tau_1,\tau_2[$, \ but on its
maximal initial subinterval of monotonicity \
$]\tau_1,\tau_{\ast}[$ \ with \ $\gamma_{\tau_{\ast}} \in
SIN_n^{<\alpha_1}$ \ it has even the \textit{unit} characteristic
stabilized on it; therefore the interval \ $[\tau_1,\tau_2[$ \ and
the corresponding interval \ $[ \gamma_{\tau_1}^{<\alpha_1},
\gamma_{\tau_2}^{<\alpha_1} [$ \ will be called further the (unit)
steps below \ $\alpha_1$; in this case the ordinal
\[
    Od \alpha S_f^{<\alpha_1}(\tau_1, \tau_{\ast})
\]
will be called the \textit{height} of this step.
\\
\quad \medskip %

\noindent 4.\quad $A_{1.1}^{Mst \vartriangleleft \alpha_1}(\tau_1,
\tau_{\ast}, \tau_2,\alpha S_f^{<\alpha_1},a_f^{<\alpha_1})$:
\[
    \qquad
    A_{1.1}^{st \vartriangleleft \alpha_1}(\tau_1, \tau_{\ast}, \tau_2,
    \alpha S_f^{<\alpha_1}, a_f^{<\alpha_1})    \wedge
    A_{1.1}^{M \vartriangleleft \alpha_1}(\tau_1, \tau_2, \alpha S_f^{<\alpha_1});
\]
\quad \\
in addition here is indicated, that the interval \ $] \tau_1,
\tau_2[$ \ is the maximal with \
\mbox{$\gamma_{\tau_2}^{<\alpha_1} \in SIN_n^{<\alpha_1}$}, \
thereby we shall call the interval \ $[ \tau_1, \tau_2 [$ \ and
the corresponding interval \ $[ \gamma_{\tau_1}^{<\alpha_1},
\gamma_{\tau_2}^{<\alpha_1} [$ \ the \textit{maximal (unit) steps}
below \ $\alpha_1$.

This survey leads to the notion of stairway:
\\
\quad \\
5.\quad $A_{8}^{\mathcal{S} t \vartriangleleft \alpha_1}
(\mathcal{S} t, \chi, \alpha S_f^{<\alpha_1}, a_f^{<\alpha_1})$:
\[
    (\mathcal{S} t \mbox{ -- is a function on } \chi^{\ast +})
    \wedge \qquad\qquad\qquad\qquad\qquad\qquad\qquad\qquad
\]
\[
    \wedge \forall \beta < \chi^{\ast +} ~ \exists \tau_1, \tau_{\ast},
    \tau_2 \big( \mathcal{S} t(\beta) = (\tau_1,\tau_{\ast},\tau_2)
    \wedge
    \qquad\qquad\qquad\qquad
\]
\[
    \qquad\qquad\qquad
    \wedge A_{1.1}^{M st \vartriangleleft \alpha_1}(\tau_1,\tau_{\ast},\tau_2,
    \alpha S_f^{<\alpha_1}, a_f^{<\alpha_1}) \wedge
\]
\[
    \wedge \forall \tau_1, \tau_{\ast}, \tau_2 \bigl(
    A_{1.1}^{M st \vartriangleleft \alpha_1}(\tau_1,\tau_{\ast},\tau_2,
    \alpha S_f^{<\alpha_1}, a_f^{<\alpha_1})
    \longrightarrow
    \qquad\qquad\qquad
\]
\[
    \qquad\qquad\qquad\qquad\qquad
    \longrightarrow
    \exists \beta < \chi^{\ast +} ~ \mathcal{S}t(\beta) =
    (\tau_1, \tau_{\ast}, \tau_2) \bigr) \wedge
\]
\[
    \wedge \forall \beta_1, \beta_2 < \chi^{\ast +}
    ~ \forall \tau_1^{\prime}, \tau_{\ast}^{\prime}, \tau_2^{\prime}
    ~ \forall \tau_1^{\prime\prime}, \tau_{\ast}^{\prime\prime},
    \tau_2^{\prime\prime} \big( \beta_1 < \beta_2
    \wedge \qquad\qquad\qquad
\]
\[
    \wedge \mathcal{S} t(\beta_1) = (\tau_1^{\prime},
    \tau_{\ast}^{\prime}, \tau_2^{\prime})
    \wedge \mathcal{S} t(\beta_2) = (\tau_1^{\prime\prime},
    \tau_{\ast}^{\prime\prime}, \tau_2^{\prime\prime} ) \rightarrow
    \tau_2^{\prime} < \tau_1^{\prime\prime} \wedge
\]
\[
    \qquad\qquad\qquad
    \wedge Od~\alpha S_f^{<\alpha_1}(\tau_1^{\prime}, \tau_{\ast}^{\prime})
    < Od~\alpha S_f^{<\alpha_1}(\tau_1^{\prime\prime}, \tau_{\ast}^{\prime\prime})
    \big) \wedge
\]
\[
    \wedge \sup \big\{ Od~\alpha S_f^{<\alpha_1}(\tau_1, \tau_{\ast}) : \exists
    \beta, \tau_2 ~ \mathcal{S} t(\beta) = (\tau_1,\tau_{\ast},\tau_2) \big\} =
    \chi^{\ast +};
\]
here is indicated, that \ $\mathcal{S} t$ \ is the function on \
$\chi^{\ast +}$, \ and its values are all triples \
$(\tau_1,\tau_{\ast},\tau_2)$ \ such that the intervals \ $[
\tau_1, \tau_2 [$ \ are maximal unit steps disposed successively
one after another. Therefore such \ $\mathcal{S} t$ \ will be
called the \textit{stairway} and the intervals \ $[ \tau_1, \tau_2
[$ \ and the corresponding intervals \ $[
\gamma_{\tau_1}^{<\alpha_1}, \gamma_{\tau_2}^{<\alpha_1} [$ -- its
steps below \ $\alpha_1$.
\\
This notion is justified by the strict increasing of their
heights; also we shall say, that the stairway \ $\mathcal{S} t$ \
consists of these steps, or contains them.
\\
Respectively, the cardinal
\[
    h(\mathcal{S} t) = \sup \big \{ Od~\alpha S_f^{<\alpha_1}(\tau_1,\tau_{\ast}):
    \exists \beta, \tau_2 ~ \mathcal{S} t(\beta) =
    (\tau_1, \tau_{\ast}, \tau_2) \big \}
\]
will be called the \textit{height} of the whole stairway \
$\mathcal{S} t$. \ So, here is required that \ $\mathcal{S} t$ \
amounts up to \ $\chi^{\ast +}$, \ that is
\[
    h(\mathcal{S} t) = \chi^{\ast +}.
\]
Also the cardinal
\[
    \upsilon = \sup \left\{ \gamma_{\tau_2}: \exists \beta, \tau_1,
    \tau_{\ast} ~ \mathcal{S} t(\beta) = (\tau_1, \tau_{\ast}, \tau_2)
    \right\}
\]
will be called the termination cardinal of \ $\mathcal{S} t$ \ and
will be denoted through
\[
    \upsilon(\mathcal{S} t);
\]
so, we shall say, that the stairway \ $\mathcal{S} t$ \ terminates
in this cardinal \ $\upsilon(\mathcal{S} t)$.

\noindent If such stairway \ $\mathcal{S} t$ \ exist below \
$\alpha_1$, \ then we shall say, that \ $\alpha_1$ \ is
\textit{provided} by this stairway \ $\mathcal{S} t$.
\\
When \ $\alpha > \chi^{\ast}$ \ is a carrier of the matrix \ $S$ \
and its prejump cardinal \ $\alpha_1 = \alpha_{\chi}^{\Downarrow}$
\ after \ $\chi^{\ast}$ \ is provided by some stairway \
$\mathcal{S} t$, \ then we shall say, that this \ $S$ \ on \
$\alpha$ \ is provided by this stairway.

And here is quite important example of the inner phenomena: the
\textit{inner} property of providing the matrix \ $S$ \ by some
stairway.
\\
This property for \ $S$ \ on its carrier \ $\alpha$ \ is definable
by the formula
\[
    \exists \mathcal{S}t \vartriangleleft \alpha^{\Downarrow +} ~
    A_{8}^{\mathcal{S} t \vartriangleleft \alpha^{\Downarrow} }
    (\mathcal{S} t, \alpha S_f^{< \alpha^{\Downarrow}},
    a_f^{< \alpha^{\Downarrow}})
\]
which can be bounded by the jump cardinal \ $\alpha^{\downarrow}$
\ of the carrier \ $\alpha$ \ after \ $\chi^{\ast}$. \ Therefore
by lemma 5.11~\cite{Kiselev9} the same property holds for \ $S$ \
on any other carrier \ $\alpha^1 > \chi^{\ast}$:
\[
    \exists \mathcal{S}t^1 \vartriangleleft \alpha^{1 \Downarrow +} ~
    A_{8}^{\mathcal{S} t \vartriangleleft \alpha^{ 1 \Downarrow} }
    (\mathcal{S} t^1, \alpha S_f^{< \alpha^{1 \Downarrow}},
    a_f^{< \alpha^{1 \Downarrow}} )
\]
being bounded by the jump cardinal \ $\alpha^{1 \downarrow}$ \ of
\ $\alpha^{1}$ \ after \ $\chi^{\ast}$, \ and, so, \ $S$ \ on \
$\alpha^1$ \ is again provided by some stairway \
${\mathcal{S}t}^1$ \ as well.

\begin{lemma}
\label{9.5.} \emph{(About stairway cut-off from above)}  \\
\hspace*{1em} Let
\medskip

(i) \quad $A_{1}^{1 \vartriangleleft \alpha_1}(\tau _{1},\tau
_{2})$;
\newline

(ii) \quad $\tau_2 \le \tau_3$ \ and \ $S^3$ \ be a matrix of
characteristic \ $a^3$ \ on a carrier
\[
    \alpha_3 \in \; ]\gamma_{\tau_3}^{<\alpha_1}, \alpha_1 [
\]
with disseminator \ $\widetilde{\delta}^3$ \ and data base \
$\rho^3$ \ admissible for \ $\gamma_{\tau_3}^{<\alpha_1}$ \ below
\ $\alpha_1$ \ and with the generating eigendisseminator \
$\check{\delta}^{S^3}$ \ on \ $\alpha^3$;
\newline

(iii) \quad $\forall \tau <\tau _{2}\quad \exists \tau ^{\prime
}\in \left[ \tau ,\tau _{2}\right[ \quad a_{\tau ^{\prime
}}^{<\alpha _{1}} = a^3$. \medskip

\noindent Then \medskip

1. \ $Od\alpha S_{f}^{<\alpha _{1}}(\tau _{1},\tau _{2}) <
Od(S^3)$;
\\

2a. \ hence, if \ $a^3 = 1$, \ then there is no stairway below \
$\alpha_1$ \ terminating in some \ $SIN_n^{<\alpha_1}$-cardinal \
$\upsilon < \alpha_3^{\Downarrow}$;
\\

2b. \ therefore if there exist some unit matrix \ $S^0$ \ on its
carriers admissible below \ $\alpha_1$ \ and disposed cofinally to
\ $\alpha_1$: \vspace{-6pt}
\begin{multline*}
    \forall \gamma < \alpha_1 ~ \exists \gamma^1 \in \; ]\gamma,
    \alpha_1[ ~ \exists \delta, \alpha, \rho < \alpha_1 ~ \big(
    SIN_{n-1}^{<\alpha_1}(\gamma^1) \wedge
\\
    \wedge \alpha \mathbf{K}^{<\alpha_1}(1, \delta, \gamma^1, \alpha, \rho,
    S^0) \big),
\end{multline*}
then \ $\alpha_1$ \ is not provided by any
stairway;
\\

3. \ if \ $S^3$ \ is the \ $\underline{\lessdot}$-minimal of all
matrices of the same characteristic \ $a^3$ \ on carriers \ $\in
\; ] \gamma_{\tau_3}^{<\alpha_1}, \alpha_1[$ \ admissible for \
$\gamma_{\tau_3}^{<\alpha_1}$, \ then
\[
    \gamma_{\tau _{2}}^{<\alpha _{1}}<\check{\delta}^{S^3} \leq
    \widetilde{\delta }^3 < \gamma _{\tau_{3}}.
\]
\end{lemma}

\noindent \textit{Proof.} \ Let us demonstrate the reasoning for
the case \ $a^3=1$, \ used in what follows; in this important case
\ $\widetilde{\delta}^3 = \check{\delta}^{S^3}$ \ and condition
$(iii) $ should be weakened up to \ $a^3=1$ \ by lemma
3.2~\cite{Kiselev9}. In this case the nonsuppression condition \
$\neg A_5^{S,0}$ \ for the unit matrix \ $S$ \ on its carriers can
be dropped, because such \ $S$ \ is always nonsuppressed and the
formulas \ $\alpha \mathbf{K}^{\ast < \alpha_1}$, $\alpha
\mathbf{K}_{n+1}^{\ast < \alpha_1}$ \ are equivalent to the
formulas \ $\alpha \mathbf{K}^{< \alpha_1}$, $\alpha
\mathbf{K}_{n+1}^{< \alpha_1}$; \ the upper indices \ $<
\alpha_{1} $, $\vartriangleleft \alpha_{1}$ \ will be omitted.
\\
By this lemma the matrix \ $S^3$ receives unit characteristic on
its admissible carriers disposed cofinally to \ $
\gamma_{\tau_{2}}$, \ as it was in part 1a. in lemma~\ref{9.4.}
proof, where \ $\gamma_{\tau^2}$, \ $S^0$ \ on \ $\alpha_{\tau^1}$
\ should be replaced with \ $\gamma_{\tau_2}$, \ $S^3$ \ on
$\alpha_3$. \ So $(i)$ and lemma~\ref{9.4.} imply
\[
    Od \alpha S_{f} (\tau_{1}, \tau_{2}) \leq Od (S^3)
    \mathrm{\ and\ } a_{\tau} \equiv 1 \mathrm{\ on\ }
    ]\tau_{1}, \tau_{2} [ ~.
\]
Now let us assume that the function \ $\alpha S_f$ \ stabilizes on
\ $[ \tau_1, \tau_2[$ \ and let \ $\tau_{2}^s $ \ be the
stabilization ordinal of \ $\alpha S_{f} $ \ on \ $[\tau_{1},
\tau_{2} [ $ \ so that there exists \ $S^0 $ \ such that
\[
    \alpha S_{\tau} \equiv \alpha S_{\tau_2^s} = S^{0}
    \mathrm{\ on\ } [ \tau_2^s, \tau_2[ ~.
\]
We apply now the restriction-and-extension argument mode, that was
used several times above. Let us turn to the matrix \ $S^{0} $ \
on the carrier \ $\alpha_{\tau_2^s+1} $ \ with the prejump
cardinal \ $\alpha^{1} = \alpha_{\tau_2^s+1}^{\Downarrow} $ \ and
the disseminator \ $\check{\delta}^{1} =
\check{\delta}_{\tau_2^s+1}$. \ The same matrix \ $S^{0} $ \ on
the carrier \ $\alpha_{\tau_2^s} $ \ of unit characteristic by
lemma 3.2~\cite{Kiselev9} about restriction receives unit
characteristic also on its admissible carriers disposed cofinally
to \ $\check{\delta}^{1} $ \ and hence below \ $\check{\delta}^{1}
$ \ the following \ $\Pi_{n+1} $-proposition holds for \
$S=S^{0}$:
\begin{equation} \label{e9.11}
    \qquad\qquad
    \forall \gamma \exists \gamma^{\prime} ~ \bigl( \gamma <
    \gamma^{\prime} \wedge SIN_{n-1}(\gamma^{\prime}) \wedge
    \qquad\qquad\qquad\qquad\qquad\qquad\qquad
\end{equation}
\[
    \qquad\qquad\qquad
    \wedge \exists \delta, \alpha, \rho ~~ \alpha\mathbf{K} (1,
    \delta, \gamma^{\prime}, \alpha, \rho, S) \bigr ) ~.
\]
\vspace{-6pt}

\noindent Disseminator \ $\check{\delta}^{1} $ \ extends it up to
\ $ \alpha^{1} $ \ and so the matrix \ $S^{0} $ \ receives unit
characteristic on its admissible carriers disposed cofinally to \
$\alpha^{1} $, \ that is (\ref{e9.11}) is fulfilled by  \ $S =
S^{0} $ \ under the boundary \ $ \vartriangleleft \alpha^{1} $. \
After the minimizing such matrices \ $S $ \ we receive the matrix
\ $S=S^{1} $ \ with the property (\ref{e9.11}) below  \ $
\alpha^{1} $ \ and by lemma 4.6~\cite{Kiselev9} \ $S^{1} \lessdot
S^{0} $. \ One should point out, that statement (\ref{e9.11}) \
$\vartriangleleft$-bounded by \ $\alpha^1$ \ with \ $S = S^1$ \ is
the \textit{inner property} of the matrix \ $S^0$.
\\
If now

\[
    Od \alpha S_{f} (\tau_{1},\tau_{2}) = Od ( S^3 ), \mathrm{\
    that\ is\ } S^{0} = S^3,
\]
then matrix \ $S^{1} $ \ by lemma 5.11~\cite{Kiselev9} about
informativeness receives its admissible carriers with the same
unit characteristic disposed cofinally to the prejump cardinal \
$\alpha^{3}=\alpha_3^{\Downarrow} $, \ because \ $S^{0} $ \ on the
carrier \ $\alpha_{\tau_2^s+1} $ \ has the same property.
\\
After that again by lemma 3.2~\cite{Kiselev9} such carriers appear
disposed cofinally to \ $\gamma _{\tau_{2}} $. \ Hence, at last,
by $(i)$ it comes the contradiction:

\begin{equation}  \label{e9.12}
Od \alpha S_{f} (\tau_{1},\tau_{2})\leq Od (S^{1})< Od (S^{0}).
\end{equation}
\vspace{0pt}

\noindent If the function \ $\alpha S_f$ \ does not stabilize on \ $[\tau_1,
\tau_2[$ then the ordinal
\[
    \rho=Od\alpha S_f(\tau_1, \tau_2)
\]
is limit. But here to finish the proof of 1. one should remind,
that G\"{o}del function \ $F$ \ has values \ $F (\alpha ) = F |
\alpha$ \ for limit ordinals \ $\alpha$ \ (see G\"{o}del
\cite{Godel}). One can see that \ $Od(S)$ \ can not be limit and
thus \ $\rho \leq Od(S) $ \ implies \ $\rho < Od(S)$.
\\

\noindent Turning to 3. let us suppose that it is wrong and
\[
    \check{\delta}^{3} = \check{\delta}^{S^3} \leq \gamma_{\tau_{2}}
\]
and standing on \ $\alpha^3 = \alpha_3^{\Downarrow}$ \ let us
review the situation below \ $ \alpha^{3} $ \ obtained. Two cases
here should be considered:
\\
1. \quad $\gamma_{\tau_{1}} < \check{\delta}^{3} \leq
\gamma_{\tau_{2}} $.
\\
Since \ $\gamma_{\tau_{2}} \in SIN_{n} $ \ and \
$\check{\delta}^{3} \in SIN_{n}^{< \alpha^{3}} $ \
lemmas~3.8~\cite{Kiselev9}, 8.5~1), 8.7 imply that \ $\check{
\delta}^{3} \in SIN_{n} $ \ and
\[
    \alpha S_{\tau}^{< \alpha^{3}} \equiv \alpha S_{\tau}
    \mathrm{\quad on \quad}  \{
    \tau : \gamma_{\tau_{1}} < \gamma_{\tau} < \check{\delta}^{3} \}.
\]
The disseminator \ $\check{\delta}^{3} $ \ extends up to \
$\alpha^{3} $ \ the \ $\Pi_{n+1} $-proposition stating the
definiteness of the function \ $ \alpha S_{f} $ \ of \textit{unit}
characteristic and with values \ $\lessdot S^3$ \ due to part 1.:
\begin{equation}  \label{e9.13}
    \qquad
    \forall \gamma^{\prime} ~ \Big( \gamma_{\tau_1} < \gamma^{\prime} \wedge
    SIN_{n-1}(\gamma^{\prime}) ~ \longrightarrow
    \qquad \qquad \qquad \qquad \qquad \qquad \qquad
\end{equation}
\[
    \qquad
    \longrightarrow \exists \delta, \alpha, \rho, S ~~
    \big ( S \lessdot S^3 \wedge
    \alpha \mathbf{K}(1, \delta, \gamma^{\prime}, \alpha, \rho, S) \big) \Big)
\]
\vspace{0pt}

\noindent and, hence, there exists the matrix \ $\alpha
S_{\tau_{3}}^{< \alpha^{3}} $ \ on the carrier
\[
    \alpha_{\tau_3}^{<\alpha^3} \in \; ] \gamma_{\tau_3}, \alpha^3 [
\]
of unit characteristic and \ $\alpha S_{\tau_{3}}^{< \alpha^{3}}
\lessdot S^3 $ \ contrary to the \ $\underline{\lessdot}
$-minimality of \ $ S^3$ \ on \  $\alpha_3$. \ It remains to
consider the case:
\\

\noindent 2. \quad $\check{\delta}^{3} \leq \gamma_{\tau_{1}} $.
\\
It should be pointed out that the condition of \
$\underline{\lessdot} $-minimality of \ $S^3$ \ is not used in
this case. Here the matrix \ $S^0 = \alpha S_{\tau_1 + 2}$ \ of
unit characteristic should be considered on its carrier \
$\alpha_{\tau_1 + 2}$ \ with the prejump cardinal \ $\alpha^1 =
\alpha_{\tau_1 + 2}^{\Downarrow}$ \ and the disseminator \
$\check{\delta}^1 = \check{\delta}_{\tau_1 + 2}$, \ just as it was
done above for \ $S^0 = \alpha S_{\tau_2^s + 1}$, \ $\alpha^1 =
\alpha_{\tau_2^s + 1}^{\Downarrow}$, \ $\check{\delta}^1 =
\check{\delta} _{\tau_2^s + 1}$ \ in the proof of 1. (let us
preserve the previous notations for convenience). And again the
matrix \ $\alpha S_{\tau_{1} + 1} $ \ receives the unit
characteristic on its admissible carriers disposed cofinally to
the disseminator \ $\check{\delta}^1 = \gamma_{\tau_1}$ \ and it
extends proposition (\ref{e9.11}) for \ $S = \alpha S_{\tau_1 +
1}$ \ up to \ $\alpha^1$; \ so, it brings by the same way the
minimal matrix \ $S^1 \lessdot S^0$ \ with the previous
properties: it receives unit characteristic on its admissible
carriers disposed cofinally to \ $\alpha^1$.
\\
By lemma 3.2~\cite{Kiselev9} there appear carriers of \ $S^{1} $ \
of unit characteristic disposed cofinally to \ $\check{\delta}^3$,
\ that is (\ref{e9.11}) is fulfilled by \ $S=S^{1} $ \ under the \
$\vartriangleleft$-bounding by \ $\check{\delta }^3$; \ hence, the
disseminator \ $\check{\delta}^{3} $ \ extends this proposition up
to \ $\alpha^{3} $. \ After that the cardinal \ $
\gamma_{\tau_{2}} \in SIN_{n}$ \ restricts this proposition with \
$\gamma $ \ substituted for an arbitrary constant \ $\gamma^1 <
\gamma_{\tau_{2}} $. \ As a result the matrix \ $S^{1} $ \
receives the unit characteristic on its admissible carriers
disposed cofinally to \ $\gamma_{\tau_{2}} $ \ and so again we
come to contradiction (\ref{e9.12}).

Turning to \ $\alpha^3=0$ \ one should repeat all this proof but
for \textit{zero} matrices \ $S$ \ on their carriers \ $\alpha$ \
admissible for cardinals \ $\gamma_{\tau}$ \ under consideration,
but only for \ $\gamma_{\tau} \notin SIN_n$. \ In all cases of
this kind such matrices \ $S$ \ are nonsuppressed by definition
and again the nonsuppression condition \ $\neg A_5^{S,0}$ \ can be
dropped, and the formulas \ $\alpha \mathbf{K}^{\ast < \alpha_1}$,
$\alpha \mathbf{K}_{n+1}^{\ast < \alpha_1}$ \ can be replaced with
the formulas \ $\alpha \mathbf{K}^{< \alpha_1}$, $\alpha
\mathbf{K}_{n+1}^{< \alpha_1}$. \ Precisely such matrices \ $S$ \
on their carriers \ $\alpha$ \ can be used in the
restriction-and-extension reasoning above, that provides the proof
for zero characteristic \ $a_3=0$.

At last proposition 2a. comes from 1. almost obviously. Suppose it
fails, that is there exist some stairway \ $\mathcal{S}t$ \ below
\ $\alpha_1$ \ terminating in the \ $SIN_n$-cardinal \
$\upsilon(\mathcal{S}t) < \alpha_3^{\Downarrow}$; \ it implies
that
\[
    \upsilon(\mathcal{S}t) < \gamma_{\tau_3}^{<\alpha_1}.
\]
By definition this stairway consists of \textit{unit} steps below
\ $\alpha_1$
\[
    \mathcal{S}t(\tau^{\prime}) = (\tau_1^{\prime},
    \tau_2^{\prime\ast}, \tau_2^{\prime})
\]
and each of them possesses the property
\[
    A_1^{1 \vartriangleleft \alpha_1} (\tau_1^{\prime},
    \tau_2^{\prime\ast})
\]
with the unit characteristic stabilizing on \ $[ \tau_1^{\prime},
\tau_2^{\prime\ast} [$ \ (see definition of the stairway notion
before lemma \ref{9.5.}). Hence 1. provides the height \
$h(\mathcal{S}t)$ \ of the whole \ $\mathcal{S}t$ \ bounded by the
ordinal
\[
    Od(S^3) < \chi^{\ast +}
\]
though \ $h(\mathcal{S}t) = \chi^{\ast +}$, \ that is \
$\mathcal{S}t$ \ amounts up to \ $\chi^{\ast +}$ \ by definition.
\\
From here it follows 2b. when the matrix \ $S^0$ \ is used instead
of the matrix \ $S^3$. \hspace*{\fill} $\dashv$
\\
One should notice, that for \ $a^3 = 1$ \ the condition of the
matrix \ $S^3$ \ minimality in point 3. of this lemma~\ref{9.5.}
can be dropped by means of the reasoning repeating the arguments
above in case 1. slightly changed.
\\

Next obvious corollary shows that such steps heights always
increase strictly:

\begin{corollary}
\label{9.6.} \hfill {} \newline
\hspace*{1em} Let \medskip

(i) \quad $A_{1}^{1 \vartriangleleft \alpha_1} ( \tau _{1},\tau
_{2})$, \ $A_{1}^{1 \vartriangleleft \alpha_1} ( \tau _{3},\tau
_{4}) $, \ $\tau _{2}<\tau _{4}$; \newline

(ii)\quad $\forall \tau <\tau _{4} \quad \exists \tau ^{\prime }\in \left[
\tau ;\tau _{4}\right[ \quad a_{\tau ^{\prime }}^{<\alpha _{1}}=1$. \medskip

\noindent Then \medskip

1) \ $\forall \tau \in \; ]\tau_1, \tau_2[ \; \cup \; ]\tau_3, \tau_4[ \quad
a_\tau^{<\alpha_1} = 1$; \newline

2) \ $Od\alpha S_{f}^{<\alpha _{1}} ( \tau _{1},\tau _{2} ) <Od\alpha
S_{f}^{<\alpha _{1}} ( \tau _{3},\tau _{4} )$; \newline

3) \ $\forall \tau \in \; ]\tau_3, \tau_4[ \quad \gamma _{\tau
_{2}}^{<\alpha_{1}} < \check{\delta}_{\tau}^{S} = \widetilde{
\delta}_{\tau}^{ < \alpha _{1}} < \gamma^{ < \alpha _{1}}_{\tau}$,
\\
where \ $\check{\delta}_{\tau}^{S}$ \ is the generating
eigendisseminator of \ $\alpha S_{\tau}^{< \alpha_1}$ \ on \
$\alpha_{\tau}^{< \alpha_1}$; \ hence
\[
    \gamma_{\tau_2}^{<\alpha_1} < \gamma_{\tau_3}^{<\alpha_1}.
\]
\end{corollary}

\noindent \textit{Proof.} \ From conditions $(i)$, $(ii)$ and
lemmas 3.2~\cite{Kiselev9}, \ref{9.4.} there follows that \
$a_{\tau}^{<\alpha_1} \equiv 1$ \ on the intervals \ $]\tau_1,
\tau_2[$, \ $]\tau_3, \tau_4[$\;. \ Therefore lemma~\ref{9.5.}
(where \ $\tau_3$ \ plays the role of any \ $\tau \in \; ]\tau_3,
\tau_4[$\;) implies statements~2),~3). For \ $\tau = \tau_3 + 1$ \
here it comes \ $\gamma_{\tau_2}^{<\alpha_1} <
\check{\delta}_{\tau}^{<\alpha_1}$ \ and at the same time by
lemmas~\ref{8.5.}~7)~$(ii)$, \ref{8.8.}~2) \ -- \
$\check{\delta}_{\tau}^{S} = \gamma_{\tau_3}^{<\alpha_1}$; \ so, \
$\tau_2 < \tau_3$.
\\
\hspace*{\fill} $\dashv$

\begin{corollary}
\label{9.7.} \hfill {} \newline
\hspace*{1em} Let \medskip

(i) \quad $A_{1}^{1 \vartriangleleft \alpha_1}(\tau _{1},\tau
_{2})$; \newline

(ii) \quad $\tau _{3}\in dom(\alpha S_{f}^{<\alpha _{1}})$, \
$\tau _{3}\geq \tau _{2}$;
\\

(iii) \quad a matrix \ $\alpha S_{\tau_3}^{< \alpha_1}$ \ on \
$\alpha_{\tau_3}^{< \alpha_1}$ \ has generating eigendisseminator
\[
    \check{\delta}_{\tau _{3}}^{S}\leq
    \gamma _{\tau _{2}}^{<\alpha _{1}}
\]
below \ $\alpha_1$.
\\
Then \medskip

1) \ $a_{\tau }^{<\alpha _{1}}\equiv 1$ \ on \ $\left] \tau
_{1},\tau _{2} \right[ $, \quad $a_{\tau _{3}}^{<\alpha _{1}}=0$;
\newline

2) \ $\check{\delta}_{\tau _{3}}^S \leq \gamma _{\tau
_{1}}^{<\alpha _{1}}$ \quad and \quad \newline

3) \ $Od\alpha S_{f}^{<\alpha _{1}}(\tau _{1},\tau _{2})>Od(\alpha
S_{\tau _{3}}^{<\alpha _{1}})$.
\\

\noindent Analogously for the generating disseminator \
$\check{\delta}_{\tau _{3}}^{<\alpha_1} $ \ of \ $\alpha S_{\tau
_{3}}^{<\alpha_1}$ \ on \ $\alpha_{\tau _{3}}^{<\alpha_1}$.
\end{corollary}

\noindent \textit{Proof.} \ We shall omit the upper indices \ $<
\alpha_{1} $, $\vartriangleleft \alpha_{1}$. \ By lemma~\ref{9.5.}
for \ $S^3 = \alpha S_{\tau_3}^{< \alpha_1}$ \ condition $(iii)$
implies that for some \ $\tau < \tau_{2}$
\begin{equation*}
\forall \tau^{\prime} \in [\tau, \tau_{2} [ \quad  a_{\tau^{\prime}} \neq
a_{\tau_{3}} ;
\end{equation*}
due to lemma 3.2~\cite{Kiselev9} it is possible only when
\begin{equation*}
\forall \tau^{\prime} \in [\tau, \tau_{2} [ \quad  a_{\tau^{\prime}} =1,
\quad a_{\tau_{3}} = 0
\end{equation*}
and then by lemma 9.4 \ $a_{\tau} \equiv 1 $ \ on \ $]\tau_{1},
\tau_{2} [ $\;.
\\
If \ $\check{\delta}_{\tau_{3}}^S \in \; ] \gamma_{\tau_{1}},
\gamma_{\tau_{2}}[\;$, \ then one can obtain \ $a_{\tau_{3}} =1 $
\ again using the reasoning from the proof of lemma~\ref{9.5.},
and extending proposition~(\ref{e9.13}) without its subformula \
$S \lessdot S^3$ \ by the disseminator \
$\check{\delta}_{\tau_{3}}^S$ \ up to \
$\alpha_{\tau_{3}}^{\Downarrow}$, \ that provides \ $a_{\tau_3}^{<
\alpha_1} = 1$.
\\
After that it is enough to conduct the reasoning from the end of
this proof (case 2.) repeated literally by means of
restriction-and-extension method.
\\
\hspace*{\fill} $\dashv$
\\
\quad \\
The immediate consequence of this lemma for \ $\tau_{2}=
\tau_{3}$ \ is the following
\newline

\noindent \textbf{Theorem 1.}\quad \newline
\emph{\hspace*{1em} Let \medskip }

\emph{(i) \quad $\alpha S_f^{<\alpha_1}$ \ be monotone on \ $[\tau_1, \tau_2[
$ \ below \ $\alpha_1$; \newline
}

\emph{(ii) \quad $\tau_1 = \min \{ \tau: \; ]\tau, \tau_2[ \; \subseteq dom
(\alpha S_f^{<\alpha_1}) \} $. \medskip }

\emph{\noindent Then
\begin{equation*}
]\gamma_{\tau_1}^{< \alpha_1}, \gamma_{\tau_2}^{< \alpha_1} [  ~
\cap ~ SIN_n^{<\alpha_1} = \varnothing.
\end{equation*}
\vspace{0pt} }

\noindent \textit{Proof.} \ Let us suppose that, on the contrary,
there exists \ $SIN_n^{<\alpha_1}$-cardinal \
$\gamma_{\tau_2^{\prime}}^{< \alpha_1} \in ]\gamma_{\tau_1}^{<
\alpha_1}, \gamma_{\tau_2}^{< \alpha_1} [ $~. \newline Then \
$\gamma_{\tau_1}^{< \alpha_1}$ \ belongs to \ $SIN_n^{< \alpha_1}$
\ as well; one can see it repeating literally the argument from
the proof of lemma~\ref{8.10.}. So, the statement \ $A_1^{1
\vartriangleleft \alpha_1}(\tau_1, \tau_2^\prime)$ \ holds on; it
remains to apply corollary~\ref{9.7.} using \ $\tau_2^{\prime}$ \
as \ $\tau_2 = \tau_3$, \ since \
$\check{\delta}_{\tau_2^\prime}^{S} <
\gamma_{\tau_2^\prime}^{<\alpha_1}$ \ by definition.
\\
\hspace*{\fill} $\dashv$

\newpage

\section{Analysis of \ $\protect\alpha $\,-Function Nonmonotonicity}

\setcounter{equation}{0}

\hspace*{1em} So, any interval of the \ $\alpha$-function
monotonicity can not be ``too long'' by theorem 1.
\\
However, such function can be defined on ``rather long''
intervals; for example, the function \ $\alpha S_f^{<\alpha_1}$ \
is defined on the final segment \ $T^{\alpha_1}$ \ of any
sufficiently great \ $SIN_n$-cardinal \ $\alpha_1 < k$ \ (lemma
8.9\;). Hence, its monotonicity on this segment is violated on
some ordinals.
\newline How does this phenomenon happen? In this section all
substantial violations of this kind are analyzed. To this end the
formula \ $A_2^{\vartriangleleft \alpha_1}(\tau_1, \tau_2,
\tau_3)$ \ should be recalled (see definition~\ref{8.1.}~1.4 for \
$X_1 = \alpha S_f^{<\alpha_1}$): \vspace{-6pt}
\begin{multline*}
    A_1^{\vartriangleleft \alpha_1}(\tau_1, \tau_3, \alpha
    S_f^{<\alpha_1}) \wedge \tau_1 + 1 < \tau_2 <  \tau_3 \wedge
\\
    \wedge \tau_2 = \sup \bigl \{ \tau < \tau_3 :  \forall
    \tau^\prime, \tau^{\prime\prime} (\tau_1 < \tau^\prime <
    \tau^{\prime\prime} < \tau \rightarrow  \alpha
    S_{\tau^\prime}^{<\alpha_1} \underline{\lessdot} \alpha
    S_{\tau^{\prime\prime}}^{<\alpha_1} \bigr ) \}.
\end{multline*}

\noindent So, here \ $\tau_2$ \ is the \textit{minimal} ordinal
breaking the monotonicity of the function \ $\alpha
S_f^{<\alpha_1}$ \ on the interval \ $ [ \tau_1, \tau_3 [$\;. \
Thus, in all reasonings of this paragraph some nonmonotonicity \
$A_2^{\vartriangleleft \alpha_1}(\tau_1, \tau_2, \tau_3)$ \ on the
intervals \ $[\tau_1, \tau_3[$ \ is treated in different
situations (but the condition \ $SIN_n^{\vartriangleleft
\alpha_1}(\gamma_{\tau_3})$ \ can be dropped everywhere except the
last lemma~\ref{10.5.}\;).

\begin{lemma}
\label{10.1.} \hfill {}
\\
Let \medskip

(i) \ $A_{2}^{\vartriangleleft \alpha_1}(\tau _{1},\tau _{2},\tau
_{3})$;
\newline

(ii) \ $SIN_{n}^{<\alpha _{\tau _{2}}^{\Downarrow }}\cap \gamma _{\tau
_{2}}^{<\alpha _{1}}\subseteq SIN_{n}^{<\alpha _{1}}$. \medskip

\noindent Then \medskip

1) \ $a_{\tau }^{<\alpha _{1}}\equiv 1$ \ on \ $] \tau _{1},
\tau_{2} [ $, \ $a_{\tau _{2}}^{<\alpha _{1}}=0$ \quad and
\newline

2) \ $\widetilde{\delta }_{\tau _{2}}^{<\alpha _{1}} \le \gamma
_{\tau _{1}}^{<\alpha _{1}}$. \label{c13}
\endnote{
\ p. \pageref{c13}. \ It can be proved, that here \
$\widetilde{\delta}_{\tau_2}^{<\alpha_1} =
\gamma_{\tau_1}^{<\alpha_1}$.
\\
\quad \\
} %
\end{lemma}

\noindent \textit{Proof.} \ The upper indices \ $< \alpha_{1} $,
$\vartriangleleft \alpha_{1} $ \ will be dropped. Since the
function \ $\alpha S_f$ \ is monotone on \ $]\tau_1,\tau_2[$, \
from theorem~1 it comes
\begin{equation} \label{e10.1}
    ]\gamma_{\tau_1}, \gamma_{\tau_2}[ \;  \cap \; SIN_n =
    \varnothing.
\end{equation}
Standing on \ $\alpha^2 = \alpha_{\tau_2}^{\Downarrow}$, \ let us
consider below \ $\alpha^{2}$ \ the function \ $\alpha
S_{f}^{<\alpha^{2}}$. \ By $(ii)$ and lemma~\ref{8.7.} about
absoluteness it coincides with \ $\alpha S_{f}$ \ on \ $] \tau
_{1}, \tau _{2} [ $ \ and is monotone on this interval.
\\
That is why \ $\widetilde{\delta}_{\tau_2} \le \gamma_{\tau_1}$, \
otherwise \ $\widetilde{\delta}_{\tau_2} \in \; ] \gamma_{\tau_1},
\gamma_{\tau_2}[$ \ contrary to $(ii)$, (\ref{e10.1}).
\\
If now \ $a_{\tau _{2}}=1$, \ then by lemma~\ref{9.2.}
\begin{equation*}
Od\alpha S_{f} ( \tau _{1},\tau _{2} ) \leq  Od ( \alpha S_{\tau _{2}} )
\end{equation*}
in spite of (\textit{i}) and so \ $a_{\tau_2} = 0$. \ The same
happens if
\[
    \forall \tau <\tau _{2} \quad \exists \tau ^{\prime}
    \in \left[ \tau ;\tau _{2}\right[ \quad a_{\tau ^{\prime}}=0
\]
because in this case due to condition $(i)$ one can consider
\mbox{$\tau_1^2 \in \; ]\tau_1, \tau_2[$} \ such that for \ $S^2 =
\alpha S_{\tau_2}$
\begin{equation} \label{e10.2}
    a_{\tau_1^2} = 0, \quad \alpha S_{\tau_1^2} \gtrdot S^2.
\end{equation}
Due to consequence 2) and lemma 3.2~\cite{Kiselev9} about
restriction \textit{zero} matrix \ $S^2$ \ receives some
admissible carrier \ $\alpha \in \; ]\gamma_{\tau_1^2},
\gamma_{\tau_1^2+1}[$ \ as a result of restricting by \
$SIN_{n-1}$-cardinal \ $\gamma_{\tau_1^2+1}$ \ of the following \
$\Sigma_{n-1}$-proposition
\begin{equation} \label{e10.3}
    \quad
    \exists \alpha \Big( \gamma_{\tau_1^2} < \alpha \wedge
    \exists \delta, \alpha, \rho ~ \bigl(
    \delta \le \gamma_{\tau_1} \wedge
    SIN_n^{<\alpha^{\Downarrow}}(\gamma_{\tau_1}) \wedge
    \qquad \qquad \qquad \qquad
\end{equation}
\begin{equation*}
    \qquad\qquad\qquad\qquad\qquad
    \wedge \alpha \mathbf{K}_{n+1}^{\exists}(0, \delta, \gamma_{\tau_1^2}, \alpha, \rho, S^2)
    \big) \Big)
\end{equation*}
which holds below \ $\gamma_{\tau_1^2 + 1}$, \ since it holds for
\ $\alpha=\alpha_{\tau_2}$ \ below \ $\alpha_1$.
\\
Hence, (\ref{e10.2}), (\ref{e10.3}) imply that \ $S^2$ \ is
rejected at the defining of the matrix value \ $\alpha
S_{\tau_1^2}$ \ by definition~\ref{8.3.}~2).
\\
But it can happen only if \ $S^2$ \ on \ $\alpha$ \ is
\textit{suppressed} for \ $\gamma_{\tau_1^2}$, \ that implies \
$SIN_n(\gamma_{\tau_1^2})$ \ contrary to (\ref{e10.1}).
\\
Hence
\[
    \exists \tau < \tau_2 ~ \forall \tau^{\prime} \in [\tau,\tau_2[
    ~~ a_{\tau^{\prime}}=1,
\]
and from lemma \ref{9.4.} it follows \ $a_{\tau} \equiv 1$ \ on \
$]\tau_1, \tau_2[\;$.
\\
\hspace*{\fill} $\dashv$
\\

\noindent From here and theorem~1 it comes directly

\begin{corollary}
\label{10.2.} \hfill {}
\\
Let \medskip

(i) \ $A_{2}^{\vartriangleleft \alpha_1}(\tau _{1},\tau _{2},\tau
_{3})$;
\newline

(ii) \ $\left] \gamma _{\tau _{1}}^{<\alpha _{1}}, \gamma _{\tau
_{2}}^{<\alpha _{1}}\right] \cap SIN_{n}^{<\alpha _{1}} \neq
\varnothing $.
\medskip

\noindent Then \medskip

1) \ $\gamma _{\tau _{2}}^{<\alpha _{1}}$ \ is the successor of \ $\gamma
_{\tau _{1}}^{<\alpha _{1}}$ \ in \ $SIN_{n}^{<\alpha _{1}}$; \medskip

2) \ $a_{\tau }^{<\alpha _{1}}\equiv 1$ on \ $\left] \tau
_{1},\tau _{2} \right[ $, \ $a_{\tau _{2}}^{<\alpha _{1}}=0$ \quad
and \medskip

3) \ $\widetilde{\delta }_{\tau _{2}}^{<\alpha _{1}} \le \gamma
_{\tau _{1}}^{<\alpha _{1}}$.  \label{c14}
\endnote{
\ p. \pageref{c14}. \ Again actually here \
$\widetilde{\delta}_{\tau_2}^{<\alpha_1} =
\gamma_{\tau_1}^{<\alpha_1}$.
\\
\quad \\
} %
\\
\hspace*{\fill} $\dashv$
\end{corollary}

\begin{lemma}
\label{10.3.} \hfill {}
\\
Let \medskip

(i) \ $A_{2}^{\vartriangleleft \alpha_1}(\tau _{1},\tau _{2},\tau
_{3})$; \newline

(ii) \ $a_{\tau _{2}}^{<\alpha _{1}}=1$. \medskip

\noindent Then for the prejump cardinal \ $\alpha ^{2}=\alpha
_{\tau _{2}}^{<\alpha _{1}\Downarrow }$ \ there exists an ordinal
\begin{equation*}
\tau _{\ast }=\min \{\tau \in \; ]\tau_1, \tau_2[ \; : \gamma
_{\tau }^{<\alpha _{1}}\in SIN_{n}^{<\alpha ^{2}}\}
\end{equation*}
such that \medskip

1) \ $\gamma _{\tau _{\ast }}^{<\alpha _{1}}<\widetilde{\delta }_{\tau
_{2}}^{<\alpha _{1}}$, \ $\gamma _{\tau _{\ast }}^{<\alpha _{1}}\notin
SIN_{n}^{<\alpha _{1}}$; \medskip

2) \ $a_{\tau }^{<\alpha ^{2}} \equiv a_{\tau }^{<\alpha_1} \equiv
1$ \ on \ $\left] \tau _{1},\tau _{\ast }\right[ $;
\medskip

3) \ $\alpha S_{f}^{<\alpha ^{2}}$ \ is monotone on \ $[\tau _{1},\tau
_{\ast }[$ \quad and \medskip

4) \ $Od\alpha S_{f}^{<\alpha _{1}}(\tau _{1},\tau _{\ast })>Od(\alpha
S_{\tau _{2}}^{<\alpha _{1}})$. \medskip
\end{lemma}

\noindent \textit{Proof.} \ The upper indices \ $< \alpha_{1} $,
$\vartriangleleft \alpha_{1} $ \ will be omitted. First one should
see that \ $\gamma_{\tau_{1}} < \widetilde{\delta}_{\tau_{2}} $, \
otherwise $(ii)$ and lemma~\ref{9.2.} break $(i)$.

\noindent Then from lemma~\ref{10.1.} and (\textit{i}),
(\textit{ii}) it comes the existence of the following ordinal
below the prejump cardinal \ $\alpha^2 =
\alpha_{\tau_2}^{\Downarrow}$:

\vspace{0pt}
\begin{equation*}
    \tau_{\ast} = \min \bigl \{ \tau >\tau _{1}:  \gamma _{\tau }\in
    \bigl ( SIN_{n}^{<\alpha^2}  - SIN_{n} \bigr ) \bigr \} .
\end{equation*}
\vspace{0pt}

\noindent On \ $[ \tau _{1},\tau_{\ast} [ $\quad $\alpha S_{f}$ \
is monotone and by theorem~1 (for \ $\alpha^2 $ \ instead of \
$\alpha_{1} $) \ $\gamma _{\tau_{\ast}}$ \ is the successor of \
$\gamma _{\tau _{1}}$ \ in \ $SIN_{n}^{<\alpha^2}$. \ Due to
lemma~\ref{9.3.} the function \ $\alpha S_{f} $ \ stabilizes on \
$[ \tau _{1},\tau_{\ast} [ $, \ so that for some \ $\tau _{0}\in
\left] \tau _{1},\tau_{\ast} \right[$, \quad $S^{0}$ \ the
proposition \ $\forall \tau \geq \tau _{0} \ \alpha S_{\tau
}=S^{0}$ \ is true below \ $ \gamma _{\tau_{\ast}}$. \ Hence,
below \ $\gamma_{\tau_{\ast}}$ \ the weaker proposition is true:
\begin{equation*}
    \forall \tau \Bigl( \tau_0 < \tau ~ \rightarrow \exists S  ( S = \alpha
    S_{\tau} \wedge S \underline{\gtrdot} S^0 )\Bigr).
\end{equation*}
It can be formulated in the \ $\Pi_n$-form, just as it was done in
the proof of lemma~\ref{9.3.}, by means of proposition
(\ref{e9.4}), where \ $\tau_{1,3}^n$, \ $\alpha S_{\tau_2}^1$ \
should be replaced with \ $\tau_0$, \ $S^0$ \ respectively.
\newline
The cardinal \ $\gamma _{\tau_{\ast}}$ \ extends this proposition
up to \ $ \alpha^2$ \ and therefore by $(i)$
\begin{equation*}
\rho =Od\alpha S_{f} ( \tau _{1},\tau_{\ast} ) >  Od ( \alpha S_{\tau _{2}}
) .
\end{equation*}
\vspace{0pt}

\noindent Next, let us discuss the unit characteristic. From lemma
3.2~\cite{Kiselev9} and (\textit{i}), (\textit{ii}) it comes that
there exist some admissible carriers of \ $\alpha $-matrices of
unit characteristic disposed cofinally to \ $\gamma _{\tau _{1}}$
\ as it was several times above. Hence such carriers must be also
in \ $] \gamma_{\tau_{\ast}}, \alpha^2 [\;$, \ otherwise \ $\gamma
_{\tau _{\ast}}$ \ would be defined below \ $\alpha^2$ \ along
with the ordinal \ $\rho $ \ and then by lemma 4.6~\cite{Kiselev9}
\[
    \rho <Od (\alpha S_{\tau _{2}} ).
\]
It remains to apply lemmas 3.2~\cite{Kiselev9}, \ref{9.4.} (where
\ $\tau_{\ast}$, \ $ \alpha^2$ \ play the role of \ $\tau_2$, \
$\alpha_1$ \ respectively), since \ $\alpha S_{f}^{<\alpha^1}$, \
$a_{f}^{<\alpha^1}$ \ coincide with \ $ \alpha S_{f}^{<\alpha^2}$,
\ $a_{f}^{<\alpha^2}$ \ on \ $[\tau_1, \tau_\ast [ $ \ due to
lemma~\ref{8.7.} about absoluteness.
\\
\hspace*{\fill} $\dashv$
\\
\quad \\
With the help of reasoning analogous to the proofs of
lemmas~\ref{10.1.}-\ref {10.3.} it is not hard to obtain

\begin{lemma}
\label{10.4.} \hfill {}
\\
Let \medskip

(i) \ $A_{2}^{\vartriangleleft \alpha_1}(\tau _{1},\tau _{2},\tau
_{3})$;
\medskip

(ii) \ $\forall \gamma <\gamma _{\tau _{2}}^{<\alpha _{1}}\exists \tau
~(\gamma <\gamma _{\tau }^{<\alpha ^{2}}\wedge a_{\tau }^{<\alpha ^{2}}=1)$
\ for \ $\alpha ^{2}=\alpha _{\tau _{2}}^{<\alpha _{1}\Downarrow }$ ;
\medskip

(iii) \ $a_{\tau _{2}}^{<\alpha _{1}}=0$. \medskip

\noindent Then \medskip

1) \ $\widetilde{\delta }_{\tau _{2}}^{<\alpha _{1}} \le \gamma
_{\tau _{1}}^{<\alpha _{1}}$ \quad and \medskip

2) \ $\exists \tau \in \left] \tau _{1},\tau _{2}\right[ \quad
(a_{\tau }^{<\alpha _{1}}=1\wedge \alpha S_{\tau }^{<\alpha
_{1}}\gtrdot \alpha S_{\tau _{2}}^{<\alpha _{1}})$. \label{c15}
\endnote{
\ p. \pageref{c15}. \ In fact here \ $a_{\tau}^{<\alpha_1} \equiv
1$ \ on \ $]\tau_1, \tau_2[$\;; \ and again \
$\widetilde{\delta}_{\tau_2}^{<\alpha_1} =
\gamma_{\tau_1}^{<\alpha_1}$.
\\
\quad \\
} %
\end{lemma}

\noindent \textit{Proof.} \ The upper indices \ $< \alpha_{1} $,
$\vartriangleleft \alpha_{1}$ \ will be omitted. Let us reveal the
situation below, standing as usual on \ $\alpha^{2} =
\alpha_{\tau_{2}}^{\Downarrow} $. \ Suppose that \
$\widetilde{\delta}_{\tau_{2}} \in ]\gamma_{\tau_{1}},
\gamma_{\tau_{2}} [ $; \ here the following two cases should be
considered:
\\
1. \quad $[\gamma_{\tau_1}, \gamma_{\tau_2} [ \; \cap \;
SIN_n^{<\alpha^2} \subseteq SIN_n$, \ then again (just as it was
in the proof of theorem 1\;) \ $ \alpha S_{f} $ \ is monotone on
the interval \ $[\tau_1, \tau_2[$ \ such that the interval \
$]\gamma_{\tau_1}, \gamma_{\tau_2}[$ \ contains \ $ SIN_{n}
$-cardinals, contrary to theorem 1.
\\
2. \quad $[\gamma_{\tau_1}, \gamma_{\tau_2} [ \; \cap \;
SIN_n^{<\alpha^2} \nsubseteq SIN_n$. \ In this case one should
again apply the restriction-and-extension techniques of reasoning
precisely as it was done in part~1a. of lemma~\ref{9.4.} proof.
Let us first repeat the argument from the proof of lemma 10.3
concerning the function \ $\alpha S_{f} $ \ defined on the
interval \ $]\tau_{1}, \tau_\ast [\; $, \ where
\begin{equation*}
\gamma_{\tau_\ast} = \min \left ( SIN_{n}^{ < \alpha^{2}} - SIN_{n} \right )
\mathrm{\quad and \quad}  a_{\tau} \equiv 1 \mathrm{~~on~~}
]\tau_{1},\tau_\ast [
\end{equation*}
due to $(ii)$. By theorem 1 (for \ $\tau_\ast$, $\alpha^2$, as
$\tau_2$, $\alpha_1$) and $(i)$ the cardinal \
$\gamma_{\tau_\ast}$ \ is the successor of \ $\gamma_{\tau_1}$ \
in \ $SIN_n^{<\alpha^2}$. Hence, the proposition
\begin{equation*}
\forall \tau > \tau_1 \quad a_{\tau} = 1
\end{equation*}
is true below \ $\gamma_{\tau_\ast}$; \ it is not hard to see it
with the help of condition $(ii)$ and lemmas~3.2~\cite{Kiselev9},
\ref{9.4.}. This proposition can be formulated in the \
$\Pi_n$-form for this case, just as it was done above in the proof
of lemma~\ref{9.4.} in part~1a. with \ $\tau^3$ \ as \ $\tau_1$:

\vspace{0pt}
\begin{equation*}
\forall \gamma \Bigl ( \gamma_{\tau_1} < \gamma  \; \wedge \; SIN_{n-1}
(\gamma) \longrightarrow \qquad \qquad \qquad  \qquad \qquad \qquad \qquad
\qquad \qquad
\end{equation*}
\begin{equation*}
\qquad \qquad \longrightarrow \exists \delta, \alpha,  \rho, S ~
\bigl (  SIN_n^{<\alpha^{\Downarrow}}(\gamma_{\tau_1}) \wedge
\alpha \mathbf{K}_{n+1}^{\exists}(1, \delta, \gamma, \alpha,\rho,
S)  \bigr ) \Bigr ).
\end{equation*}
\vspace{0pt}

\noindent After that the \ $SIN_n^{<\alpha^2} $-cardinal \
$\gamma_{\tau_\ast} $ \ extends this proposition up to \
$\alpha^{2} $ \ and so \ $a_{\tau_{2}}=1 $ \ contrary to $(iii)$.
\\
Thus, \ $\widetilde{\delta}_{\tau_{2}} \leq \gamma_{\tau_{1}}$; \
coming to the end of the proof one should apply lemma~9.2. Suppose
that
\begin{equation*}
\exists \tau \in \; ]\tau_{1},\tau_{2} [ ~~ \forall \tau^{\prime}
\in [\tau ,\tau_{2} [ \ \ \ a_{\tau^{\prime}} = 0,
\end{equation*}
then by this lemma the monotonicity of \ $\alpha S_{f} $ \ on \
$]\tau_{1},\tau_{2}[$ \ implies that
\[
    Od \alpha S_{f} (\tau_{1},\tau_{2}) \leq Od (\alpha S_{\tau_2})
\]
contrary to $(i)$. This contradiction along with $(i)$ provides 2)
and ends the proof.
\\
\hspace*{\fill} $\dashv$
\\

The following lemma will be used at the end of the Main theorem
proof, again relying strongly on the formula \ $A^{0
\vartriangleleft \alpha_1}(\tau)$ \ (recall definition
\ref{8.1.}~3.2\;): \vspace{-6pt}
\begin{multline*}
    \exists \gamma < \alpha_1 \Big( \gamma =
    \gamma_{\tau}^{<\alpha_1} \wedge \neg \exists a, \delta,
    \alpha, \rho < \alpha_1 \exists S \vartriangleleft \rho ~ \big(
    \mathbf{K}_n^{\forall <\alpha_1}(\gamma,\alpha_{\chi}^{\Downarrow})
    \wedge
\\
    \wedge \alpha \mathbf{K}_{n+1}^{\exists \vartriangleleft
    \alpha_1}(a, \delta, \gamma, \alpha, \rho, S) \big) \Big)
\end{multline*}
meaning, that there is no \ $\alpha$-matrix on some carrier \
$\alpha > \gamma_{\tau}^{<\alpha_1}$ \ admissible for \
$\gamma_{\tau}^{<\alpha_1}$ \ below \ $\alpha_1$.
\\
Accordingly, through \ $A_1^{0 \vartriangleleft \alpha_1}(\tau_1,
\tau_2, \alpha S_f^{<\alpha_1})$ \ is denoted the formula (remind
definition \ref{8.1.}~1.1 for \ $X_1 = \alpha S_f^{<\alpha_1}$):
\[
    A^{0 \vartriangleleft \alpha_1}(\tau_1) \wedge
    A_1^{\vartriangleleft \alpha_1}(\tau_1, \tau_2, \alpha
    S_f^{<\alpha_1}),
\]
and also should be used the formula $A_2^{0 \vartriangleleft
\alpha_1}(\tau_1, \tau_2^{\prime}, \tau_3, \alpha
S_f^{\vartriangleleft \alpha_1})$ \ (recall definition
\ref{8.2.}~3.3\;):
\[
    A^{0 \vartriangleleft \alpha_1}(\tau_1) \wedge
    A_2^{\vartriangleleft \alpha_1}(\tau_1, \tau_2^{\prime}, \tau_3, \alpha
    S_f^{<\alpha_1}).
\]

\begin{lemma}
\label{10.5.} \hfill {}
\\
Let \medskip

(i)\quad $A_{1}^{0 \vartriangleleft \alpha_1}(\tau _{1},\tau _{2},
\alpha S_f^{<\alpha_1})$;
\medskip

(ii)\quad $\tau_2 \le \tau_3$ \ and \ $S^3$ \ be a matrix of
characteristic \ $a^3$ \ on a carrier
\[
    \alpha_3 \in \; ] \gamma_{\tau_3}^{<\alpha_1}, \alpha_1 [
\]
preserving \ $SIN_n^{<\alpha_1}$-cardinals \ $\le
\gamma_{\tau_2}^{<\alpha_1}$ \ below \ $\alpha_1$ \ and with
generating eigendisseminator \ $\check{\delta}^{S^3}$;
\medskip

(iii)\quad $\check{\delta}^{S^3} \leq \gamma _{\tau _{1}}^{<\alpha
_{1}}$. \medskip

\noindent Then \quad $a^3=0$.
\\
Analogously for any disseminator \ $\widetilde{\delta}$ \ of \
$S^3$ \ on \ $\alpha^3$ \ with any base \ $\rho \ge \rho^{S^3}$.
\end{lemma}

\noindent \textit{Proof.} \ As usual, we shall look over the
situation below, standing on \ $a^3 = \alpha_3^{\Downarrow}$ \ and
considering the disseminator \ $\check{ \delta}^{S^3} $ \ with the
data base \ $\rho^3 = \rho^{S^3} = \widehat{\rho_1}$, $\rho_1 =
Od(S^3)$; \ the upper indices \ $<\alpha _{1}$, $\vartriangleleft
\alpha_1$ \ will be dropped for some convenience.
\\
Suppose that this lemma fails and \ $a^3 = 1$, \ then
$\check{\delta}^{S^3}$ \ is admissible and nonsuppressed
disseminator of \ $S^3$ \ on \ $\alpha_3$ \ for any \
$\gamma_{\tau} \in \; ]\tau_1, \tau_3[$ \ and by lemma~\ref{9.2.}
\ $a_{\tau } \equiv 1$ \ on \ $ \left] \tau _{1}, \tau _{2}\right[
$ \ and
\begin{equation*}
Od \alpha S_{f}(\tau_{1},\tau_{2})\leq Od(\alpha  S_{\tau_{3}}).
\end{equation*}
From lemma~\ref{9.5.} it follows that due to $(iii)$ \ $\alpha
S_{f}$ \ is nonmonotone on \ $[ \tau _{1},\tau _{2} [ $ \ (remind
case~2. in the proof of point 3. of lemma~\ref{9.5.}\;) and that
is why there exists \ $\tau _{2}^{\prime }$ \ for which there
holds
\begin{equation}  \label{e10.4}
A_2^0 ( \tau _{1},\tau _{2}^{\prime }, \tau_{2} ), \quad a_{\tau
_{2}^{\prime }}=1
\end{equation}
Now one should repeat precisely the reasoning from part~1b. of
lemma~\ref{9.4.} proof. Below \ $\alpha^3$ \ the following \
$\Sigma _{n+1}$-proposition is true by lemma~\ref{8.7.} about
absoluteness (remind (\ref{e9.8})\;):
\[
    \exists \gamma^0 \exists \tau _{1}^{\prime },  \tau _{2}^{\prime
    \prime },\tau _{3}^{\prime } < \gamma^0 ~  \Bigl( SIN_{n} ( \gamma^0 ) \wedge
    \gamma_{\tau_1^{\prime}}^{<\gamma^0} <
    \gamma_{\tau_2^{\prime\prime}}^{<\gamma^0} <
    \gamma_{\tau_3^{\prime}}^{<\gamma^0} < \gamma^0 \wedge
    \qquad \qquad
\]
\begin{equation}  \label{e10.5}
    \qquad
    \wedge A_2^{0 \vartriangleleft \gamma^0}
    ( \tau _{1}^{\prime },\tau _{2}^{\prime \prime },
    \tau_{3}^{\prime },\alpha S_{f}^{<\gamma^0 } )
    \wedge \forall \tau^{\prime\prime\prime} \in \;
    ]\tau_1^{\prime}, \tau_2^{\prime\prime}] ~
    a_{\tau^{\prime\prime\prime}}^{<\gamma^0}  = 1 \wedge
\end{equation}
\[
    \qquad \qquad \qquad \qquad \qquad \qquad \qquad \qquad \wedge
    \alpha S_{\tau _{2}^{\prime \prime }}^{<\gamma^0 } =  \alpha
    S_{\tau_{2}^{\prime }} \Bigr).
\]
\vspace{0pt}

\noindent It contains the constants \ $< \rho^3 $ \ and \ $\alpha
S_{\tau _{2}^{\prime }} \vartriangleleft \rho^3 $, \ therefore the
disseminator \ $\check{\delta}^{S^3}$ \ restricts this proposition
and it fulfills below \ $\check{ \delta}^{S^3}$.
\\
Now let us reveal the situation below the prejump cardinal
\\
$\alpha^2 = \alpha_{\tau_2^\prime}^\Downarrow$.
\\
By $(i)$ \ $\gamma_{\tau_1} \in SIN_n$, \ therefore by
lemma~\ref{8.5.}~1) \ $\gamma_{\tau_1} \in SIN_n^{<\alpha^3}$. \
Since \ $\check{\delta}^{S^3} \le \gamma_{\tau_1}$ \ and \
$\check{\delta}^{S^3} \in SIN_n^{<\alpha^3}$, \ lemma
3.8~\cite{Kiselev9} (for \ $\alpha^3$, \ $\gamma_{\tau_1}$ \ as \
$\alpha_1$, \ $\alpha_2$) \ implies \ $\check{\delta}^{S^3} =
\gamma_{\tau_1}$ \ or \ $\check{\delta}^{S^3} \in
SIN_n^{<\gamma_{\tau_1}}$; \ then by the same lemma (for \
$\gamma_{\tau_1}$ \ as \ $\alpha_2$) \ $\check{\delta}^{S^3} \in
SIN_n$.
\\
From here and lemma~\ref{8.5.}~1) it comes out \
$\check{\delta}^{S^3} \in SIN_n^{<\alpha^2}$; \ hence, in
proposition (\ref{e10.5}) we can replace \ $\gamma^0 $ \ with \
$\alpha^{2}$ \ by lemma~\ref{8.7.} about absoluteness and due to
lemma \ref{8.5.}~5) (for \ $\alpha_{\tau_2^\prime}$ \ as \
$\alpha$) \ it comes out \ $a_{\tau _{2}^{\prime }}=0$ \ contrary
to~(\ref{e10.4}).
\\
\hspace*{\fill} $\dashv$
\\

Now the special theory of matrix functions is developed enough to
start the proof of the Main theorem.

\newpage

\chapter{Applications of Special Theory}

\setcounter{section}{10}

\section{Proof of Main Theorem}

\setcounter{equation}{0}

\hspace*{1em} The contradiction, which proves the Main theorem, is
the following:
\\
On one hand, by lemma 8.9 the function \ $\alpha S_f^{<\alpha_1}$
\ is defined on the nonempty set
\begin{equation*}
T^{\alpha_1} = \{\tau:\alpha\delta^{\ast}<\gamma_{\tau}\ <
\alpha_1 \}
\end{equation*}
for every sufficiently great cardinal \ $\alpha_1 \in SIN_n$.
\\
Its monotonicity on this set is excluded by theorem 1.
\\
But on the other hand, this monotonicity is ensured by the
following theorem for every \ $SIN_n$-cardinal \
$\alpha_1>\alpha\delta^{\ast}$ \ of sufficiently great cofinality.
Remind, that bounding cardinals \ $\alpha_1$ \ are always assumed
to be equiformative with \ $\chi^{\ast}$, \ that is there holds \
$A_6^e(\chi^{\ast}, \alpha_1)$ \ (recall definition \ref{8.1.}~5.1
for \ $\chi=\chi^{\ast}$, $\alpha^0=\alpha_1$).
\\

\noindent \textbf{Theorem 2.}\quad \newline \emph{\hspace*{1em}
Let the function \ $\alpha S_f^{<\alpha_1}$ \ be defined on
nonempty set
\begin{equation*}
T^{\alpha_1} =\{\tau:\gamma_{\tau_1}^{<\alpha_1} <
\gamma_{\tau}^{<\alpha_1} < \alpha_1 \}
\end{equation*}
such that \ $\alpha_1 < k $ \ and: \medskip }

\emph{(i) \ $\tau_1=\min\{\tau:\forall\tau^{\prime}
(\gamma_{\tau}^{< \alpha_1} <
\gamma_{\tau^{\prime}}^{<\alpha_1}\longrightarrow\tau^{\prime}\in
dom(\alpha S_f^{<\alpha_1}))\};$ \newline }

\emph{(ii) \ $\sup SIN_n^{<\alpha_1}=\alpha_1;$ \newline
}

\emph{(iii) \ $cf(\alpha_1)\geq\chi^{\ast +}.$ \newline
}

\emph{\noindent Then \ $\alpha S_f^{<\alpha_1}$ \ is monotone on
this set:
} %

\[
    \forall \tau_1, \tau_2 \in T^{\alpha_1} \bigl( \tau_1 <
    \tau_2 \rightarrow \alpha S_{\tau_1}^{<\alpha_1}
    \underline{\lessdot} \alpha S_{\tau_2}^{<\alpha_1} \bigr).
\]

\quad \\
\noindent \textit{Proof.} \ The scheme of the proof is the
following in outline.
\\
The reasoning will be carried out by the induction on the cardinal
\ $\alpha_1$.

Let us suppose, that this theorem fails and the cardinal \
$\alpha_1^\ast$ \ is \textit{minimal} breaking this theorem, that
is the function \ $\alpha S^{<\alpha_1^\ast}_f$ \ is nonmonotone
on \ the set
\begin{equation*}
T^{\alpha_1^\ast} = \bigl \{ \tau:
\gamma_{\tau_1^\ast}^{<\alpha_1^\ast} <
\gamma_{\tau}^{<\alpha_1^\ast} < \alpha_1^{\ast} \bigr \}
\end{equation*}
with specified properties $(i)$--$(iii)$ for some \
$\tau_1^{\ast}$, \ so that the \textit{first inductive hypothesis}
holds:
\\
for every \ $\alpha_1 < \alpha_1^{\ast}$ \ the function \ $\alpha
S_f^{<\alpha_1}$ \ is monotone on the set \ $T^{\alpha_1}$ \ with
properties $(i)$--$(iii)$.
\\
It follows straight from theorem 1, that this \ $\alpha_1^{\ast}$
\ is simply the \textit{minimal} cardinal \ $\alpha_1$, \ for
which the set \ $T^{\alpha_1}$ \ exist, because for every such \
$\alpha_1 < \alpha_1^{\ast}$ \ the function \ $\alpha
S_f^{<\alpha_1}$ \ on \ $T^{\alpha_1}$ \ is \textit{nonmonotone}
by theorem~1 and at the same time is \textit{monotone} by the
minimality of \ $\alpha_1^{\ast}$.

The reasoning will be conducted below \ $\alpha_1^\ast$ \ (and all
variables will be bounded by \ $\alpha_1^{\ast}$), \ or below
bounding cardinals \ $\alpha_1 \le \alpha_1^{\ast}$, \ so the
upper indices \ $< \alpha_1^\ast$, \ $\vartriangleleft
\alpha_1^\ast$ \ will be omitted for some shortness up to the end
of theorem 2 proof.

First, let us notice that in conditions of this theorem~2 there
holds
\begin{equation*}
\gamma_{\tau_1}^{<\alpha_1}\in SIN_n^{<\alpha_1};
\end{equation*}
to see it one should repeat once more the argument applied before
several times (first in the proofs of lemmas 7.7, 8.10\;).
Therefore it is not hard to see that for \textit{every}
sufficiently great \ $ \tau_3^\ast \in T^{\alpha_1^\ast}$ \ the
interval \ $[\gamma_{\tau_1^\ast},\ \gamma_{\tau_3^\ast}[$ \ can
be considered \textit{as the block}, that is there exist some
ordinals \ $\tau_1^{\ast\prime},\ \tau_2^\ast,\ \eta^{\ast 3}$ \
which fulfill the statement (remind definition \ref{8.1.}~1.6 for
\ $X_1 = \alpha S_f$, $X_2 = a_f$):

\begin{equation*}
A_4^b (\tau_1^\ast,\ \tau_1^{\ast\prime},\ \tau_2^\ast,\
\tau_3^\ast,\ \eta^{\ast 3}, \alpha S_f, a_f).
\end{equation*}

\noindent Here (due to this definition~\ref{8.1.}\;) \
$\tau^{\ast\prime}_1$ \ is the index of the matrix \ $\alpha
S_{\tau_1^{\ast \prime}}$ \ of \textit{unit characteristic} \
$a_{\tau^{\ast\prime}_1}=1$ \ on its carrier \ $
\alpha_{\tau^{\ast\prime}_1}$ \ and \ $\eta^{\ast 3}$ \ --- \ the
type of this interval.
\\
Next, due to condition \ ($iii$) \ of this theorem 2 we can use
the index \ $ \tau_3^\ast \in T^{\alpha_1^\ast}$ \ such, that the
interval \ $[\gamma_{\tau_1^\ast}$, \ $\gamma_{\tau_3^\ast}[$ \
has just the type

\begin{equation*}
\eta^{\ast 3} > Od(\alpha S_{\tau_1^{\ast\prime}}), \quad
\eta^{\ast 3} < \chi^{\ast +}.
\end{equation*}

\noindent Now the formula \ $\mathbf{K}^0$ \ starts to work and
\textit{closes the diagonal reasoning}:
\\
There arises the matrix \ $\alpha S_{\tau_3^\ast}$ \ on the
carrier \ $\alpha_{\tau_3^\ast}$ \ along with its disseminator \
$\widetilde{\delta} ^{\ast 3}=\widetilde{\delta}_{\tau_3^\ast}$ \
and data base \ $\rho^{\ast 3}=\rho_{\tau_3^\ast}$, \ and we shall
see, that by lemma~\ref{10.5.} it has \textit{zero} characteristic
on this carrier.
\\
Standing on the prejump cardinal \ $\alpha^{\ast
3}=\alpha_{\tau_3^\ast}^{\Downarrow}$ \ one should reveal the
following situation below \ $\alpha^{\ast 3}$:
\\
One shall see soon, that by lemma \ref{8.8.} disseminator \
$\widetilde{\delta}^{\ast 3}$ \ falls into some maximal block \
$[\gamma_{\tau_1^\ast},\ \gamma^{\ast 3}[$ \ below \ $\alpha^{\ast
3}$ \ of a type \ $\eta^{\ast 3 \prime}< \chi^{\ast +}$, \ where \
$\gamma^{\ast 3}$ \ is some \ $
\gamma_{\tau^{\ast\prime}_3}^{<\alpha^{\ast 3}}$. \ It is not hard
to see, that
\begin{equation*}
\gamma_{\tau_3^\ast} \leq \gamma^{\ast 3}\ \ \wedge\ \  \eta^{\ast
3} \leq \eta^{\ast 3 \prime};
\end{equation*}
so, there holds
\\
\begin{equation*}
A^{M b \vartriangleleft \alpha^{\ast 3} } _4(\tau_1^\ast,
\tau_1^{\ast \prime}, \tau_2^\ast, \tau_3^{\ast\prime}, \eta^{\ast
3 \prime}, \alpha S_f^{<\alpha^{\ast 3}}, a_f^{<\alpha^{\ast 3}}).
\end{equation*}
\newline
All these facts together constitute the premise of lemma 8.5~6):
\\

\[
    a_{\tau_3^\ast}=0\wedge\gamma_{\tau_1^\ast}^{<\alpha^{\ast 3}}
    \leq \widetilde{\delta}^{\ast 3}<
    \gamma_{\tau^{\ast\prime}_3}^{<\alpha^{\ast3}} \wedge
    \qquad\qquad\qquad\qquad
\]
\[
    \qquad\qquad\qquad\qquad
    \wedge A_4^{M b \vartriangleleft \alpha^{\ast 3}} (\tau_1^\ast,
    \tau_1^{\ast \prime}, \tau_2^\ast, \tau_3^{\ast\prime}, \eta^{\ast
    3 \prime}, \alpha S_f^{<\alpha^{\ast 3}}, a_f^{<\alpha^{\ast 3}}).
\]

\noindent Therefore this lemma implies

\begin{equation*}
\eta^{\ast 3 \prime}<\rho^{\ast 3}\vee\rho^{\ast 3} =\chi^{\ast
+};
\end{equation*}

\noindent thus, at any rate, \newline
\begin{equation*}
Od(\alpha S_{\tau_1^{\ast\prime}} ) < \eta^{\ast 3} \leq
\eta^{\ast 3 \prime} < \rho^{\ast 3}.
\end{equation*}

\noindent But we shall see soon, that it is impossible, because by
lemma~\ref{9.5.} (about stairway cut-off from above) and lemma
\ref{11.3.} below there holds:
\begin{equation*}
\rho^{\ast 3} \leq Od(\alpha S_{\tau_1^{\ast\prime}} ).
\end{equation*}
\vspace{0pt}

\noindent This contradiction will end the proof of theorem~2.
\quad \\
\\

\noindent To realize this scheme some more information is needed.

The reasoning sketched above relies on the following easy
auxiliary lemmas 11.1, 11.3, which are coming as its carrying
construction and are describing some important properties of zero
matrices behavior; they were not presented earlier because of
their rather special character. For this purpose one should remind
the formula (recall definition \ref{8.1.}~1.1 for \ $X_1 = \alpha
S_f^{<\alpha_1}$)
\[
    A_1^{\vartriangleleft \alpha_1} (\tau_1, \tau_2, \alpha
    S_f^{<\alpha_1}):
    \qquad \qquad \qquad \qquad \qquad \qquad \qquad \qquad
\]
\[
    \tau_1+1 < \tau_2 \wedge \tau_1 = \min \bigl\{ \tau:
    \; ]\tau, \tau_2[ \; \subseteq dom(\alpha S_f^{<\alpha_1} )\bigr\}
\]
\[
    \qquad \qquad \qquad \qquad \qquad \wedge
    \gamma_{\tau_1}^{<\alpha_1} \in SIN_n^{<\alpha_1} \wedge
    \gamma_{\tau_2}^{<\alpha_1} \in SIN_n^{<\alpha_1};
\]
\vspace{0pt}

Remind also that we often omit the functions \ $\alpha
S_f^{<\alpha_1}$, $a_f^{<\alpha_1}$ \ in notations of formulas
below \ $\alpha_1$; \ recall that the type of interval \ $[
\gamma_{\tau_1}^{<\alpha_1}, \gamma_{\tau_2}^{<\alpha_1}[$ \ below
\ $\alpha_1$ \ is the order type of the set (see definition
\ref{8.1.}~1.3\;):
\[
    \bigl \{ \gamma: \gamma_{\tau_1}^{<\alpha_1} < \gamma <
    \gamma_{\tau_2}^{<\alpha_1} \wedge SIN_n^{<\alpha_1}(\gamma)
    \bigr \}.
\]

\noindent Beforehand one should delay on the following auxiliary
arguments, suitable to shorten the succeeding reasonings; with
this aim one should introduce the following notion:
\\
an interval \ $[\tau_1, \tau_2[$ \ and the corresponding interval
\ $[\gamma_{\tau_1}^{<\alpha_1}, \gamma_{\tau_2}^{<\alpha_1}[$ \
will be called the intervals of matrix admissibility, or
\textit{admissibility intervals}, below \ $\alpha_1$, \ if for
every \ $\tau^{\prime} \in \; ]\tau_1,\tau_2[$ \ there exist some
\ $\alpha$-matrix \ $S$ \ on some carrier \ $>
\gamma_{\tau^{\prime}}^{<\alpha_1}$ \ \textit{admissible} for \
$\gamma_{\tau^{\prime}}^{<\alpha_1}$ \ below \ $\alpha_1$:
\[
    \forall \tau^{\prime} \in \; ]\tau_1, \tau_2[ ~ \exists
    a^{\prime}, \delta^{\prime}, \alpha^{\prime}, \rho^{\prime},
    S^{\prime} ~ \alpha \mathbf{K}^{<\alpha_1} (a^{\prime},
    \delta^{\prime}, \gamma_{\tau^{\prime}}^{<\alpha_1},
    \alpha^{\prime}, \rho^{\prime}, S^{\prime}),
\]
and \ \mbox{$\gamma_{\tau_1}^{<\alpha_1} \in SIN_n$},
\mbox{$\gamma_{\tau_2}^{<\alpha_1} \in SIN_n$} \ and \ $\tau_1$ \
is the minimal ordinal with these attributes.
\\
Next, the following properties of any zero matrix \ $S$ on its
carrier \ $\alpha$ \ admissible for \ $\gamma_{\tau}^{<\alpha_1}$
\ along with its \textit{minimal} disseminator \
$\widetilde{\delta}$ \ with base \ $\rho$ \ below \ $\alpha_1$ \
should be treated for \ $\alpha_1 \le \alpha_1^{\ast}$:

(1a.) if \ $\gamma_{\tau_1}^{<\alpha_1} <
\gamma_{\tau_2}^{<\alpha_1} \le  \gamma_{\tau}^{<\alpha_1}$ \ and
\ $\widetilde{\delta}$ \ falls in the admissibility interval \
$[\gamma_{\tau_1}^{<\alpha_1}, \gamma_{\tau_2}^{<\alpha_1}[\;$, \
that is \ $ \gamma_{\tau_1}^{<\alpha_1} \le \widetilde{\delta} <
\gamma_{\tau_2}^{<\alpha_1}$, \ then \
$\gamma_{\tau_1}^{<\alpha_1} = \widetilde{\delta}$;

(1b.) if there exist some zero matrix \ $S^1$ \ on some
\textit{another} carrier \ $\alpha^1 \neq \alpha$, \ admissible
for the same \ $\gamma_{\tau}^{<\alpha_1}$ \ along with its
minimal disseminator \ $\widetilde{\delta}^1$ \ with base \
$\rho^1$, \ then \ $S$ \ on \ $\alpha$ \ is \textit{nonsuppressed}
for \ $\gamma_{\tau}^{<\alpha_1}$ \ along with \
$\widetilde{\delta}$, $\rho$ \ below \ $\alpha_1$.

\noindent The testing of these properties will be conducted by the
induction on triples \ $(\alpha_1, \alpha, \tau)$ \ ordered
canonically as usual (with \ $\alpha_1$ \ as the first component,
\ $\alpha$ \ as the second and \ $\tau$ \ as the third).
\\
Suppose the triple $(\alpha_1^0, \alpha^0, \tau^0)$ \ is the
minimal violating (1a.) or (1b.); so the \textit{second inductive
hypothesis} is accepted:
\\
for every lesser triple $(\alpha_1, \alpha, \tau)$ there hold
(1a.) and (1b.).
\\
We shall see, that it provides contradictions; the reasoning
forthcoming will be conducted below \ $\alpha_1^0$, \ so the upper
indices \ $<\alpha_1^0$, $\vartriangleleft \alpha_1^0$ \ will be
dropped as usual (up to the special remark, if the context will
not point out to another case clearly).

1. \ Let us begin with (1a.); suppose it is wrong, that is there
exist some zero matrix \ $S^0$ \ on its carrier \ $\alpha^0 >
\gamma_{\tau^0}$, \ its minimal disseminator \
$\widetilde{\delta^0}$ \ with base \ $\rho^0$, \ all admissible
for \ $\gamma_{\tau^0}$, \ and \ $\widetilde{\delta^0}$ \ falls
into the \textit{admissibility interval} \ $[\gamma_{\tau_1^0},
\gamma_{\tau_2^0}[$, \ but
\begin{equation} \label{e11.1}
    \gamma_{\tau_1^0} < \widetilde{\delta}^0 < \gamma_{\tau_2^0}
    \le \gamma_{\tau^0}\;,
    \mbox{~~ that is ~} \widetilde{\delta}^0 = \gamma_{\tau_3^0}\;,
    ~ \tau_1^0 < \tau_3^0 < \tau_2^0.
\end{equation}
From here and lemma~3.8 \cite{Kiselev9} it follows immediately,
that
\[
    \widetilde{\delta}^0 \in SIN_n
\]
since \ $\widetilde{\delta}^0 < \gamma_{\tau_2^0}$, \
$\widetilde{\delta}^0 \in SIN_n^{<\alpha^{0 \Downarrow}}$, \
$\gamma_{\tau_2^0} \in SIN_n$.
\\
By definition of admissibility interval there exist the matrix \
$\alpha S_{\tau_3^0}$ \ on its carrier \ $\alpha_{\tau_3^0}$ \
admissible for \ $\gamma_{\tau_3^0}$ \ along with its minimal
disseminator \ $\widetilde{\delta}_{\tau_3^0}$ \ with base \
$\rho_{\tau_3^0}$ \ (all below \ $\alpha_1^0$).
\\
From the second inductive hypothesis it comes \
$\widetilde{\delta}_{\tau_3^0} = \gamma_{\tau_1^0}$; \ hence lemma
3.2~\cite{Kiselev9} implies, that for every \ \mbox{$\gamma_{\tau}
\in \; ]\gamma_{\tau_1^0}, \gamma_{\tau_3^0}[$} \ the matrix \
$\alpha S_{\tau_3^0}$ \ possesses many carriers \ $\alpha \in \;
]\gamma_{\tau}, \gamma_{\tau+1}[$ \ admissible for \
$\gamma_{\tau}$, \ which are nonsuppressed for this \
$\gamma_{\tau}$ \ due to the same inductive hypothesis and so
there holds \ $A_1^0(\tau_1^0, \tau_3^0, \alpha
S_f^{<\alpha_1^0})$:
\[
    A^0(\tau_1^0) \wedge
    A_1(\tau_1^0, \tau_3^0, \alpha S_f^{<\alpha_1^0})
\]
below \ $\alpha_1^0$. \ The same arguments work below the prejump
cardinal \ $\alpha^{0 \Downarrow}$, \ so there holds below \
$\alpha^{0 \Downarrow}$ \ as well:
\[
    A_1^{0 < \alpha^{0 \Downarrow}}(\tau_1^0, \tau_3^0, \alpha S_f^{<\alpha^{0 \Downarrow}}).
\]
It implies for the generating disseminator \ $\check{\delta}^0$ \
of \ $S^0$ \ on \ $\alpha^0$ \ with the base \ $\rho^0$:
\begin{equation} \label{e11.2}
    \check{\delta}^0 \le \gamma_{\tau_1^0},
\end{equation}
because in the opposite case \ $\check{\delta}^0$ \ falls strictly
in the admissibility interval \ $]\gamma_{\tau_1^0},
\gamma_{\tau_3^0}]$:
\begin{equation} \label{e11.3}
    \gamma_{\tau_1^0} < \check{\delta}^0 \le
    \gamma_{\tau_3^0},
\end{equation}
and then \ $\check{\delta^0}$ \ extends up to \ $\alpha^{0
\Downarrow}$ \ the \ $\Pi_{n+1}$-proposition about admissibility
of some matrices for every \ $\gamma_{\tau}^{<\alpha^{0
\Downarrow}} > \gamma_{\tau_1^0}$, \ that become even
nonsuppressed for all such \ $\gamma_{\tau}^{<\alpha^{0
\Downarrow}}$ \ by the second inductive hypothesis (all it below
$\alpha^{0 \Downarrow}$); \ hence, there arises the set \
$T^{\alpha^{0 \Downarrow}}$ \ with properties $(i)$--$(iii)$
specified in theorem~2, contrary to the first inductive hypothesis
and theorem~1, that is to the minimality of \ $\alpha_1^{\ast}$.

\noindent From this place the reasoning passes to the matter below
\ $\alpha^{0 \Downarrow}$, \ and the upper indices \ $<\alpha^{0
\Downarrow}$, $\vartriangleleft \alpha^{0 \Downarrow}$ \ will be
dropped.
\\
Below \ $\alpha^{0 \Downarrow}$ \ the function \ $\alpha S_f$ \ is
defined on the interval \ $]\tau_3^0, \tau_3^1[$ \ for \
$\gamma_{\tau_3^1}$ \ the successor of \ $\widetilde{\delta}^0$ \
in \ $SIN_n$, \ by lemma~\ref{8.7.} about absoluteness. From here
and (\ref{e11.2}) it follows
\begin{equation} \label{e11.4}
\tau_3^0 \notin dom(\alpha S_f),
\end{equation}
otherwise it again comes (\ref{e11.3}) or \ $\widetilde{\delta}^0
= \gamma_{\tau_1^0}$ \ as the result of the minimizing of \
$\widetilde{\delta}^0$ \ within \ $[\gamma_{\tau_1^0},
\gamma_{\tau_3^0}[$ \ contrary to the supposition (all it below \
$\alpha^{0 \Downarrow}$).
\\
But (\ref{e11.4}) can be carried out only when the admissible
matrix \ $\alpha S_{\tau_3^0}$ \ is \textit{suppressed} for \
$\gamma_{\tau_3^0}$, \ that is  when there holds the suppression
condition \ $A_5^{S,0}$ \ for \ $\alpha S_{\tau_3^0}$ \ on \
$\alpha_{\tau_3^0}$ \ of the characteristic \ $a_{\tau_3^0}$ \
with the base \ $\rho_{\tau_3^0}$ (see definition
\ref{8.1.}~2.6\;) below \ $\alpha^{0 \Downarrow}$, --- and now all
boundings should be pointed out clearly:
\[
    a_{\tau_3^0} = 0 \wedge SIN_n^{<\alpha^{0 \Downarrow}}(\gamma_{\tau_3^0}) \wedge
    \rho_{\tau_3^0} < \chi^{\ast +} \wedge \sigma(\chi^{\ast},
    \alpha_{\tau_3^0}, S_{\tau_3^0}) \wedge \qquad
\]
\[
    \wedge \exists \eta^{\ast} < \gamma_{\tau_3^0} \Big(
    A_{5.4}^{sc \vartriangleleft \alpha^{0 \Downarrow}}(\gamma_{\tau_3^0},
    \eta^{\ast}, \alpha S_f^{< \alpha^{0 \Downarrow}}|\tau_3^0,
    a_f^{< \alpha^{0 \Downarrow}}|\tau_3^0) \wedge \qquad
\]
\begin{equation} \label{e11.5}
    \wedge \forall \tau^{\prime} \big( \gamma_{\tau_3^0} <
    \gamma_{\tau^{\prime}}^{< \alpha^{0 \Downarrow}} \wedge
    SIN_n^{<\alpha^{0 \Downarrow}} (\gamma_{\tau^{\prime}}^{< \alpha^{0 \Downarrow}})
    \rightarrow \qquad \qquad
\end{equation}
\[
    \rightarrow \exists \alpha^{\prime}, S^{\prime} \big[
    \gamma_{\tau^{\prime}}^{< \alpha^{0 \Downarrow}} < \alpha^{\prime} <
    \gamma_{\tau^{\prime}+1}^{< \alpha^{0 \Downarrow}} \wedge
    SIN_n^{<\alpha^{\prime \Downarrow}} (\gamma_{\tau^{\prime}}^{< \alpha^{0 \Downarrow}} )
    \wedge \sigma(\chi^{\ast}, \alpha^{\prime}, S^{\prime}) \wedge
    \qquad \qquad
\]
\[
    \qquad \qquad \qquad \qquad \wedge
    A_{5.5}^{sc \vartriangleleft \alpha^{0 \Downarrow}} (\gamma_{\tau_3^0}, \eta^{\ast},
    \alpha^{\prime \Downarrow},
    \alpha S_f^{< \alpha^{\prime \Downarrow}},
    a_f^{< \alpha^{\prime \Downarrow}}) \big] \big) \Big).
\]

\noindent Hence, there exist the cardinals
\[
    \gamma^m < \gamma^{\ast} \le \gamma_{\tau_1^0} <
    \gamma_{\tau_3^0} \mbox{~ and the limit type ~} \eta^{\ast}
\]
which carry out all its constituents \ $A_{5.1}^{sc} -
A_{5.5}^{sc}$ \ below \ $\alpha^{0 \Downarrow}$ \ (see definition
\ref{8.1.}~2.1--2.5\;); in particular the interval \
$[\gamma_{\tau_1^0}, \gamma_{\tau_3^0}[$ \ is the block of the
type \ $\eta^{\ast}$ \ due to the condition
\[
    A_4^{b < \alpha^{0 \Downarrow}}(\tau_1^0, \tau_3^0,
    \eta^{\ast}, \alpha S_f^{< \alpha^{0 \Downarrow}} | \tau_3^0,
    a_f^{< \alpha^{0 \Downarrow}} | \tau_3^0)
\]
from the condition \ $A_{5.4}^{sc}$ \ (see definition
\ref{8.1.}~2.4,~2.3\;). Moreover, there exist some its succeeding
maximal block
\begin{equation} \label{e11.6}
    [\gamma_{\tau_3^0}, \gamma_{\tau_2^2}^{< \alpha^{0
    \Downarrow}}[ \mathrm{~~of~the~type~~} \ge \eta^{\ast}
    \mathrm{~~below~~} \alpha^{0 \Downarrow}.
\end{equation}
Really, let us take any cardinal
\[
    \gamma^{\prime} = \gamma_{\tau_2^{\prime}}^{< \alpha^{0
    \Downarrow}} \in SIN_n^{< \alpha^{0 \Downarrow}},
    \tau_2^{\prime} > \tau_2
\]
such that below \ $\alpha^{0 \Downarrow}$
\begin{equation} \label{e11.7}
    \tau_2^{\prime} \not \in  dom(\alpha S_f^{< \alpha^{0
    \Downarrow}}).
\end{equation}
Then by (\ref{e11.5}) there exist some singular matrix \
$S^{\prime}$ \ on its indicated carrier \ $\alpha^{\prime} >
\gamma^{\prime}$ \ with prejump cardinal \ $\alpha^{\prime
\Downarrow}$ \ preserving all \ $SIN_n^{< \alpha^{0
\Downarrow}}$-cardinals \ $\le \gamma^{\prime}$ \ and carrying out
the condition below \ $\alpha^{\prime \Downarrow}$ (recall
definition \ref{8.1.}~2.5\;):
\[
    A_{5.5}^{sc}(\gamma_{\tau_3^0}, \eta^{\ast}, \alpha^{\prime
    \Downarrow}, \alpha S_f^{< \alpha^{\prime \Downarrow}},
    a_f^{< \alpha^{\prime \Downarrow}});
\]
it means, that the whole interval \ $[\gamma_{\tau_3^0},
\alpha^{\prime \Downarrow}[$ \ is covered by blocks below \
$\alpha^{\prime \Downarrow}$ \ of types \ $\ge \eta^{\ast}$. \
Among them there exist the \textit{succeeding} block
\[
    [\gamma_{\tau_3^0}, \gamma_{\tau_2^3}^{< \alpha^{\prime
    \Downarrow}}[ \mbox{~ of the type ~} \ge \eta^{\ast},
\]
so one can treat its subblock \ $[\gamma_{\tau_3^0},
\gamma_{\tau_2^4}^{< \alpha^{\prime \Downarrow}}[$ \ of the type
exactly \ $\eta^{\ast}$.
\\
Remind, the type \ $\eta^{\ast}$ \ is limit, thereafter for every
\ $\gamma_{\tau}$ \ from this subblock there exist many different
matrix carriers admissible for such \ $\gamma_{\tau}$ \ due to
lemma 3.2~\cite{Kiselev9} about restriction; after that due to the
second inductive hypothesis all of them are nonsuppressed for all
such corresponding \ $\gamma_{\tau}$ \ -- and all it below \
$\alpha^{\prime \Downarrow}$.
\\
The same argument works below \ $\alpha^{0 \Downarrow}$ \ and we
return to the situation below this cardinal. From (\ref{e11.7}) it
follows
\[
    \gamma_{\tau_2^4}^{< \alpha^{\prime \Downarrow}} =
    \gamma_{\tau_2^4}^{< \alpha^{0 \Downarrow}}, \quad
    \gamma_{\tau_2^4}^{< \alpha^{0 \Downarrow}}
    \in SIN_n^{< \alpha^{0 \Downarrow}}
\]
and the interval \ $[\gamma_{\tau_3^0}, \gamma_{\tau_2^4}^{<
\alpha^{0 \Downarrow}}[$ \ is really the block of the type \
$\eta^{\ast}$, \ but already below \ $\alpha^{0 \Downarrow}$, \
which contains the admissible disseminator \
$\widetilde{\delta}^0$ \ of the matrix \ $S^0$ \ on \ $\alpha^0$.

\noindent But it provides the contradiction. On one hand, \ $S^0$
\ is admissible for \ $\gamma_{\tau^0}$ \ and then by the closing
condition \ $\mathbf{K}^0$ \ it has the disseminator \
$\widetilde{\delta}^0$ \ with base \ $\rho^0 > \eta^{\ast}$. \ But
on the other hand, the preceding block \ $[\gamma_{\tau_1^0},
\gamma_{\tau_3^0}[$ \ below \ $\alpha^{0 \Downarrow}$ \ has the
same type \ $\eta^{\ast}$ \ and by (\ref{e11.2}) its left end
$\gamma_{\tau_1^0}$ can serve as the admissible disseminator for \
$S^0$ \ on \ $\alpha^0$ \ with the same base \ $\rho^0$, \ and
thereby \ $\widetilde{\delta^0} \le \gamma_{\tau_1^0}$ \ due to
the minimality of \ $\widetilde{\delta}^0$, \ contrary to
supposition (\ref{e11.1}).

2. \ So, (1a.) holds for \ $(\alpha_1^0, \alpha^0, \tau^0)$ \ and
it remains to suppose that (1b.) is wrong for this triple, and we
return to the matter below \ $\alpha_1^0$; \ it means:
\\
there exist some zero matrix \ $S^{0 1}$ \ on the carrier \
$\alpha^{0 1} \neq \alpha^0$ \ admissible for \ $\gamma^0 =
\gamma_{\tau^0}$ \ along with its minimal disseminator \
$\widetilde{\delta}^{0 1}$ \ and the generating disseminator \
$\check{\delta}^{0 1}$ \ with the base \ $\rho^{0 1}$,
\\
but still \ $S^0$ \ on \ $\alpha^0$ \ is suppressed for \
$\gamma^0 = \gamma_{\tau^0}$ \ (below \ $\alpha_1^0$); \ we shall
consider the \textit{minimal} \ $\alpha^{0 1}$ \ with this
property.
\\
Since zero matrix \ $S^0$ \ on \ $\alpha^0$ \ is admissible for \
$\gamma^0$, \ this suppression means, that there holds the
suppression condition (\ref{e11.5}) below \ $\alpha_1^0$, \ that
is for \ $\alpha^{0 \Downarrow}$, \ $\gamma_{\tau_3^0}$ \ replaced
with \ $\alpha_1^0$, \ $\gamma^0$ \ respectively everywhere in
(\ref{e11.5}).

\noindent From here it follows that
\[
    \alpha^0 < \alpha^{0 1},
\]
because if \ $\alpha^0 > \alpha^{0 1}$, \ then the second
inductive hypothesis states, that \ $S^{0 1}$ \ on \ $\alpha^{0
1}$ is nonsuppressed for \ $\gamma^0$ \ below \ $\alpha_1^0$, \
and at the same time it is suppressed by the same suppression
condition. Besides that \ $S^0$ \ on \ $\alpha^0$ \ is the only
matrix admissible for \ $\gamma^0$ \ with the carrier \ $\alpha^0
\in \; ]\gamma^0, \alpha^{0 1}[$ \ due to the minimality of \
$\alpha^{0 1}$.
\\
Now this condition (\ref{e11.5}) with \ $\alpha_1^0$, $\gamma^0$ \
instead of \ $\alpha^{0 \Downarrow}$, $\gamma_{\tau_3^0}$ \
respectively states the existence of the cardinals (we preserve
the previous notations to stress the analogy with the reasoning in
part~1.):
\[
    \gamma^m < \gamma^{\ast} \le \gamma_{\tau_1^0} < \gamma^0
    \mbox{~ and the limit type ~} \eta^{\ast},
\]
holding all the constituents \ $A_{5.1}^{sc} - A_{5.5}^{sc}$; \ in
particular the interval \ $[\gamma^m, \gamma^{\ast}[$ \ is covered
by the maximal blocks of types nondecreasing substantially up to
the limit ordinal \ $\eta^{\ast}$; \ $[\gamma^{\ast},
\gamma_{\tau_1^0}[$ \ is covered by the maximal blocks of the type
exactly \ $\eta^{\ast}$; \ $[\gamma_{\tau_1^0}, \gamma^0[$ \ is
also the block of the same type \ $\eta^{\ast}$ \ -- and so on.

All these conditions define \ $\gamma^m$, $\gamma^{\ast}$,
$\gamma_{\tau_1^0}$, \ $\eta^{\ast}$ \ uniquely through \
$\gamma^0$ \ below \ $\alpha_1^0$ \ and provide the very special
kind of this covering; to operate with it one should use the
following auxiliary \ $\Sigma_n$-formulas treating only the notion
of admissibility (remind definition \ref{8.2.}~5\;):
\\
\quad \\
$
    \alpha \mathbf{K}^1(\gamma): \quad \exists \alpha^{\prime},
    S^{\prime} ~ \alpha \mathbf{K}(\gamma, \alpha^{\prime},
    S^{\prime});
$
\quad \\
\quad \\
$
    \alpha \mathbf{K}^2(\gamma): \quad \exists \alpha^{\prime},
    S^{\prime} ~ \exists \alpha^{\prime\prime}, S^{\prime\prime}
    \big( \alpha^{\prime} \neq \alpha^{\prime\prime} \wedge
$
\[
    \qquad\qquad\qquad
    \wedge \alpha \mathbf{K}(\gamma, \alpha^{\prime},
    S^{\prime}) \wedge \alpha \mathbf{K}(\gamma, \alpha^{\prime\prime},
    S^{\prime\prime}) \big).
\]
The first of them means, that there exist at least \textit{one}
matrix carrier \ $\alpha^{\prime}$ \ admissible for \ $\gamma$; \
the second -- that there exist \textit{more then one} such
carriers \ $\alpha^{\prime} \neq \alpha^{\prime\prime}$; \ thus \
$\neg \alpha \mathbf{K}^1(\gamma)$ \ means, that there is no such
carriers at all.
\\
Since the type \ $\eta^{\ast}$ \ is limit, every maximal block \
$[\gamma_{\tau_1}, \gamma_{\tau_2}[$ \ from the covering of \
$[\gamma^{\ast}, \gamma_{\tau_1^0}[$ \ possesses two properties:
\\
(i) \ if \ $\gamma_{\tau}$ \ is inner in \ $[\gamma_{\tau_1},
\gamma_{\tau_2}[$, $\tau_1 < \tau < \tau_2$, \ then \ $\alpha
\mathbf{K}^2(\gamma_{\tau})$ \ holds; it follows from the second
inductive hypothesis and lemma 3.2~\cite{Kiselev9} about
restriction;
\\
(ii) \ if \ $\gamma_{\tau}$ \ is the end of this block, then \
$\alpha \mathbf{K}^1(\gamma_{\tau})$ \ fails.
\\
It can be verified in the following way. Suppose that \
$\gamma_{\tau}$ \ is the right end, \ $\gamma_{\tau} =
\gamma_{\tau_2}$, \ then the existence of some \ $S^{\prime}$ \ on
\ $\alpha^{\prime}$ \ admissible for \ $\gamma_{\tau_2}$ \
provides the union of this block and of the succeeding block \
$[\gamma_{\tau_2}, \gamma_{\tau_3}[$ \ in the common admissibility
interval \ $[\gamma_{\tau_1}, \gamma_{\tau_3}[$ \ of the type \ $2
\eta^{\ast}$. \ And again by the second inductive hypothesis and
lemma 3.2~\cite{Kiselev9} there exist several matrix carriers \
$\alpha^{\prime}$ \ admissible for \ $\gamma_{\tau_2}$ \ which
become nonsuppressed for \ $\gamma_{\tau_2}$ \ and, so, the
function \ $\alpha S_f$ \ is defined on the whole interval \
$[\tau_1, \tau_3[$, \ though \ $[\gamma_{\tau_1},
\gamma_{\tau_2}[$ \ is the \textit{maximal} block (all it below \
$\alpha_1^0$). \ The left end \ $\gamma_{\tau} = \gamma_{\tau_1}$
\ should be treated in the analogous way.

\noindent Hence for every \ $\gamma_{\tau} \in [\gamma^{\ast},
\gamma^0[$ \ there holds the \ $\Delta_{n+1}$-formula:
\begin{equation} \label{e11.8}
    \alpha \mathbf{K}^2(\gamma_{\tau}) \vee \neg \alpha
    \mathbf{K^1}(\gamma_{\tau});
\end{equation}
it is not hard to see, that the same situation holds below \
$\alpha^{0 1 \Downarrow}$ \ by the same reasons.
\\
Now the generating disseminator \ $\check{\delta}^{0 1}$ \ of \
$S^{0 1}$ \ on \ $\alpha^{0 1}$ \ starts to work and realizes the
restriction-and-extension method.
\\
First,
\[
    \widetilde{\delta}^{0 1} \le \gamma_{\tau_1^0};
\]
in the opposite case
\[
    \gamma_{\tau_1^0} < \check{\delta}^{0 1} =
    \widetilde{\delta}^{0 1} < \gamma^0
\]
and \ $\check{\delta}^{0 1}$ \ extends up to \ $\alpha^{0 1
\Downarrow}$ \ the \ $\Sigma_{n+1}$-proposition
\[
    \forall \gamma_{\tau} > \gamma_{\tau_1^0} \quad \alpha
    \mathbf{K}^2(\gamma_{\tau}).
\]
This fact along with the second inductive hypothesis provides the
definiteness of \ $\alpha S_f^{<\alpha^{0 1 \Downarrow}}$ \ on
some nonempty set \ $T^{\alpha^{0 1 \Downarrow}}$ \ with
properties $(i)$--$(iii)$ specified in theorem~2, contrary to the
minimality of \ $\alpha_1^{\ast}$.
\\
From \ $\widetilde{\delta}^{0 1} \le \gamma_{\tau_1^0}$ \ it
follows
\begin{equation} \label{e11.9}
\gamma^{\ast} < \check{\delta}^{0 1} \le \widetilde{\delta}^{0 1}
\le \gamma_{\tau_1^0}.
\end{equation}
Really, the block \ $[\gamma_{\tau_1^0}, \gamma^0[$ \ obviously
provides the following \ $\Sigma_{n+1}$-formula \
$\varphi(\tau_1^0, \tau^0, \eta^{\ast})$ \ below \ $\alpha^{0 1}$:
\[
    \exists \gamma \big( \gamma_{\tau_1^0} < \gamma_{\tau^0} \le
    \gamma \wedge SIN_n(\gamma) \wedge A_{1.2}^{\vartriangleleft
    \gamma} (\tau_1^0, \tau^0, \eta^{\ast}) \wedge
\]
\[
    \wedge \forall \tau \in \; ]\tau_1, \tau^0[ \quad \alpha
    \mathbf{K}^{2 \vartriangleleft \gamma} (\gamma_{\tau}) \big);
\]
remind, here \ $A_{1.2}(\tau_1, \tau^0, \eta^{\ast})$ \ means,
that the interval \ $[\gamma_{\tau_1^0}, \gamma_{\tau^0}[$ \ has
the type \ $\eta^{\ast}$.
\\
The disseminator \ $\widetilde{\delta}^{0 1}$ \ falls into \
$[\gamma_{\tau_1^0}, \gamma^0[$ \ and hence \
$\widetilde{\delta}^{0 1} = \gamma_{\tau_1^0}$, \ otherwise \
$\widetilde{\delta}^{0 1} < \gamma_{\tau_1^0}$ \ and by lemma
3.2~\cite{Kiselev9} there appear many carriers of \ $S^{0 1}$ \
admissible for this \ $\gamma_{\tau_1^0}$; \ then by the second
inductive hypothesis all of them are nonsuppressed for \
$\gamma_{\tau_1^0}$; \ hence the matrix function \ $\alpha S_f$ \
becomes defined for \ $\tau_1^0$ \ contrary to the minimality of
the left end \ $\gamma_{\tau_1^0}$ \ by definition of the block
notion.

\noindent Since \ $\widetilde{\delta}^{0 1} = \gamma_{\tau_1^0}$,
\ the closing condition \ $\mathbf{K}^0$ \ for \ $S^{0 1}$ \ on \
$\alpha^{0 1}$ \ implies \ $\eta^{\ast} < \rho^{0 1}$ \ for the
base \ $\rho^{0 1}$ \ of \ $\widetilde{\delta}^{0 1}$.
\\
But then the generating disseminator \ $\check{\delta}^{0 1}$ \
with this base restricts the \ $\Sigma_{n+1}$-proposition
\[
    \exists \tau_1^{\prime}, \tau^{\prime} ~ \varphi(\tau_1^{\prime},
    \tau^{\prime}, \eta^{\ast}),
\]
because it contains only constants, bounded by \ $\rho^{0 1}$.
\\
Therefore below \ $\check{\delta}^{0 1}$ \ there appear blocks of
types \ $\ge \eta^{\ast}$ \ (again due to second inductive
hypothesis).
\\
Now if \ $\check{\delta}^{0 1} \le \gamma^{\ast}$, \ then it
violates the condition \ $A_{5.2}^{sc}$ \ about nondecreasing of
covering types of \ $[\gamma^m, \gamma^{\ast}[$ \ up to \
$\eta^{\ast}$ \ substantially.

Thus (\ref{e11.9}) holds. Due to (\ref{e11.8}) below \
$\check{\delta}^{0 1}$ \ there holds the \ $\Pi_{n+1}$-proposition
\[
    \forall \tau \big( \gamma^{\ast} < \gamma_{\tau} \rightarrow
    (\alpha \mathbf{K}^2(\gamma_{\tau}) \vee \neg \alpha
    \mathbf{K}^1(\gamma_{\tau} ) ) \big)
\]
and the disseminator \ $\check{\delta}^{0 1}$ \ extends it up to \
$\alpha^{0 1 \Downarrow}$ \ by lemma 6.6~\cite{Kiselev9} (for \
$m=n+1$, $\delta=\check{\delta}^{0 1}$, $\alpha_0 =
\gamma^{\ast}$, $\alpha_1 = \alpha^{0 1 \Downarrow}$). \ But it
provides the contradiction: (\ref{e11.8}) holds for \
$\gamma_{\tau} = \gamma^0$, \ though there is exactly one matrix \
$S^0$ \ on \ $\alpha^0$ \ below \ $\alpha^{0 1 \Downarrow}$ \
admissible for \ $\gamma^0$.
\\
So, properties (1a.) and (1b.) are carried out. Now one can return
to lemmas \ref{8.5.}~8), \ref{8.8.}, \ref{8.10.} (for \ $\alpha_1
\le \alpha_1^{\ast}$):
\\
(2) First one should dwell on lemma~\ref{8.8.}~1); now it is not
hard to receive (\ref{e8.5}). Toward this end let us compare two
intervals
\[
    ]\gamma_{\tau_1}, \gamma_{\tau_2}[\;, \quad
    ]\gamma_{\tau_1^{\prime}}, \gamma_{\tau_2}[\;.
\]
Due to (\ref{e8.3}), (\ref{e8.4}) \ $\widetilde{\delta}^3$ \
contains in both of them, that provides
\[
    \gamma_{\tau_1^{\prime}} \le \gamma_{\tau_1},
\]
otherwise \ $\gamma_{\tau_1} < \gamma_{\tau_1^{\prime}}$ \ and
(\ref{e8.3}) causes the existence of some matrix \ $S^1$ \ on its
carrier \ $\alpha^1$ \ with the disseminator \
$\widetilde{\delta}^1$ \ and base \ $\rho^1$ \ admissible for \
$\gamma_{\tau_1^{\prime}}$. \ By lemma 3.2~\cite{Kiselev9} \ $S^1$
\ receives its carriers admissible for every \ $\gamma_{\tau} \in
\; ]\widetilde{\delta}^1, \gamma_{\tau_1^{\prime}}]$ \ along with
the same \ $\widetilde{\delta}^1$, $\rho^1$. \ So, there arises
some admissibility interval \ $[\gamma_{\tau_1^{\prime\prime}},
\gamma_{\tau_3}[$ \ with \ $\gamma_{\tau_1^{\prime\prime}} <
\gamma_{\tau_1^{\prime}}$ \ and by (1a.) \ $\widetilde{\delta}^3 =
\gamma_{\tau_1^{\prime\prime}}$ \ contrary to (\ref{e8.4}). Thus
there holds \ $\gamma_{\tau_1^{\prime}} \le \gamma_{\tau_1}$ \ and
along with (\ref{e8.3}), (\ref{e8.4}) it implies (\ref{e8.5}),
that provides the rest part of lemma \ref{8.8.}~1) proof.
\\
Turning to lemma \ref{8.8.}~2), let \ $\alpha$-matrix \ $S$ \ of
characteristic \ $a$ \ on a carrier \ $\alpha$ \ be admissible for
\ $\gamma_{\tau}^{<\alpha_1}$ \ along with its disseminator \
$\widetilde{\delta}$ \ and base \ $\rho$ \ below \ $\alpha_1$; \
one can prove
\[
    \{ \tau^{\prime}: \widetilde{\delta} <
    \gamma_{\tau^{\prime}}^{<\alpha_1} < \gamma_{\tau}^{<\alpha_1}
    \} \subseteq dom(\alpha S_f^{<\alpha_1})
\]
by the reasoning already used above:
\\
for every \ $\gamma_{\tau^{\prime}}^{<\alpha_1} \in \;
]\widetilde{\delta}, \gamma_{\tau}^{<\alpha_1}[$ \ there exist
many admissible carriers of the matrix \ $S$ \ by lemma
3.2~\cite{Kiselev9}, therefore all of them are nonsuppressed due
to (1b.), hence \ $\tau^{\prime} \in dom(\alpha S_f^{<\alpha_1})$.
\\
By the similar reasons in lemma \ref{8.10.} the function \ $\alpha
S_f^{<\alpha_1}$ \ is defined on the whole interval \ $]\alpha
\tau_1^{\ast}, \alpha \tau^{\ast 1}[$ \ for any \ $\alpha_1 >
\alpha \delta^{\ast 1}$, $\alpha_1 \in SIN_n$, \ and \ $\alpha
\delta^{\ast} = \gamma_{\alpha \tau_1^{\ast}}$ \ is the
disseminator of the matrix \ $\alpha S_{\alpha \tau^{\ast
1}}^{<\alpha_1}$ \ on its carrier \ $\alpha_{\alpha \tau^{\ast
1}}^{<\alpha_1}$ \ with the base \ $\alpha \rho^{\ast 1} =
\rho_{\alpha \tau^{\ast 1}}^{<\alpha_1}$.
\\
The similar reasoning in lemma \ref{8.5.}~8) proof should be used.
To finish this proof for the nonsuppressibility one should notice,
that if \ $S$ \ along with \ $\delta$, $\rho$ \ has a carrier \
$\alpha$ \ admissible and nonsuppressed for \
$\gamma_{\tau}^{<\alpha_1}$ \ only in \
$[\gamma_{\tau+1}^{<\alpha_1}, \alpha_1[\;$, \ then \
$\gamma_{\tau+1}^{<\alpha_1}$ \ restricts the \
$SIN_{n-1}^{<\alpha_1}$-proposition
\[
    \exists \alpha^{\prime} \big( \gamma < \alpha^{\prime} \wedge
    \alpha \mathbf{K}_{n-2}(\delta,
    \gamma_{\tau^n}, \gamma_{\tau}^{<\alpha_1},
    \alpha^{\prime}, \rho, S) \big)
\]
for any \ $\gamma \in \; ]\gamma_{\tau}^{<\alpha_1},
\gamma_{\tau+1}^{<\alpha_1}[\;$, \ therefore \ $S$ \ receives many
carriers in \ $]\gamma_{\tau}^{<\alpha_1},
\gamma_{\tau+1}^{<\alpha_1}[$ \ admissible for \
$\gamma_{\tau}^{<\alpha_1}$, \ that also become nonsuppressed for
\ $\gamma_{\tau}^{<\alpha_1}$ \ below \ $\alpha_1$ \ by (1b.).

\begin{sloppypar}
Next lemma shows that intervals \ $[\gamma_{\tau_1}^{<\alpha_1},
\gamma_{\tau_3}^{<\alpha_1}[$ \ of the matrix function \ $\alpha
S_f^{<\alpha_1}$ \ definiteness with the minimal left end \
$\gamma_{\tau_1}^{<\alpha_1} \in SIN_n^{<\alpha_1}$ \ are composed
in a special way: for every \ $SIN_n^{<\alpha_1}$-cardinal \
\mbox{$\gamma_{\tau}^{<\alpha_1} \in \;
]\gamma_{\tau_1}^{<\alpha_1}, \gamma_{\tau_3}^{<\alpha_1}[$} \ the
matrix \ $\alpha S_{\tau}^{<\alpha_1}$ \ has zero characteristic
and disseminators \ $\check{\delta}_{\tau} \le
\widetilde{\delta}_{\tau} = \gamma_{\tau_1}^{<\alpha_1}$ \ below \
$\alpha_1$:
\end{sloppypar}

\begin{lemma}
\label{11.1.} \hfill {}
\\
Let \medskip

(i) \ $A_{1}^{< \alpha_1}(\tau _{1},\tau _{2})$;
\newline

(ii) \ $S^2$ \ be \ $\alpha$-matrix of characteristic \ $a^2$ \ on
a carrier
\[
    \alpha^2 \in \; ] \gamma_{\tau_2}^{<\alpha_1}, \alpha_1 [
\]
admissible for \ $\gamma_{\tau_2}^{<\alpha_1}$ \ below \
$\alpha_1$ \ along with its minimal disseminator \
$\widetilde{\delta}^2$ \ with a base \ $\rho^2$ \ and with
generating eigendisseminator \ $\check{\delta}^{S^2}$;
\\
Then
\[
    \check{\delta}^{S^2} \le \widetilde{\delta}^2 = \gamma_{\tau_1}^{<\alpha_1}
    \mbox{~ and ~} a^3 = 0.
\]
\end{lemma}

\noindent \textit{Proof.} \ The upper indices \ $< \alpha_1,
\vartriangleleft \alpha_1$ \ will be  dropped for shortness as
usual.
\\
Let us consider the \ $\alpha$-matrix \ $S^2$ \ on its carrier \
$\alpha^2$ \ admissible for \ $\gamma_{\tau_2}$ \ along with the
minimal disseminator \ $\widetilde{\delta}^2$ \ of the base \
$\rho^2$ \ and with the generating eigendisseminator \
$\check{\delta}^2 = \check{\delta}^{S^2}$, \ and examine the
situation below the prejump cardinal \ $\alpha^{2 \Downarrow}$.
\\
1. \ Suppose that, on the contrary, this lemma is wrong and \
$\check{\delta^2} \nleq \gamma_{\tau_1}$, \ so it comes
\[
    \gamma_{\tau_1} < \check{\delta}^2 < \gamma_{\tau_2}.
\]
By definition \ $\gamma_{\tau_2} \in SIN_n$ \ and, hence, \
$\gamma_{\tau_2} \in SIN_n^{<\alpha^{2 \Downarrow}}$. \ Due to
this fact and lemma \ref{8.7.} the admissibility below \
$\alpha_1$ \ is equivalent to the admissibility below \ $\alpha^{2
\Downarrow}$ \ for every \ $\gamma_{\tau} \in \; ]\chi^{\ast},
\gamma_{\tau_2}[$.
\\
Then the generating disseminator \ $\check{\delta}^2$ \ extends up
to \ $\alpha^{2 \Downarrow}$ \ the \ $\Pi_{n+1}$-proposition
stating the definiteness of the function \ $\alpha S_f^{<\alpha^{2
\Downarrow}}$, \ as it was several times above, for instance, in
the form:
\begin{equation*}
    \forall \gamma^{\prime}
    \Bigl( \gamma_{\tau_1} < \gamma^{\prime} \wedge
    SIN_{n-1}(\gamma^{\prime}) \rightarrow
    \exists \alpha, S ~
    \alpha \mathbf{K} (\gamma^{\prime}, \alpha, S) \Bigr).
\end{equation*}
After that there appears the function \ $\alpha S_f^{<\alpha^{2
\Downarrow}}$ \ defined on the set with properties $(i)$--$(iii)$
from theorem~2:
\[
    T^{\alpha^{2 \Downarrow}} = \bigl \{ \tau: \gamma_{\tau_1} <
    \gamma_{\tau}^{<\alpha^{2 \Downarrow}} < \alpha^{2 \Downarrow}
    \bigl \},
\]
because for every \ $\tau \in T^{\alpha^{2 \Downarrow}}$ \ there
appear some matrices on many carriers admissible for \
$\gamma_{\tau}^{< \alpha^{2 \Downarrow}}$ \ which are
nonsuppressed for \ $\gamma_{\tau}^{< \alpha^{2 \Downarrow}}$ \
due to arguments (1b.), (1a.); but it contradicts the minimality
of \ $\alpha_1^{\ast}$.
\\
2. \ So, \ $\check{\delta}^2 \le \gamma_{\tau_1}$; \ moreover,
there holds \ $A^0(\tau_1)$. \ Suppose it fails and there exist
some \ $\alpha$-matrix \ $S^1$ \ on its carrier \ $\alpha^1$ \
admissible for \ $\gamma_{\tau_1}$ \ along with its minimal
disseminator \ $\widetilde{\delta}^1 = \gamma_{\tau_1^1}$ \ with
base \ $\rho^1$; \ hence \ $\tau_1^1 < \tau_1$. \ By lemma
\ref{8.8.} \ $]\tau_1^1, \tau_1[\; \subseteq dom(\alpha S_f)$ \
and, so, for every \ $\tau \in \; ]\tau_1^1, \tau_2[$ \ there
exist some \ $\alpha$-matrix on a carrier admissible for \
$\gamma_{\tau}$. \ Therefore there arises some admissibility
interval \ $]\tau_1^{1 \prime}, \tau_2[$ \ with \ $\tau_1^{1
\prime} \le \tau_1^1 < \tau_2$. \ Due to (1a.) \
$\widetilde{\delta}^2 = \gamma_{\tau_1^{1 \prime}}$ \ and by lemma
\ref{8.8.} \ $]\tau_1^{1 \prime}, \tau_2[ \; \subseteq dom(\alpha
S_f)$ \ contrary to the minimality of \ $\tau_1$ \ stated here in
$(i)$.
\\
3. \ Thus \ $A^0(\tau_1)$ \ and, hence, \ $A_1^0(\tau_1, \tau_2)$
\ hold on; thereby lemma~\ref{10.5.} implies \ $a^2=0$. \ At last
again by (1a.) there comes \ $\widetilde{\delta}^2 =
\gamma_{\tau_1}$.
\\
\hspace*{\fill} $\dashv$
\\

Let us turn now to the following suitable notion that already was
used above several times, but further it will play the key role;
thereby it should be emphasized in the following

\begin{definition} \label{11.2.}
\hfill {} \newline \hspace*{1em} Let \ $S$ \ be a matrix on some
carrier \ $\alpha$ \ along with its disseminator
$\widetilde{\delta} < \gamma_{\tau }^{<\alpha _{1}}$ \ with a base
\ $\rho$.
\\

{\em 1)}\quad We say that \ $S$ \ leans on \ $\widetilde{\delta}$
\ on this carrier \ $\alpha $ \ below \ $\alpha_1$ \ if  \
$\widetilde{\delta}$ \ falls in some block \
$[\gamma_{\tau_1}^{<\alpha_1},
\gamma_{\tau_3}^{<\alpha^{\Downarrow}}[ $ \ of type \ $\eta$, \
that is if there exist ordinals \ $\tau _{1}$, $\tau _{1}^{\prime
}$, $\tau _{2}$, $\tau _{3}, \eta$ \ such that
\[
    \gamma _{\tau _{1}}^{<\alpha _{1}} \leq \widetilde{\delta} <
    \gamma _{\tau _{3}}^{<\alpha^{\Downarrow }} \wedge
    A_{4}^{b \vartriangleleft \alpha^{\Downarrow}}
    (\tau _{1},\tau _{1}^{\prime },\tau_{2}, \tau_{3}, \eta).
\]
If in addition this block is the maximal below \
$\alpha^{\Downarrow}$ \ and \ $\eta < \rho $ \ then we say that \
$S$ \ leans on \ $\widetilde{\delta}$ \ very strongly.
\\

{\em 2)}\quad Let a cardinal \ $\upsilon \in \; ]\chi^{\ast},
\alpha^{\Downarrow}]$ \ be covered by blocks; this covering is
called \ $\eta$-bounded (below \ $\alpha^{\Downarrow}$) \ if all
types of its blocks in \ $]\chi^{\ast}, \alpha^{\Downarrow}]$ are
bounded by some constant ordinal \ $\eta < \chi^{\ast +}$:
\\
\quad \\
\hspace*{1em} $\forall \gamma^{\prime} < \upsilon ~ \forall
\tau_1^{\prime}, \tau_2^{\prime}, \eta^{\prime} \big( \chi^{\ast}
< \gamma_{\tau_1^{\prime}}^{<\alpha^{\Downarrow}} <
\gamma_{\tau_2^{\prime}}^{<\alpha^{\Downarrow}} < \upsilon \wedge$
\\
\quad \\ \hspace*{12em} $\wedge A_4^{M b\vartriangleleft
\alpha^{\Downarrow}}(\tau_1^{\prime}, \tau_2^{\prime},
\eta^{\prime}) \rightarrow \eta^{\prime} \le \eta \big).$ \
\\
\hspace*{\fill} $\dashv$
\end{definition}

Part (I) of the following lemma comes out as the carrying
construction of the further reasoning; part (II) will be used at
the endpoint of the proof of theorem 2 strongly. Here one should
remind the notion of stairway and its various attributes, that
were introduced just before lemma \ref{9.5.} by means of formulas
1.--8.; such stairway, being defined below the prejump cardinal \
\mbox{$\alpha_1 = \alpha^{\Downarrow}$} \ of matrix \ $S$ \
carrier \ $\alpha$ \ by means of the formula \ $A_8^{\mathcal{S}t
\vartriangleleft \alpha^{\Downarrow}}(\mathcal{S}t, \alpha
S_f^{<\alpha^{\Downarrow}}, a_f^{<\alpha^{\Downarrow}})$, \ should
be used as the function on \ $\chi^{\ast +}$:
\[
    \mathcal{S}t = \bigl( (\tau_1^{\beta}, \tau_{\ast}^{\beta},
    \tau_2^{\beta}) \bigr) _{\beta < \chi^{\ast +}} ,
\]
so, that for every \ $\beta$, $\beta_1$, $\beta_2$: \\
\quad \\
(i) \ $ \beta < \chi^{\ast +} \rightarrow \tau_{1}^{\beta} <
\tau_{\ast }^{\beta } \leq \tau_{2}^{\beta} \wedge A_{1.1}^{Mst
\vartriangleleft \alpha^{\Downarrow}} (\tau_1^{\beta},
\tau_{\ast}^{\beta}, \tau_2^{\beta}, \alpha
S_f^{<\alpha^{\Downarrow}}, a_f^{<\alpha^{\Downarrow}}),$
\\
\quad \\
that is \ $[ \gamma_{\tau_{1}^{\beta}}^{< \alpha^{\Downarrow}},
\gamma_{\tau_{2}^{\beta}}^{< \alpha^{ \Downarrow}} [$ \ is the
maximal unit step below \ $\alpha^{\Downarrow}$:
\[
    A_{1.1}^{st \vartriangleleft
    \alpha^{\Downarrow}} (\tau_{1}^{\beta}, \tau_{\ast}^{\beta},
    \tau_{2}^{\beta}, \alpha S_f^{<\alpha^{\Downarrow}},
    a_f^{<\alpha^{\Downarrow}}) \wedge
    A_{1.1}^{M \vartriangleleft \alpha^{\Downarrow}}
    (\tau_{1}^{\beta}, \tau_{2}^{\beta}, \alpha S_f^{<\alpha^{\Downarrow}});
\]
(ii) \ $\beta _{1}<\beta _{2}<\chi ^{\ast +}\longrightarrow$
\[
    \longrightarrow \tau _{2}^{\beta _{1}}<\tau _{1}^{\beta
    _{2}}\wedge Od\alpha S_{f}^{<\alpha ^{\Downarrow }} ( \tau
    _{1}^{\beta _{1}},\tau _{\ast }^{\beta _{1}} ) <Od\alpha
    S_{f}^{<\alpha ^{\Downarrow }} ( \tau _{1}^{\beta
    _{2}},\tau _{\ast }^{\beta _{2}} ) ,
\]
that is such steps are disposed successively one after another and
their heights are strictly increasing;
\\
\quad \\
(iii) \ $\sup_{\beta } Od\alpha S_{f}^{<\alpha ^{\Downarrow }} (
\tau _{1}^{\beta },\tau _{\ast }^{\beta } ) =\chi ^{\ast +}$,
\\
\quad \\
that is \ $h(\mathcal{S}t) = \chi^{\ast +}$ \ and heights of these
steps amounts strictly up to \ $\chi^{\ast +}$;
\\
\quad \\
(iv) \ for every maximal unit step \ $[ \gamma_{\tau_{1}}^{<
\alpha^{\Downarrow}}, \gamma_{\tau_{2}}^{< \alpha^{ \Downarrow}}
[$ \ below \ $\alpha^{\Downarrow}$ \ the corresponding triple \
$(\tau_1, \tau_{\ast}, \tau_2)$ \ is the value of this
function.\label{c16}
\endnote{
\ p. \pageref{c16}. \ This condition is not necessary in what
follows, but it is still accepted to make such stairway be single
for some convenience.
\\
\quad \\
} %

\quad \\
Respectively, this stairway \ $\mathcal{S}t$ \ terminates in \
$\upsilon(\mathcal{S}t) = \alpha^{\Downarrow}$, \ if its steps are
disposed cofinally to \ $\alpha^{\Downarrow}$, \ that is if there
holds the property \ $H(\alpha^{\Downarrow})$: \vspace{-6pt}
\begin{multline*}
    \forall \gamma < \alpha^{\Downarrow} \exists \beta < \chi^{\ast
    +} \exists \tau_1^{\beta}, \tau_{\ast}^{\beta}, \tau_1^{\beta},
    \bigl(\gamma < \gamma_{\tau_1^{\beta}}^{<\alpha^{\Downarrow}} <
    \gamma_{\tau_2^{\beta}}^{<\alpha^{\Downarrow}} < \alpha^{\Downarrow}
    \wedge
\\
    \wedge \mathcal{S}t(\beta) = (\tau_1^{\beta},
    \tau_{\ast}^{\beta}, \tau_2^{\beta}) \bigr).
\end{multline*}

\begin{lemma}
\label{11.3.} \hfill {} \newline \hspace*{1em} For every matrix \
$S$ \ of zero characteristic on a carrier \ $\alpha >
\chi^{\ast}$:
\\
\quad \\
(I) \ $S$ \ on \ $\alpha$ \ is provided by some stairway \
$\mathcal{S}t$.
\\
\quad \\
(II) This stairway \ $\mathcal{S}t$ \ terminates in \
$\alpha^{\Downarrow}$, \ that is
\[
    \upsilon(\mathcal{S}t) =
    \alpha^{\Downarrow} = \sup \big\{ \gamma_{\tau_2}^{<\alpha^{\Downarrow}} :
    \exists \beta, \tau_1, \tau_{\ast} ~~ \mathcal{S} t(\beta) = (\tau_1, \tau_{\ast}, \tau_2)
    \big\}.
\]
\end{lemma}

\noindent \textit{Proof} I. \ Let us consider any carrier \
$\alpha_{0}>\chi ^{\ast }$ \ of the matrix \ $S^0$ \ of zero
characteristic on \ $\alpha_{0} $ \ and \
$\alpha^{0}=\alpha_{0}^{\Downarrow} $.
\\
By lemma~\ref{8.5.}~5) there exist \ $\tau _{1}^{\prime }$, $\tau
_{2}^{\prime }$, $\tau _{3}^{\prime }$ \ such that below \
$\alpha^0$
\[
    A_{2}^{0 \vartriangleleft \alpha^0}
    ( \tau _{1}^{\prime }, \tau _{2}^{\prime }, \tau _{3}^{\prime},
    \alpha S_f^{< \alpha^0} ) \wedge
    \forall \tau^{\prime\prime} \in \; ]\tau_1^{\prime},
    \tau_2^{\prime}] ~ a_{\tau^{\prime\prime}}^{<\alpha^0} = 1
    \wedge \alpha S_{\tau _{2}^{\prime }}^{< \alpha
    ^{0}} = S^0
\]
and from here it follows \ $a_{\tau_{2}^{\prime }}^{<\alpha
^{0}}=1$. \ Now let us consider \ $S^0$ \ on the carrier \ $\alpha
_{2}=\alpha _{\tau _{2}^{\prime }}^{<\alpha ^{0}}$ \ and \
$\alpha^{2}=\alpha_{2}^{\Downarrow}$. \ By lemma~\ref{10.3.} there
are \ $\tau _{\ast }^{\prime }$, $\tau _{2}^{\prime \prime }$ \
such, that \ $\tau _{1}^{\prime }<\tau _{\ast }^{\prime }\leq \tau
_{2}^{\prime \prime }$ \ and

\begin{equation} \label{e11.10}
    A_{1.1}^{st \vartriangleleft \alpha^2} (\tau_{1}^{\prime },
    \tau_{\ast}^{\prime}, \tau_{2}^{\prime \prime }, \alpha S_f^{< \alpha^2},
    a_f^{< \alpha^2})
    \wedge A_{1.1}^{M \vartriangleleft \alpha^2}
    (\tau _{1}^{\prime},\tau _{2}^{\prime \prime }, \alpha S_f^{< \alpha^2} );
\end{equation}

\begin{equation} \label{e11.11}
    Od\alpha S_{f}^{<\alpha ^{2}} ( \tau _{1}^{\prime},
    \tau _{\ast }^{\prime } ) > Od ( S^0 ).
\end{equation}
\vspace{0pt}

\noindent Let us enumerate all the triples of ordinals \ $( \tau
_{1}^{\prime },\tau _{\ast }^{\prime },\tau _{2}^{\prime \prime} )
$ \ possessing property (\ref{e11.10}) without omission in the
order of increasing of their first components, that is let us
define the function
\[
    \mathcal{S}t = \bigl( ( \tau _{1}^{\beta },\tau _{\ast }^{\beta },
    \tau_{2}^{\beta } ) \bigr) _{\beta }
\]
with property (i) of the stairway presentation above for \
\mbox{$\alpha^{\Downarrow} = \alpha^2$}; \ statement (ii) comes
after that from corollary~\ref{9.6.} for \ $\alpha_1 = \alpha^2$.
\\
From here and (\ref{e11.11}) it follows that the ordinal \
$Od\alpha S_{f}^{<\alpha ^{2}} ( \tau _{1}^{\beta },\tau _{\ast
}^{\beta } ) $ \ is strictly increasing along with \ $\beta $ \ up
to \ $\chi ^{\ast +}$; \ in the opposite case it is possible to
define below \ $\alpha ^{2}$ \ the upper bound of the set of these
ordinals
\[
    \rho \in \bigl[ Od\alpha S_{f}^{<\alpha ^{2}} ( \tau _{1}^{\prime
    },\tau _{\ast }^{\prime } ) ;\chi ^{\ast +} \bigr[
\]
and then by lemma 4.6~\cite{Kiselev9} about spectrum type
\[
    \rho <Od ( S^0 )
\]
in spite of (\ref{e11.11}). Thus \ $dom(\mathcal{S}t) = \chi^{\ast
+}$; \ statement (iv) is obvious due to the construction of \
$\mathcal{S} t$.
\\
So, statements \ (i)--(iv) are proved for the carrier \ $\alpha
_{2}$ \ of the matrix \ $S^0$ \ and on \ $\alpha_2$ \ it is
provided with this stairway. Then by lemma 5.11~\cite{Kiselev9}
about informativeness the matrix \ $S^0$ \ is provided by some
stairway \ $\mathcal{S}t^0$ \ on its carrier \ $\alpha_0$ \ as
well, because this property is the \textit{inner} property of \
$S^0$ (see comments before lemma~\ref{9.5.}\;).

\noindent II. Turning to statement (II) let us suppose that it is
wrong and this \ $\mathcal{S}t^0$ \ terminates in some cardinal \
$\upsilon^{0} < \alpha^0 = \alpha^{0 \Downarrow}$:
\[
    \upsilon^0 = \sup \big\{ \gamma_{\tau_2}^{<\alpha^0} :
    \exists \beta, \tau_1, \tau_{\ast} ~~ \mathcal{S} t^0(\beta) = (\tau_1, \tau_{\ast}, \tau_2)
    \big\};
\]
evidently, \ $\upsilon^{0}$ \ belongs to \ $SIN_n^{<\alpha^0}$ \
and has the cofinality \ $\chi^{\ast +}$.

\noindent The rest part of the reasoning relies on the method
which may be called \textit{sewing method}; here is its outline
(below \ $\alpha^0$):
\\
Considering some cardinal \ $\upsilon$ \ one can face with the
situation when there are cardinals \ $\gamma_{\tau} < \upsilon$ \
disposed nearby this \ $ \upsilon$ \ and such that the function \
$\alpha S_f$ \ is assumed to be not defined for the corresponding
\ $\tau$; \ therefore such cardinals \ $\gamma_{\tau} < \upsilon$
\ may be called ``holes'' in the set
\[
    \upsilon \cap \big \{ \gamma_{\tau}: \tau \in dom (\alpha S_f) \big \}.
\]
In order to get over this situation and, nevertheless, to see \
$\alpha S_f$ \ be defined for such holes, one should perform the
following two steps:
\\
It should be discovered some \ $\alpha$-matrix \ $S$ \ on a
carrier \ $ \alpha \geq \upsilon$ \ of some characteristic \ $a$ \
along with its generating disseminator \ $\check{\delta}^{\rho} <
\upsilon$ \ and its base \ $\rho$ \ so, that the interval \
$]\check{\delta}^{\rho}, \upsilon[$ \ contains such holes.
\\
At the same time it should be discovered some cardinal
\begin{equation*}
    \gamma^{\delta} \in [ \check{\delta}^{\rho}, \upsilon [ \;  \cap \; SIN_n
\end{equation*}
which along with \ $S$, \ $\rho$ \ \textit{destroys} the premise
of the closing condition
\begin{equation*}
    \mathbf{K}^0 (a, \gamma^\delta, \alpha, \rho)
\end{equation*}
or fulfils its conclusion and therefore holds it on, and cause of
that by lemma 6.8~1)~\cite{Kiselev9} \ $\gamma^\delta$ \ becomes
also the admissible disseminator of \ $ S$ \ on \ $\alpha$ \ with
the same base. Moreover, one should see that this new disseminator
is \textit{admissible and nonsuppressed for every}
\begin{equation*}
    \gamma_{\tau} \in \{ \gamma_{\tau}: \gamma^\delta <  \gamma_{\tau} <
    \upsilon \},
\end{equation*}
because this condition trivially fulfils for \ $\gamma_{\tau}$ \
and, hence, the whole statement of admissibility
\begin{equation*}
    \alpha \mathbf{K} (a, \gamma^{\delta}, \gamma_\tau, \alpha,  \rho, S)
\end{equation*}
holds for many carriers \ $\alpha > \gamma_{\tau}$ \ of the matrix
\ $S$ \ too. Therefore due to argument (1b.) the function \
$\alpha S_f$ \ is found to be defined on the whole set
\begin{equation*}
    \{ \tau: \gamma^\delta < \gamma_\tau < \upsilon \},
\end{equation*}
and thus it happens the ``sewing'' of the interval \ $[
\gamma^\delta, \upsilon [$ \ -- it means, that this set includes
in \ $dom(\alpha S_f)$ \ and this interval does not contain any
holes in spite of the assumption.
\\
The contradiction of this kind will help to advance on the proof
of lemma~\ref{11.3.}, and, hence, the proof of theorem 2 at each
critical stage.
\\

So, let us consider as such \ $\upsilon$ \ the cardinal \
$\upsilon^1 \in SIN_n^{<\alpha^0}$ \ which is \ $\chi^{\ast +}$ \
by order in \ $SIN_n^{<\alpha^0}$, \ that is the set
\[
    \upsilon^1 \cap SIN_n^{<\alpha^0}
\]
has the order type \ $\chi^{\ast +}$; \ this cardinal \
$\upsilon^1 \le \upsilon^0$ \ really exist due to \ $\upsilon^0 <
\alpha^0$, $cf(\upsilon^0) = \chi^{\ast +}$.
\\
Since \ $\upsilon^1 \in SIN_n^{<\alpha^0}$ \ and \ $cf(\upsilon^1)
= \chi^{\ast +}$, \ there exist the \ $\delta $-matrix \ $S^1$ \
of the characteristic \ $a^1$ \ reduced to \ $\chi^{\ast}$ \ and
produced by the cardinal \ $\upsilon ^{1}$ \ on the carrier \
$\alpha^1 < \alpha^0$ \ with the prejump cardinal \ $\alpha^{1
\Downarrow} = \upsilon^1$ \ and the generating eigendisseminator \
$\check{\delta}^1 = \check{\delta}^{S^1} < \upsilon^{1}$ \ with
the base \ $\rho^1 = \rho^{S^1}$ \ by lemma~6.13~\cite{Kiselev9}
(used for \ $m=n+1$, $\alpha_{0}=\upsilon^1$, $\alpha_1 =
\alpha^0$) \ and the function
\[
    f(\beta)= OT(\beta \cap SIN_n^{<\upsilon^1});
\]
we shall consider the minimal \ $\alpha^1$ \ for some
definiteness.

\noindent We shall see, that it provides the contradiction: there
arises certain set
\[
    T^{\upsilon^1} = \{ \tau: \gamma < \gamma_{\tau}^{\upsilon^1} <
    \upsilon^1 \} \subseteq dom(\alpha S_f^{<\upsilon^1})
\]
meeting all the conditions of theorem~2 (for \ $\upsilon^1$ \
instead of \ $\alpha_1$), \ contrary to the minimality of \
$\alpha_1^{\ast}$; \ this contradiction shows, that in fact the
stairway \ $\mathcal{S}t^0$ \ terminates in the cardinal \
$\upsilon^0 = \alpha^0 = \alpha_0^{\Downarrow}$. \ This effect
will be achieved by the sewing method applied to \ $\upsilon^1$.
\\
First, there arises the covering of the interval \
$[\check{\delta}^1, \upsilon^1[$ \ by maximal blocks \ (below \
$\upsilon^1$). \ Suppose it is wrong, then there exist some
cardinal
\[
    \gamma^1 \in SIN_n^{<\upsilon^1} \cap \; ]\check{\delta}^1, \upsilon^1[
\]
which does not belong to any block (below \ $\upsilon^1$). \
Therefore this \ $\gamma^1$ \ can serve as the disseminator
\[
    \widetilde{\delta}^1 = \gamma^1
\]
with the same base \ $\rho^1$ \ by lemma 6.8~\cite{Kiselev9} (for
\ $m=n+1$) \ admissible for every \ $\gamma_{\tau} \in \;
]\widetilde{\delta}^1, \upsilon^1[\;$, \ since there holds the
closing \ $\Delta_1$-condition \ $\mathbf{K}^0(\alpha^1,
\widetilde{\delta}^1, \alpha^1, \rho^1)$

\begin{equation} \label{e11.12}
    \big( a^1 = 0 \rightarrow \forall \tau_1^{\prime},
    \tau_1^{\prime\prime}, \tau_2^{\prime}, \tau_3^{\prime},
    \eta{\prime} < \alpha^{1 \Downarrow} \big[
    \gamma_{\tau_1^{\prime}}^{<\alpha^{1 \Downarrow}} \le
    \widetilde{\delta}^1 < \gamma_{\tau_3^{\prime}}^{<\alpha^{1
    \Downarrow}} \wedge
    \qquad\qquad\qquad\qquad
\end{equation}
\[
    \wedge A_4^{M \vartriangleleft \alpha^{1 \Downarrow}}
    (\tau_1^{\prime}, \tau_1^{\prime\prime},
    \tau_2^{\prime}, \tau_3^{\prime}, \eta{\prime}, \alpha
    S_f^{<\alpha^{1 \Downarrow}}, a_f^{<\alpha^{1 \Downarrow}})
    \rightarrow \eta^{\prime} < \rho^1 \vee \rho^1 = \chi^{\ast +}
    \big] \big)
\]
\vspace{0pt}

\noindent due to the failure of its premise \ $A_4^{M
\vartriangleleft \alpha^{1 \Downarrow}}$. \ Now the sewing method
works: for every \ $\gamma_{\tau} \in \; ]\widetilde{\delta}^1,
\upsilon^1[$ \ and for
\[
    \gamma_{\tau^n} = \sup\{\gamma \le \gamma_{\tau}: \gamma \in
    SIN_n^{<\upsilon^1} \}
\]
there holds the \ $\Pi_{n-2}$-proposition \ $\varphi(a^1,
\widetilde{\delta}^1, \gamma_{\tau^n}, \gamma_{\tau}, \alpha^1,
\rho^1, S^1)$:
\[
    \gamma_{\tau} < \alpha^1 \wedge SIN_n^{<\alpha^{1
    \Downarrow}}(\gamma_{\tau^n}) \wedge
    \alpha \mathbf{K}_{n+1}^{\exists}(a^1, \widetilde{\delta}^1,
    \gamma_{\tau}, \alpha^1, \rho^1, S^1)
\]
stating, that \ $S^1$ \ on \ $\alpha^1$ \ is admissible for \
$\gamma_{\tau}$ \ along with the same \ $a^1$,
$\widetilde{\delta}^1$, $\rho^1$. \ Then the \
$SIN_{n-1}$-cardinal \ $\gamma_{\tau+1}$ \ restricts the \
$\Sigma_{n-1}$-proposition
\begin{equation} \label{e11.13}
    \exists \alpha \big( \gamma_{\tau} < \alpha \wedge \varphi(a^1,
    \widetilde{\delta}^1, \gamma_{\tau^n}, \gamma_{\tau},
    \alpha, \rho^1, S^1) \big)
\end{equation}
and therefore in \ $]\gamma_{\tau}, \gamma_{\tau+1}[$ \ there
appear many admissible for \ $\gamma_{\tau}$ \ carriers \ $\alpha$
\ with this property (\ref{e11.13}) and it provides the
contradiction below \ $\upsilon^1$:
\\
all of them are nonsuppressed for \ $\gamma_{\tau}$ \ due to
благодаря (1b.) and the function \ $\alpha S_f^{< \upsilon^1}$ \
becomes defined for \ $\gamma_{\tau}$ \ and, so, becomes defined
for all interval \ $[\check{\delta}^1, \upsilon^1[$ \ (that is
there happened the sewing of this interval); but it just present
the contradiction -- there appears some set \ $T^{\upsilon^1}$ \
of the function определённости этой функции \ $\alpha S_f^{<
\upsilon^1}$ \ definability with properties со свойствами
$(i)$-$(iii)$ from theorem 2 (for для \ $\alpha_1=\upsilon^1$), \
contrary to the minimality of в противоречии с минимальностью \
$\alpha_1^{\ast}$.
\\So, the
interval \ $[\check{\delta}^1, \upsilon^1[$ \ is covered by
maximal blocks below \ $\upsilon^1$ \ and there holds
\[
    A_{5.1}^{sc \vartriangleleft \upsilon^1} (\gamma^m, \alpha
    S_f^{< \upsilon^1}, a_f^{< \upsilon^1})
\]
stating the covering of \ $[\gamma^m, \upsilon^1[$ \ by the
maximal blocks below \ $\upsilon^1$ \ and the minimality of the
cardinal \ $\gamma^m$ \ with this property (remind definition
\ref{8.1.}~2.1a.,~2.1b.).

\noindent The rest part of this lemma~\ref{11.3.} proof is
conducted below \ $\upsilon^1$ \ and the upper indices \ $<
\upsilon^1$, \ $\vartriangleleft \upsilon^1$ \ and notations of
the functions \ $\alpha S_f^{< \upsilon^1}$, $a_f^{< \upsilon^1}$
\ will be dropped as usual (when the context will point out them
clearly).
\\
Here the final contradiction of this proof comes:
\\
this covering cannot be \ $\eta$-bounded, and at the same time it
have to be \ $\eta$-bounded below \ $\upsilon^1$ \ (remind
definition \ref{11.2.}~2) for \ $\alpha^{\Downarrow} =
\upsilon^1$).

Really, this covering cannot be \ $\eta$-bounded, because in the
opposite case there exist some constant type of its maximal
blocks, disposed cofinally to \ $\upsilon^1$. \ The minimal type \
$\eta^1$ \ of such types is obviously defined below \ $\upsilon^1
= \alpha^{1 \Downarrow}$ \ and by lemma 4.6~\cite{Kiselev9} about
spectrum type there holds
\[
    \eta^1 < Od(S^1) < \rho^1.
\]
Here again the sewing method works. Let \ $[\gamma_{\tau_1},
\gamma_{\tau_2}[$ \ be the maximal block in \
$[\widetilde{\delta}^1, \upsilon^1[$ \ of this type \ $\eta^1$ \
with the minimal left end \ $\gamma_{\tau_1}$, \ then \
$SIN_n$-cardinal \ $\gamma_{\tau_1}$ \ again can serve as the
disseminator \ $\widetilde{\delta}^{1 \prime} = \gamma_{\tau_1}$ \
for \ $S^1$ \ on \ $\alpha^1$ \ with the same base \ $\rho^1$.
\\
And again proposition~(\ref{e11.12}) holds (where \
$\widetilde{\delta}^1$ \ should be replaced with
$\widetilde{\delta}^{1 \prime}$), \ but now because \ $S^1$ \ on \
$\alpha^1$ \ leans on \ $\widetilde{\delta}^{1 \prime}$ \ very
strongly: there is the single maximal block \
$[\gamma_{\tau_1^{\prime}}^{<\alpha^{1 \Downarrow}},
\gamma_{\tau_3^{\prime}}^{<\alpha^{1 \Downarrow}}[\;$, \ that is
precisely \ $[\gamma_{\tau_1}, \gamma_{\tau_2}[$ \ of the type \
$\eta^1$, \ which contains \ $\widetilde{\delta}^{1 \prime}$ \ and
which meet its conclusion \ $\eta^1 < \rho^1$. \ And again there
holds (\ref{e11.13}) for every \ $\gamma_{\tau} \in \;
]\widetilde{\delta}^{1 \prime}, \upsilon^1[$ \ and there arises
some set \ $T^{\upsilon^1}$ \ in the same way, contrary to the
minimality of \ $\alpha_1^{\ast}$.
\\
Thus the covering of \ $[\gamma^m, \upsilon^1[$ \ must not be \
$\eta$-bounded, that is the types of its maximal blocks have to
nondecrease up to \ $\chi^{\ast +}$ \ substantially (below \
$\upsilon^1$):
\begin{equation} \label{e11.14}
    \forall \eta < \chi^{\ast +} ~ \exists \gamma^{\prime} <
    \upsilon^1 ~ \forall \tau_1^{\prime}, \tau_2^{\prime},
    \eta^{\prime} ~ \big( \gamma^{\prime} < \gamma_{\tau_2^{\prime}} \le
    \upsilon^1 \wedge \qquad \qquad
\end{equation}
\[
    \qquad \qquad \wedge A_4^{M b}(\tau_1^{\prime}, \tau_2^{\prime},
    \eta^{\prime}, \alpha S_f^{<\upsilon^1}, a_f^{<\upsilon^1})
    \rightarrow \eta < \eta^{\prime} \big )
\]
otherwise there appears в some constant type of its maximal blocks
disposed cofinally to \ $\upsilon^1$ and leading to the same
contradiction .
\\But it also provides the contradiction: infinitely many values of the
matrix function \ $\alpha S_f^{<\upsilon^1}$ \ become suppressed
(below \ $\upsilon^1$), though they are nonsuppressed by
definition~\ref{8.3.} of this matrix function.

\noindent To see it one should apply the reasoning mode, that may
be called the ``blocks cut-off from the right'' and that consists
in the ``shortening'' of the blocks, which are ``too long'', from
their right ends.
\\
This method works here quite well because this covering is not \
$\eta$-bounded and therefore it is working on these blocks as the
suppressing covering -- more precisely, holding the condition (see
definition \ref{8.1.}~2.4 for \ $X_1 = \alpha
S_f^{<\upsilon^1}|\tau$, $X_2 = a_f^{<\upsilon^1}|\tau$):
\[
    A_{5.4}^{sc}(\gamma_{\tau}, \eta^{\ast}, \alpha S_f^{<\upsilon^1}|\tau,
    a_f^{<\upsilon^1}|\tau)
\]
for infinitely many cardinals \ $\gamma_{\tau}$ \ disposed
cofinally to \ $\upsilon^1$ and for some \ $\gamma^m$,
$\gamma^{\ast}$, $\gamma^1$, $\eta^{\ast}$ .\

\noindent The cardinal \ $\gamma^m$ \ here is already defined
above as the minimal of blocks ends covering of the cardinal \
$\upsilon^1$.
\\
Next, the ordinals \ $\gamma^{\ast}$, $\eta^{\ast}$ \ can be
defined here in different ways, for instance, as limit points of
the following sequences due to (\ref{e11.14}) below \ $\upsilon^1$
\ (the notations \ $\alpha S_f^{<\upsilon^1}$, $a_f^{<\upsilon^1}$
\ will be dropped):
\[
    \gamma_0 = \gamma^m; \qquad\qquad\qquad\qquad\qquad\qquad\qquad\qquad\qquad\qquad
\]
\[
    \eta_i = \sup \big \{ \eta: ~ \exists \tau_1, \tau_2 \big (
    \gamma_{\tau_1} < \gamma_{\tau_2} < \gamma_i \wedge A_4^{M
    b}(\tau_1, \tau_2, \eta) \big) \big \};
\]
\[
    \gamma_{i+1} = \min \big \{ \gamma: ~ \exists \tau_1, \tau_2, \eta \big (
    \gamma_i < \gamma_{\tau_1} < \gamma_{\tau_2} = \gamma \wedge
    \eta_i < \eta \wedge
\]
\begin{equation} \label{e11.15}
    \qquad \qquad \wedge A_4^{M b}(\tau_1, \tau_2, \eta) \wedge
    \forall \tau_1^{\prime}, \tau_2^{\prime}, \eta^{\prime} \big(
    \gamma \le \gamma_{\tau_1^{\prime}} < \gamma_{\tau_2^{\prime}}
    \wedge
\end{equation}
\[
    \qquad \qquad \qquad \wedge A_4^{M b}(\tau_1^{\prime}, \tau_2^{\prime}, \eta^{\prime})
    \rightarrow \eta \le \eta^{\prime} \big ) \big) \big \};
\]
\[
    \qquad\qquad\qquad\qquad
    \eta^{\ast} = \sup_{i \in \omega_0} \eta_i; \quad \gamma^{\ast} =
    \sup_{i \in \omega_0} \gamma_i.
\]
Since the types of maximal blocks of this covering have to
nondecrease up to \ $\chi^{\ast +}$ \  sustantially, there exist
the maximal block in \ $[\gamma^{\ast}, \upsilon^1[$
\[
    [\gamma_{\tau_1^{\ast}}, \gamma_{\tau_2^{\ast}}[ \mbox{ of the greater type }
    \eta^{\ast 1} > \eta^{\ast}
\]
and one should take such block just with the \textit{minimal} left
end \ $\gamma_{\tau_1^{\ast}} > \gamma^{\ast}$. \ It is clear,
that this block includes the initial subinterval
\[
    [\gamma_{\tau_1^{\ast}}, \gamma^{\ast 1}[ \mbox{ exactly of the type }
    \eta^{\ast},
\]
which is also the block (not maximal) with right end \
$\gamma^{\ast 1} \in SIN_n$, \mbox{$\gamma^{\ast 1} =
\gamma_{\tau_2^{\ast 1}}$}.
\\
Hence, there exist the matrix
\[
    S^{\ast 1} = \alpha S_{\tau_2^{\ast 1}}
\]
on its carrier \ $\alpha^{\ast 1} = \alpha_{\tau_2^{\ast 1}}$ \ of
the characteristic \ $a^{\ast 1} = a_{\tau_2^{\ast 1}}$, \ which
is admissible and \textit{nonsuppressed} for \ $\gamma^{\ast 1}$ \
along with its disseminator \ $\widetilde{\delta}^{\ast 1} =
\widetilde{\delta}_{\tau_2^{\ast 1}}$ \ with the base \
$\rho^{\ast 1} = \rho_{\tau_2^{\ast 1}}$ \ by definition (all it
below \ $\upsilon^1$).
\\
But at the same time this matrix \ $S^{\ast 1}$ \ and all its
accessories are, vice versa, suppressed for \ $\gamma^{\ast 1}$ \
because there holds the suppression condition \ $A_5^{S,0}$ \ for
them below \ $\upsilon^1$ \ (remind definition \ref{8.1.}~2.6 and
also (\ref{e11.5})\;), which has now the form:
\[
    a^{\ast 1} = 0 \wedge SIN_n^{<\upsilon^1}(\gamma^{\ast 1})
    \wedge \rho^{\ast 1} < \chi^{\ast +} \wedge \sigma(\chi^{\ast},
    \alpha^{\ast 1}, S^{\ast 1}) \wedge
\]
\[
    \wedge \exists \eta^{\ast}, \tau < \gamma^{\ast 1} \Big(
    \gamma^{\ast 1} = \gamma_{\tau}^{<\upsilon^1} \wedge
    A_{5.4}^{sc \vartriangleleft \upsilon^1} (\gamma^{\ast 1},
    \eta^{\ast}, \alpha S_f^{<\upsilon^1}|\tau, a_f^{<\upsilon^1}|\tau)
    \wedge
\]
\begin{equation} \label{e11.16}
    \wedge \forall \tau^{\prime} \Big( \tau < \tau^{\prime} \wedge
    SIN_n^{< \upsilon^1} (\gamma_{\tau^{\prime}}^{< \upsilon^1}) \rightarrow
    \qquad \qquad \qquad
\end{equation}
\[
    \rightarrow \exists \alpha^{\prime}, S^{\prime} \big [
    \gamma_{\tau^{\prime}}^{<\upsilon^1} < \alpha^{\prime} <
    \gamma_{\tau^{\prime} + 1}^{<\upsilon^1} \wedge
    SIN_n^{<\alpha^{\prime \Downarrow}} (\gamma_{\tau^{\prime}}^{< \upsilon^1})
    \wedge \sigma(\chi^{\ast}, \alpha^{\prime}, S^{\prime}) \wedge
\]
\[
    \qquad\qquad\qquad\qquad\qquad\qquad \wedge A_{5.5}^{sc} (\gamma^{\ast 1}, \eta^{\ast},
    \alpha^{\prime \Downarrow}, \alpha S_f^{<\alpha^{\prime \Downarrow}},
    a_f^{<\alpha^{\prime \Downarrow}} ) \big] \Big) \Big).
\]
Here actually \ $a^{\ast 1} = 0$ \ by lemma \ref{11.1.}; \
$SIN_n^{<\upsilon^1}({\gamma^{\ast 1}})$ \ by definition; \
$\rho^{\ast 1} < \chi^{\ast +}$ \ since \ $\upsilon^1$ \ is \
$\chi^{\ast +}$ \ by order in $SIN_n^{<\upsilon^1}$; \
$\sigma(\chi^{\ast}, \alpha^{\ast 1}, S^{\ast 1})$ \ due to the
admissibility of \ $S^{\ast 1}$ \ on \ $\alpha^{\ast 1}$ \ for \
$\gamma^{\ast 1}$; \ $A_{5.4}^{sc \vartriangleleft \upsilon^1}$ \
holds since types of \ $\gamma^{\ast}$ \ covering are
nondecreasing up to \ $\eta^{\ast}$ \ substantially by
(\ref{e11.15}); and the maximal blocks from the interval \
$[\gamma^{\ast}, \gamma_{\tau_1^{\ast}}[$ \ have the constant type
\ $\eta^{\ast}$ \ due to the minimality of \
$\gamma_{\tau_1^{\ast}}$ -- and it remains only to check the
condition \ $A_{5.5}^{sc}$ \ from (\ref{e11.16}). To this end one
should apply the usual restriction argument:
\\
Every maximal block \ $[\gamma_{\tau_1}, \gamma_{\tau_2}[$ \ in \
$[\gamma^{\ast}, \upsilon^1[$ \ has the type \ $\eta \ge
\eta^{\ast}$ \ due to (\ref{e11.15}) and that is why there holds
the following \ $\Pi_{n-2}$-proposition \ $\psi(\gamma^{\ast 1},
\eta^{\ast}, \alpha^1, S^1, \alpha S_f^{<\alpha^{1 \Downarrow}},
a_f^{<\alpha^{1 \Downarrow}})$:
\[
    \sigma(\chi^{\ast}, \alpha^1, S^1) \wedge A_{5.5}^{sc}
    (\gamma^{\ast 1}, \eta^{\ast}, \alpha^{1 \Downarrow}, \alpha
    S_f^{<\alpha^{1 \Downarrow}}, a_f^{<\alpha^{1 \Downarrow}}),
\]
where \ $A_{5.5}^{sc}$ \ is the \ $\Delta_1$-formula (see
definition \ref{8.1.}~2.5\;):
\begin{multline*}
    \forall \gamma^{\prime} \Big( \gamma^{\ast 1} \le
    \gamma^{\prime} < \alpha^{1 \Downarrow} \rightarrow
    \exists \tau_1^{\prime}, \tau_2^{\prime}, \eta^{\prime} \big(
    \gamma_{\tau_1^{\prime}}^{<\alpha^{1 \Downarrow}} \le
    \gamma^{\prime} < \gamma_{\tau_2^{\prime}}^{<\alpha^{1
    \Downarrow}}\wedge
\\
    \wedge A_4^{M \vartriangleleft \alpha^{1
    \Downarrow}} (\tau_1^{\prime}, \tau_2^{\prime}, \eta^{\prime},
    \alpha S_f^{<\alpha^{1 \Downarrow}}, a_f^{<\alpha^{1 \Downarrow}}
    ) \wedge \eta^{\prime} \ge \eta^{\ast} \big ) \Big).
\end{multline*}
\vspace{0pt}

\noindent Now one should take any \ $\gamma_{\tau^{\prime}}
> \gamma^{\ast 1}$, $\gamma_{\tau^{\prime}} \in
SIN_n^{<\upsilon^1}$; \ by lemma 3.2~\cite{Kiselev9} about
restriction \ $SIN_{n-1}$-cardinal \ $\gamma_{\tau^{\prime}+1}$ \
restricts the \ $\Sigma_{n-1}$-proposition \ $\exists
\alpha^{\prime} ~ \psi_1(\gamma^{\ast 1}, \eta^{\ast},
\alpha^{\prime}, \gamma_{\tau^{\prime}})$, where \ $\psi_1$ is the
formula: \vspace{-6pt}
\begin{multline*}
    \exists S^{\prime} \big [
    \gamma_{\tau^{\prime}} < \alpha^{\prime} \wedge SIN_n^{<\alpha^{\prime
    \Downarrow}}(\gamma_{\tau^{\prime}}) \wedge
\\
    \wedge \psi(\gamma^{\ast 1},
    \eta^{\ast}, \alpha^{\prime}, S^{\prime}, \alpha S_f^{<\alpha^{\prime \Downarrow}},
    a_f^{<\alpha^{\prime \Downarrow}}) \big ],
\end{multline*}
\vspace{0pt} and some carrier \ $\alpha^{\prime}$ \ of matrix \
$S^{\prime}$ \ with this property appears in \
$]\gamma_{\tau^{\prime}}, \gamma_{\tau^{\prime}+1}[$. \ Thus below
\ $\upsilon^1$ \ there holds the proposition:

\[
    \forall \tau^{\prime} \big ( \tau_2^{\ast 1} < \tau^{\prime}
    \wedge SIN_n^{<\upsilon^1}(\gamma_{\tau^{\prime}}) \
    \rightarrow \exists \alpha^{\prime} < \gamma_{\tau^{\prime}+1}
    ~ \psi_1^{\vartriangleleft \upsilon^1} (\gamma^{\ast 1},
    \eta^{\ast}, \alpha^{\prime}, \gamma_{\tau^{\prime}}) \big);
\]
\vspace{-6pt}

\noindent as the result the whole suppression condition
(\ref{e11.16}) is carried out for the matrix \ $S^{\ast 1}$ \ on
its carrier \ $\alpha^{\ast 1}$ \ and it cannot be the value of
the matrix function \ $\alpha S_f^{<\upsilon^1}$ \ below \
$\upsilon^1$ \ contrary to the assumption.
\\
\hspace*{\fill} $\dashv$
\quad \\

Here  the final part of theorem~2 proof comes to its close. Again
all reasonings will be relativized to \ $\alpha _{1}^{\ast }$ \
and therefore the upper indices \ $<\alpha _{1}^{\ast }$,
$\vartriangleleft \alpha_{1}^{\ast }$ \ and the denotations of the
functions \ $\alpha S_f^{<\alpha_1^{\ast}}$,
$a_f^{<\alpha_1^{\ast}}$ \ will be dropped.
\\
By the supposition this theorem fails for the minimal cardinal \
$\alpha_{1}^{\ast }$, \  so there exist \ $\tau _{2}^{\ast }$,
$\tau _{3}^{\ast }$ \ such that there holds
\[
    A_{2} (\tau_{1}^{\ast }, \tau_{2}^{\ast }, \tau _{3}^{\ast } )
\]
where, remind, \ $\tau_1^{\ast}$ \ is the minimal ordinal in \
$T^{\alpha_1^{\ast}}$ \ and where \ $\tau _{2}^{\ast } $ \ is the
minimal ordinal at which the monotonicity of \ $\alpha S_{f}$ \ on
\ $T^{\alpha_1^{\ast}}$ \ is broken; let us consider any arbitrary
great ordinal \ $\tau_3^{\ast}$ \ from the set
\[
    Z^{\ast }=\{ \tau : \gamma _{2}^{\ast } < \gamma _{\tau } <
    \alpha_{1}^{\ast } \wedge \gamma _{\tau } \in SIN_{n} \},
\]
and consider the corresponding cardinals
\[
    \gamma _{i}^{\ast } = \gamma_{\tau _{i}^{\ast }}, ~
    i=\overline{1,3}\;,
\]
and the matrix \ $S^{\ast 2}=\alpha S_{\tau _{2}^{\ast }}$ \ on
the carrier \ $\alpha_{\tau_{2}^{\ast }}$ \ of the characteristic
\ $a^{\ast 2} = a_{\tau_{2}^{\ast}}$ \ with the prejump cardinal \
$\alpha^{\ast 2} = \alpha_{\tau_2^{\ast}}^{\Downarrow}$ \ and its
generating eigendisseminator \ $\check{\delta}^{\ast 2} =
\check{\delta}_{\tau_{2}^{\ast }}^{S}$.
\\
But the main role will be played by the matrix
\[
    S^{\ast 3} = \alpha S_{\tau_3^{\ast}} \mathrm{\ on\ the\ carrier\
    } \alpha_{\tau_3^{\ast}}
\]
for this \ $\tau_3^{\ast} \in Z^{\ast}$ \ with the prejump
cardinal \ $\alpha^{\ast 3} = a_{\tau_3^{\ast}}^{\Downarrow}$, \
and the generating and floating disseminators
\[
    \check{\delta}^{\ast 3} = \check{\delta}_{\tau_3^{\ast}},~
    \widetilde{\delta}^{\ast 3} =
    \widetilde{\delta}_{\tau_3^{\ast}}
    \mathrm{\ with\ the\ base\ } \rho^{\ast 3} =
    \rho_{\tau_3^{\ast}}.
\]

\noindent From lemma \ref{11.1.} (for \ $\tau_1^{\ast}, \tau,
\alpha S_{\tau}, \alpha_{\tau}, \alpha_1^{\ast}$ \ as \ $\tau_1,
\tau_2, S^2, \alpha^2, \alpha_1$) \ it follows:
\begin{equation}
\label{e11.6o} \forall \tau \in Z^{\ast }\  ( a_{\tau }=0\wedge
\widetilde{\delta }_{\tau }=\gamma _{1}^{\ast } ).
\end{equation}

\noindent Now the following cases come:

Case 1.\quad $a^{\ast 2}=1$.\quad Then by lemma~\ref{10.3.} there
exist \ $\tau _{1}^{\ast \prime }$ \ such that
\begin{equation}
\label{e11.7o} A_{3} (\tau _{1}^{\ast }, \tau_{1}^{\ast \prime },
\tau _{2}^{\ast }, \tau_{3}^{\ast } )
\end{equation}
where, remind, the matrix \ $\alpha S_{\tau_1^{\ast \prime}}$ \
has unit characteristic on its carrier \ $\alpha_{\tau_1^{\ast
\prime}}$.
\\
From this moment one should use only the ordinal \ $\tau_3^{\ast}$
\ such that the interval
\[
    [ \gamma_1^{\ast}, \gamma_{\tau_3^{\ast}}[
\]
has the type
\[
    \eta^{\ast 3} > Od(\alpha S_{\tau_{1}^{\ast \prime}});
\]
the existence of such ordinal \ $\tau_3^{\ast}$ \ follows from
conditions $(i)$, $(iii)$ of this theorem~2.

\noindent Let us consider the following subcases:
\\
Case 1a. \quad Suppose that
\[
    \alpha S_{\tau_1^{\ast \prime}} \vartriangleleft \rho^{\ast 3}.
\]
But it is excluded by the restriction-and-extension, providing the
following argument that may be called the ``stairway cut-off from
above'' and that consists in cutting of heights of stairway steps.
\\
The matrix \ $\alpha S_{\tau_1^{\ast \prime}}$ \ of unit
characteristic has the admissible carrier
\[
    \alpha_{\tau_1^{\ast \prime}} \in \; ]
    \widetilde{\delta}_3^{\ast}, \gamma_3^{\ast}[ \; ,
\]
since \ $\widetilde{\delta}_3^{\ast} = \gamma_1^{\ast}$. \ By
lemma 3.2~\cite{Kiselev9} this matrix receives its carriers of
unit characteristic disposed cofinally to \ $\check{\delta}^{\ast
3} \le \widetilde{\delta}^{\ast 3}$ \ and therefore there holds
the \ $\Pi_{n+1}$-proposition below \ $\check{\delta}^{\ast 3}$:
\[
    \forall \gamma ~ \exists \gamma^{1} > \gamma ~
    \exists \delta, \alpha, \rho ~ \big( SIN_{n-1}(\gamma^1) \wedge \alpha
    \mathbf{K}(1, \delta, \gamma^{1}, \alpha, \rho,
    \alpha S_{\tau_1^{\ast \prime}}) \big),
\]
and due to lemma 6.6~\cite{Kiselev9} (for \ $m=n+1$) \ this
disseminator extends this proposition up to the prejump cardinal \
$\alpha^{\ast 3}$ \ and, so, \ $\alpha S_{\tau_1^{\ast \prime}}$ \
receives its admissible carriers of unit characteristic below \
$\alpha^{\ast 3}$ \ disposed cofinally up to \ $\alpha^{\ast 3}$.
\\
By lemma \ref{11.3.} there exist the stairway \ $\mathcal{S}t$ \
terminating in \ $\alpha^{\ast 3}$, \ but by lemma~\ref{9.5.}~2b.
(about stairway cut-off from above, where \ $\alpha_1$, $S^0$ \
are replaced with \ $\alpha^{\ast 3}$, \ $\alpha S_{\tau_1^{\ast
\prime}}$) \ this stairway \ $\mathcal{S} t$ \ is impossible,
because all steps \ $\mathcal{S}t(\beta)$ \ of this stairway
receives heights less then the ordinal
\[
    Od( \alpha S_{\tau_1^{\ast \prime}}) < \chi^{\ast
    +},
\]
though they amount up to \ $\chi^{\ast +}$ \ by definition.
\\
Case 1b. \quad So,
\[
    \rho^{\ast 3} \le Od( \alpha S_{\tau_1^{\ast \prime}} ).
\]
But recall that here is used the type
\[
    \eta^{\ast 3} > Od(\alpha S_{\tau_1^{\ast\prime}})
\]
of the block \ $[\gamma_1^{\ast}, \gamma_3^{\ast}[\;$.
\\
Due to (\ref{e11.7o}) there exist ordinals \ $\tau_{3}^{\prime }$,
$\eta _{3}^{\prime }$, \ such that for \ $\alpha^{\ast 3} =
\alpha_{\tau_{3}^{\ast}}^{\Downarrow}$ \ there holds
\begin{equation} \label{e11.8o}
    A_{4}^{M b \vartriangleleft \alpha^{\ast 3}}
    (\tau_{1}^{\ast}, \tau_{1}^{\ast \prime},
    \tau_{2}^{\ast}, \tau_{3}^{\prime },
    \eta_{3}^{\prime} ).
\end{equation}
These ordinals \ $\tau_{1}^{\ast}, \tau_{1}^{\ast \prime},
\tau_{2}^{\ast}, \tau_{3}^{\prime}, \eta_{3}^{\prime} $ \ are
uniquely defined through \ $\chi^{\ast}$, $\gamma_{1}^{\ast} =
\widetilde{\delta }^{\ast 3}$ \ below \ $\alpha^{\ast 3} $ \ and
it is not hard to see, that
\[
    \tau_3^{\ast} \le \tau_3^{\prime}, ~~ \eta^{\ast 3} \le
    \eta_3^{\prime}.
\]
From the admissibility of \ $S^{\ast 3}$ \ on \
$\alpha_{\tau_3^{\ast}}$ \ and lemma \ref{8.5.}~6) it comes now
the proposition \ $\mathbf{K}^{0} (a^{\ast 3},
\widetilde{\delta}^{\ast 3}, \alpha_{\tau_{3}^{\ast}}, \rho^{\ast
3})$:
\begin{eqnarray*}
    \qquad a^{\ast 3} =0 \longrightarrow \forall \tau _{1}^{\prime },
    \tau_{1}^{\prime \prime },\tau _{2}^{\prime },
    \tau _{3}^{\prime },\eta ^{\prime }
    \left[ \gamma _{\tau _{1}^{\prime }}^{<\alpha^{ \ast 3}}
    \leq \widetilde{\delta}^{\ast 3} <
    \gamma _{\tau _{3}^{\prime}}^{<\alpha ^{\ast 3}}
    \wedge \right. \qquad \qquad
    \\
    \left. \wedge A_{4}^{M b \vartriangleleft \alpha^{\ast 3}}
    (\tau _{1}^{\prime }, \tau_{1}^{\prime \prime },
    \tau _{2}^{\prime },\tau _{3}^{\prime },
    \eta^{\prime }, \alpha S_f^{<\alpha^{\ast 3}}, a_f^{<\alpha^{\ast 3}})
    \rightarrow \eta^{\prime} < \rho^{\ast 3} \vee \rho^{\ast 3} = \chi^{\ast +} \right].
\end{eqnarray*}
Here again the block \ $[ \gamma_{\tau_1^{\prime}}^{<\alpha^{\ast
3}}, \gamma_{\tau_3^{\prime}}^{<\alpha^{\ast 3}}[$ \ is defined
through \ $\widetilde{\delta}^{\ast 3} = \gamma_1^{\ast}$ \
uniquely and therefore
\[
    \widetilde{\delta}^{\ast 3} = \gamma_1^{\ast} =
    \gamma_{\tau_1^{\prime}}^{<\alpha^{\ast 3}}, \quad
    \eta^{\ast 3} \le \eta_3^{\prime} = \eta^{\prime}.
\]
From here, from (\ref{e11.6o}) and (\ref{e11.8o}) it follows that
\ $S^{\ast 3}$ \ on \ $\alpha_{\tau_3^{\ast}}$ \ leans on \
$\widetilde{\delta}^{\ast 3}$ \ very strongly and
\[
    \alpha S_{\tau _{1}^{\ast \prime }} \vartriangleleft
    \rho^{\ast 3}
\]
contrary to the condition of this subcase.
\\
Case 2. \quad $a^{\ast 2}=0$. \ In this case below the prejump
cardinal \ $\alpha^{\ast 2} = \alpha_{\tau_2^{\ast}}^{\Downarrow}$
\ there holds
\begin{equation*}
    \forall \gamma < \gamma _{2}^{\ast } ~ \exists \tau ~
    ( \gamma < \gamma_{\tau}^{<\alpha^{\ast 2}} \wedge
    a_{\tau }^{<\alpha ^{\ast 2}}=1 )  .
\end{equation*}

\noindent This proposition obviously follows from lemma
\ref{11.3.} because there exist some stairway \ $\mathcal{S}t$ \
of unit steps disposed cofinally to \ $\alpha^{\ast 2}$. \ From
here and lemma \ref{10.4.} comes the existence of the ordinal \
$\tau _{1}^{\ast \prime }$ \ for which (\ref{e11.7o}) holds again
(we preserve here the notation from case~1. for some convenience).
It remains to repeat literally the reasoning coming after
(\ref{e11.7o}). The proof of theorem~2 came to its close.
\\
\hspace*{\fill} $\dashv$
\\

\noindent Now let us sum up.
\\
All the reasonings were conducted in the system
\[
    ZF + \exists k ~ (k \mbox{ is weakly inaccessible cardinal});
\]
there was considered the countable standard model
\[
    \mathfrak{M} = (L_{\chi^0}, \in, =)
\]
of the theory
\[
    ZF + V = L + \exists k ~ (k \mbox{ is weakly inaccessible cardinal}),
\]
where any weakly inaccessible cardinal becomes strongly
inaccessible.
\\
In this model the matrix functions were considered; such function
\ $\alpha S_f^{<\alpha_1}$ \ is defined on any nonempty set \
$T^{\alpha_1}$, \ which exist for any sufficiently great cardinal
\ $\alpha_1 < k$, $\alpha_1 \in SIN_n$ \ due to lemma~\ref{8.9.}.
\\
It provides the final contradiction: let us take any \
$SIN_n$-cardinal \ $\alpha_1 > \alpha \delta^{\ast}$ \ limit for \
$SIN_n \cap \alpha_1$ \ of the cofinality \ $cf(\alpha_1) \ge
\chi^{\ast +}$ \ providing such nonempty set \ $T^{\alpha_1}$ \
with properties $(i)$--$(iii)$ from theorem~2, then the function \
$\alpha S_f^{<\alpha_1}$ \ is nonmonotone on this \ $T^{\alpha_1}$
\ by theorem~1 and at the same time is monotone on this set by
theorem~2.
\\
This contradiction ends the Main theorem proof.
\\
\hspace*{\fill} $\dashv$

\newpage

\section{Some Consequences}

\setcounter{equation}{0}

Let us return to the beginning in the Introduction~\cite{Kiselev9}
where the diverse well-known interrelations between Large Cardinal
Hypotheses, Axiom of Determinacy, the regular set properties, etc.
were mentioned (see Drake~\cite{Drake}, Kanamori~\cite{Kanamori}).
Here we shall set forth some easy consequences of such results and
the Main theorem.

\paragraph{I. Large Cardinal Hypotheses}

\hfill {} \newline \hspace*{1em} The hierarchy of large cardinals
arranges them ``by the extent of inaccessibility'' and bases on
the (weakly) inaccessible cardinals. The existence of some of them
(Mahlo cardinals, weakly compact, etc.) is directly forbidden by
the Main theorem. Hence, there are no cardinals possessing
stronger partition properties, for example, indescribable, Ramsey,
Erd\"{o}s cardinals and others; the measurable cardinals do not
exist also since they are Ramsey cardinals. In some cases in the
proof of the inconsistency of Large Cardinal Hypotheses \ $AC$ \
can be used, but it is possible to avoid it, deriving from such
hypotheses the existence of model of \ $ZFC+\exists ~
$inaccessible cardinal (see Silver~\cite{Silver} for example).

We omit the reformulation of this results in terms of filters, trees,
infinitary languages etc.

\paragraph{II. Singular Cardinals. Sharps}

\hfill {} \newline \hspace*{1em} By the nonexistence of
inaccessibles every uncountable limit cardinal is singular. It is
known that in \ $ZFC$ \ every successive cardinal is regular.
Hence any cardinal is singular \textit{iff} it is uncountable and
limit.

It is well-known the remarkable result of Jensen: the negation of
the Inaccessible Cardinal Hypothesis implies the Singular Cardinal
Hypothesis (see also Stern~\cite{Stern}). The non-existence of
inner model with a measurable cardinal implies the Covering Lemma
for the core model \ $K $: \ for every uncountable \ $X \subseteq
On $ \ there is \ $Y \in K $ \ such that \ $X \subseteq Y $ \ and
\ $|X | = |Y |$.

\noindent From here the Singular Cardinal Hypothesis comes out
(Dodd, Jensen~ \cite{Dodd,Dodd1}). So, the nonexistence of
inaccessibles provides this Covering Lemma and the Singular
Cardinal Hypothesis. \newline This situation sheds a new light on
the problem of sharps. \newline It is well-known that the
existence of \ $0^{\sharp} $ \ implies the existence of cardinals
inaccessible in \ $L $ \ (Gitik, Magidor, Woodin~\cite {Gitik}).
Hence \ $0^{\sharp} $ \ does not exist; this result implies the
Covering Lemma for \ $L $ \ by the celebrated theorem of Jensen
(see Devlin, Jensen~\cite{Devlin}).
\newline From here and the famous result of Kunen, establishing
equivalence of the existence of elementary embedding \ $L \prec L
$ \ and the existence of \ $0^{\sharp} $, \ it follows that there
are no elementary embeddings \ $ L \prec L $ \ and, next, there
are no elementary embeddings \ $L_{\alpha} \prec L_{\beta} $ \
with a critical point less than \ $|\alpha | $.

\paragraph{III. Axiom of Determinacy}

\hfill {} \newline \hspace*{1em} It is known that \ $AD$ \ implies
certain Large Cardinal Hypotheses. For example, Solovay had
proved, that \ $AD$ \ implies the measurability of the cardinal \
$\omega _{1}$; \ cardinals \ $\omega _{2}$, \ $\omega _{\omega
+1}$, \ $\omega _{\omega +2} $ \ also are measurable (see also
Kleinberg~\cite{Kleinberg}, Mignone~\cite{Mignone}). Besides that,
\ $AD $ \ implies that cardinals \ $\omega _{1}$, \ $\omega _{2} $
\ are \ $\delta $-supercompact for inaccessible cardinal \ $\delta
$ \ (Becker~\cite {Becker}). Mycielski~\cite{Mycielsky} had proved
that the consistency of
\begin{equation*}
ZF+AD
\end{equation*}
implies the consistency of
\begin{equation*}
ZFC + \exists \mathrm{\ inaccessible\ cardinal.}
\end{equation*}
Hence, \ $AD$ \ is inconsistent and it is possible to define more
exactly its inconsistent part:
\\
Namely, \ $AD ( \Sigma _{2}^{1} ) $ \ is equiconsistent with the
Measurable Cardinal Hypothesis (Louveau~\cite{Louveau}). So, there
exist undetermined \ \mbox{$\Sigma_{2}^{1}$-games}.

\noindent Also, using \ $AC_{\omega}({}^{\omega} \omega) $, \ one
can prove \ $\neg Det (\Pi_{1}^{1}) $ \ (see
Kanamori~\cite{Kanamori}). Therefore in
\begin{equation*}
ZF+AC_{\omega}({}^{\omega} \omega)
\end{equation*}
there exist undetermined \ $\Pi_{1}^{1} $-games, because otherwise
there exist inaccessibles. This result can be hardly improved
because every \ $\bigtriangleup_{1}^{1}$-game is determined
(Martin~ \cite{Martin}).

\paragraph{IV. Continual Sets}

\hfill {} \newline \hspace*{1em} There are well-known series of
outstanding results linking the regular set properties with large
cardinals (by the relative consistency). For example,
Shelah~\cite{Shelar} established the necessity of inaccessible
cardinals for the proposition that every set of reals is Lebesgue
measurable; the measurability of \ $\Sigma _{3}^{1}$-sets implies
the inaccessibility of \ $\omega _{1} $ \ in \ $L $ \ (see also
Raisonnier~\cite {Raisonnier}). From here it follows the existence
of unmeasurable \ $\Sigma_{3}^{1} $-set of reals. Analogously, the
consistence of
\[
    ZF+DC+ \qquad\qquad\qquad\qquad\qquad\qquad\qquad\qquad\qquad\qquad\qquad
\]
\vspace{-18pt}
\[
    \qquad +~(\mathrm{every\ uncountable\ set\ }  \subseteq {}^{\omega}\omega
    \mathrm{\ has\ the\ perfect\ set\ property})
\]
cannot be proved without the Measurable Cardinal Hypothesis (see
Mycielsky~ \cite{Mycielsky}). Besides that, the following
hypotheses are equiconsistent over $ZF$:
\\
\quad \\
1) \hspace{\stretch{0.4}} $AC+\exists \mathrm{\ inaccessible\
cardinal}$; \hspace{\stretch{0.6}}
\\
\quad \\
2) \hspace{\stretch{0.4}} $DC + ~ \mathrm{every\ uncountable\ set\
of\ reals\ has\ the\ perfect\ set\ property}$;
\hspace{\stretch{0.6}}
\\
\quad \\
3) \hspace{\stretch{0.4}} $\omega_{1} \mathrm{\ is\ regular\ } + \
\forall a \in {}^{\omega} \omega \left ( \omega_{1}^{L[a]}<
\omega_{1} \right )$ \hspace{\stretch{0.6}}
\\
\quad \\
(see Solovay~\cite{Solovay}, Specker~\cite{Specer},
Levy~\cite{Levy}). So, $DC$ implies the existence of a set of
reals without the perfect set property and the regularity of \
$\omega_{1} $ \ implies

\begin{equation*}
  \omega_{1}^{L[a]} \nless  \omega_{1}
\end{equation*}
\vspace{0pt}

\noindent for some \ $a \in {}^{\omega} \omega$. \ Also it is
well-known that the Weakly Compact Cardinal Hypothesis is
equiconsistent (over \ $ZFC~$) to the statement about the regular
properties of all continual sets in \ $ZF+MA$ \ (Harrington,
Shelah~\cite{Harrington}). Hence, the Main theorem implies the
existence of uncountable sets \ $\subseteq \omega_{\omega}$ \
without regular properties.

\noindent These results should be made more precise; for example,
Solovay~ \cite{Solovay1} proved that for any \ $a \in
{}^{\omega}\omega$ \quad $ \omega_{1}^{L[a]} < \omega_{1}$ \ is
equivalent to the perfect set property of every \ $\Pi_{1}^{1}(a)
$-set of reals. From here it comes now that the regularity of \
$\omega_{1} $ \ implies the existence of \ $\Pi_{1}^{1}(a) $-set
of reals without the perfect set property for some \ $a \in
{}^{\omega}\omega $.

\paragraph{V. Martin Axiom}

\hfill {} \newline \hspace*{1em} The consequences of the results
of Harrington, Shelah~\cite {Harrington} and the Main theorem
mentioned above should be formulated more exactly, for example, in
the following way:

\noindent $MA$ \ implies the existence of \ $\Delta _{3}^{1}$-sets
without the Baire property and unmeasurable \ $\Sigma
_{3}^{1}$-sets \ $\subseteq {}^{\omega}\omega $.

\noindent Kanovei~\cite{Kanovei} had proved in
\begin{equation*}
ZFC+MA+\left| R\right| >\omega _{1}+\forall x\subseteq  \omega _{1}~~\omega
_{1}^{L\left[ x\right] }=\omega _{0}
\end{equation*}
the consistency of
\begin{equation*}
ZFC+\exists ~ \mathrm{\ inaccessible\ cardinal}.
\end{equation*}
Hence, \ $MA$ \ is inconsistent in
\begin{equation*}
ZFC+\left| R\right| >\omega _{1}+\forall x\subseteq  \omega _{1}~~\omega
_{1}^{L\left[ x\right] }=\omega _{0}.
\end{equation*}

The host of other consequences of the Main theorem is too long to
enumerate; so, the author intends to represent the more detailed
analysis of such consequences in publications forthcoming.

\theendnotes

\newpage

\thispagestyle{empty}

\begin{center}
\quad \\
\quad \\
\quad \\
\quad \\
\quad \\
\quad \\
\quad \\
{Scientific edition} \\
\quad \\
{\bf Kiselev} Alexander\\
\quad \\
{\bf INACCESSIBILITY \\
AND \\
SUBINACCESSIBILITY}\\
\quad \\
{In two parts}\\
{Part II}\\
\quad \\
{\small Responsible for release {\em T. E. Yanchuk}} \\
\quad \\
{\footnotesize Signed for publication on 01.04.2010. Format
60$\times$84$^{1}/_{16}$. Paper offset. Risograph. Conditional
quires 8,6. Registration-publishing sheets 7,7. Circulation is 100
ex. Order \selectlanguage{russian} \No
\selectlanguage{english} 610.} \\
\quad \\
{\footnotesize Republican unitary enterprise} \\
{\footnotesize ``Publishing center of Belarusian State University''} \\
{\footnotesize \selectlanguage{russian} ЛИ \No
\selectlanguage{english}
02330/0494361 on 16.03.2009.} \\
{\footnotesize Krasnoarmeyskaya Str., 6, 220030, Minsk.}\\
\quad \\
{\footnotesize Printed from the author's layout in Republican
unitary enterprise} \\
{\footnotesize ``Publishing center of Belarusian State University''} \\
{\footnotesize \selectlanguage{russian} ЛП \No
\selectlanguage{english} 02330/0494178 on 03.04.2009.} \\
{\footnotesize Krasnoarmeyskaya Str., 6, 220030, Minsk.}\\
\end{center}
\label{end}

\end{document}